\documentclass[11pt,reqno]{amsart}
\usepackage[T1]{fontenc}
\usepackage[utf8]{inputenc}
\usepackage[english]{babel}
\usepackage{microtype}
\usepackage[a4paper,margin=2.6cm]{geometry}
\usepackage{graphicx}
\usepackage{tikz}
\usetikzlibrary{arrows.meta,positioning,calc}
\usepackage{cancel}
\usepackage{eurosym}
\usepackage{enumerate}
\usepackage{mathtools}
\usepackage{hyperref}
\usepackage{etoolbox}
\makeatletter
\renewcommand{\paragraph}{\@startsection{paragraph}{4}%
  \z@{.35\linespacing\@plus.5\linespacing}{-.5em}%
  {\normalfont\bfseries}}
\makeatother

\numberwithin{equation}{section}

\theoremstyle{plain}
\newtheorem{theorem}{Theorem}[section]
\newtheorem{lemma}[theorem]{Lemma}
\newtheorem{proposition}[theorem]{Proposition}
\newtheorem{corollary}[theorem]{Corollary}

\theoremstyle{definition}
\newtheorem{definition}[theorem]{Definition}
\newtheorem{assumption}[theorem]{Assumption}
\newtheorem{hypothesis}[theorem]{Hypothesis}
\newtheorem{example}[theorem]{Example}

\theoremstyle{remark}
\newtheorem{remark}[theorem]{Remark}

\newcounter{proofpart}
\newenvironment{proofpart}[1][]
  {\refstepcounter{proofpart}\par\smallskip\noindent\textit{(\roman{proofpart})%
  \if\relax\detokenize{#1}\relax\else\ #1.\fi}\ }
  {\par\smallskip}

\makeatletter
\tagsleft@false
\@fleqnfalse
\def\@setauthors{%
  \begingroup
  \def\thanks{\protect\thanks@warning}%
  \trivlist
  \centering\normalsize \@topsep30\p@\relax
  \advance\@topsep by -\baselineskip
  \item\relax
  \author@andify\authors
  \def\\{\protect\linebreak}%
  \authors
  \ifx\@empty\contribs
  \else
    ,\penalty-3 \space \@setcontribs
    \@closetoccontribs
  \fi
  \endtrivlist
  \endgroup
}
\AtBeginEnvironment{proof}{\setcounter{proofpart}{0}}
\makeatother

\hypersetup{
  pdftitle  = {Geometric foundations of mesoscopic nucleation theory: commutativity of kinetic and probabilistic reduction for the Dean--Kawasaki equation},
  pdfauthor = {Miguel A. Dur\'{a}n-Olivencia},
}

\title[Geometric foundations of mesoscopic nucleation theory]{Geometric foundations of mesoscopic nucleation theory: commutativity of kinetic and probabilistic reduction for the Dean--Kawasaki equation}
\author[M. A. Dur\'{a}n-Olivencia]{Miguel A. Dur\'{a}n-Olivencia$^{1,2}$\\
{\small $^{1}$Imperial College London, South Kensington Campus, London SW7 2AZ, United Kingdom}\\
{\small $^{2}$Research, Vortico Tech, M\'{a}laga 29100, Spain}}
\date{}

\begin{document}
\begin{abstract}
The reduction of a stochastic field equation to a few degrees of freedom is a central problem in many areas of physics, applied mathematics and engineering.
In this work we present two routes for the dimensional reduction of conservative stochastic partial differential equations of Dean--Kawasaki type.
This problem is of particular interest in mesoscopic nucleation theory, which is formulated in terms of the cumulative mass, the integrated form of the density field. 
The stochastic equation for the mass constitutes the basis for the derivation of the whole theory. 
To our knowledge, such reductions have always been carried out at the level of the trajectories, i.e.\ on Dean's equation, and never at the level of the probability distribution, i.e.\ on Kawasaki's equation. 
The first route starts from the parametrisation of the density introduced by Lutsko [J. F. Lutsko, J. Chem. Phys. 136, 034509 (2012)], who imposed it on the field equation and obtained stochastic equations for the order parameters by integrating against a suitable weight. 
Here, this parametrisation is turned into a projection of the trajectories onto the manifold of order parameters, formulated within a precise functional and geometric setting.
The second route, which reduces Kawasaki's equation, is new. 
This route averages the functional forward Kolmogorov equation over the fast transverse degrees of freedom, using Fermi tubular coordinates around the manifold.
In both cases the metric on the manifold is not postulated but derived, once from the noise and once from the Wasserstein geometry of the parent equation, and the two results agree.
Since the equations of Dean and Kawasaki are formally equivalent, one would expect both reductions to be equivalent as well. 
We show that this is not the case. 
In particular, we prove that the two reduced generators coincide if and only if the noise-induced drift corrections arising along each route, which we call anomalies, are equal, and that this condition fails in general. 
For noise regularised over a correlation length $\delta$ on a spherical domain, the mismatch diverges as $\delta^{-1}$ for the ideal gas and, under a closure of the transverse dynamics, as $\delta^{-3}$ for a fluid with a local excess free energy, the latter order being set by the closure. 
Moreover, in the limit of infinite transverse stiffness, the second route recovers, at the order retained, the generator that the full dynamics with the local covariance induces on the manifold, while the first one keeps a drift of order $\delta^{-1}$.
Finally, we show that, under the same closure, regularity hypotheses and the normalisation $\varepsilon^2 = N^{-1}$, the two reduced theories differ after a time $t$ by at most $\mathcal{O}(t/n\delta^3)$, where $n$ is the mean number density.
Thus, over any fixed time, both routes yield the same predictions when the correlation volume contains many particles, which is the condition underlying fluctuating hydrodynamics.
Overall, the present framework provides a geometric foundation for mesoscopic nucleation theory, where the metric, the orthogonality of the projection and the transformation of the noise follow from the theory instead of being fixed by convention.
\end{abstract}
\maketitle

\section{Introduction}
\label{sec:introduction}

\subsection{From equilibrium density functional theory to its deterministic dynamics}

The equilibrium statistical mechanics of a classical fluid in an external field rests on a variational principle.
Evans~\cite{evans1979nature} showed that the grand potential of an inhomogeneous fluid is a functional of the one-body density alone ($\Omega[\rho]$), and that the equilibrium profile is the one minimising it.
However, exact as it is, the result of Evans is not constructive.
What it leaves open is the excess free energy ($\mathcal{F}_{\mathrm{exc}}[\rho]$), which is only known as a series in the correlation functions of the homogeneous fluid, and since then a huge amount of work has been devoted to approximating it.
An exception to this approximation is due to Percus~\cite{percus1976equilibrium}, who obtained the exact functional for hard rods in one dimension. 
Several years later, Rosenfeld~\cite{rosenfeld1989free} constructed the fundamental measure theory of hard spheres as a generalisation of scaled particle theory, which turned out to be the natural extension of Percus' result to three dimensions, since it requires no mapping onto an effective liquid. 
Tarazona, Cuesta and Mart\'{\i}nez-Rat\'{o}n~\cite{tarazona2008density} then reviewed that construction for general hard-particle and liquid-crystalline systems (attractive interactions being treated as a perturbation of the hard-core reference). 
A thorough review of these density functionals can be found in Ref.~\cite{lutsko2010recent}.

We note, however, that density functional theory says nothing about how a fluid reaches equilibrium. 
To this end, we need to postulate a dynamics.
We adopt here the one presented by Lutsko~\cite[\S VI.A]{lutsko2010recent}, which combines mass conservation (which is exact) with a constitutive relation of linear response, in which the current follows the gradient of the local chemical potential.
Local equilibrium then fixes that chemical potential as the functional derivative of the free energy. 
For the mobility, we take a coefficient proportional to the density (as required to recover the diffusion of a dilute tracer in a bath). 
The outcome is the well-known equation of dynamical density functional theory,
\begin{equation}
    \partial_t \rho = \nabla \cdot \!\bigl(\rho\,\nabla \tfrac{\delta \mathcal{F}}{\delta \rho}\bigr).
    \label{eq:ddft}
\end{equation}

Marconi and Tarazona~\cite{marconi1999dynamic} derived the same equation from the particle dynamics.
The derivation started from the Langevin equations of the individual particles and gave an evolution equation for the averaged density. 
Unfortunately, such an equation is not closed, since it involves the two-point correlation of the fluid out of equilibrium. 
The authors closed it by requiring that this correlation keep its equilibrium form at every instant, which is the adiabatic assumption in its sharpest version, and they emphasised that the resulting equation is deterministic. 
In the following years, the framework kept evolving to include inertia, active particles and mixtures, and te~Vrugt, L\"{o}wen and Wittkowski~\cite{tevrugt2020classical} wrote a review of it.

But why should the mobility be the density? 
Interestingly, the answer came from mathematics.
Jordan, Kinderlehrer and Otto~\cite{jordan1998variational} showed that the Fokker--Planck equation can be solved by a sequence of minimisations. 
Each step trades the free energy against the Wasserstein-2 distance travelled since the previous one, and the steps converge to the Fokker--Planck dynamics. 
Ever since the pioneering work of Otto~\cite{otto2001geometry} on the porous medium equation, which showed that such equations are gradient flows for the 2-Wasserstein metric, there has been an explosion of interest in their formal study.
Within that formal Riemannian calculus, the mobility operator ($\rho\,\nabla(\cdot)$) is nothing else than the inverse metric tensor, and hence dynamical density functional theory is a gradient flow in the strict sense. 
Rigorous versions on general metric spaces came later, with the work of Ambrosio, Gigli and Savar\'{e}~\cite{ambrosio2008gradient}. 
For the transport side, the work of Villani~\cite{villani2009optimal} is the standard reference.

More recently, Bruna, Burger and Carrillo~\cite{bruna2020coarse} checked whether this structure survives coarse-graining for the particular case of the many-particle limit of hard-core interacting diffusions. Indeed, the structure survives. Moreover, the macroscopic gradient flow comes out together with its metric, order by order, from the high-dimensional Wasserstein gradient flow of the particles. In this work we ask the very same question, although for a different kind of reduction (a projection onto a few order parameters at fixed number of particles), and for the stochastic equation instead of the deterministic one.

There is nevertheless one thing which Eq.~\eqref{eq:ddft} cannot do. Since it relaxes monotonically towards a minimum of $\mathcal{F}$, it cannot describe the passage from a metastable basin to a stable one, a passage which requires the system to climb a free-energy barrier and which only fluctuations can supply. As Lutsko~\cite{lutsko2019how} observed, a theory of nucleation requires two elements (a means of computing the free energy of a cluster, and a dynamical description of the fluctuations). The classical theory of nucleation supplies each of them separately, the first through the capillarity approximation (which assigns a macroscopic surface tension to a cluster of microscopic size), and the second through Becker--D\"{o}ring kinetics. We now know how the first element relates to density functional theory, since the classical barrier was shown~\cite{lutsko2011density} to be what the density functional barrier becomes in the limit of vanishing interface width relative to the cluster radius. The same work showed also that the classical picture admits a single degree of freedom (the radius). This makes the nucleation pathway, by construction, a monotone growth of that radius, and whatever the radius does not record is not recorded at all. But how can we improve this description? As we discuss next, the fluctuations have to be restored to the field equation itself.

\subsection{Fluctuating hydrodynamics and the Dean and Kawasaki equation}

Once the fluctuations are restored, we arrive at the Dean and Kawasaki class of stochastic partial differential equations,
\begin{equation}
    d\rho = \nabla \cdot \!\bigl(\rho\,\nabla \tfrac{\delta\mathcal{F}}{\delta\rho}\bigr)\,dt \;+\; \varepsilon\,\nabla\cdot\!\bigl(\sqrt{\rho}\;\boldsymbol{\eta}\bigr).
    \label{eq:dean-kawasaki}
\end{equation}
Here $\varepsilon$ sets the noise intensity and $\boldsymbol{\eta}$ is a spatially regularised Gaussian field (which we make precise in Definition~\ref{def:mollifier}). Note that the noise enters in divergence form, which means that mass is conserved along every realisation. Moreover, its amplitude vanishes with the density (a region with no particles cannot show density fluctuations).

The two names of the equation are due to two independent derivations, which do not lead to the same object. On the one hand, Dean~\cite{dean1996langevin} took the It\^{o} differential of the empirical measure of $N$ interacting Brownian particles and found that, at finite $N$ and without any approximation, it obeys the white-noise form of Eq.~\eqref{eq:dean-kawasaki}. Note that in this derivation the multiplicative noise emerges from the calculation itself, without being imposed. Kawasaki~\cite{kawasaki1994stochastic}, on the other hand, worked with the law, and by coarse-graining the $N$-body Smoluchowski equation onto cells large enough to hold many particles obtained the functional Fokker--Planck equation for the probability distribution $P[\rho,t]$ of the density field. The same equation can also be obtained by a second route, which starts from nonlinear fluctuating hydrodynamics for the density and the momentum density, eliminates the momentum with a Zwanzig projector and treats its relaxation as fast. Both ways agree once the diffusion constant is identified as $D_0 = k_BT\tau/m$, $\tau$ being the momentum relaxation time.

Both objects are related exactly by the It\^{o} formula (Kawasaki's equation is nothing else than the forward Kolmogorov equation of Dean's). Hence, at the level of the full field, there is no loss of information in going from the trajectory to its law, and we may use either description. Does this remain true after we carry out a dimensional reduction? In the present work we give an answer to this question.

For some time, the meaning of the density field in these equations was not clear. In one line of work a deterministic equation had been obtained for a quantity called the density, while in the other the equation obtained was stochastic. The matter was settled by Archer and Rauscher~\cite{archer2004dynamical}, who noted that both are correct and that the answer depends on what is meant by the one-body density. The ensemble-averaged density obeys a deterministic equation, while a spatially or temporally coarse-grained density obeys a stochastic one. Later, Archer~\cite{archer2009dynamical} derived dynamical density functional theory from the Kramers equation for the phase-space distribution, obtaining a generalisation of the Euler equation that resolves the density down to the scale of the individual particles. Dur\'{a}n-Olivencia et~al.~\cite{duran2017framework} carried the programme to the Hamiltonian level. The bath was eliminated by Zwanzig projection, which gave a generalised Dean and Kawasaki model for the density and momentum fields, and a nonequilibrium ensemble average was then performed under a decomposition into local-equilibrium and nonequilibrium parts. The result is a fluctuating dynamical density functional theory for mesoscopic fields, i.e.\ fields that are neither microscopic nor macroscopic, since they have been averaged over an ensemble of initial conditions while the noise has been kept. The same work identifies deterministic dynamical density functional theory as the most likely path of the fluctuating equation in the weak-noise limit. This is the mesoscopic level at which we carry out the present work. For a recent survey of the Dean and Kawasaki equation, and of its uses in active matter, electrolytes and mode-coupling theory, see Illien~\cite{illien2024deankawasaki}.

When it comes to its distributional white-noise form ($\delta \to 0$), the mathematical status of Eq.~\eqref{eq:dean-kawasaki} is delicate. First, the noise coefficient $\sqrt{\rho}$ vanishes at $\rho = 0$. Second, it acts as a distribution-valued operator, so that classical It\^{o} theory cannot be applied directly. Unfortunately, the rigidity result of Konarovskyi, Lehmann and von Renesse~\cite{konarovskyi2019dean} goes the opposite way to what one might have hoped. Measure-valued solutions exist only for a restricted set of parameters, and whenever one exists it is the empirical measure of a finite system of Langevin particles with mean-field interaction. In other words, the distributional equation has no genuinely mesoscopic solutions beyond the finite particle systems it was meant to coarse-grain. Cornalba, Shardlow and Zimmer~\cite{cornalba2019regularised,cornalba2020dean} restored high-probability well-posedness by regularising the noise in space, replacing each particle by a narrow Gaussian instead of a Dirac mass. Throughout this work we also regularise the noise, at finite correlation length $\delta > 0$, although their existence theory (which is for an inertial model, and was later extended to several space dimensions~\cite{cornalba2021slender}) does not cover the equation studied here (Assumption~\ref{ass:parent-equation}).

Since then, the analysis has advanced on several fronts, and it is the boundaries drawn by this analysis that make our question well posed. M\"uller, von Renesse and Zimmer~\cite{muller2024wellposedness} extended the well-posedness framework to Dean--Kawasaki models of Vlasov--Fokker--Planck type, covering second-order Langevin dynamics with inertia, and found the same existence dichotomy in that broader class, which includes active matter and flocking models. Hao, Wu and Zimmer~\cite{hao2025kinetic} then proved strong well-posedness of the Vlasov--Fokker--Planck--Dean--Kawasaki system with regular correlated noise, by means of renormalised kinetic solutions. Fehrman and Gess~\cite{fehrman2021well,fehrman2024deankawasaki} established well-posedness for nonlinear diffusion equations driven by nonlinear conservative noise on the torus, the Dean--Kawasaki equation with correlated noise among them, and Fehrman~\cite{fehrman2025neumann} did the same on a bounded domain with a reflecting wall, in each case for the Stratonovich equation with a drift local in the density, a class that does not cover the equation studied here (Assumption~\ref{ass:parent-equation}). In addition, the variational structure of the equation has been studied in a direct way. Dirr, Stamatakis and Zimmer~\cite{dirr2016entropic} developed entropic and gradient-flow formulations for nonlinear diffusion where the noise operator of Eq.~\eqref{eq:dean-kawasaki} is formally dictated by the same Riemannian metric that defines the deterministic drift. Dur\'{a}n-Olivencia et~al.~\cite{duran2019thinfilms} studied a conservative equation of the same gradient-flow class for thin liquid films, imposing the detailed-balance condition for state-dependent multiplicative noise and developing spectral methods for thermally activated rupture. All these works deal with the well-posedness of regularised Dean--Kawasaki models. However, in none of them is the geometric structure of the equation studied after the reduction has been carried out. This is the problem that we address in the present work.

\subsection{The mesoscopic nucleation programme}

Lutsko~\cite{lutsko2012mean} formulated a dynamical theory of nucleation in the framework of fluctuating hydrodynamics, in which the starting point is the stochastic field equation (Eq.~\eqref{eq:dean-kawasaki}) instead of the deterministic gradient flow. The physical idea is that the density profiles visited during a nucleation event can be described, to good accuracy, by a small number of collective coordinates, so one writes $\rho(r;t) = \rho\bigl(r;\mathbf{X}(t)\bigr)$, substitutes it into the field equation and integrates over the spatial variable, which turns the field equation into a system of stochastic differential equations for $\mathbf{X}$. The construction was developed for an arbitrary number of order parameters and specialised to a single one, the cluster radius, only in the closing illustration of the paper.

Lutsko and Dur\'{a}n-Olivencia~\cite{lutsko2013classical} then showed that classical nucleation theory emerges as a limiting case of this dynamical framework. A subsequent work~\cite{lutsko2015two} considered a second coordinate, namely the density inside the cluster, which in the capillarity approximation is fixed to its bulk value. Allowing this density to vary led to a two-dimensional Langevin system able to resolve non-classical effects. In the case of confined geometries, the finite mass available imposes a largest supercritical cluster and even forbids nucleation once the container is small enough~\cite{duran2015mesoscopic}. Dur\'{a}n-Olivencia et~al.~\cite{duran2018general} introduced the most likely path formulation, where the dynamics reduces to a single coordinate (the arc-length along the nucleation pathway). This body of work is known as mesoscopic nucleation theory. It was shown to improve substantially on the classical theory for complex fluids~\cite{lutsko2018systematically}.

The present work formalises this reduction, and we briefly recall how it was carried out. The field equation was multiplied by a weight $W_i(r;\mathbf{X})$, integrated over the domain and gave one equation per coordinate, and these equations hold whatever weight is used. The weight finally chosen, $W_i = (4\pi r^2\rho)^{-1}\,\partial_{X^i} m$ (written in terms of the cumulative mass $m$), makes the diffusion matrix proportional to the matrix $g_{ij}$ multiplying $\dot{\mathbf{X}}$, and this is what lets a fluctuation-dissipation relation survive the reduction. The construction is exact when the parametrisation of the density is complete. For a finite number of coordinates, the same equations were kept as an approximation.

That paper also identified the correction that is the main concern of this paper. Writing the reduced equations with the ordinary chain rule gives a Stratonovich system. Unfortunately, the two readings are not equivalent for the reduced system, and the conversion produces a residual drift of order $\varepsilon^2$, the so-called \emph{anomalous force}, which, as we will see, is only one of the two pieces of the geometric anomaly $\mathfrak{s}_a$. The other piece, the projected field correction, is absent from that treatment, since there the parent equation is taken to be Stratonovich. Note that the discrepancy is not the price of an approximation, as was already realised there, since a mere change of variables is enough to break the It\^{o}--Stratonovich equivalence. We make this observation geometric (Remark~\ref{rem:flat-manifold-test}), with a manifold of identically vanishing second fundamental form where the terms that build $\mathfrak{s}_a$ still do not vanish.

On the other hand, the parent equation was treated differently, and here the difference between fixing a convention and establishing an equivalence becomes important. For the full field equation, in the general case as well as in the spherically symmetric case, both readings were reported to coincide. Following van Saarloos, the field was discretised on a lattice of spacing $\Delta$, and the drift so obtained was found to cancel term by term. Note that on that lattice the cancellation is exact. This cancellation is, however, a consequence of the stencil, the noise amplitude at each site being built from the neighbouring sites. Therefore, what was shown is a property of the discretisation, while the continuum limit, where the question is about an equation and no longer about a scheme, was not considered. Note that taking that limit would require a space where it can be taken. In this work we follow the opposite path, and mollify the noise from the very beginning. At finite correlation length $\delta$ the correction becomes a convergent trace in a stated Hilbert space, which can be computed.

Our conclusion differs, since no anomalous force was found there in the parent field equation. The correction, evaluated in Lemma~\ref{lemma:ito-strat-correction} as a trace against the covariance kernel, is generically non-zero at every finite $\delta$. As a consequence, the It\^{o} reading and the Stratonovich reading of one and the same written field equation lead to two different equations. Hence we are left with a choice of frame. The parent equation is posed in the It\^{o} frame used by Dean, with its Stratonovich form obtained by subtracting the correction. Going to the Stratonovich frame and back thus leaves the equation unchanged, as well as the bare thermodynamic force, and for that reason no residue appears at field level, although the statement says nothing about the size of the correction.

The contrast with the reduced system is the subject of this work. In the reduced system the projection sits between the two conversions and does not commute with them, so going out and back is no longer the identity, and what survives is $\mathfrak{s}_a$.

The derivations of mesoscopic nucleation theory use the idiom of theoretical physics, where manipulations are formal and the function space and the stochastic calculus are left implicit. This idiom is well suited to the purpose of those papers (i.e.\ physical insight into nucleation), and their results are the ones we set out to ground here. However, it does not provide any way to check a reduction against anything other than itself. First, no function space is fixed, so the weight determining the metric on the manifold of profiles is chosen for the structure it confers, instead of being derived from the parent dynamics. Second, and as a consequence, the residual fluctuation (i.e.\ the part of the field not described by the coordinates) has no inner product with respect to which its orthogonality could be stated. Third, the anomalous force is computed within the trajectory description only, and Kawasaki's density-level description of the same dynamics has never been reduced, so that there is no independent construction to compare the reduced equations with. Supplying the geometric and measure-theoretic foundation, and with it that comparison, is the object of the present work.

In finite dimensions, the corresponding problem is currently being studied, and the results obtained there are instructive. Ayala et~al.~\cite{ayala2025reversibility} showed, by Mosco convergence of Dirichlet forms, that reversibility and the Klimontovich noise interpretation are preserved under coarse-graining for reversible slow-fast systems with block-diagonal diffusion. Working at the same time as us, Nateghi et~al.~\cite{nateghi2025consistent} built a Zwanzig-projection scheme for underdamped Langevin dynamics, designed to preserve both thermodynamic consistency and kinetic observables such as transition timescales. These two results are positive, which makes the comparison with our case instructive. In our case the reduction starts from an infinite-dimensional state space, with a Dean--Kawasaki diffusion that is neither block-diagonal nor state-independent. We will show that in this case the structure is not preserved automatically, and we will write the obstruction down explicitly.

\subsection{The dimensional reduction problem}

As shown in the preceding sections, the parent equation comes in two descriptions that are equivalent before any reduction (the trajectory and its law), and the reduction may be applied to either of them. In both cases the target is the same, which is given by a finite-dimensional submanifold $\mathcal{M}$ of the ambient function space, made of the density profiles reachable by the collective coordinates $\mathbf{X}$ of the preceding section, and the reduced model has to keep the gradient-flow structure of the parent equation. Reducing the dynamics onto $\mathcal{M}$ means eliminating the fast degrees of freedom of the field in favour of a slow effective description, a classical problem for which Pavliotis and Stuart~\cite{pavliotis2008multiscale} collected the averaging and homogenisation theorems and E and Engquist~\cite{e2003heterogeneous} developed the numerical framework. What changes is the object being reduced, and one cannot tell a priori whether the two reductions lead to the same reduced model.

The first route (which we call the kinetic projection) acts on the stochastic vector field itself. Under a spectral-gap hypothesis, the field is projected onto the tangent bundle of $\mathcal{M}$, and the orthogonal fast modes are discarded. It is worth mentioning that this is not a new idea, since Ciccotti, Leli\`{e}vre and Vanden-Eijnden~\cite{ciccotti2008projection} projected diffusions on submanifolds when studying constrained Langevin systems. In a related way, Legoll and Leli\`{e}vre~\cite{legoll2010effective} obtained effective dynamics for a reaction coordinate by conditional expectation, and Pavliotis and Stuart~\cite{pavliotis2008multiscale} presented stochastic model reduction methods. Mesoscopic nucleation theory follows this route, working in coordinates. As a result, one obtains a system of Stratonovich Langevin equations carrying the Riemannian metric induced by the embedding. Written in It\^{o} form, the same dynamics acquires a drift correction (the geometric anomaly $\mathfrak{s}_a$). It is an effect of the state dependence of the noise amplitude, not of the curvature of $\mathcal{M}$ (Remark~\ref{rem:flat-manifold-test}).

The second route (developed in this work) never passes through the stochastic trajectories, and reduces the problem at the level of probability densities. In other words, we reduce Kawasaki's description, not Dean's. We take the forward Kolmogorov equation associated with Eq.~\eqref{eq:dean-kawasaki}, introduce Fermi tubular coordinates around $\mathcal{M}$, factorise the probability density into a marginal and a conditional part, and integrate out the transverse degrees of freedom in the adiabatic limit. In this way we obtain a covariant Fokker--Planck equation on $\mathcal{M}$. The method builds on the finite-dimensional averaging principles of Freidlin and Wentzell~\cite{freidlin2012random}, and on their infinite-dimensional extensions by Cerrai and Freidlin~\cite{cerrai2009averaging}. To the best of our knowledge, however, it has never been applied to the functional forward Kolmogorov equation of a conservative SPDE with state-dependent, spatially correlated noise. In the probabilistic route we obtain a drift correction made of three parts. On the one hand, the statistical anomaly $\Upsilon$ originates from the curvature of the deterministic drift averaged against the transverse fluctuations, and is produced by the excess free energy alone. The frame correction $\mathfrak{f}$ and the connection correction $\mathfrak{c}$, on the other hand, are geometric in origin.

In what follows, we call both corrections anomalies, in the field-theoretic sense of an identity that holds in the classical theory and fails once fluctuations are retained. If the reduction commuted with the calculus, the reduced generator would be the intrinsic operator dictated by geometry and free energy alone, i.e.\ the Laplace--Beltrami operator (Theorem~\ref{thm:geo-operator}) twisted by the gradient of the free energy ($\nabla\mathcal{V}$), and the fibre-averaged drift would be the drift evaluated on the manifold. Neither statement holds. The kinetic route leaves the geometric anomaly $\mathfrak{s}_a$ as a residue, and the fibre average of the drift leaves the statistical anomaly plus the frame correction ($\Upsilon + \mathfrak{f}$). Both residues vanish in the deterministic theory, are generated by the quadratic variation of the noise and are of order $\varepsilon^2$. Moreover, the probabilistic generator carries the connection correction ($\mathfrak{c}$), already induced on the manifold by the full dynamics. The name is the one used in mesoscopic nucleation theory~\cite{lutsko2012mean}, where the correction of the kinetic route is called the anomalous force (since its origin is not the thermodynamic driving force). In this work we keep this name. The more common term (\emph{spurious drift}) would suggest that only the value of a coefficient is wrong, while what fails is a covariance statement about the reduction.

Unfortunately, only one of the two routes has been explored, since the kinetic projection is the only one that has been carried out in the case of conservative stochastic partial differential equations, and always at the level of the trajectories. This would be of no consequence for finite-dimensional systems with additive or state-independent noise, where trajectory and density carry the same information. There the Fokker--Planck equation follows from the It\^{o} equation by a change of test function, and no anomalous drift appears. This is the case, for instance, of the classical overdamped reductions of Kramers and Smoluchowski~\cite{kramers1940brownian,smoluchowski1916drei}.

However, the Dean and Kawasaki equation is a different story. Here we have to deal with a conservative, state-dependent and spatially correlated noise, acting on an infinite-dimensional state space. The density-level reduction of such an equation was missing, even in the spherically symmetric geometry considered here, so there was nothing to compare the kinetic projection with. The difficulty is due to the way the It\^{o}-to-Stratonovich conversion interacts with the spatial regularisation of the noise. Moreover, the comparison has to be made at the level of the infinitesimal generator, where neither the order-parameter reduction~\cite{lutsko2012mean} nor the well-posedness theory~\cite{cornalba2019regularised} operates.

Already in finite dimensions, the question turns out to be relevant. It was found that the reversibility of a diffusion depends on the noise convention (It\^{o}, Stratonovich or Klimontovich)~\cite{ayala2025reversibility}. Hence, already in that case, the noise convention and the reduction need not commute, and we address the infinite-dimensional version of this for conservative equations.

Note that different generators would make the effective theory ambiguous, since its predictions would then depend on the route taken, and the theory would be well defined only up to the size of that dependence. Grmela and \"{O}ttinger~\cite{grmela1997dynamics} based their GENERIC formalism on the reversible and irreversible structure of macroscopic dynamics. Its variational counterpart, which connects gradient-flow structure and large-deviation rate functionals, was established by Mielke, Peletier and Renger~\cite{mielke2014relation}.

\begin{figure}[ht]
    \centering
    \begin{tikzpicture}[node distance=2.8cm and 4.5cm, auto, >=Latex]
        % Nodes with mathematical descriptions
        \node (SPDE) [align=center] {\textsf{Full field dynamics} \\ (SPDE) \\ $dm_t = [\mathcal{A}_m - \varepsilon^2\mathfrak{S}^\delta_m]\, dt + \varepsilon\,\mathcal{B}_m \circ d\xi^\delta_t$};
        \node (FFPE) [right=of SPDE, align=center] {\textsf{Functional probability} \\ (Forward Kolmogorov) \\ $\partial_t \mu_t = \mathcal{L}_{\varepsilon}^* \mu_t$};
        \node (SDE) [below=of SPDE, align=center] {\textsf{Reduced trajectory} \\ (Geometric Langevin) \\ generator $\mathcal{L}_A$, anomaly $\mathfrak{s}_a$};
        \node (gFPE) [below=of FFPE, align=center] {\textsf{Reduced probability} \\ (Generalised FPE) \\ generator $\mathcal{L}_B$, correction $\Upsilon + \mathfrak{f} + \mathfrak{c}$};
        % Arrows
        \draw[<->, thick] (SPDE) -- node [above, align=center, font=\small] {Stochastic duality \\ (It\^{o} formula)} (FFPE);
        \draw[->, thick] (SPDE) -- node [left, align=right, font=\small] {\emph{Kinetic}\\ \emph{projection} \\ $\Pi_{\mathbf{X}}^{\parallel} (dm_t)$ \\ (Stiffness $\lambda_0 \gg 1$)} (SDE);
        \draw[->, thick] (FFPE) -- node [right, align=left, font=\small] {\emph{Probabilistic}\\ \emph{projection} \\ $\int (\cdot) \Gamma_t J_{\Psi}^{-1} d\lambda_{\perp}$ \\ (Limit $\varepsilon \to 0$)} (gFPE);
        \draw[<->, thick, dashed] (SDE) -- node [below, align=center, font=\small] {Generator equivalence \\ $\mathcal{L}_A \overset{?}{=} \mathcal{L}_B$} (gFPE);
    \end{tikzpicture}
    \caption{The two reduction pathways compared in this work. Applied to the SPDE, the It\^{o} formula gives the forward Kolmogorov equation for the law of the field (top edge, Section~\ref{subsec:functional-fpe}). On the left, we carry the vector field of the SPDE onto the tangent bundle of the concentration manifold and discard the transverse modes under the stiffness hypothesis (\emph{kinetic projection}). On the right, we disintegrate the law $\mu_t$ in Fermi tubular coordinates and integrate over the normal fibres in the weak-noise limit $\varepsilon\to0$ (\emph{probabilistic projection}). Each route leaves a different correction behind, and in Theorem~\ref{thm:commutativity} we show that the two routes generate the same dynamics (dashed edge) if and only if the geometric anomaly and the probabilistic correction coincide as vector fields. This identity fails generically, and in Corollary~\ref{cor:asymptotic-agreement} we bound the resulting discrepancy over a fixed time horizon by $\mathcal{O}(t/n\delta^3)$.}
    \label{fig:reduction_diagram}
\end{figure}
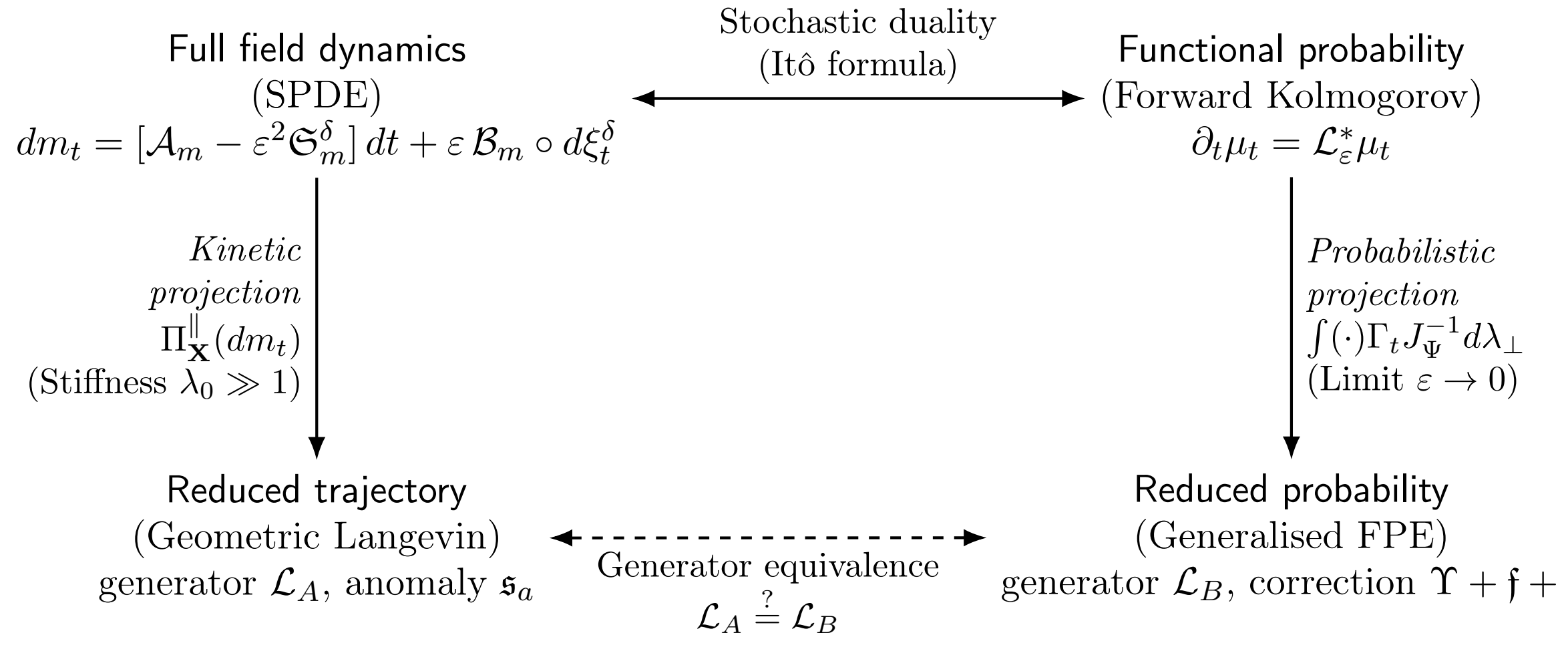
\subsection{Main result and contributions}

In this work we resolve such an ambiguity for conservative stochastic partial differential equations of the Dean and Kawasaki type. Here we consider them on a spherically symmetric domain (the ball $V = B_L(0)$ of radius $L$ in $\mathbb{R}^3$), driven by mollified noise with a smooth, finite-range covariance kernel $C_\delta(r,s)$ of correlation length $\delta > 0$. In Theorem~\ref{thm:commutativity} we compare the infinitesimal generators $\mathcal{L}_A$ and $\mathcal{L}_B$ of the two routes on smooth functions on the manifold, $C^\infty(\mathcal{M}_m)$. In particular, we prove that they agree if and only if the geometric anomaly of the first route and the probabilistic correction of the second coincide as vector fields. In order to prove it, we make use of the following three structural conditions:
\begin{enumerate}[(i)]
    \item \emph{Spectral stiffness} (Hypothesis~\ref{hyp:stiffness}): we take the transverse spectrum to have a uniform gap $\lambda_0 \gg 1$, so that the fast modes equilibrate on the timescale $\tau_\perp = \lambda_0^{-1}$. In addition, the adiabatic approximation is assumed to be valid.
    \item \emph{Adiabatic concentration} (Hypothesis~\ref{hyp:concentration}): the transverse fluctuations become Gaussian, with the conditional measure on each normal fibre converging to a Gaussian of covariance $\varepsilon^2\mathbf{C}_\perp$, with a mean of order $\varepsilon^2$ fixed by the mean curvature of the manifold.
    \item \emph{Kinematic covariance} (Hypothesis~\ref{hyp:kinematic-covariance}): we take the stationary transverse covariance to be approximately proportional to the noise diffusion tensor in the Karhunen--Lo\`{e}ve eigenbasis, which gives the same relaxation rate to all the transverse modes. In what follows, we call this approximation the scalar closure.
\end{enumerate}
The metric weight $\varpi$, on the other hand, needs no hypothesis of its own, since we obtain it twice, once from the noise covariance (Lemma~\ref{lemma:metric-selection}) and once, independently, from the Wasserstein geometry of the parent gradient flow (Theorem~\ref{thm:wasserstein-isometry}), and both give the same weight.

As a first step, we reduce the question of commutativity to a single vector identity between the geometric anomaly and the probabilistic correction. Note that we do this algebraically, for every smooth symmetric covariance kernel $C_\delta(r,s)$ and without ever taking the white-noise limit. When the reduction fixes its metric by means of the local covariance, the equivalence turns out to be exact. If the generator carries instead the exact finite-$\delta$ diffusion tensor, it differs from the former by $\mathcal{O}(\varepsilon^2\delta^2)$ under Assumption~\ref{ass:uniform-metric-expansion} (i.e.\ at higher order than the discrepancy we want to measure). We then test the identity and find that, as an exact statement, it fails generically, due to the fact that its two sides depend on the covariance kernel at incompatible rates. For a fluid with a local excess free energy, and under the scalar closure, the statistical anomaly depends on the kernel through its diagonal mixed curvature $\kappa_\delta$, while the geometric anomaly only depends on milder invariants of it ($V_\delta$ and $V_\delta'$). In the case of the ideal gas the statistical anomaly vanishes, and the geometric anomaly is not compensated by the remaining corrections, which are of order unity. In the limit of infinite transverse stiffness the probabilistic generator becomes, at the order retained, the generator induced on the manifold by the full dynamics with the local covariance. On the other hand, the kinetic generator differs from it by a drift of the order of $\delta^{-1}$ (Corollary~\ref{cor:stiff-limit}). Nevertheless, over a fixed time horizon the two reduced semigroups differ by at most a term of the order of $t/n\delta^3$, so the failure is under quantitative control. Hence the two reductions cannot be told apart in terms of observables, at the order retained, as soon as the noise correlation volume contains many particles.

The main contributions of this work are the following:
\begin{enumerate}[(i)]
    \item \emph{Construction of the probabilistic projection.} We develop a density-level reduction of the Dean and Kawasaki equation on the spherically symmetric domain, which never passes through the stochastic trajectories (Section~\ref{sec:kolmogorov-reduction}). It is based on Fermi tubular coordinates and on covariant measure disintegration on the normal fibres. To the best of our knowledge, this is the first application of such a probabilistic reduction to an infinite-dimensional conservative SPDE with state-dependent, spatially correlated noise.

    \item \emph{Identification and comparison of the drift corrections.} Each of the two routes generates a noise-induced drift correction that is absent from a na\"{\i}ve gradient-flow truncation. In the kinetic case this correction is the geometric anomaly $\mathfrak{s}_a$, which comes from the It\^{o}-to-Stratonovich conversion. In the probabilistic case it is the statistical anomaly $\Upsilon$ (from the curvature of the drift averaged against the transverse fluctuations), as well as the frame correction $\mathfrak{f}$ and the connection correction $\mathfrak{c}$. Moreover, we compute all of them in closed form, show that the difference between the two generators is carried by them alone, and find that the criterion does not hold identically. For a fluid with a local excess free energy the discrepancy is governed by the diagonal mixed curvature of the covariance kernel, and for the ideal gas by the projected field correction (Theorem~\ref{thm:commutativity}, Proposition~\ref{prop:validity-condition} and, for the ideal gas, Corollary~\ref{cor:finite-delta-comparison}).

    \item \emph{Invariant measure of the reduced dynamics.} It can be shown that the Fokker--Planck equation of the kinetic reduction is covariant. When its anomaly is conservative (i.e.\ the gradient of a potential), $\mathfrak{s}_a = -\nabla\Phi_\xi$, and the resulting exponent is integrable, it admits an explicit stationary density, $\bar{\varrho}_{\mathrm{eq}} \propto \exp(-\mathcal{V}/\varepsilon^2 - \Phi_\xi/D_0)$, with respect to the Riemannian volume measure $\mathfrak{m}_g$ (Proposition~\ref{prop:invariant-measures}). With respect to the Boltzmann--Gibbs form, the noise multiplies the density by a factor that carries no power of $\varepsilon$, which leaves the large-deviation exponent $\mathcal{V}/\varepsilon^2$ unchanged but still cannot be regarded as a small correction to the density. Note that no Jacobian factor of the Fermi coordinates enters, since the Jacobian cancels from every fibre average (Remark~\ref{rem:volume-no-drift}).

    \item \emph{A geometric formulation of mesoscopic nucleation theory.} We formalise from first principles the reductions of mesoscopic nucleation theory~\cite{lutsko2012mean,lutsko2013classical,duran2015mesoscopic}. These reductions were obtained by testing the field equation against a weight chosen so that a fluctuation-dissipation relation should survive the reduction. Those reductions presuppose three structures, which we derive instead of positing them. In particular, we obtain the metric weight from the Wasserstein geometry of the parent gradient flow, the orthogonality of the projection from that metric, and the noise transformation from the It\^{o}-to-Stratonovich calculus. The framework rests on infinite-dimensional stochastic analysis, optimal transport and covariant measure disintegration.
    \item \emph{A consistency criterion for mesoscopic reductions.} The commutativity of the two routes is what makes the effective theory independent of the description from which it was obtained. Here we provide the criterion that decides it in terms of an explicit functional, which vanishes if and only if the two routes commute, and we evaluate it in closed form for the spherically symmetric setting. Passing this test is necessary, although not sufficient, for a reduction to be correct. A reduction that fails it is only well defined up to the size of the discrepancy, which is bounded in Corollary~\ref{cor:asymptotic-agreement}.
    \item \emph{The stiff limit.} We show that in the limit of infinite transverse stiffness the probabilistic generator becomes, at the order retained, the generator induced on the manifold by the full dynamics with the local covariance, while the kinetic generator keeps a drift of order $\delta^{-1}$ (Corollary~\ref{cor:stiff-limit}). Mesoscopic nucleation theory obtains its reduced equations by projecting the trajectories, and this projection is not the stiff limit of the dynamics it reduces. On the other hand, the averaging of the law (which reduces the generator of the observable itself) does recover that limit at the order retained, with the local covariance.
\end{enumerate}

The result also complements the analytical framework of Zimmer and collaborators~\cite{cornalba2019regularised,cornalba2020dean,muller2024wellposedness}, who established the well-posedness of related equations. Here we ask what happens to the geometric structure of the equation once it has been reduced. Unfortunately, at the exact level the answer is negative, because the structure is not preserved automatically. The obstruction is the diagonal mixed curvature of the covariance kernel for a fluid with a local excess free energy under the scalar closure, and the projected field correction for the ideal gas. At the approximate level, on the other hand, the answer is quantitative, as the residual discrepancy between the two reduced semigroups is at most of the order of $t/n\delta^3$ over a fixed time horizon under the same closure. The result may also be seen as the infinite-dimensional counterpart of the finite-dimensional work of Ayala et~al.~\cite{ayala2025reversibility} on the preservation of noise conventions under coarse-graining. That work established that the Klimontovich convention survives averaging for reversible diffusions in finite dimensions. For a conservative SPDE on a function space over a spherically symmetric domain, we identify the exact condition for the full infinitesimal generator to be invariant under two distinct reduction routes. Moreover, we quantify the departure from that condition at finite correlation length.

We separate what is structural in the comparison from what is tied to the symmetry of the domain. Both routes produce a generator with the same thermodynamic force and the same diffusion term, and in this sense the comparison itself is structural. As a consequence, the whole of their difference is carried by the geometric anomaly and the probabilistic correction, so agreement is nothing else than the coincidence of two vector fields, which makes no reference to the geometry of $V$. To compute the two anomalies, on the other hand, we rely on machinery that is specific to the spherically symmetric setting, i.e.\ the cumulative mass variable and the weight $\varpi = (4\pi r^2\rho)^{-1}$ that it induces, the Dirichlet condition at the two ends of the radial interval, the reduction of the covariance to a one-dimensional chord integral and the spherical Bessel spectrum of the transverse operator. Within that setting we obtain the closed forms of the two anomalies, $\mathfrak{s}_a$ and $\Upsilon$, and with them the verdict that the criterion fails. The same setting also carries the spherical-mean closure of Assumption~\ref{ass:spherical-closure}, which mesoscopic nucleation theory adopts without estimating the error it introduces, and the mass equation (i.e.\ the common starting point of both routes) holds under it. Note that we did not adopt this setting for convenience, since mesoscopic nucleation theory is itself formulated in it, and this makes it the natural first case. The extension to domains of lower symmetry is taken up in Section~\ref{sec:conclusions}.

\subsection{Organisation of the paper}

In Section~\ref{sec:preliminaries} we introduce the prerequisites that we need from functional analysis, operator theory, stochastic analysis and differential geometry. Section~\ref{sec:reduction} is devoted to the dimensional reduction of the stochastic equation in terms of cumulative mass coordinates. There we state the concentration manifold hypothesis, project the dynamics onto the tangent bundle, and derive the geometric Langevin equations as well as the Fokker--Planck equation that they induce. In Section~\ref{sec:kolmogorov-reduction} we follow the other route, and reduce the functional forward Kolmogorov equation. To do so, we construct Fermi coordinates, factorise the probability measure and take the adiabatic limit on the normal fibres. In Section~\ref{sec:commutativity} we bring the two routes together, and compare their generators. We then establish the criterion that decides whether they agree, and quantify the departure from it at finite correlation length. Finally, in Section~\ref{sec:conclusions} we summarise the main results of this work and point out some directions for future work.

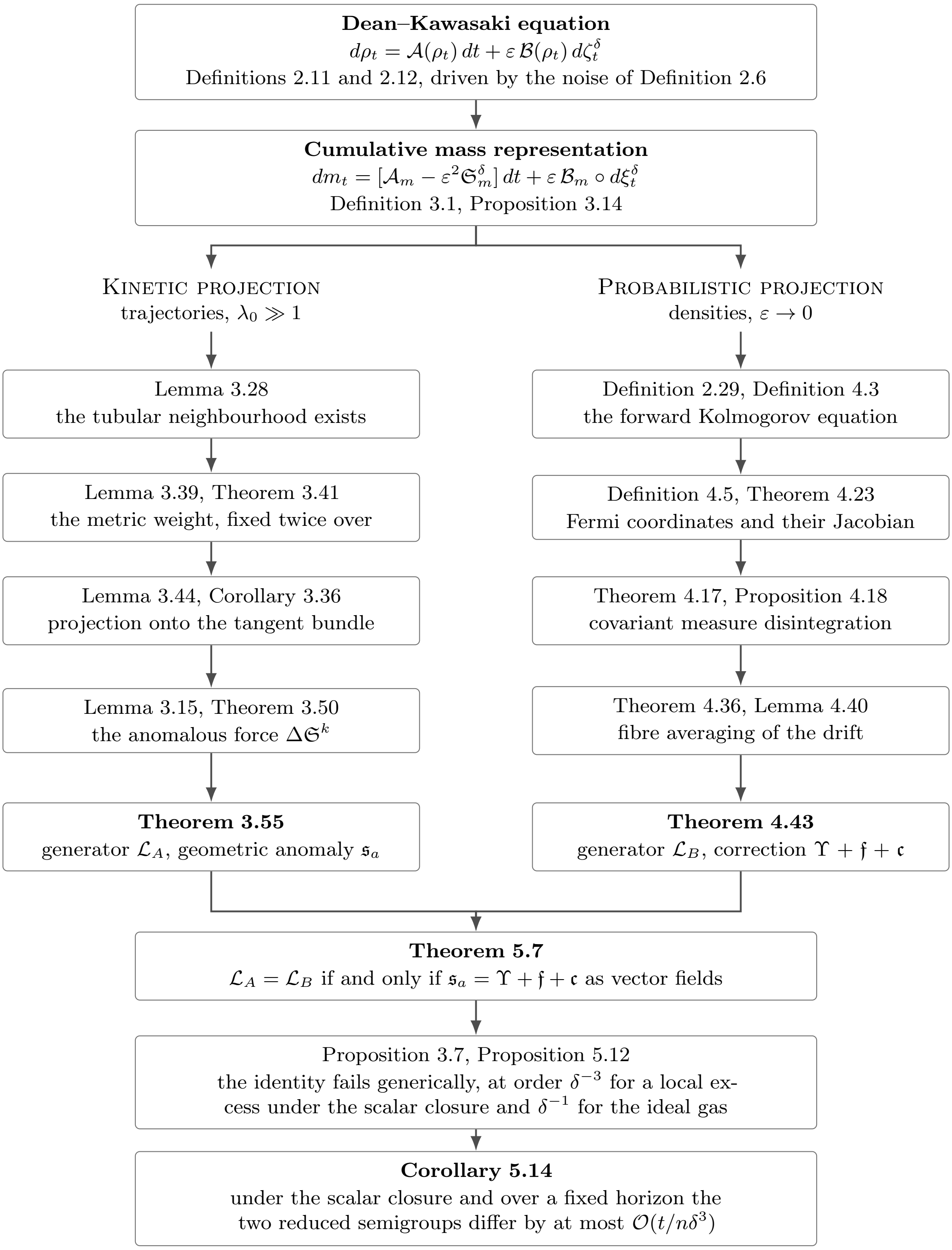
\begin{figure}[htp]
    \centering
    \begin{tikzpicture}[
        >=Latex,
        every node/.style={font=\footnotesize},
        bx/.style={draw=black!55, rounded corners=2pt, align=center,
                   inner sep=5pt, text width=#1},
        hd/.style={align=center, font=\small\scshape},
        ar/.style={->, thick, draw=black!70}]

    \def\lx{-4.1}
    \def\rx{4.1}

    \node[bx=10.2cm] (spde) at (0,0)
      {\textbf{Dean--Kawasaki equation}\\[1pt]
       \mbox{$d\rho_t = \mathcal{A}(\rho_t)\,dt + \varepsilon\,\mathcal{B}(\rho_t)\,d\zeta^\delta_t$}\\[1pt]
       Definitions~\ref{def:drift-op} and~\ref{def:diff-op}, driven by the noise of Definition~\ref{def:mollifier}};

    \node[bx=10.2cm] (mass) at (0,-1.95)
      {\textbf{Cumulative mass representation}\\[1pt]
       \mbox{$dm_t = [\mathcal{A}_m - \varepsilon^2\mathfrak{S}^\delta_m]\,dt + \varepsilon\,\mathcal{B}_m\circ d\xi^\delta_t$}\\[1pt]
       Definition~\ref{def:mass-operator}, Proposition~\ref{prop:mass-dynamics}};

    \node[hd] (lh) at (\lx,-3.85) {Kinetic projection\\[-1pt] \normalfont\footnotesize trajectories, $\lambda_0\gg1$};
    \node[hd] (rh) at (\rx,-3.85) {Probabilistic projection\\[-1pt] \normalfont\footnotesize densities, $\varepsilon\to0$};

    \node[bx=6.1cm] (l1) at (\lx,-5.45)
      {Lemma~\ref{lemma:tubular-existence}\\[1pt] the tubular neighbourhood exists};
    \node[bx=6.1cm] (l2) at (\lx,-7.05)
      {Lemma~\ref{lemma:metric-selection}, Theorem~\ref{thm:wasserstein-isometry}\\[1pt]
       the metric weight, fixed twice over};
    \node[bx=6.1cm] (l3) at (\lx,-8.65)
      {Lemma~\ref{lemma:projected-dynamics}, Corollary~\ref{cor:adiabatic-elimination}\\[1pt]
       projection onto the tangent bundle};
    \node[bx=6.1cm] (l4) at (\lx,-10.35)
      {Lemma~\ref{lemma:ito-strat-correction}, Theorem~\ref{thm:equivalent-fluctuations}\\[1pt]
       the anomalous force $\Delta\mathfrak{S}^k$};
    \node[bx=6.1cm] (l5) at (\lx,-12.15)
      {\textbf{Theorem~\ref{thm:geo-operator}}\\[1pt]
       generator $\mathcal{L}_A$, geometric anomaly $\mathfrak{s}_a$};

    \node[bx=6.1cm] (r1) at (\rx,-5.45)
      {Definition~\ref{def:ffpe}, Definition~\ref{def:adjoint-evolution}\\[1pt]
       the forward Kolmogorov equation};
    \node[bx=6.1cm] (r2) at (\rx,-7.05)
      {Definition~\ref{def:fermi-map}, Theorem~\ref{thm:jacobian-expansion}\\[1pt]
       Fermi coordinates and their Jacobian};
    \node[bx=6.1cm] (r3) at (\rx,-8.65)
      {Theorem~\ref{thm:disintegration}, Proposition~\ref{prop:measure-decomposition}\\[1pt]
       covariant measure disintegration};
    \node[bx=6.1cm] (r4) at (\rx,-10.35)
      {Theorem~\ref{thm:marginal-equation}, Lemma~\ref{lemma:mechanical-drift-explicit}\\[1pt]
       fibre averaging of the drift};
    \node[bx=6.1cm] (r5) at (\rx,-12.15)
      {\textbf{Theorem~\ref{thm:reduced-fpe}}\\[1pt]
       generator $\mathcal{L}_B$, correction $\Upsilon + \mathfrak{f} + \mathfrak{c}$};

    \node[bx=10.2cm] (cmp) at (0,-14.15)
      {\textbf{Theorem~\ref{thm:commutativity}}\\[1pt]
       \mbox{$\mathcal{L}_A = \mathcal{L}_B$} if and only if \mbox{$\mathfrak{s}_a = \Upsilon + \mathfrak{f} + \mathfrak{c}$} as vector fields};
    \node[bx=10.2cm] (obs) at (0,-15.95)
      {Proposition~\ref{prop:kernel-asymptotics}, Proposition~\ref{prop:validity-condition}\\[1pt]
       the identity fails generically, at order \mbox{$\delta^{-3}$} for a local excess under the scalar closure and \mbox{$\delta^{-1}$} for the ideal gas};
    \node[bx=10.2cm] (cns) at (0,-17.75)
      {\textbf{Corollary~\ref{cor:asymptotic-agreement}}\\[1pt]
       under the scalar closure and over a fixed horizon the two reduced semigroups differ by at most $\mathcal{O}(t/n\delta^3)$};

    \draw[ar] (spde) -- (mass);
    \draw[ar] (mass.south) -- ++(0,-0.30) -| (lh.north);
    \draw[ar] (mass.south) -- ++(0,-0.30) -| (rh.north);
    \draw[ar] (lh) -- (l1); \draw[ar] (rh) -- (r1);
    \draw[ar] (l1) -- (l2); \draw[ar] (l2) -- (l3);
    \draw[ar] (l3) -- (l4); \draw[ar] (l4) -- (l5);
    \draw[ar] (r1) -- (r2); \draw[ar] (r2) -- (r3);
    \draw[ar] (r3) -- (r4); \draw[ar] (r4) -- (r5);
    \draw[ar] (l5.south) |- ++(0,-0.62) -| (cmp.north);
    \draw[ar] (r5.south) |- ++(0,-0.62) -| (cmp.north);
    \draw[ar] (cmp) -- (obs);
    \draw[ar] (obs) -- (cns);
    \end{tikzpicture}
    \caption{The route from the parent equation to the comparison, together with the result
    that establishes each step. The left branch is the route of
    Figure~\ref{fig:reduction_diagram} that reduces the trajectories, and the right branch is
    the one that reduces their law. On the left, Theorem~\ref{thm:geo-operator} computes the
    correction of the first route, and on the right,
    Lemma~\ref{lemma:mechanical-drift-explicit} and Definition~\ref{def:connection-correction}
    compute that of the second, and both branches end at Theorem~\ref{thm:commutativity},
    where the question of agreement is reduced to the coincidence of the geometric anomaly
    and the probabilistic correction as vector fields. The two boxes at the bottom show the
    failure of the coincidence, through the diagonal mixed curvature of the covariance kernel
    in Proposition~\ref{prop:kernel-asymptotics} and Proposition~\ref{prop:validity-condition},
    as well as the bound on the distance between the two reduced semigroups in
    Corollary~\ref{cor:asymptotic-agreement}.}
    \label{fig:argument_spine}
\end{figure}

Figure~\ref{fig:argument_spine} shows the same course in the form of a chain of results, where each step is accompanied by the lemma, proposition or theorem that establishes it.

\subsection*{Notation and conventions}
\begingroup
\small
\noindent Two ambient spaces appear throughout this work. These are the unweighted Hilbert space $\mathcal{H}$ of the stochastic description and the weighted one ($\mathcal{H}_\varpi$) of the projected description. In cumulative-mass variables, the reduced state space is the manifold $\mathcal{M}_m$, parametrised by the embedding $M(\mathbf{X})$ (with tangent basis $\mathbf{e}_k^m = \partial_{X^k}M$).\par
\smallskip\noindent On $\mathcal{M}_m$, $g_{ij}$ is the induced Riemannian metric and $g^{ij}$ its inverse, while the weight entering the projected inner product is $\varpi = (\partial_r M)^{-1}$. Here $\Pi_{\mathbf{X}}^{\parallel}$ and $\Pi_{\mathbf{X}}^{\perp}$ denote the associated tangential and normal projectors.\par
\smallskip\noindent As for the free energy $\mathcal{F}[\rho]$, it splits into ideal-gas and excess contributions. The corresponding effective potential on the reduced manifold is $\mathcal{V}(\mathbf{X}) = \beta\,\mathcal{F}_m[M(\mathbf{X})]$ (in units of $k_BT$), and hence the thermodynamic force reads $\mathbf{F}_0 = -D_0 \nabla \mathcal{V}$, where $D_0$ is the bare diffusion coefficient (and the noise intensity is $\varepsilon$).\par
\smallskip\noindent For the noise, $C_\delta(r,s)$ denotes the mollified covariance kernel (with correlation length $\delta$), while in the mass representation $\mathcal{B}_m$ and $\mathcal{Q}_m$ are the noise operator and its covariance, and the transverse covariance on the normal fibre is $\mathbf{C}_\perp$.\par
\smallskip\noindent We denote the geometric and statistical anomaly vectors by $\mathfrak{s}_a$ and $\Upsilon$ (note that neither of them is set in bold), and the frame and connection corrections by $\mathfrak{f}$ and $\mathfrak{c}$. The kinetic and probabilistic reductions produce the generators $\mathcal{L}_A$ and $\mathcal{L}_B$, respectively. Finally, $J_\Psi$ is the Jacobian of the Fermi map, $A_{\mathbf{n}}$ the shape operator, and $\overline{\mathbf{F}}$ the effective fibre-averaged force.\par
\endgroup

\section{Preliminaries}
\label{sec:preliminaries}
In this section we collect only the material that the reductions of Sections~\ref{sec:reduction} and~\ref{sec:kolmogorov-reduction} presuppose. We first fix the ambient function spaces and the operator classes that we need to state that the noise is trace class. This condition is secured by the mollification, and it is what makes the stochastic calculus available at all. Then we record the regularised forcing itself and the Stratonovich convention, since with it the projection obeys the ordinary chain rule. After that comes the manifold-theoretic apparatus (the immersion, its tangent and normal bundles and the second fundamental form) in the generality that the later sections need. The section closes with the transport structure that fixes the metric and with the weak formulation of the functional forward Kolmogorov equation, which is the object that the probabilistic pathway reduces.
\subsection{Functional analytic setting}
\label{subsec:functional-analysis}
We fix the inverse temperature $\beta \coloneqq (k_B T)^{-1} > 0$, the bare diffusivity $D_0 > 0$, and the fluctuation amplitude $\varepsilon > 0$. The physical domain is a bounded connected set $V \subset \mathbb{R}^3$ with a boundary of class $C^2$. The ambient Hilbert space is
\begin{equation}
    \mathcal{H} \coloneqq L^2(V), \qquad \langle u, v \rangle_{\mathcal{H}} \coloneqq \int_V u(\mathbf{r}) v(\mathbf{r}) \, d\mathbf{r},
    \label{eq:l2-inner-product}
\end{equation}
and the corresponding Sobolev space is $H^1(V) = \{u \in \mathcal{H} : \nabla u \in L^2(V;\mathbb{R}^3)\}$. We write $(H^1(V))^*$ for its topological dual.

We follow the usual conventions for classes of differentiable functions. A set or a map is of class $C^k$ when it is $k$ times continuously differentiable, the derivative being the Fr\'{e}chet one wherever the target is infinite-dimensional. On a domain or a manifold, $C^k$, $C^\infty$ and $C_c^\infty$ denote the correspondingly differentiable functions, the smooth functions and the smooth functions of compact support, while $C_b$ denotes the bounded continuous ones. On a Hilbert space, $C_b^2$ denotes the bounded real functionals that are twice continuously Fr\'{e}chet differentiable with first and second derivatives bounded, and it is the boundedness of the second derivative that makes the trace term of the infinitesimal generator well defined.

The Neumann Laplacian $\Delta_N$ is the self-adjoint operator associated with the quadratic form $\mathfrak{e}(u,v) = \int_V \nabla u \cdot \nabla v \, d\mathbf{x}$ on the domain
\begin{equation}
    D(\Delta_N) = \left\{u \in H^2(V) : \frac{\partial u}{\partial \mathbf{n}} = 0 \text{ on } \partial V\right\},
\end{equation}
so that $-\Delta_N$ is non-negative with discrete spectrum, the identification of that domain with the form being Theorem~2.2.2.5 of Grisvard~\cite{grisvard2011elliptic}. Mass-conserving fluctuations live in the closed subspace
\begin{equation}
    \mathcal{H}_0 \coloneqq \left\{u \in \mathcal{H} : \int_V u(\mathbf{r}) \, d\mathbf{r} = 0\right\},
\end{equation}
and on $\mathcal{H}_0$ the inverse Neumann Laplacian $-\Delta_N^{-1}$ is well defined, compact, and self-adjoint.

\begin{definition}[Configuration space of densities]
\label{def:density-space}
The state of the mesoscopic system is described by a density field $\rho$. We define the configuration space $\mathcal{D}_N$ as the convex set of non-negative functions in $\mathcal{H}$ with fixed total mass $N > 0$:
\begin{equation}
    \mathcal{D}_N \coloneqq \left\{ \rho \in \mathcal{H} : \rho(\mathbf{r}) \ge 0 \text{ a.e.}, \quad \int_V \rho(\mathbf{r}) \, d\mathbf{r} = N \right\}.
\end{equation}
We denote by $\mathcal{D}_N^+$ the subset of strictly positive densities. Displacements which preserve the total mass are exactly the elements of $\mathcal{H}_0$.
\end{definition}
\subsection{Stochastic and variational structure}
\label{subsec:stochastic-analysis}
The reduction uses only a small portion of the infinite-dimensional stochastic toolkit, but it uses that portion essentially. Here we record the operator classes that control the covariance of the regularised noise, the notion of differentiability needed for the Hilbert-space It\^{o} formula to apply, and the driving process itself. Throughout, we fix a filtered probability space $(\Omega, \mathfrak{F}, (\mathfrak{F}_t)_{t \ge 0}, \mathbb{P})$ satisfying the usual conditions. The set $\{h_k\}_{k \in \mathbb{N}}$ denotes a complete orthonormal basis of the separable Hilbert space that carries the operator in question, and both quantities below are independent of that choice. The operator classes are treated in Appendix~C of Da~Prato and Zabczyk~\cite{daprato2014stochastic} and the $\mathcal{Q}$-Wiener process, the stochastic integral and the It\^{o} formula in Chapter~4 of the same work, while the conversion between the two stochastic conventions is developed at the level of the equations in Chapter~4 of Duan and Wang~\cite{duan2014effective}.
\begin{definition}[Hilbert--Schmidt and trace-class operators]
\label{def:hilbert-schmidt}
Let $\mathcal{U}$ and $\mathcal{K}$ be separable Hilbert spaces and let $\{h_k\}_{k\ge1}$ be a complete orthonormal basis of $\mathcal{U}$. A bounded linear operator $T: \mathcal{U} \to \mathcal{K}$ is \emph{Hilbert--Schmidt} if
\begin{equation}
    \|T\|_{\mathrm{HS}}^2 \coloneqq \sum_{k=1}^{\infty}\|T h_k\|_{\mathcal{K}}^2 < \infty ,
\end{equation}
and we write $T \in \mathcal{L}_2(\mathcal{U},\mathcal{K})$, with $\mathcal{L}_2(\mathcal{H}) \coloneqq \mathcal{L}_2(\mathcal{H},\mathcal{H})$. Neither this sum nor the trace below depends upon the basis chosen. A non-negative self-adjoint operator $Q$ upon a separable Hilbert space $\mathcal{K}$ is \emph{trace class}, written $Q \in \mathcal{L}_1(\mathcal{K})$, if
\begin{equation}
    \operatorname{Tr}(Q) \coloneqq \sum_{k=1}^{\infty}\langle Q f_k, f_k\rangle_{\mathcal{K}} < \infty
\end{equation}
for one, and then every, complete orthonormal basis $\{f_k\}$ of $\mathcal{K}$. If $T \in \mathcal{L}_2(\mathcal{U},\mathcal{K})$ then $Q = TT^*$ is non-negative, self-adjoint and trace class upon $\mathcal{K}$, with $\operatorname{Tr}(TT^*) = \|T\|_{\mathrm{HS}}^2$. The two spaces are kept apart because the diffusion operator of this work carries the space the noise lives upon into the dual $(H^1(V))^*$, and not into $\mathcal{H}$.
\end{definition}
\noindent{}Trace-class covariance is exactly the condition under which a Gaussian measure exists on $\mathcal{H}$, and it is what the mollification (Definition~\ref{def:mollifier}) is designed to secure. The spectral theorem for compact self-adjoint operators then supplies the eigenbasis used for every mode sum in this work.
\begin{definition}[Cylindrical Wiener process]
\label{def:cylindrical-wiener}
A \emph{cylindrical Wiener process} upon a separable Hilbert space $\mathcal{U}$ is the formal series $\mathcal{W}_t \coloneqq \sum_{k=1}^{\infty}\beta_k(t)\,h_k$, where $\{\beta_k\}_{k \ge 1}$ are independent standard scalar Brownian motions, and their increments after any time $s$ are jointly independent of $\mathfrak{F}_s$. The series does not converge in $\mathcal{U}$, but converges in any larger Hilbert space into which $\mathcal{U}$ embeds by a Hilbert--Schmidt map. For a separable Hilbert space $\mathcal{K}$, the integral $\int_0^t B_s\,d\mathcal{W}_s$ is a genuine $\mathcal{K}$-valued It\^{o} integral whenever $B$ is adapted with left-continuous trajectories and $B_s \in \mathcal{L}_2(\mathcal{U},\mathcal{K})$ with $\mathbb{E}\int_0^T\|B_s\|_{\mathrm{HS}}^2\,ds < \infty$, and it then obeys the It\^{o} isometry $\mathbb{E}\|\int_0^t B_s\,d\mathcal{W}_s\|_{\mathcal{K}}^2 = \mathbb{E}\int_0^t\|B_s\|_{\mathrm{HS}}^2\,ds$.
\end{definition}
\begin{definition}[Fr\'{e}chet derivative]
\label{def:frechet}
A functional $F:\mathcal{H}\to\mathbb{R}$ is \emph{Fr\'{e}chet differentiable} at $u$ if there is a vector $DF[u]\in\mathcal{H}$ for which the linear approximation is uniform in the direction:
\begin{equation}
    \lim_{\|h\|_{\mathcal{H}}\to 0}\frac{\bigl|F(u+h)-F(u)-\langle DF[u],h\rangle_{\mathcal{H}}\bigr|}{\|h\|_{\mathcal{H}}} = 0 .
\end{equation}
The weaker G\^{a}teaux derivative, obtained as the directional limit $\lim_{s\to0}s^{-1}[F(u+sv)-F(u)]$, suffices to define the chemical potential $\mu = \delta\mathcal{F}/\delta\rho$ as a critical-point condition, but not to run the It\^{o} calculus. The infinite-dimensional It\^{o} formula requires the stronger notion, since its correction term is quadratic in the diffusion and is therefore built from the second Fr\'{e}chet derivative. That term, $\operatorname{Tr}(BB^*D^2F)$, is finite because $B$ is Hilbert--Schmidt, which makes $BB^*$ trace class, and because $D^2F$ is bounded. The trace class property comes from the diffusion and not from $F$, of which nothing beyond a bounded second derivative is asked, and it is the point at which Definition~\ref{def:hilbert-schmidt} enters the dynamics rather than merely the construction of the noise.
\end{definition}

\begin{theorem}[Infinite-dimensional It\^{o} formula]
\label{thm:ito-formula}
Let $\mathcal{K}$ be a separable Hilbert space, let $\rho_t$ be a $\mathcal{K}$-valued It\^{o} process satisfying $d\rho_t = A_t\,dt + B_t\,d\mathcal{W}_t$, with $A$ adapted with left-continuous trajectories and pathwise integrable upon $[0,T]$ and with $B$ adapted, $\mathcal{L}_2(\mathcal{U},\mathcal{K})$-valued and square integrable in the sense of Definition~\ref{def:cylindrical-wiener}, and let $F:\mathcal{K}\to\mathbb{R}$ be twice Fr\'{e}chet differentiable with $F$, $DF$ and $D^2F$ uniformly continuous upon the bounded subsets of $\mathcal{K}$. Then
\begin{align}
    F(\rho_t) &= F(\rho_0) + \int_0^t \langle DF[\rho_s], A_s \rangle_{\mathcal{K}}\,ds + \int_0^t \langle DF[\rho_s], B_s\,d\mathcal{W}_s \rangle_{\mathcal{K}} \nonumber\\
    &\quad + \frac{1}{2}\int_0^t \operatorname{Tr}\!\left(B_sB_s^*D^2F[\rho_s]\right)\,ds.
\end{align}
\end{theorem}
\noindent{}The proof of this classical result may be found in the reference works on stochastic analysis in infinite dimensions, among them Da~Prato and Zabczyk~\cite{daprato2014stochastic}, Theorem~4.32, stated for processes with a predictable drift, a class which contains the adapted processes with left-continuous trajectories used here. Note that the trace term is well defined precisely because $B_sB_s^*$ is trace class on $\mathcal{K}$ (in the sense of Definition~\ref{def:hilbert-schmidt}) and $D^2F$ is bounded. For the parent equation we read the formula with $\mathcal{K} = (H^1(V))^*$ (the target space of the diffusion operator of Definition~\ref{def:diff-op}), and for the reduced equations of Section~\ref{sec:reduction} with $\mathcal{K} = \mathbb{R}$. We do not prove that the solution of the parent equation is an It\^{o} process in that space, and we assume it.

To obtain a trace-class covariance, and with it the It\^{o} calculus of Theorem~\ref{thm:ito-formula}, we do not work with bare space-time white noise. Instead we regularise the driving field by convolution with an even compactly supported mollifier, which is the only noise construction used later.
\begin{definition}[Mollified noise field]
\label{def:mollifier}
Let $\psi \in C_c^\infty(\mathbb{R}^3)$ be radial, which is to say that $\psi(\mathbf{x}) = \psi(\mathbf{y})$ whenever $|\mathbf{x}| = |\mathbf{y}|$, and therefore even, compactly supported in $B_1(0)$, and normalised by $\int_{\mathbb{R}^3}\psi(\mathbf{x})\,d\mathbf{x}=1$. For $\delta>0$, define
\begin{equation}
    \psi_\delta(\mathbf{x}) \coloneqq \delta^{-3}\psi(\mathbf{x}/\delta).
\end{equation}
Write $\mathcal{U}_0 \coloneqq L^2(\mathbb{R}^3;\mathbb{R}^3)$ and $\mathcal{U} \coloneqq L^2(V;\mathbb{R}^3)$, let $\mathcal{W}$ be a cylindrical Wiener process upon $\mathcal{U}_0$ in the sense of Definition~\ref{def:cylindrical-wiener}, and let
\begin{equation}
    J_\delta : \mathcal{U}_0 \to \mathcal{U} , \qquad J_\delta h \coloneqq (\psi_\delta * h)\big|_V
\end{equation}
be the mollification operator. The regularised conservative noise is
\begin{equation}
    \zeta_t^\delta \coloneqq \sum_{k=1}^{\infty} \beta_k(t) \, J_\delta h_k ,
\end{equation}
and the integral against it is $\int_0^t \Phi_s \, d\zeta^\delta_s \coloneqq \int_0^t \Phi_s J_\delta \, d\mathcal{W}_s$. The field is the restriction to $V$ of the mollified field upon the whole of $\mathbb{R}^3$, so that its covariance is the translation-invariant kernel $g_\delta(|\mathbf{x}-\mathbf{y}|)$ with $g_\delta = \psi_\delta * \psi_\delta$. That is the covariance used by the chord representation of Lemma~\ref{lemma:radial-closed-form} and by every kernel estimate which follows. No boundary condition is imposed upon $\zeta^\delta$ itself, since a convolution of a cylindrical Wiener field with a compactly supported mollifier does not in general satisfy one, and imposing one would destroy the translation invariance the chord formula requires. Conservation of total mass is secured instead by the divergence form of the flux in Eq.~\eqref{eq:full_spde_ito}, read with no flux through the wall, which is a property of the equation and not of the driving field. In the cumulative-mass representation that conservation appears as $m(L)=N$, while $m(0)=0$ holds by construction, from which the Dirichlet condition of Lemma~\ref{lemma:boundary-vanishing} upon the tangent fields follows in turn.
\end{definition}

\begin{lemma}[The mollified noise and its covariance]
\label{lemma:mollified-noise}
The operator $J_\delta$ is Hilbert--Schmidt from $\mathcal{U}_0$ into $\mathcal{U}$, with $\|J_\delta\|_{\mathrm{HS}}^2 = 3\,\lambda^3(V)\,g_\delta(0)$ and $g_\delta(0) = \delta^{-3}\|\psi\|_{L^2}^2$, the series defining $\zeta^\delta_t$ converges in $L^2(\Omega;\mathcal{U})$, and the covariance operator $\mathcal{Q}_\delta \coloneqq J_\delta J_\delta^*$ is non-negative, self-adjoint and trace class upon $\mathcal{U}$, acting as the integral operator with kernel $g_\delta(|\mathbf{x}-\mathbf{y}|)$. Writing $\{e_k\}$ and $\{\lambda_k\}$ for its eigenvectors and eigenvalues, and $\beta^\delta_k \coloneqq \lambda_k^{-1/2}\langle \zeta^\delta, e_k\rangle_{\mathcal{U}}$ for those $k$ with $\lambda_k > 0$, the field admits the Karhunen--Lo\`{e}ve expansion
\begin{equation}
    \zeta^\delta(t,\mathbf{x}) = \sum_{k=1}^{\infty}\sqrt{\lambda_k}\,\beta^\delta_k(t)\,e_k(\mathbf{x}),
    \label{eq:noise_expansion}
\end{equation}
convergent in $L^2(\Omega;\mathcal{U})$, with coefficients satisfying $\mathbb{E}[\beta^\delta_k(t)\beta^\delta_j(s)] = \delta_{kj}\,(t\wedge s)$.
\end{lemma}
\begin{proof}
Parseval's identity in $\mathcal{U}_0$ gives $\sum_k |(J_\delta h_k)(\mathbf{x})|^2 = 3\,\|\psi_\delta\|_{L^2}^2$ for every $\mathbf{x}$, and integration over $V$ gives the stated Hilbert--Schmidt norm, the evenness of $\psi$ identifying $\|\psi_\delta\|_{L^2}^2$ with $g_\delta(0)$. The increments being independent and centred, the partial sums of the series obey $\mathbb{E}\|\sum_{k=N+1}^{M}\beta_k(t)J_\delta h_k\|_{\mathcal{U}}^2 = t\sum_{k=N+1}^{M}\|J_\delta h_k\|_{\mathcal{U}}^2$, which is a tail of a convergent series, so the sums are Cauchy in $L^2(\Omega;\mathcal{U})$. The operator $\mathcal{Q}_\delta = J_\delta J_\delta^*$ is therefore trace class by Definition~\ref{def:hilbert-schmidt}, with $\operatorname{Tr}(\mathcal{Q}_\delta) = \|J_\delta\|_{\mathrm{HS}}^2$, and the adjoint $(J_\delta^* f)(\mathbf{y}) = \int_V \psi_\delta(\mathbf{x}-\mathbf{y})f(\mathbf{x})\,d\mathbf{x}$ turns $\langle \mathcal{Q}_\delta f,\varphi\rangle_{\mathcal{U}} = \langle J_\delta^* f, J_\delta^*\varphi\rangle_{\mathcal{U}_0}$ into the double integral of $g_\delta$ against $f$ and $\varphi$, and the substitution $\mathbf{y} = R\mathbf{w}$ in the convolution gives $g_\delta(R\mathbf{v}) = g_\delta(\mathbf{v})$ for every orthogonal $R$, the mollifier being radial, so that $g_\delta$ depends upon $|\mathbf{v}|$ alone. The stated covariance of the coefficients follows from $\langle \mathcal{Q}_\delta e_k, e_j\rangle_{\mathcal{U}} = \lambda_k\delta_{kj}$, and the expansion is that of $\zeta^\delta_t$ in the basis $\{e_k\}$.
\end{proof}

The geometric reduction projects a stochastic flow onto a submanifold, and this projection is a nonlinear change of variables. For this reason we work in the Stratonovich interpretation, which obeys the ordinary chain rule. This choice fixes the calculus of the reduction, but on its own it does not put the bare thermodynamic drift in that calculus. Indeed, the parent equation is posed in the It\^{o} frame, and its Stratonovich form differs from it by the correction given below. Since the symbol $\circ$ carries the whole of that convention, we define it through the exact It\^{o} correction that it generates, which is the seed of every anomalous drift that we compute later in this work.
\begin{definition}[Stratonovich integral against the noise]
\label{def:stratonovich}
Let $\zeta^\delta_t$ be the noise of~\eqref{eq:noise_expansion}, with eigenvectors $\{e_k\}$ and eigenvalues $\{\lambda_k\}$, and consider an It\^{o} equation $d\rho_t = a(\rho_t)\,dt + \Phi(\rho_t)\,d\zeta^\delta_t$ with a diffusion $\Phi(\rho)$ linear upon the noise and differentiable in $\rho$, with derivative $D\Phi[\rho]$. Write $\Phi_t \coloneqq \Phi(\rho_t)$, and let its \emph{sensitivity} $\Psi_t$ be the map $v \mapsto D\Phi[\rho_t]\bigl(\Phi(\rho_t)v\bigr)$. The \emph{It\^{o}--Stratonovich correction} is the series
\begin{equation}
    \mathfrak{S}(\Psi_s) \coloneqq \sum_{k=1}^{\infty} \lambda_k (\Psi_s e_k) e_k ,
    \label{eq:strat-correction-series}
\end{equation}
wherever it converges in $\mathcal{H}$, and the \emph{Stratonovich integral} of $\Phi$ against the noise is
\begin{equation}
    \int_0^t \Phi_s \circ d\zeta^\delta_s \coloneqq \int_0^t \Phi_s \, d\zeta^\delta_s + \frac{1}{2} \int_0^t \mathfrak{S}(\Psi_s) \, ds .
    \label{eq:ito-stratonovich-exact}
\end{equation}
The It\^{o} equation is accordingly the same equation as its \emph{Stratonovich form}
\begin{equation}
    d\rho_t = \Bigl[ a(\rho_t) - \frac{1}{2} \mathfrak{S}(\Psi_t) \Bigr] dt + \Phi(\rho_t) \circ d\zeta^\delta_t .
    \label{eq:stratonovich-form}
\end{equation}
\end{definition}
\noindent{}For real continuous semimartingales $X$ and $Y$ the same definition, $\int_0^t X \circ dY \coloneqq \int_0^t X \, dY + \tfrac{1}{2}[X,Y]_t$, is the one of Chapter~15 of Kallenberg~\cite{kallenberg2002foundations}, and there it is the limit in probability, as the mesh tends to zero, of the sums $\sum_i \tfrac{1}{2}(X_{t_i} + X_{t_{i+1}})(Y_{t_{i+1}} - Y_{t_i})$, which follows from Corollary~15.14 and Proposition~15.18 of that work and is the origin of the name. For the diffusion of this work, which contains a derivative, the convergence of those sums in $\mathcal{H}$ is not established, and so we take the correction as the definition. The conversion is carried out formally, from Riemann sums evaluated at the midpoints, in Section~4.5 of Duan and Wang~\cite{duan2014effective}, and Fehrman and Gess~\cite{fehrman2024deankawasaki}, who establish the well-posedness of the Dean--Kawasaki equation with correlated noise upon the torus, pass between the two forms in the same formal way.

\noindent{}Dean's derivation computes the It\^{o} differential of the empirical measure, so the regularised Dean--Kawasaki dynamics carry the thermodynamic drift of Definition~\ref{def:drift-op} and the diffusion operator of Definition~\ref{def:diff-op} in the It\^{o} frame,
\begin{equation}
    d\rho_t = \mathcal{A}(\rho_t) \, dt + \varepsilon \mathcal{B}(\rho_t) \, d\zeta^\delta_t .
    \label{eq:full_spde_ito}
\end{equation}
Since the geometric reduction is a nonlinear change of variables, it requires the Stratonovich form of the same dynamics, and we display this form in Remark~\ref{rem:correction-trace} below, once we have named the correction that it carries.
\begin{remark}[The correction for the density SPDE, and the $\mathcal{Q}$-trace]
\label{rem:correction-trace}
In the case of~\eqref{eq:full_spde_ito}, the integrand is the diffusion operator itself ($\Phi_t = \mathcal{B}(\rho_t)$), and so the sensitivity is given by the composition $\Psi_t = D\mathcal{B}[\rho_t]\cdot\mathcal{B}(\rho_t)$. Writing the \emph{$\mathcal{Q}$-trace} of an operator-valued bilinear form $T$ as
\begin{equation}
    \operatorname{Tr}_{\mathcal{Q}}(T) \coloneqq \sum_{k=1}^{\infty} T\bigl(\sqrt{\mathcal{Q}}e_k,\sqrt{\mathcal{Q}}e_k\bigr) = \sum_{k=1}^{\infty}\lambda_k\,T(e_k,e_k),
\end{equation}
Definition~\ref{def:stratonovich} gives the Stratonovich drift correction compactly as
\begin{equation}
    \mathfrak{S}_\rho = \tfrac{1}{2}\operatorname{Tr}_{\mathcal{Q}}\bigl(D\mathcal{B}[\rho]\cdot\mathcal{B}(\rho)\bigr) = \tfrac{1}{2}\sum_{k=1}^{\infty}\lambda_k\,D\mathcal{B}[\rho]\bigl(\mathcal{B}(\rho)e_k\bigr)e_k ,
\end{equation}
the factor $\tfrac12$ being the one carried by~\eqref{eq:ito-stratonovich-exact}. This means that the correction is intrinsic to the geometry of the noise (through $\mathcal{Q}_\delta$) and to the state dependence of the diffusion coefficient $\mathcal{B}$. It vanishes identically for additive noise, and it does not vanish here precisely because $\mathcal{B}(\rho)\propto\sqrt{\rho}$ is state dependent.
Since~\eqref{eq:full_spde_ito} is posed in the It\^{o} frame, the Stratonovich form of the same dynamics carries this correction with the opposite sign,
\begin{equation}
    d\rho_t = \bigl[\mathcal{A}(\rho_t) - \varepsilon^2\,\mathfrak{S}_\rho(\rho_t)\bigr]\,dt + \varepsilon \mathcal{B}(\rho_t) \circ d\zeta^\delta_t ,
    \label{eq:full_spde}
\end{equation}
which is~\eqref{eq:stratonovich-form} for this equation. In Section~\ref{sec:reduction} we start from~\eqref{eq:full_spde_ito} and form the Stratonovich form of the reduced equation in the same way.
The symbol $\mathfrak{S}$ denotes the It\^{o}--Stratonovich correction at every level of the reduction, and each decoration refers to a distinct object of the hierarchy. At the level of the field, $\mathfrak{S}(\Psi)$ is the abstract correction of Definition~\ref{def:stratonovich}, $\mathfrak{S}_\rho$ its density-level specialisation, and $\mathfrak{S}_m^\delta$ the cumulative-mass correction of Lemma~\ref{lemma:ito-strat-correction}. On the reduced manifold, $\mathfrak{S}^k$ denotes the projected It\^{o} drift (Theorem~\ref{thm:equivalent-fluctuations}), and $\mathfrak{s}_a^k$ (in lowercase) the residual geometric anomaly left after the standard Stratonovich--It\^{o} back-correction is subtracted.
\end{remark}

\smallskip

The deterministic part is generated by the Helmholtz free energy $\mathcal{F}[\rho] = \mathcal{F}_0[\rho] + \mathcal{F}_{\mathrm{exc}}[\rho]$, the ideal-gas part $\mathcal{F}_0[\rho] = \beta^{-1}\int_V \rho(\ln\rho - 1)\,d\mathbf{r}$ together with an excess part. We assume that the excess term is twice continuously Fr\'{e}chet differentiable, and that the associated chemical potential $\mu[\rho] \coloneqq \delta\mathcal{F}/\delta\rho$ lies in $H^2(V)$ for strictly positive profiles. In the case $V = B_L(0) \subset \mathbb{R}^3$, we also assume that the free energy is invariant under rotations, $\mathcal{F}[\rho \circ R] = \mathcal{F}[\rho]$ for every orthogonal $3 \times 3$ matrix $R$.
\begin{lemma}[The chemical potential of a spherically symmetric density]
\label{lemma:radial-chemical-potential}
Let $V = B_L(0) \subset \mathbb{R}^3$, with the free energy and the chemical potential above. Let $\rho \in \mathcal{D}_N^+ \cap H^1(V)$ be spherically symmetric, which is to say that $\rho(R\mathbf{x}) = \rho(\mathbf{x})$ for almost every $\mathbf{x} \in V$ and every orthogonal $3 \times 3$ matrix $R$, with $\operatorname{ess\,inf}_V \rho > 0$, and suppose that $\mu = \mu[\rho]$ has a representative continuous upon $V$. Then $\mu(R\mathbf{x}) = \mu(\mathbf{x})$ for every $\mathbf{x} \in V$ and every orthogonal $R$.
\end{lemma}
\begin{proof}
\begin{proofpart}[Rotations preserve Lebesgue measure]
\label{part:rotation-measure}
Let $R$ be an orthogonal $3 \times 3$ matrix and let $\nu_R(A) \coloneqq \lambda^3(R^{-1}A)$. For $\mathbf{t} \in \mathbb{R}^3$ the linearity of $R^{-1}$ gives $R^{-1}(A + \mathbf{t}) = R^{-1}A + R^{-1}\mathbf{t}$, a set of measure $\lambda^3(R^{-1}A)$, so $\nu_R$ is translation invariant. It gives a finite value to $[0,1]^3$, which lies in a ball that $R^{-1}$ carries onto itself. Lemma~1.29 of Kallenberg~\cite{kallenberg2002foundations}, taken with a space $S$ of one point, therefore gives $\nu_R = c_R\lambda^3$ for a constant $c_R$. Since $R^{-1}\bar B_1 = \bar B_1$, evaluating upon $\bar B_1$ gives $c_R = 1$. Lemma~1.22 of the same work then gives, for every integrable $f$ upon $V$, which $R$ carries onto itself,
\begin{equation}
    \int_V f(R\mathbf{x}) \, d\mathbf{x} = \int_V f(\mathbf{x}) \, d\mathbf{x} .
    \label{eq:rotation-invariance}
\end{equation}
\end{proofpart}
\begin{proofpart}[The chemical potential is invariant up to a constant]
\label{part:mu-invariant}
Let $h \in \mathcal{H}_0$ be bounded. Since $\rho$ is spherically symmetric, $\rho(R\mathbf{x}) = \rho(\mathbf{x})$ for almost every $\mathbf{x}$, and~\eqref{eq:rotation-invariance} shows that $h \circ R$ is bounded with integral zero. For $|s| \, \sup|h| < \operatorname{ess\,inf}\rho$ both $\rho + s h$ and $\rho + s\, h\circ R = (\rho + s h)\circ R$ lie in $\mathcal{D}_N^+$, and the rotation invariance of $\mathcal{F}$ gives
\begin{equation}
    \mathcal{F}[\rho + s \, h \circ R] = \mathcal{F}\bigl[(\rho + s h) \circ R\bigr] = \mathcal{F}[\rho + s h] .
\end{equation}
Differentiating at $s = 0$, with the derivative $\mu[\rho]$, and then using~\eqref{eq:rotation-invariance} with $R^{\mathsf T} = R^{-1}$ upon the left,
\begin{equation}
    \int_V \mu(R^{\mathsf T}\mathbf{x}) \, h(\mathbf{x}) \, d\mathbf{x}
      = \int_V \mu(\mathbf{x}) \, h(R\mathbf{x}) \, d\mathbf{x}
      = \int_V \mu(\mathbf{x}) \, h(\mathbf{x}) \, d\mathbf{x} .
\end{equation}
The bounded elements of $\mathcal{H}_0$ are dense in $\mathcal{H}_0$. For $f \in \mathcal{H}_0$ the truncations $f_n \coloneqq \max(-n,\min(f,n))$ tend to $f$ in $\mathcal{H}$ by dominated convergence, so their integrals tend to $\int_V f = 0$, and the bounded functions $f_n - \lambda^3(V)^{-1}\int_V f_n$, which lie in $\mathcal{H}_0$, tend to $f$ as well. The space $\mathcal{H}_0$ is the orthogonal complement of the constant function $1$, so $\mu \circ R^{\mathsf T} - \mu$ is orthogonal to $\mathcal{H}_0$ and is therefore a constant $c$ almost everywhere.
\end{proofpart}
\begin{proofpart}[The chemical potential is spherically symmetric]
The function $\mu \circ R^{\mathsf T} - \mu - c$ is continuous upon the open set $V$ and vanishes almost everywhere, so the open set upon which it does not vanish is null and therefore empty. Evaluating at $\mathbf{x} = 0$, which $R^{\mathsf T}$ fixes, gives $c = 0$, so $\mu(R^{\mathsf T}\mathbf{x}) = \mu(\mathbf{x})$ for every $\mathbf{x} \in V$ and every orthogonal $R$, which is the spherical symmetry of $\mu$.
\end{proofpart}
\end{proof}
\begin{definition}[The drift operator]
\label{def:drift-op}
For a strictly positive $\rho \in H^1(V)$ with $\rho\nabla\mu[\rho] \in L^2(V;\mathbb{R}^3)$, the drift operator gives the element $\mathcal{A}(\rho)$ of $(H^1(V))^*$ defined by
\begin{equation}
    \langle \mathcal{A}(\rho), \varphi \rangle \coloneqq - \int_V D_0 \, \rho(\mathbf{r}) \, \nabla\bigl(\beta\mu[\rho](\mathbf{r})\bigr) \cdot \nabla\varphi(\mathbf{r}) \, d\mathbf{r} ,
    \qquad \varphi \in H^1(V) .
    \label{eq:drift_explicit}
\end{equation}
Integration by parts writes this, for smooth arguments, as the integral of $\varphi$ against $\nabla \cdot (D_0\rho\nabla\beta\mu)$ together with the boundary integral of $\varphi$ against the outward flux $-D_0\rho\nabla\beta\mu\cdot\hat{\mathbf{n}}$, so the weak form carries the divergence together with the no-flux condition $\hat{\mathbf{n}} \cdot (\rho\nabla\mu) = 0$ upon $\partial V$. The choice $\varphi = 1$ gives $\langle \mathcal{A}(\rho), 1\rangle = 0$, so the drift conserves the total mass.
\end{definition}

The stochastic part is a conservative flux driven by the regularised field $\zeta^\delta_t$ (Definition~\ref{def:mollifier}). Here $\delta > 0$ denotes the microscopic correlation length of that field, which is distinct from the macroscopic amplitude parameter $\varepsilon$.
\begin{definition}[The diffusion operator]
\label{def:diff-op}
For a vector field $\mathbf{u} \in L^2(V;\mathbb{R}^3)$ with $\sqrt{\rho}\,\mathbf{u}$ square integrable, the diffusion operator gives the element $\mathcal{B}(\rho)\mathbf{u}$ of $(H^1(V))^*$, the dual of $H^1(V)$, defined by
\begin{equation}
    \langle \mathcal{B}(\rho)\mathbf{u}, \varphi \rangle \coloneqq -\int_V \sqrt{2 D_0 \rho(\mathbf{r})} \, \mathbf{u}(\mathbf{r}) \cdot \nabla\varphi(\mathbf{r}) \, d\mathbf{r} ,
    \qquad \varphi \in H^1(V) .
    \label{eq:diffusion_explicit}
\end{equation}
When $\rho$ and $\mathbf{u}$ are continuously differentiable, integration by parts writes~\eqref{eq:diffusion_explicit} as the integral of $\varphi$ against the divergence $\nabla\cdot(\sqrt{2D_0\rho}\,\mathbf{u})$ together with the integral over $\partial V$ of $\varphi$ against the outward flux $-\sqrt{2D_0\rho}\,\mathbf{u}\cdot\hat{\mathbf{n}}$, the second term returning through the wall what the divergence lets out. The wall is therefore reflecting for the stochastic flux, as the no-flux condition in the domain of the drift makes it for the deterministic one, and the choice $\varphi = 1$ gives $\langle \mathcal{B}(\rho)\mathbf{u}, 1 \rangle = 0$, so that the noise does not change the total mass.
\end{definition}

\begin{assumption}[Strict positivity of the density]
\label{ass:strict-positivity}
Throughout this work, we restrict the stochastic dynamics to density profiles in the subset $\mathcal{D}_N^+$ of strictly positive densities (Definition~\ref{def:density-space}). Formally, the diffusion coefficient $\sqrt{\rho}$ is not Lipschitz continuous at the boundary $\rho = 0$, and the bare Dean--Kawasaki equation ($\delta = 0$) is ill-posed as a classical SPDE. This difficulty is circumvented by two features of our formulation:
\begin{enumerate}[(i)]
    \item The noise is spatially mollified at scale $\delta > 0$ (Definition~\ref{def:mollifier}), which yields a smooth trace-class covariance operator and a well-defined regularised stochastic flux for strictly positive initial data.
    \item The stochastic dynamics are understood to hold up to the stopping time
    \begin{equation}
        \tau_\alpha \coloneqq \inf\bigl\{t > 0 : \inf_{\mathbf{r} \in V} \rho_t(\mathbf{r}) \le \alpha \bigr\},
    \end{equation}
    for an arbitrarily small threshold $\alpha > 0$. All subsequently derived formulae are valid on the event $\{t < \tau_\alpha\}$, and the law considered below is that of the process stopped at $\tau_\alpha$.
\end{enumerate}
For the nucleation problems considered in this work, the density profiles of physical interest remain uniformly bounded away from zero (cf.\ Remark~\ref{rem:gaussian-ansatz} below), and the restriction $t < \tau_\alpha$ is not expected to bind, although establishing non-attainment of the positivity boundary for a conservative equation with multiplicative noise lies beyond the scope of the present work, so that every calculation below is conditional upon $t < \tau_\alpha$.
\end{assumption}
\begin{assumption}[The parent equation]
\label{ass:parent-equation}
Upon the filtered probability space fixed above there is a process $\rho = \{\rho_t\}$ which is an It\^{o} process in the sense of Theorem~\ref{thm:ito-formula}, taken with $\mathcal{K} = (H^1(V))^*$, with initial value in $\mathcal{D}_N^+$, drift $\mathcal{A}(\rho_t)$ of Definition~\ref{def:drift-op} and diffusion $\varepsilon\mathcal{B}(\rho_t)J_\delta$, so that
\begin{equation}
    \rho_t = \rho_0 + \int_0^t \mathcal{A}(\rho_s)\,ds + \varepsilon \int_0^t \mathcal{B}(\rho_s)\,d\zeta^\delta_s
    \label{eq:parent-equation}
\end{equation}
holds in $(H^1(V))^*$ for $t < \tau_\alpha$. This is Eq.~\eqref{eq:full_spde_ito} read in that space, and it is assumed rather than proved: Konarovskyi, Lehmann and von Renesse~\cite{konarovskyi2019dean} solve the It\^{o} equation only for initial data which are finite sums of point masses, Cornalba, Shardlow and Zimmer~\cite{cornalba2019regularised,cornalba2020dean} solve a second-order model upon a periodic interval with high probability, as they do upon the torus in several dimensions~\cite{cornalba2021slender}, Fehrman and Gess~\cite{fehrman2021well,fehrman2024deankawasaki} solve the Stratonovich equation upon the torus with kinetic solutions, and Fehrman~\cite{fehrman2025neumann} solves it upon a bounded domain with a reflecting wall, in each case with a drift local in the density. None of these covers the present equation.
\end{assumption}
\begin{remark}[Geometric singularity at the origin in the mass formulation]
\label{rem:origin-singularity}
Even under Assumption~\ref{ass:strict-positivity}, the diffusion coefficient in the cumulative mass representation $\sigma[m](r) = \sqrt{2D_0 \partial_r m(r)} = \sqrt{2D_0 \cdot 4\pi r^2 \rho(r)}$ vanishes at $r = 0$ at rate $\sigma(r) = \mathcal{O}(r)$. This is due to the fact that the spherical volume element $4\pi r^2$ degenerates at the origin. This geometric zero is an artefact of the radial coordinate system, not a failure of the physical diffusion.
In order to ensure that the stochastic trace (Lemma~\ref{lemma:ito-strat-correction} below) and the subsequent projections are well defined, we work in the weighted Hilbert space $\mathcal{H}_\varpi$ (Definition~\ref{def:weighted-space}) with weight $\varpi(r) = (4\pi r^2 \rho)^{-1} = (\partial_r m)^{-1}$. The relevant projected quantity is
\begin{equation}
    \langle \mathfrak{S}_m^\delta, \mathbf{e}_k^m \rangle_\varpi = \int_0^L \mathfrak{S}_m^\delta(r) \, \mathbf{e}_k^m(r) \, \varpi(r) \, dr,
\end{equation}
and we now compute its rate of vanishing instead of asserting it. By Lemma~\ref{lemma:ito-strat-correction} of Section~\ref{sec:reduction} below, $\mathfrak{S}_m^\delta = \frac{D_0}{2\sigma}[\sigma' V_\delta + \tfrac12\sigma V_\delta']$. The local variance $V_\delta(r) = C_\delta(r,r)$ does \emph{not} tend to a finite non-zero constant at the origin. Indeed, Eq.~\eqref{eq:variance-origin} below gives $V_\delta(r) = \mathcal{O}(r^{4})$ as $r \to 0^+$ (under the normalisation of Definition~\ref{def:radial-noise-norm-ref}), since the two concentric spheres degenerate. With $\sigma = \mathcal{O}(r)$ and $\sigma' = \mathcal{O}(1)$ this gives
\begin{equation}
    \mathfrak{S}_m^\delta(r) = \mathcal{O}(r^{3}) \qquad (r \to 0^+),
    \label{eq:origin-rate}
\end{equation}
and with $\mathbf{e}_k^m = \mathcal{O}(r^{3})$ (Lemma~\ref{lemma:boundary-vanishing} below) the integrand is $\mathcal{O}(r^{3})\cdot\mathcal{O}(r^{3})\cdot\mathcal{O}(r^{-2}) = \mathcal{O}(r^{4})$ and the pairing converges absolutely. This means that there is in fact no singularity at the origin at all.
Another consequence of $\sigma(0)=0$ is that the Nemytskii map $m \mapsto \sigma[m]$ is not pointwise Fr\'{e}chet differentiable at the origin. For this reason we never apply Theorem~\ref{thm:ito-formula} to that map. Instead, we apply it to the field, with $\mathcal{K} = (H^1(V))^*$, and to the finitely many coordinates of the reduced dynamics, where the space is $\mathbb{R}^d$.
\end{remark}
\subsection{Differential geometry of submanifolds in \texorpdfstring{$\mathcal{H}$}{H}}
\label{subsec:embedding-geometry}
The reduction hypothesis states that the stochastic dynamics concentrate near a smooth finite-dimensional family of profiles. The reduced states then form a submanifold of the ambient Hilbert space, and the reduction depends on the extrinsic geometry of that submanifold. That geometry is classical. The second fundamental form, the Gauss formula and the Weingarten map appear in the form used here in Chapter~8 of Lee~\cite{lee2018riemannian}, the tubular neighbourhood on which the Fermi coordinates of Section~\ref{sec:kolmogorov-reduction} rest is Theorem~5.25 of the same work, and the Fermi construction itself is developed in Chapters~2 and~3 of Gray~\cite{gray2004tubes}, and Section~9.3 of that work gives the power-series expansion of the associated volume Jacobian.

In both of these accounts the ambient manifold is taken to be finite-dimensional, but in the present case the ambient space is $\mathcal{H}$. The curvature constructions carry over unchanged to $\mathcal{H}$, since each of them makes use of only the ambient inner product, the orthogonal splitting into tangential and normal parts that a Hilbert space supplies, and the Fr\'{e}chet derivative in place of the classical one. The case of the tubular neighbourhood is different. It needs an argument of its own, which is given in Lemma~\ref{lemma:tubular-existence}. Moreover, the Weingarten duality fails outright under the weighted metric (Section~\ref{sec:reduction}), and we discuss its corrected form in Remark~\ref{rem:weingarten-duality-fails} below. Hence we fix the dimension $d \in \mathbb{N}$ of the reduced model and an open connected parameter set $\mathcal{X} \subseteq \mathbb{R}^d$ with compact closure $\overline{\mathcal{X}}$ and coordinates $\mathbf{X} = (X^1,\dots,X^d)$. The embedding is required to extend smoothly to $\overline{\mathcal{X}}$, its density profiles to be bounded above and below by positive constants uniformly on $\overline{\mathcal{X}}$, and the parameter derivatives $\partial_{X^k}\rho(\cdot\,;\mathbf{X})$ of its density profiles to be continuous and bounded near the centre $r = 0$ of the ball, locally uniformly in $\mathbf{X}$. We then regard the reduced states as an embedded submanifold of $\mathcal{H}$.
\begin{definition}[Finite-dimensional immersion]
\label{def:immersion}
Let $\iota: \mathcal{X} \to \mathcal{H}$ be a map that assigns a field profile to every parameter configuration, denoted by $\iota(\mathbf{X}) \equiv f(\cdot; \mathbf{X})$. We say $\iota$ is a \textit{smooth immersion} if:
\begin{enumerate}[(i)]
    \item The map is of class $C^2$ in the Fr\'{e}chet sense. That is, for every $\mathbf{X}$, there exists a bounded linear operator $D\iota(\mathbf{X}): \mathbb{R}^d \to \mathcal{H}$ (the total derivative) describing the local linear approximation of the field variation.
    \item For every $\mathbf{X} \in \mathcal{X}$, the derivative $D\iota(\mathbf{X})$ is injective.
\end{enumerate}
\end{definition}
Although $\mathcal{X}$ is finite-dimensional, we retain the Fr\'{e}chet notation in order to emphasise that all differentiability statements are taken in the topology of $\mathcal{H}$.

\begin{definition}[The physical manifold]
\label{def:manifold}
We define the \textit{physical manifold} $\mathcal{M}$ as the image of the immersion:
\begin{equation}
    \mathcal{M} \coloneqq \iota(\mathcal{X}) = \left\{ \rho \in \mathcal{H} : \rho = f(\cdot; \mathbf{X}) \text{ for some } \mathbf{X} \in \mathcal{X} \right\}.
\end{equation}
We assume $\iota$ is a topological embedding, such that $\mathcal{M}$ is a submanifold of $\mathcal{H}$ inheriting the subspace topology.
\end{definition}

The local linearisation of $\mathcal{M}$ is encoded by the tangent vectors obtained by differentiating the embedding with respect to the collective coordinates.
\begin{definition}[Tangent space]
\label{def:tangent-space}
For any point $\rho = \iota(\mathbf{X}) \in \mathcal{M}$, the \textit{tangent space} $T_{\rho}\mathcal{M}$ is the $d$-dimensional linear subspace of $\mathcal{H}$ spanned by the image of the differential $D\iota(\mathbf{X})$. We define the natural coordinate basis vectors $\{ \mathbf{e}_k \}_{k=1}^d \subset \mathcal{H}$ as the partial derivatives of the embedding:
\begin{equation}
    \mathbf{e}_k(\cdot; \mathbf{X}) \coloneqq \frac{\partial f(\cdot; \mathbf{X})}{\partial X^k}.
\end{equation}
Since $\iota$ is an immersion, the set $\{ \mathbf{e}_1, \dots, \mathbf{e}_d \}$ is linearly independent in $\mathcal{H}$.
\end{definition}

To project the dynamics, it is convenient to package the tangent basis into the Jacobian of the embedding.
\begin{definition}[The Jacobian / pushforward]
\label{def:jacobian}
We identify the tangent space of the parameter domain $T_{\mathbf{X}}\mathcal{X}$ with $\mathbb{R}^d$. The \textit{Jacobian operator} (or differential pushforward) $\mathcal{J}_{\iota}(\mathbf{X}): \mathbb{R}^d \to \mathcal{H}$ is defined by the action of the Fr\'{e}chet derivative on a parameter velocity vector $\mathbf{v} \in \mathbb{R}^d$:
\begin{equation}
    \mathcal{J}_{\iota}(\mathbf{X}) \mathbf{v} \coloneqq D\iota(\mathbf{X}) \cdot \mathbf{v} = \sum_{k=1}^d v^k \mathbf{e}_k(\cdot; \mathbf{X}).
\end{equation}
The image of $\mathcal{J}_{\iota}(\mathbf{X})$ is precisely $T_{\rho}\mathcal{M}$.
\end{definition}
If the parameters evolve along a path $\mathbf{X}_t$, then the induced variation of the field is simply $\frac{d\rho}{dt} = \mathcal{J}_{\iota}(\mathbf{X}_t)\frac{d\mathbf{X}}{dt}$.

\begin{definition}[The adjoint Jacobian / pullback]
\label{def:adjoint-jacobian}
The \textit{adjoint Jacobian} $\mathcal{J}_{\iota}^*(\mathbf{X}): \mathcal{H} \to \mathbb{R}^d$ is defined as the unique operator satisfying the duality relation with respect to the standard $L^2$ inner product:
\begin{equation}
    \langle u, \mathcal{J}_{\iota}(\mathbf{X}) \mathbf{v} \rangle_{\mathcal{H}} = \langle \mathcal{J}_{\iota}^*(\mathbf{X}) u, \mathbf{v} \rangle_{\mathbb{R}^d}, \quad \forall u \in \mathcal{H}, \mathbf{v} \in \mathbb{R}^d.
\end{equation}
The components of the adjoint action are given by the projection of the field $u$ onto the coordinate basis vectors:
\begin{equation}
    \left[ \mathcal{J}_{\iota}^*(\mathbf{X}) u \right]_k = \langle u, \mathbf{e}_k(\cdot; \mathbf{X}) \rangle_{\mathcal{H}} = \int_V u(\mathbf{r}) \frac{\partial f(\mathbf{r}; \mathbf{X})}{\partial X^k} \, d\mathbf{r}.
\end{equation}
\end{definition}
\begin{remark}[Metric dependence]
\label{rem:metric-dependence}
Note that the adjoint depends on the chosen inner product. Here we use the flat $L^2$ pairing, while in the projection theory of Section~\ref{sec:projection} the relevant adjoint is the one induced by the weighted metric that the noise selects.
\end{remark}
Finally, we record the extrinsic curvature data of the embedding, which we use from Assumption~\ref{ass:uniform-curvature} onwards and again in the Fermi-coordinate expansion of Section~\ref{sec:kolmogorov-reduction}. Moreover, we write $\langle\cdot,\cdot\rangle$ for the ambient inner product in force at the relevant stage, and $\Pi_{\mathbf{X}}^{\parallel}$ and $\Pi_{\mathbf{X}}^{\perp}$ for the orthogonal projections onto the tangent space and its complement.
\begin{definition}[Normal bundle, second fundamental form and shape operator]
\label{def:second-fundamental-form}
The \emph{normal space} at $\mathbf{X}$ is the orthogonal complement of the tangent space within the ambient space,
\begin{equation}
    \mathcal{N}_{\mathbf{X}} \coloneqq \bigl(T_{\mathbf{X}}\mathcal{M}\bigr)^{\perp},
\end{equation}
and the \emph{normal bundle} is the disjoint union of those spaces over $\mathbf{X}\in\mathcal{X}$. Writing $D_{\mathbf{u}}\mathbf{w} \coloneqq \sum_i U^i \partial_i \mathbf{w}$ for the ambient derivative of a field along a tangent field $\mathbf{u} = \sum_i U^i\mathbf{e}_i$, the \emph{second fundamental form} is
\begin{equation}
    \mathbf{I\!I}(\mathbf{u},\mathbf{v}) \coloneqq \Pi_{\mathbf{X}}^{\perp}\bigl(D_{\mathbf{u}}\mathbf{v}\bigr),
\end{equation}
and for a normal field $\tilde{\mathbf{n}}$ of class $C^1$ with value $\mathbf{n}$ at $\mathbf{X}$ the \emph{shape operator} is $A_{\mathbf{n}}\mathbf{v} \coloneqq -\Pi^{\parallel}_{\mathbf{X}}(D_{\mathbf{v}}\tilde{\mathbf{n}})$.
\end{definition}
\begin{lemma}[Properties of the second fundamental form]
\label{lemma:second-fundamental-form}
The second fundamental form takes its values in $\mathcal{N}_{\mathbf{X}}$ and is symmetric, it vanishes identically when the immersion is affine, and whenever the ambient inner product is independent of the base point the shape operator obeys the classical Weingarten duality
\begin{equation}
    \langle A_{\mathbf{n}}\mathbf{u},\mathbf{v}\rangle = \langle \mathbf{I\!I}(\mathbf{u},\mathbf{v}),\mathbf{n}\rangle ,
    \label{eq:classical-weingarten-duality}
\end{equation}
so that $A_{\mathbf{n}}$ depends upon the normal field through its value $\mathbf{n}$ alone.
\end{lemma}
\begin{proof}
The normal projector takes its values in the normal space. Expanding the ambient derivative, $D_{\mathbf{u}}\mathbf{v} - D_{\mathbf{v}}\mathbf{u}$ is a combination of the frame vectors together with $\sum_{i,j}U^iV^j(\partial_i\mathbf{e}_j - \partial_j\mathbf{e}_i)$, and the second sum vanishes because the second derivatives of a map of class $C^2$ commute, which is Proposition~2.4.14 of Abraham, Marsden and Ratiu~\cite{abraham1988manifolds}, so the difference is tangent and the form is symmetric. For an affine immersion the frame is constant and $D_{\mathbf{u}}\mathbf{v}$ is tangent, so the form vanishes. For the duality, the function $\langle\tilde{\mathbf{n}},\mathbf{v}\rangle$ vanishes identically along the manifold, and differentiating it along $\mathbf{u}$ gives $\langle D_{\mathbf{u}}\tilde{\mathbf{n}},\mathbf{v}\rangle + \langle\tilde{\mathbf{n}},D_{\mathbf{u}}\mathbf{v}\rangle = 0$, the two terms of which are $-\langle A_{\mathbf{n}}\mathbf{u},\mathbf{v}\rangle$ and $\langle\mathbf{n},\mathbf{I\!I}(\mathbf{u},\mathbf{v})\rangle$ by the self-adjointness of the projectors.
\end{proof}
\begin{remark}[Failure of the Weingarten duality under a state-dependent weight]
\label{rem:weingarten-duality-fails}
Identity~\eqref{eq:classical-weingarten-duality} does \emph{not} carry over to the weighted pairing $\langle\cdot,\cdot\rangle_{\varpi}$ of Definition~\ref{def:weighted-space}, used from Section~\ref{sec:reduction} onwards. Lemma~\ref{lemma:second-fundamental-form} follows the proof of Proposition~8.4 of Lee~\cite{lee2018riemannian}. There we differentiated the orthogonality relation $\langle\mathbf{n},\mathbf{e}_k^m\rangle = 0$ along the manifold and kept only the two terms where $\mathbf{n}$ and $\mathbf{e}_k^m$ are differentiated. If the weight $\varpi(\cdot\,;\mathbf{X})$ depends on the base point, the product rule also produces a third term, proportional to $\partial_{X^l}\varpi$, which does not vanish in general. Hence we prove a corrected relation that includes this term (Lemma~\ref{lemma:weingarten-weighted}), and in the Fermi-coordinate expansion of Section~\ref{sec:kolmogorov-reduction} we use it in place of~\eqref{eq:classical-weingarten-duality}.
\end{remark}
\subsection{The Wasserstein-2 geometry}
\label{subsec:wasserstein-geometry}
The Hilbert space $\mathcal{H}=L^2(V)$ carries the stochastic calculus, but the flat $L^2$ inner product is the wrong instrument for a mass-conserving system, since it measures pointwise discrepancy between densities rather than the cost of transporting mass, so that spatially disjoint fluctuations are orthogonal regardless of how far apart they lie. Hence the transport structure we use is that of Otto~\cite{otto2001geometry}, Equations~(8) and~(9), and Chapter~15 of Villani~\cite{villani2009optimal}, and only in the minimal form that Theorem~\ref{thm:wasserstein-isometry} needs in order to fix the weight $\varpi$ from the transport geometry alone, without reference to the noise. It is worth mentioning that in Section~1.3 of Otto's work a tangent vector is identified with the potential that generates it through $-\nabla\cdot(\rho\nabla\psi)$ (his Equation~(8)), and the space of densities is equipped with the metric $\int_V \rho\,\nabla\psi_u\cdot\nabla\psi_v\,d\mathbf{r}$ (his Equation~(9)). The calculus built on that identification is the subject of Chapter~15 of Villani.

Note that both Otto and Villani point out that the construction is formal, not a Riemannian manifold in the strict sense. This is the status that the definitions below inherit. Little is lost by this here. Indeed, we use the transport structure to fix $\varpi$ and for nothing else. From then on we carry out every computation in the weighted Hilbert space of Definition~\ref{def:weighted-space}, or on the finite-dimensional manifold of Section~\ref{subsec:embedding-geometry}. Ambrosio, Gigli and Savar\'{e}~\cite{ambrosio2008gradient} set out the rigorous infinite-dimensional theory in Chapters~8 and~11, for a reader who wants the statements themselves instead of the calculus. Our Theorem~\ref{thm:wasserstein-isometry} is the radial counterpart of Theorem~6.0.2 of that work, where the Wasserstein distance between measures on the line becomes an $L^2$ distance between inverse cumulative distribution functions. We also show in Theorem~\ref{thm:identification} that the drift of~\eqref{eq:full_spde_ito} carries the gradient-flow structure that Jordan, Kinderlehrer and Otto~\cite{jordan1998variational} established for the Fokker--Planck equation in Theorem~5.1 of their work.
\begin{definition}[The elliptic transport operator]
\label{def:transport-op}
For a profile $\rho \in \mathcal{D}_N^+\cap H^1(V)\cap L^\infty(V)$ with $\operatorname{ess\,inf}_V \rho > 0$, and with $\dot{H}^1(V) \coloneqq H^1(V)/\mathbb{R}$ the potentials modulo constants, the transport operator $\mathcal{K}_\rho: \dot{H}^1(V)\to(\dot{H}^1(V))^*$ is
\begin{equation}
    \langle \mathcal{K}_\rho\psi, \phi \rangle \coloneqq \int_V \rho(\mathbf{r})\,\nabla\psi(\mathbf{r})\cdot\nabla\phi(\mathbf{r})\,d\mathbf{r} ,
    \qquad \phi \in \dot{H}^1(V) ,
    \label{eq:transport-op}
\end{equation}
which is the weak form of $\mathcal{K}_\rho\psi = -\nabla\cdot(\rho\nabla\psi)$ together with the homogeneous Neumann condition $\hat{\mathbf{n}}\cdot\nabla\psi = 0$ on $\partial V$, that condition removing the boundary term of the integration by parts.
\end{definition}
\begin{lemma}[Isomorphism property]
\label{lemma:transport-isomorphism}
On the quotient $\dot{H}^1(V) \coloneqq H^1(V)/\mathbb{R}$ of potentials modulo constants, $\mathcal{K}_\rho$ is an isomorphism onto the dual space,
\begin{equation}
    \mathcal{K}_\rho : \dot{H}^1(V) \longrightarrow \bigl(\dot{H}^1(V)\bigr)^{*} \eqqcolon \dot{H}^{-1}(V) ,
    \label{eq:transport-weak}
\end{equation}
understood weakly through $\langle \mathcal{K}_\rho\psi, \phi\rangle = \int_V \rho\,\nabla\psi\cdot\nabla\phi\,d\mathbf{r}$ for every $\phi \in \dot{H}^1(V)$. The image is the dual space and not $\mathcal{H}_0$, and an $L^2$-valued realisation would require elliptic regularity of $\rho$ and of $\partial V$ which we neither assume nor use. Since $\mathcal{H}_0$ embeds continuously in $\dot{H}^{-1}(V)$, the inverse restricts to $\mathcal{H}_0$, and the composition
\begin{equation}
    \mathcal{K}_\rho^{-1} : \mathcal{H}_0 \longrightarrow \dot{H}^1(V) \hookrightarrow \mathcal{H}_0
    \label{eq:transport-compact}
\end{equation}
is a compact, positive, self-adjoint operator on $\mathcal{H}_0$. It is this realisation that the Wasserstein metric below employs.
\end{lemma}
\begin{proof}
Since $\rho$ is bounded below by a positive constant and above by Definition~\ref{def:transport-op}, the bilinear form $a(\psi,\phi)=\int_V\rho\,\nabla\psi\cdot\nabla\phi\,d\mathbf{r}$ is bounded and coercive on $\dot{H}^1(V)$ by the Poincar\'{e}--Wirtinger inequality. The Lax--Milgram theorem identifies $\mathcal{K}_\rho$ with the isomorphism $\dot{H}^1(V)\to\dot{H}^{-1}(V)$ induced by $a$, which is~\eqref{eq:transport-weak}. For $u\in\mathcal{H}_0$ the functional $\phi\mapsto\langle u,\phi\rangle_{\mathcal{H}}$ lies in $\dot{H}^{-1}(V)$, so $\psi = \mathcal{K}_\rho^{-1}u$ exists and is unique in $\dot{H}^1(V)$. Composing with the embedding $\dot{H}^1(V)\hookrightarrow\mathcal{H}_0$, compact by the Rellich--Kondrachov theorem, gives a compact operator on $\mathcal{H}_0$. It is symmetric because $a$ is symmetric, and positive because $\langle\mathcal{K}_\rho^{-1}u,u\rangle_{\mathcal{H}} = a(\psi,\psi) > 0$ for $u \neq 0$, which is~\eqref{eq:transport-compact}.
\end{proof}
\begin{definition}[Wasserstein-2 metric tensor]
\label{def:wasserstein-metric}
For tangent vectors $u,v \in T_\rho\mathcal{D}_N^+\cong\mathcal{H}_0$, the \emph{Wasserstein-2 inner product} is
\begin{equation}
    g^{\mathcal{W}_2}_\rho(u,v) \coloneqq \langle u, \mathcal{K}_\rho^{-1}v\rangle_{\mathcal{H}}
    = \int_V \rho(\mathbf{r})\,\nabla\psi_u(\mathbf{r})\cdot\nabla\psi_v(\mathbf{r})\,d\mathbf{r},
\end{equation}
where $\psi_u = \mathcal{K}_\rho^{-1}u$ and $\psi_v = \mathcal{K}_\rho^{-1}v$ are the generating potentials, which exist and are unique by Lemma~\ref{lemma:transport-isomorphism}. The second expression exhibits the metric as penalising the flux $\mathbf{j}=\rho\nabla\psi$ rather than the density change itself.
\end{definition}
\noindent{}As is well known, the Riemannian gradient $\operatorname{grad}_g\mathcal{F}[\rho]$ is the unique tangent vector with $g_\rho(\operatorname{grad}_g\mathcal{F},w) = \langle D\mathcal{F}[\rho],w\rangle_{\mathcal{H}}$ for all $w$. With this structure the drift of~\eqref{eq:full_spde_ito} becomes a steepest descent (Theorem~\ref{thm:identification} below), instead of an arbitrary transport term. In this sense the parent equation carries a gradient-flow structure.
\begin{theorem}[Identification of the drift]
\label{thm:identification}
For every $\rho\in\mathcal{D}_N^+\cap H^1(V)\cap L^\infty(V)$ with $\operatorname{ess\,inf}_V \rho > 0$ in the domain of $\mathcal{A}$, and every $w \in \mathcal{H}_0$,
\begin{equation}
    g^{\mathcal{W}_2}_\rho\bigl(\mathcal{A}(\rho), w\bigr) = -D_0\,\beta\,\langle D\mathcal{F}[\rho], w\rangle_{\mathcal{H}} ,
    \label{eq:identification}
\end{equation}
the pairing on the left being the duality pairing against the potential $\mathcal{K}_\rho^{-1}w$, so that the drift is the Wasserstein gradient flow of the free energy, $\mathcal{A}(\rho) = -D_0\,\beta\,\operatorname{grad}_{\mathcal{W}_2}\mathcal{F}[\rho]$.
\end{theorem}
\begin{proof}
Let $w\in\mathcal{H}_0$ and let $\psi_w = \mathcal{K}_\rho^{-1}w \in \dot{H}^1(V)$ be its potential. The weak form of Definition~\ref{def:drift-op}, evaluated at that potential, gives
\begin{equation}
    g^{\mathcal{W}_2}_\rho\bigl(\mathcal{A}(\rho), w\bigr)
      = \langle \mathcal{A}(\rho), \psi_w \rangle
      = -D_0\beta \int_V \rho\,\nabla\mu[\rho]\cdot\nabla\psi_w\,d\mathbf{r} ,
\end{equation}
with no integration by parts. The chemical potential lies in $H^2(V)$, so it is an admissible test function in~\eqref{eq:transport-op}, which turns the integral into $\langle \mathcal{K}_\rho\psi_w, \mu[\rho]\rangle = \langle w, \mu[\rho]\rangle_{\mathcal{H}}$, and $D\mathcal{F}[\rho] = \mu[\rho]$. As $w$ was arbitrary, the defining property of the Riemannian gradient recorded above identifies $-\mathcal{A}(\rho)/(D_0\beta)$ with $\operatorname{grad}_{\mathcal{W}_2}\mathcal{F}[\rho]$, the gradient being understood in the dual sense of the pairing above, which is the claim.
\end{proof}
\subsection{Weak formulation of the forward Kolmogorov equation}
\label{subsec:functional-fpe}
For the probabilistic reduction of Sections~\ref{sec:kolmogorov-reduction} and~\ref{sec:commutativity} we only need the weak generator description of the law of the regularised SPDE. This description is classical. Da~Prato and Zabczyk~\cite{daprato2014stochastic} set it out in Section~14.2.2, and the backward equation on the transition semigroup is treated in Chapter~9 of the same work, and Cerrai~\cite{cerrai2001second} treated the same equations in Hilbert spaces by probabilistic means in Chapter~5 of the book. Let us consider a cylindrical test functional $\Phi(\rho) = f(\langle\rho,\varphi_1\rangle,\dots,\langle\rho,\varphi_n\rangle)$, with $\varphi_i \in H^1(V)$, with $f$ of class $C^2$ and bounded together with its first two derivatives, and with the pairings those of $(H^1(V))^*$ against $H^1(V)$. Hence its derivative $D\Phi[\rho] = \sum_i \partial_i f\,\varphi_i$ lies in $H^1(V)$, and the drift of Definition~\ref{def:drift-op} can be paired against it. For such a functional, we read the It\^{o} formula of Theorem~\ref{thm:ito-formula} in $(H^1(V))^*$ for the process of Assumption~\ref{ass:parent-equation} stopped at $\tau_\alpha$, with the inner product of that space read through the Riesz map as the duality pairing below, and this gives
\begin{equation}
    \begin{aligned}
    \mathbb{E}[\Phi(\rho_{t\wedge\tau_\alpha})] &= \mathbb{E}[\Phi(\rho_0)] + \mathbb{E}\int_0^{t\wedge\tau_\alpha} \mathcal{L}\Phi(\rho_s)\,ds ,
    \\
    \mathcal{L}\Phi(\rho) &\coloneqq \langle \mathcal{A}(\rho), D\Phi[\rho] \rangle + \frac{\varepsilon^2}{2} \operatorname{Tr}\!\left(Q_{\rho,\delta} D^2\Phi[\rho]\right),
    \end{aligned}
    \label{eq:generator-representation}
\end{equation}
where $Q_{\rho,\delta} = \mathcal{B}(\rho)J_\delta(\mathcal{B}(\rho)J_\delta)^*$ is the covariance of the diffusion, and the stochastic term has zero mean. We recognise here the Kolmogorov operator of Section~14.2.2 of Da~Prato and Zabczyk, and the identity recorded in Definition~\ref{def:ffpe} below is Equation~(14.8) of that work. In the same section, Da~Prato and Zabczyk observed that in infinite dimensions the class $C_b^2(\mathcal{H})$ is not an adequate space of test functionals, since the Kolmogorov operator is not defined on the whole of it, and so we use the cylindrical class here.

\begin{definition}[Weak solution to the functional forward Kolmogorov equation]
\label{def:ffpe}
Let $\mu_t$ denote the law of the process $\rho_{t\wedge\tau_\alpha}$ stopped at the time of Assumption~\ref{ass:strict-positivity}, and let $\mathcal{D}_\alpha$ be the set of profiles $\rho$ with $\inf_V \rho > \alpha$. We say that $\mu_t$ is a \textit{weak solution} of the functional forward Kolmogorov equation if, for every cylindrical test functional $\Phi$,
\begin{equation}
    \frac{d}{dt}\int \Phi(\rho) \, d\mu_t(\rho) = \int_{\mathcal{D}_\alpha} \mathcal{L}\Phi(\rho) \, d\mu_t(\rho).
    \label{eq:weak-ffpe}
\end{equation}
\end{definition}
\noindent{}Moreover, the restriction to $\mathcal{D}_\alpha$ is exact, since before $\tau_\alpha$ the stopped process lies in $\mathcal{D}_\alpha$, and from $\tau_\alpha$ onwards it is frozen at $\rho_{\tau_\alpha}$, which lies outside $\mathcal{D}_\alpha$ whenever $t \mapsto \inf_V \rho_t$ is continuous. As a result, $\mathbb{E}[\mathbf{1}_{\{s<\tau_\alpha\}}\,\mathcal{L}\Phi(\rho_s)] = \int_{\mathcal{D}_\alpha} \mathcal{L}\Phi \, d\mu_s$, and Equation~\eqref{eq:generator-representation} is the integrated form of this identity. In what follows, we will apply the generator only to functionals that factor through a finite-dimensional projection, which is the restriction imposed at Definition~\ref{def:adjoint-evolution}.%%
\section{Dimensional reduction via cumulative mass}
\label{sec:reduction}
The governing SPDE posed in Section~\ref{sec:preliminaries} describes the stochastic evolution of the density field in $\mathcal{H} = L^2(V)$, driven by the three-dimensional mollified noise field $\zeta^\delta_t$ (Definition~\ref{def:mollifier}), and it holds as an identity in the dual $(H^1(V))^*$. In this section we project the SPDE onto a finite-dimensional submanifold, making use of the spherical symmetry of the domain and of a change of dependent variable from the density to the cumulative mass. We reduce the SPDE in three steps. Since the domain is a ball, $V = B_L(0) \subset \mathbb{R}^3$, we can restrict to the closed subspace of spherically symmetric functions $\mathcal{H}_{\text{rad}} \coloneqq L^2_{\text{rad}}(V)$, provided the three-dimensional noise $\zeta^\delta_t$ is projected onto its radial component by angular integration over concentric spheres, which gives a one-dimensional stochastic forcing, and the departure from spherical symmetry that the noise produces is neglected under the spherical-mean closure of Assumption~\ref{ass:spherical-closure}. In that subspace we change the dependent variable from the density to the cumulative mass $m(r,t)$, which absorbs the factor $4\pi r^2$ and turns mass conservation into a Dirichlet boundary condition. The mass SPDE is then projected onto the tangent bundle of a $d$-dimensional submanifold, in the weighted inner product consistent with the noise covariance structure, and the Stratonovich drift correction comes out in closed form. This resolves the It\^{o}--Stratonovich ambiguity, and the comparison of Section~\ref{sec:commutativity} rests on it.
\subsection{The cumulative mass representation}
\label{subsec:mass-rep}
Throughout this section the domain is the ball $V = B_L(0) \subset \mathbb{R}^3$, and the state is restricted to the subspace $\mathcal{H}_{\text{rad}} \coloneqq L^2_{\text{rad}}(V)$ of spherically symmetric functions (this restriction is the spherical-mean closure of Assumption~\ref{ass:spherical-closure} below).
\subsubsection{The mass operator}
\begin{definition}[The cumulative mass operator]
\label{def:mass-operator}
We define the operator $\mathcal{T}$ upon $\mathcal{H}$ via integration over the balls $B_r$, which for a spherically symmetric density is radial integration:
\begin{equation}
    m(r) = \mathcal{T}[\rho](r) \coloneqq \int_{B_r} \rho \, d\mathbf{x} ,
    \qquad\text{so that}\qquad
    \mathcal{T}[\rho](r) = \int_0^r \rho(s) \, 4\pi s^2 \, ds \quad \text{for } \rho \in \mathcal{H}_{\text{rad}} .
\end{equation}
Restricted to $\mathcal{H}_{\text{rad}}$, this is a bounded and injective linear operator into $H^1([0,L])$, with image the subspace $\mathcal{M}_{\mathrm{mass}}$ of Eq.~\eqref{eq:mass-image} rather than the whole of $H^1([0,L])$. Its inverse, $\mathcal{T}^{-1}$, is defined upon that image by
\begin{equation}
    \rho[m](r) \coloneqq \mathcal{T}^{-1}[m](r) = \frac{1}{4\pi r^2} \frac{\partial m(r)}{\partial r}.
\end{equation}
\end{definition}
\begin{remark}[The inverse density relation]
\label{rem:inverse-density}
For any $m \in H^1([0, L])$ in the image of $\mathcal{T}$, the corresponding symmetric density $\hat{\rho} = \mathcal{T}^{-1}[m]$ is recovered by the Fundamental Theorem of Calculus. Differentiating the integral definition with respect to the upper limit $r$:
\begin{equation}
    \frac{\partial m}{\partial r} = \frac{\partial}{\partial r} \int_0^r \rho(s) 4\pi s^2 \, ds = \rho(r) 4\pi r^2.
\end{equation}
Solving for density yields the reconstruction formula: $\hat{\rho}[m](r) = \frac{1}{4\pi r^2} \frac{\partial m(r)}{\partial r}.$ This relation allows us to express functionals of $\rho$ in terms of the gradient $\partial_r m$ alone.
\end{remark}
\subsubsection{The mollified radial noise}
\label{subsubsec:radial-noise}
To resolve the geometrical scaling factors and avoid ultraviolet divergences in the subsequent calculus, we must project the three-dimensional regularised noise onto the radial coordinate. Instead of the formal abstraction of a one-dimensional white noise, we work strictly with the radial projection of the mollified noise field $\zeta^\delta_t$ (Definition~\ref{def:mollifier}). We obtain this projection by integrating the normal stochastic flux over each sphere. Note that its variance scales with the area $4\pi r^2$ of the sphere, and the normalised field divides it by the square root $\sqrt{4\pi r^2}$ of that area.
\begin{definition}[Mollified radial noise]
\label{def:radial-noise-norm-ref}
Let $\zeta^\delta_t(\mathbf{r})$ be the three-dimensional mollified noise field constructed in Definition~\ref{def:mollifier} and restricted to $V = B_L(0)$, with ambient covariance $g_\delta(|\mathbf{x}-\mathbf{y}|)$. We define the \textit{integrated radial noise process} $\tilde{W}^\delta(r, t)$ by the surface integral of the normal stochastic flux:
\begin{equation}
    d\tilde{W}^\delta(r, t) \coloneqq \oint_{\partial B_r} d\zeta^\delta_t(\mathbf{s}) \cdot \hat{\mathbf{n}} \, dS.
\end{equation}
We define the normalised scalar radial noise field $\xi^\delta(r, t)$ by
\begin{equation}
    d\tilde{W}^\delta(r, t) = \sqrt{4\pi r^2} \, d\xi^\delta(r, t).
\end{equation}
\end{definition}
\begin{lemma}[Smooth radial covariance]
\label{lemma:radial-covariance}
Let $C_\delta(r, s) \, dt \coloneqq \mathbb{E}[d\xi^\delta(r, t) d\xi^\delta(s, t)]$ be the spatial covariance kernel of the radial noise field:
\begin{enumerate}[(i)]
    \item $C_\delta(r, s)$ is a symmetric function in $r$ and $s$.
    \item The partial derivative of the covariance kernel evaluated at the diagonal $s=r$ is related to the radial gradient of the variance profile $V_\delta(r) \coloneqq C_\delta(r, r)$ by the identity:
    \[
    \left. \frac{\partial C_\delta(r, s)}{\partial s} \right|_{s=r} = \frac{1}{2} \frac{d V_\delta(r)}{dr}.
    \]
\end{enumerate}
\end{lemma}
\begin{proof}
\begin{proofpart}[Symmetry of the kernel]
By Definition~\ref{def:radial-noise-norm-ref} the radial field is normalised by the square root of the area of the sphere, so that the kernel is obtained from the ambient covariance $\mathcal{Q}_\delta(\mathbf{x},\mathbf{y}) = \mathbb{E}[d\zeta^\delta_t(\mathbf{x}) \otimes d\zeta^\delta_t(\mathbf{y})]/dt$ by division with $\sqrt{4\pi r^2}\,\sqrt{4\pi s^2} = 4\pi r s$:
\begin{equation}
    C_\delta(r,s) = \frac{1}{4\pi r s} \oint_{\partial B_r}\!\oint_{\partial B_s} \hat{\mathbf{n}}(\mathbf{x})^{T} \mathcal{Q}_\delta(\mathbf{x},\mathbf{y})\, \hat{\mathbf{n}}(\mathbf{y}) \, dS_x \, dS_y .
    \label{eq:kernel-surface}
\end{equation}
We evaluate $C_\delta(s,r)$ and transform it back to $C_\delta(r,s)$. Writing~\eqref{eq:kernel-surface} with the arguments interchanged, relabelling the integration variables $\mathbf{x} \leftrightarrow \mathbf{y}$, using the fact that the integrand is a scalar and therefore equal to its own transpose, and finally invoking the symmetry $\mathcal{Q}_\delta(\mathbf{y},\mathbf{x})^{T} = \mathcal{Q}_\delta(\mathbf{x},\mathbf{y})$ of the ambient covariance:
\begin{align}
    C_\delta(s,r)
    &= \frac{1}{4\pi s r} \oint_{\partial B_s}\!\oint_{\partial B_r} \hat{\mathbf{n}}(\mathbf{x})^{T} \mathcal{Q}_\delta(\mathbf{x},\mathbf{y})\, \hat{\mathbf{n}}(\mathbf{y}) \, dS_x \, dS_y \nonumber \\
    &= \frac{1}{4\pi r s} \oint_{\partial B_r}\!\oint_{\partial B_s} \hat{\mathbf{n}}(\mathbf{y})^{T} \mathcal{Q}_\delta(\mathbf{y},\mathbf{x})\, \hat{\mathbf{n}}(\mathbf{x}) \, dS_y \, dS_x \nonumber \\
    &= \frac{1}{4\pi r s} \oint_{\partial B_r}\!\oint_{\partial B_s} \Bigl[\hat{\mathbf{n}}(\mathbf{y})^{T} \mathcal{Q}_\delta(\mathbf{y},\mathbf{x})\, \hat{\mathbf{n}}(\mathbf{x})\Bigr]^{T} dS_x \, dS_y \nonumber \\
    &= \frac{1}{4\pi r s} \oint_{\partial B_r}\!\oint_{\partial B_s} \hat{\mathbf{n}}(\mathbf{x})^{T} \mathcal{Q}_\delta(\mathbf{y},\mathbf{x})^{T} \hat{\mathbf{n}}(\mathbf{y}) \, dS_x \, dS_y \nonumber \\
    &= \frac{1}{4\pi r s} \oint_{\partial B_r}\!\oint_{\partial B_s} \hat{\mathbf{n}}(\mathbf{x})^{T} \mathcal{Q}_\delta(\mathbf{x},\mathbf{y})\, \hat{\mathbf{n}}(\mathbf{y}) \, dS_x \, dS_y
     = C_\delta(r,s),
\end{align}
where the second equality is the relabelling, the third uses scalar invariance under transposition, the fourth is the transposition itself, and the fifth is the symmetry of $\mathcal{Q}_\delta$. This establishes part~(i).
\end{proofpart}
\begin{proofpart}[Differentiation on the diagonal]
The variance profile is the restriction of the kernel to the diagonal, $V_\delta(r) = C_\delta(r,r)$. Since $\zeta^\delta$ is a mollified field, the kernel $g_\delta$ of $\mathcal{Q}_\delta$ is smooth by Definition~\ref{def:mollifier}, and~\eqref{eq:kernel-surface} therefore defines a $C^\infty$ function of $(r,s)$ on $(0,L]^2$, so the chain rule applies to the composition $r \mapsto (r,r) \mapsto C_\delta(r,r)$:
\begin{align}
    \frac{d V_\delta}{dr}(r)
    &= \frac{d}{dr}\,C_\delta(r,r) \nonumber \\
    &= \left.\frac{\partial C_\delta(r,s)}{\partial r}\right|_{s=r} \cdot \frac{dr}{dr}
     + \left.\frac{\partial C_\delta(r,s)}{\partial s}\right|_{s=r} \cdot \frac{dr}{dr} \nonumber \\
    &= \left.\frac{\partial C_\delta(r,s)}{\partial r}\right|_{s=r}
     + \left.\frac{\partial C_\delta(r,s)}{\partial s}\right|_{s=r} .
    \label{eq:diagonal-chain}
\end{align}
The two partial derivatives in~\eqref{eq:diagonal-chain} coincide. Differentiating the symmetry relation $C_\delta(r,s) = C_\delta(s,r)$ of part~(i) with respect to its second argument, and then setting $s = r$:
\begin{align}
    \left.\frac{\partial C_\delta(r,s)}{\partial s}\right|_{s=r}
    &= \left.\frac{\partial}{\partial s}\,C_\delta(s,r)\right|_{s=r}
     = \left.\frac{\partial C_\delta(a,b)}{\partial a}\right|_{a=b=r}
     = \left.\frac{\partial C_\delta(r,s)}{\partial r}\right|_{s=r},
    \label{eq:diagonal-partials-equal}
\end{align}
in which the middle equality renames the differentiation variable as the first slot of the kernel. Substituting~\eqref{eq:diagonal-partials-equal} into~\eqref{eq:diagonal-chain}:
\begin{equation}
    \frac{d V_\delta}{dr}(r)
    = 2\left.\frac{\partial C_\delta(r,s)}{\partial s}\right|_{s=r} .
\end{equation}
Dividing by two and solving for the diagonal derivative, therefore,
\begin{equation}
    \left.\frac{\partial C_\delta(r,s)}{\partial s}\right|_{s=r} = \frac{1}{2}\,\frac{d V_\delta(r)}{dr},
\end{equation}
which is part~(ii).
\end{proofpart}
\end{proof}
\begin{lemma}[Closed form of the radial covariance]
\label{lemma:radial-closed-form}
Let $g_\delta \coloneqq \psi_\delta * \psi_\delta$ be the auto-correlation of the mollifier of Definition~\ref{def:mollifier}, so that $g_\delta$ is radial, non-negative definite, and normalised by $\int_{\mathbb{R}^3} g_\delta = 1$. Then for all $r, s \in (0, L]$,
\begin{equation}
    C_\delta(r,s) \;=\; \frac{\pi}{r s}\int_{|r-s|}^{\,r+s} g_\delta(u)\,\bigl(r^2+s^2-u^2\bigr)\,u \; du .
    \label{eq:radial-closed-form}
\end{equation}
\end{lemma}
\begin{proof}
\begin{proofpart}[Reduction of the angular integral to a chord integral]
By Eq.~\eqref{eq:kernel-surface} the kernel is the double surface integral of the ambient covariance, which for a radial $g_\delta$ depends on the two points only through the chord length $u$ separating them. Writing $\gamma$ for the angle subtended at the origin, $u^2 = r^2+s^2-2rs\cos\gamma$, and the double surface integral reduces to
\begin{equation}
    C_\delta(r,s) = 2\pi r s \int_0^\pi g_\delta(u)\,\cos\gamma\,\sin\gamma\;d\gamma .
    \label{eq:kernel-angular}
\end{equation}
Substitute $x = \cos\gamma$, so that $dx = -\sin\gamma\,d\gamma$, and the limits $\gamma = 0,\pi$ map to $x = 1,-1$. Differentiating $u^2 = r^2+s^2-2rsx$ at fixed $r,s$ gives $2u\,du = -2rs\,dx$, hence $dx = -u\,du/(rs)$, while inverting the same relation gives $x = (r^2+s^2-u^2)/(2rs)$. The endpoints transform as $x = 1 \mapsto u = |r-s|$ and $x = -1 \mapsto u = r+s$. Carrying out the substitution:
\begin{align}
    \int_0^\pi g_\delta(u)\cos\gamma\,\sin\gamma\;d\gamma
    &= \int_{-1}^{1} g_\delta(u)\,x\,dx \nonumber \\
    &= \int_{r+s}^{|r-s|} g_\delta(u)\,\frac{r^2+s^2-u^2}{2rs}\left(-\frac{u\,du}{rs}\right) \nonumber \\
    &= \frac{1}{2r^2s^2}\int_{|r-s|}^{r+s} g_\delta(u)\bigl(r^2+s^2-u^2\bigr)\,u\,du ,
    \label{eq:kernel-chord}
\end{align}
where the last line reverses the orientation of the interval, which cancels the minus sign.
\end{proofpart}
\begin{proofpart}[Restoration of the prefactor]
Multiplying~\eqref{eq:kernel-chord} by the prefactor $2\pi rs$ of~\eqref{eq:kernel-angular}:
\begin{align}
    C_\delta(r,s)
    &= 2\pi r s \cdot \frac{1}{2r^2s^2}\int_{|r-s|}^{r+s} g_\delta(u)\bigl(r^2+s^2-u^2\bigr)\,u\,du \nonumber \\
    &= \frac{\pi}{rs}\int_{|r-s|}^{r+s} g_\delta(u)\bigl(r^2+s^2-u^2\bigr)\,u\,du ,
\end{align}
which is~\eqref{eq:radial-closed-form}. The reduction is exact, since no approximation on $g_\delta$, $r$ or $s$ has been used.
\end{proofpart}
\end{proof}
Three scalar functions, all built from the covariance kernel, enter every subsequent estimate, and we define them here.

\begin{definition}[Kernel invariants]
\label{def:kernel-invariants}
Let $C_\delta(r,s)$ be the covariance kernel of Lemma~\ref{lemma:radial-covariance}. The \emph{variance profile}, its radial derivative, and the \emph{diagonal mixed curvature} are defined respectively by
\begin{equation}
    V_\delta(r) \coloneqq C_\delta(r,r) ,
    \qquad
    V_\delta'(r) \coloneqq \frac{d V_\delta}{dr}(r) ,
    \qquad
    \kappa_\delta(r) \coloneqq -\left.\frac{\partial^2 C_\delta(r,s)}{\partial r\,\partial s}\right|_{s=r} .
    \label{eq:kappa-convention}
\end{equation}
\end{definition}

\begin{proposition}[Small-$\delta$ behaviour of the kernel invariants]
\label{prop:kernel-asymptotics}
Write $g_\delta(u) = \delta^{-3}g(u/\delta)$ and set the moments
\begin{equation}
    A \coloneqq \int_0^\infty g(v)\,v\,dv , \qquad B \coloneqq \int_0^\infty g(v)\,v^3\,dv .
\end{equation}
Then, for fixed $r \in (0,L]$ and $\delta \ll r$,
\begin{align}
    V_\delta(r) &= \frac{2\pi A}{\delta} - \frac{\pi B}{r^2}\,\delta + \mathcal{O}(\delta^3),
    \label{eq:V-asymptotic}\\[2pt]
    V_\delta'(r) &= \frac{2\pi B}{r^{3}}\,\delta + \mathcal{O}(\delta^3),
    \label{eq:Vprime-asymptotic}\\[2pt]
    \kappa_\delta(r) &= -\,2\pi\, g_\delta(0) + \mathcal{O}(\delta^{-1})
    \;=\; -\,\frac{2\pi\, g(0)}{\delta^{3}} + \mathcal{O}(\delta^{-1}).
    \label{eq:kappa-asymptotic}
\end{align}
\end{proposition}
\begin{proof}
\begin{proofpart}[The diagonal value and its first derivative]
Setting $s = r$ in the chord representation of Lemma~\ref{lemma:radial-closed-form}, Eq.~\eqref{eq:radial-closed-form}, the lower limit becomes $|r-s| = 0$ and the upper limit becomes $2r$, while $r^2+s^2 = 2r^2$:
\begin{align}
    V_\delta(r) = C_\delta(r,r)
    &= \frac{\pi}{r^2}\int_0^{2r} g_\delta(u)\bigl(2r^2-u^2\bigr)\,u\,du \nonumber \\
    &= 2\pi\int_0^{2r} g_\delta(u)\,u\,du \;-\; \frac{\pi}{r^2}\int_0^{2r} g_\delta(u)\,u^3\,du .
    \label{eq:V-exact}
\end{align}
Since $\psi_\delta$ vanishes outside $B_\delta(0)$, its auto-correlation $g_\delta$ vanishes for $u \ge 2\delta$, so that for $r \ge \delta$ extending the upper limit to infinity leaves each integral unchanged. Rescaling $u = \delta v$ in the extended integrals, and using $g_\delta(u) = \delta^{-3}g(u/\delta)$ so that $g_\delta(\delta v)\,\delta\,dv = \delta^{-2}g(v)\,dv$:
\begin{align}
    \int_0^\infty g_\delta(u)\,u\,du &= \delta^{-1}\int_0^\infty g(v)\,v\,dv = \frac{A}{\delta}, \nonumber \\
    \int_0^\infty g_\delta(u)\,u^3\,du &= \delta\int_0^\infty g(v)\,v^3\,dv = B\,\delta .
    \label{eq:moments}
\end{align}
Substituting~\eqref{eq:moments} into~\eqref{eq:V-exact} gives~\eqref{eq:V-asymptotic}, and for $r \ge \delta$ the remainder there vanishes. Differentiating that expression in $r$, in which only the factor $r^{-2}$ carries $r$ dependence:
\begin{equation}
    V_\delta'(r) = \frac{d}{dr}\left[\frac{2\pi A}{\delta} - \frac{\pi B}{r^2}\,\delta\right] = \frac{2\pi B}{r^{3}}\,\delta ,
\end{equation}
which is~\eqref{eq:Vprime-asymptotic}.
\end{proofpart}
\begin{proofpart}[Differentiation of the chord integral off the diagonal]
Write $C_\delta = (\pi/rs)\,I$ with $I(r,s) = \int_{a}^{b} g_\delta(u)(r^2+s^2-u^2)u\,du$, where $a = |r-s|$ and $b = r+s$. Take $r > s$, so that $a = r-s$ with $\partial_s a = -1$ and $\partial_r a = +1$, while $\partial_s b = \partial_r b = +1$. At the two limits the polynomial factor takes the values $r^2+s^2-b^2 = -2rs$ and $r^2+s^2-a^2 = +2rs$. Applying the Leibniz rule for a parameter-dependent integral with moving limits, and writing $J \coloneqq \int_a^b g_\delta(u)\,u\,du$:
\begin{align}
    \partial_s I &= \int_a^b g_\delta(u)\,(2s)\,u\,du \;+\; g_\delta(b)\bigl(r^2+s^2-b^2\bigr)b\,\partial_s b \;-\; g_\delta(a)\bigl(r^2+s^2-a^2\bigr)a\,\partial_s a \nonumber \\
    &= 2sJ \;-\; 2rs\,b\,g_\delta(b) \;+\; 2rs\,a\,g_\delta(a) .
    \label{eq:dsI}
\end{align}
Differentiating~\eqref{eq:dsI} with respect to $r$, using $\partial_r J = b\,g_\delta(b) - a\,g_\delta(a)$ and $\partial_r[a g_\delta(a)] = g_\delta(a) + a g_\delta'(a)$ together with the corresponding identity at $b$:
\begin{align}
    \partial_r\partial_s I
    &= 2s\bigl[b g_\delta(b) - a g_\delta(a)\bigr]
     - 2s\,b\,g_\delta(b) - 2rs\bigl[g_\delta(b) + b g_\delta'(b)\bigr] \nonumber \\
    &\qquad + 2s\,a\,g_\delta(a) + 2rs\bigl[g_\delta(a) + a g_\delta'(a)\bigr] \nonumber \\
    &= 2rs\bigl[g_\delta(a) - g_\delta(b) + a g_\delta'(a) - b g_\delta'(b)\bigr] ,
\end{align}
in which every term carrying the factor $s$ alone has cancelled. Setting $s = r$, so that $a = 0$ and $b = 2r$, and noting that the term $a\,g_\delta'(a)$ vanishes at $a = 0$:
\begin{equation}
    \partial_r\partial_s I\big|_{s=r} = 2r^2\bigl[\,g_\delta(0) - g_\delta(2r) - 2r\,g_\delta'(2r)\,\bigr].
    \label{eq:drdsI}
\end{equation}
\end{proofpart}
\begin{proofpart}[The mixed curvature]
Applying the product rule to $C_\delta = (\pi/rs)I$, and using $\partial_r(1/rs) = -1/(r^2s)$, $\partial_s(1/rs) = -1/(rs^2)$ and $\partial_r\partial_s(1/rs) = 1/(r^2s^2)$, then evaluating at $s = r$ where $\partial_r I = \partial_s I$ by the symmetry of $C_\delta$:
\begin{align}
    \partial_r\partial_s C_\delta\big|_{s=r}
    &= \pi\left[\frac{I}{r^4} - \frac{\partial_s I}{r^3} - \frac{\partial_r I}{r^3} + \frac{\partial_r\partial_s I}{r^2}\right]_{s=r} \nonumber \\
    &= \pi\left[\frac{I}{r^4} - \frac{2}{r^3}\bigl(2rJ - 4r^3 g_\delta(2r)\bigr) + 2\bigl(g_\delta(0) - g_\delta(2r) - 2r g_\delta'(2r)\bigr)\right] .
    \label{eq:drdsC}
\end{align}
For $r \ge \delta$ the point $2r$ lies where $g_\delta$ vanishes identically, so that $g_\delta(2r) = g_\delta'(2r) = 0$. At $s = r$, Eqs.~\eqref{eq:V-exact} and~\eqref{eq:moments} give $I = 2r^2A/\delta - B\delta$ and $J = A/\delta$, and the first two terms of~\eqref{eq:drdsC} combine as
\begin{equation}
    \frac{I}{r^4} - \frac{4J}{r^2} = \frac{2A}{r^2\delta} - \frac{B\delta}{r^4} - \frac{4A}{r^2\delta} = -\frac{2A}{r^2\delta} - \frac{B\delta}{r^4},
\end{equation}
which is $\mathcal{O}(\delta^{-1})$ and therefore subordinate to the term $2\pi g_\delta(0) = 2\pi g(0)\delta^{-3}$. Collecting, and using the definition of $\kappa_\delta$ given in Definition~\ref{def:kernel-invariants}:
\begin{equation}
    \kappa_\delta(r) = -2\pi g_\delta(0) + \frac{2\pi A}{r^2\delta} + \frac{\pi B\delta}{r^4}
    = -\frac{2\pi g(0)}{\delta^3} + \mathcal{O}(\delta^{-1}),
\end{equation}
which is~\eqref{eq:kappa-asymptotic}.
\end{proofpart}
\end{proof}
\begin{remark}[Uniformity of the expansion near the origin]
\label{rem:boundary-layer}
Note that the expansions of Proposition~\ref{prop:kernel-asymptotics} are stated for fixed $r \in (0,L]$ and $\delta \ll r$, and they fail in the boundary layer $r = \mathcal{O}(\delta)$. Since we integrate the invariants over the whole of $[0,L]$, that layer needs its own estimate.

There we use the exact kernel, which degenerates in the layer. Let us write $p = r^2+s^2$ and $q = 2rs$. The limits of Eq.~\eqref{eq:radial-closed-form} then satisfy $b^2 - a^2 = 2q$, with $a = |r-s|$ and $b = r+s$. We expand $g_\delta$ about the origin. The term of first order is absent, since the auto-correlation of a radial mollifier is an even function of $u$ and hence $g_\delta'(0) = 0$, and the term carried by $g_\delta(0)$ vanishes identically, since
\begin{equation}
    \int_a^b (r^2+s^2-u^2)\,u\,du = \tfrac{1}{2}p\,(b^2-a^2) - \tfrac{1}{4}(b^4-a^4) = pq - pq = 0 ,
    \label{eq:chord-cancellation}
\end{equation}
and the first surviving contribution comes from the curvature of the mollifier. Using $b^4-a^4 = 4pq$ and $b^6-a^6 = 2q(3p^2+q^2)$ in the same way gives $\tfrac12 g_\delta''(0)\bigl[\tfrac14 p(b^4-a^4) - \tfrac16(b^6-a^6)\bigr] = -\tfrac16 g_\delta''(0)\,q^3$, and therefore
\begin{equation}
    C_\delta(r,s) = -\tfrac{4}{3}\pi\, g_\delta''(0)\, r^{2}s^{2} + \mathcal{O}\bigl(\delta^{-7}r^{2}s^{2}(r^{2}+s^{2})\bigr) , \qquad r,s \ll \delta .
    \label{eq:kernel-origin}
\end{equation}
Differentiating Eq.~\eqref{eq:kernel-origin} yields the two invariants required on the layer,
\begin{equation}
    V_\delta(r) = -\tfrac{4}{3}\pi\, g_\delta''(0)\, r^{4} ,
    \qquad
    \kappa_\delta(r) = -\,\partial_r\partial_s C_\delta\big|_{s=r} = \tfrac{16}{3}\pi\, g_\delta''(0)\, r^{2} ,
    \label{eq:variance-origin}
\end{equation}
the first positive and the second negative, since $g_\delta''(0) < 0$ at the maximum of $g_\delta$, consistently with $\kappa_\delta \le 0$. Hence the variance profile vanishes quartically at the centre, instead of growing like $\delta^{-1}$, and the mixed curvature vanishes quadratically, instead of growing like $\delta^{-3}$. Since $g_\delta''(0) = \mathcal{O}(\delta^{-5})$, Eq.~\eqref{eq:variance-origin} gives $|\kappa_\delta(r)| = \mathcal{O}(\delta^{-5}r^2)$ on the layer, which matches the bulk value $\mathcal{O}(\delta^{-3})$ at $r = \mathcal{O}(\delta)$. In other words, the two regimes join without a jump.

At this point the layer contribution can be bounded, instead of estimated. In the case of a local excess, the statistical anomaly of Eq.~\eqref{eq:upsilon-final} meets the kernel through the diagonal $\Phi_{\mathrm{d}}$ of the mixed derivative of the projected covariance, and the dominant part of this diagonal is $\sigma^2\lvert\kappa_\delta\rvert$, with $\sigma^2 = 2D_0\,\partial_r M = \mathcal{O}(r^2)$ (Proposition~\ref{prop:validity-condition}). We split the integration range as $[0,L] = [0,c\delta] \cup [c\delta,L]$ for a fixed $c > 0$. Making use of $\varpi\,\mathbf{e}_k^m = \mathcal{O}(r)$ together with $\partial_r\mathbf{e}_k^m = \mathcal{O}(r^{2})$ (Lemma~\ref{lemma:boundary-vanishing}), the term that $\Phi_{\mathrm{d}}$ produces, divided by $r^4$, obeys
\begin{equation}
    \int_0^{c\delta} \bigl| \partial_r\mathbf{e}_k^m\,\sigma^2\kappa_\delta\,r^{-4} \bigr| \, dr
    \;\le\; K \delta^{-5}\!\int_0^{c\delta}\! r^{2}\,dr
    \;=\; \tfrac{1}{3}Kc^{3}\,\delta^{-2} ,
    \label{eq:layer-bound}
\end{equation}
with $K$ depending on the profile but not on $\delta$, and the term divided by $r^3$, which carries the factor $\varpi\,\mathbf{e}_k^m$, obeys the same bound. On the layer $V_\delta = \mathcal{O}(\delta^{-5}r^{4})$ and $V_\delta' = \mathcal{O}(\delta^{-5}r^{3})$ by Eq.~\eqref{eq:variance-origin}, and hence the parts of $\Phi_{\mathrm{d}}$ carrying them are of the same order $\delta^{-5}r^{4}$ as $\sigma^2\kappa_\delta$, and they obey the same bound. Note that the singular geometric prefactors remove the separation that the invariants display in the bulk. The layer contributes $\mathcal{O}(\delta^{-2})$, against a bulk contribution of order $\delta^{-3}$ (Proposition~\ref{prop:validity-condition}), and since this is one power of the correlation length smaller, we can insert the expansions of Proposition~\ref{prop:kernel-asymptotics} under the integral sign at the retained order.
\end{remark}
\begin{remark}[Consequences of the kernel asymptotics]
\label{rem:kernel-quantitative}
Let us now discuss three consequences of Proposition~\ref{prop:kernel-asymptotics}, which bear directly on the analysis of the subsequent sections. First, $V_\delta' \not\equiv 0$ whenever the third moment $B$ is non-zero, and hence the geometric anomaly does not vanish identically. The objection that a translation-invariant kernel would force $V_\delta$ constant does not apply here, since the radial kernel is not translation invariant. Indeed, Eq.~\eqref{eq:Vprime-asymptotic} exhibits the induced $r$-dependence as $V_\delta' \propto \delta\,r^{-3}$. The anomaly is largest at distances of order $\delta$ from the origin, where the spherical geometry varies most rapidly, and below them it vanishes by Remark~\ref{rem:boundary-layer}. Second, the two diagonal invariants scale oppositely, since $V_\delta'$ vanishes with the correlation length while $V_\delta$ grows. The regime of interest is thus finite $\delta$, of the order of the microscopic correlation length, not the white-noise limit. Third, and decisively for Section~\ref{sec:commutativity}, the diagonal mixed curvature satisfies $\kappa_\delta = -2\pi g_\delta(0) + \mathcal{O}(\delta^{-1})$, and so it is of order $\delta^{-3}$, two orders more singular than $V_\delta$ and four more than $V_\delta'$, which identifies $\kappa_\delta$ as the source of the obstruction that Proposition~\ref{prop:validity-condition} obtains there for a fluid with a local excess free energy. Note that nothing diverges at fixed $\delta > 0$, since the mollification (Definition~\ref{def:mollifier}) renders the covariance smooth and trace class, and Eq.~\eqref{eq:kappa-asymptotic} records only how fast the limit $\delta \to 0$ would be approached, were it taken, and we never take it.
\end{remark}
\subsubsection{The transformed stochastic dynamics}
\begin{lemma}[Duality of chemical potentials]
\label{lemma:mu-duality}
Let $\mathcal{F}[\rho]$ be the free energy functional in the density representation, and let $\mathcal{F}_m[m]$ be its equivalent in the cumulative mass representation, such that $\mathcal{F}_m[\mathcal{T}\rho] = \mathcal{F}[\rho]$. The radial gradient of the chemical potential $\mu = \frac{\delta \mathcal{F}}{\delta \rho}$ is related to the functional derivative with respect to mass by the negative identity:
\begin{equation}
    \frac{\partial \mu(r)}{\partial r} = - \frac{\delta \mathcal{F}_m}{\delta m(r)}.
\end{equation}
\end{lemma}
\begin{proof}
Consider a variation of the mass profile $\delta m(r)$ compatible with the boundary conditions (that is, $\delta m(0) = 0$ by construction and $\delta m(L) = 0$ by mass conservation). This induces a variation in the density $\delta \rho$.
The variation of the free energy can be written in two equivalent ways.
\begin{proofpart}[The variation in the density representation]
First, in the density representation:
\begin{equation}
    \delta \mathcal{F} = \int_V \frac{\delta \mathcal{F}}{\delta \rho} \delta \rho(\mathbf{r}) \, d\mathbf{r} = \int_0^L \mu(r) \delta \rho(r) \, 4\pi r^2 \, dr.
\end{equation}
From the definition $m(r) = \int_0^r \rho(s) 4\pi s^2 ds$, we have the differential relation $\partial_r m = 4\pi r^2 \rho$. Taking variations, $\partial_r (\delta m) = 4\pi r^2 \delta \rho$. Substituting this into the integral:
\begin{equation}
    \delta \mathcal{F} = \int_0^L \mu(r) \frac{\partial (\delta m)}{\partial r} \, dr = \left[ \mu(r) \delta m(r) \right]_0^L - \int_0^L \frac{\partial \mu}{\partial r} \delta m(r) \, dr.
\end{equation}
The boundary terms vanish because the total mass is fixed ($\delta m(L)=0$) and the mass at the origin is zero ($\delta m(0)=0$).
\end{proofpart}
\begin{proofpart}[The variation in the mass representation]
Second, by definition of the functional derivative with respect to $m$:
\begin{equation}
    \delta \mathcal{F} = \int_0^L \frac{\delta \mathcal{F}_m}{\delta m(r)} \delta m(r) \, dr.
\end{equation}
Comparing the two integral expressions for arbitrary variations $\delta m(r)$ implies pointwise equality of the kernels $\frac{\delta \mathcal{F}_m}{\delta m(r)} = - \frac{\partial \mu}{\partial r}.$
\end{proofpart}
\end{proof}
To pass from the density to the cumulative mass, mesoscopic nucleation theory integrates the fluctuating equation over a spherical volume and then seeks a spherically symmetric solution. As Section~II.B of Lutsko~\cite{lutsko2012mean} observes, a spherically symmetric profile cannot simply be assumed in the fluctuating equation, since the noise breaks the symmetry. In that theory the symmetry is therefore imposed after the integration, and the error it introduces is neither stated nor estimated. Since the integration over the ball holds for every solution, we integrate first without the symmetry and measure the departure from it on each sphere, keeping apart the terms that a spherically symmetric density retains and those produced by the departure. We measure it by the following quantities.
\begin{definition}[Spherical means and asymmetries]
\label{def:spherical-asymmetry}
For a function $f$ upon $V$, we define its mean over the sphere of radius $r$ and its radial derivative by
\begin{equation}
    \langle f \rangle_r \coloneqq \frac{1}{4\pi r^2} \oint_{\partial B_r} f \, dS ,
    \qquad
    \partial_r f \coloneqq \nabla f \cdot \hat{\mathbf{x}} , \quad \hat{\mathbf{x}} \coloneqq \mathbf{x}/|\mathbf{x}| .
\end{equation}
For a density $\rho$ with chemical potential $\mu$, we define the asymmetries of $\rho$ and of $\partial_r\mu$ upon that sphere by
\begin{equation}
    a_\rho(r) \coloneqq \frac{\langle (\rho - \langle \rho \rangle_r)^2 \rangle_r^{1/2}}{\langle \rho \rangle_r} ,
    \qquad
    a_\mu(r) \coloneqq \frac{\bigl\langle (\partial_r\mu - \langle \partial_r\mu \rangle_r)^2 \bigr\rangle_r^{1/2}}{|\langle \partial_r\mu \rangle_r|} ,
    \label{eq:asymmetries}
\end{equation}
the second wherever $\langle \partial_r\mu \rangle_r \neq 0$.
\end{definition}
\begin{proposition}[The error of the spherical approximation]
\label{prop:shell-integration}
Let $\rho_t$ satisfy Eq.~\eqref{eq:full_spde_ito}, with $\rho_t \ge \alpha > 0$, with $\rho_t$ continuously differentiable and $\mu[\rho_t]$ twice continuously differentiable upon $V$, and let $m_t = \mathcal{T}[\rho_t]$. Then $\partial_r m_t = 4\pi r^2 \langle \rho_t \rangle_r$, and
\begin{equation}
    dm_t(r) = \Bigl[ D_0\beta \, \frac{\partial m_t}{\partial r} \, \frac{\partial \langle \mu \rangle_r}{\partial r} + \mathcal{E}_{\mathcal{A}}(r) \Bigr] dt
      + \varepsilon \sqrt{\frac{2D_0 \, \partial_r m_t}{4\pi r^2}} \oint_{\partial B_r} d\zeta^\delta_t \cdot \hat{\mathbf{n}} \, dS
      + \varepsilon \, d\mathcal{E}_{\mathcal{B}}(r)_t ,
    \label{eq:closure-error}
\end{equation}
with the remainders
\begin{align}
    \mathcal{E}_{\mathcal{A}}(r) &\coloneqq 4\pi r^2 D_0\beta \, \bigl\langle (\rho_t - \langle \rho_t \rangle_r)(\partial_r\mu - \langle \partial_r\mu \rangle_r) \bigr\rangle_r ,
    \nonumber \\
    \mathcal{E}_{\mathcal{B}}(r)_t &\coloneqq \int_0^t \oint_{\partial B_r} \bigl( \sqrt{2D_0\rho_s} - \sqrt{2D_0\langle \rho_s \rangle_r} \bigr) \, d\zeta^\delta_s \cdot \hat{\mathbf{n}} \, dS .
    \label{eq:closure-remainders}
\end{align}
The drift remainder satisfies
\begin{equation}
    |\mathcal{E}_{\mathcal{A}}(r)| \le a_\rho \, a_\mu \, \Bigl| D_0\beta \, \frac{\partial m_t}{\partial r} \, \frac{\partial \langle \mu \rangle_r}{\partial r} \Bigr| .
\end{equation}
The quadratic variation of the noise remainder is $[\mathcal{E}_{\mathcal{B}}(r)]_t = \int_0^t q_s(r) \, ds$, with
\begin{equation}
    q_s(r) \le \Gamma_\delta(r) \, a_\rho^2 \; 2D_0 \, \langle \rho_s \rangle_r \, 4\pi r^2 V_\delta(r) ,
    \qquad
    \Gamma_\delta(r) \coloneqq \frac{2\pi \int_0^{2r} |g_\delta(u)| \, u \, du}{V_\delta(r)} ,
\end{equation}
where $2D_0 \langle \rho_s \rangle_r \, 4\pi r^2 V_\delta(r)$ is the density of the quadratic variation of the retained noise term. Both remainders vanish when $\rho_t$ is spherically symmetric.
\end{proposition}
\begin{proof}
Throughout the proof we write $\rho$ and $\mu$ for $\rho_t$ and $\mu[\rho_t]$.
\begin{proofpart}[The integral over the ball]
Integrating both sides of Eq.~\eqref{eq:full_spde_ito} over $B_r$ and applying the divergence theorem to the two divergences, both of continuously differentiable fields, the outward normal of $\partial B_r$ being $\hat{\mathbf{n}} = \hat{\mathbf{x}}$, gives
\begin{equation}
    dm_t(r) = \Bigl[ \oint_{\partial B_r} D_0 \rho \, \nabla(\beta\mu) \cdot \hat{\mathbf{n}} \, dS \Bigr] dt + \varepsilon \oint_{\partial B_r} \sqrt{2D_0\rho} \; d\zeta^\delta_t \cdot \hat{\mathbf{n}} \, dS .
    \label{eq:shell-integrated}
\end{equation}
The derivative of the integral of a continuous and summable function upon $\mathbb{R}^3$ over $B_r$ with respect to $r$ is its integral over $\partial B_r$, which is Theorem~4(ii) of Section~C.3 of Evans~\cite{evans2010pde}. Applying it, for $r < r' < r'' < L$, to the product of $\rho_t$ with a continuous function equal to one upon $B_{r'}$ and to zero outside $B_{r''}$, which leaves the integrals over $B_s$ unchanged for $s < r'$, gives
\begin{equation}
    \frac{\partial m_t}{\partial r}(r) = \oint_{\partial B_r} \rho_t \, dS = 4\pi r^2 \, \langle \rho_t \rangle_r .
\end{equation}
\end{proofpart}
\begin{proofpart}[The drift]
Upon $\partial B_r$ we have $\nabla(\beta\mu) \cdot \hat{\mathbf{n}} = \beta \, \partial_r\mu$. Writing $\rho = \langle \rho \rangle_r + (\rho - \langle \rho \rangle_r)$ and $\partial_r\mu = \langle \partial_r\mu \rangle_r + (\partial_r\mu - \langle \partial_r\mu \rangle_r)$, and using that $\langle \rho \rangle_r$ and $\langle \partial_r\mu \rangle_r$ are constant upon the sphere while the two fluctuations have mean zero there, the two cross terms vanish and
\begin{align}
    \oint_{\partial B_r} D_0 \rho \, \nabla(\beta\mu) \cdot \hat{\mathbf{n}} \, dS
      &= 4\pi r^2 D_0\beta \, \langle \rho \, \partial_r\mu \rangle_r \nonumber \\
      &= 4\pi r^2 D_0\beta \, \Bigl[ \langle \rho \rangle_r \, \langle \partial_r\mu \rangle_r + \bigl\langle (\rho - \langle \rho \rangle_r)(\partial_r\mu - \langle \partial_r\mu \rangle_r) \bigr\rangle_r \Bigr] \nonumber \\
      &= 4\pi r^2 \langle \rho \rangle_r \, D_0\beta \, \langle \partial_r\mu \rangle_r + \mathcal{E}_{\mathcal{A}}(r) .
\end{align}
The mean over the sphere commutes with the radial derivative. Writing the sphere as $\{r\omega : \omega \in S^2\}$, with $d\omega$ the surface measure upon the unit sphere, and differentiating under the integral sign, $\mu$ being continuously differentiable,
\begin{equation}
    \frac{\partial \langle \mu \rangle_r}{\partial r}
      = \frac{\partial}{\partial r} \frac{1}{4\pi} \int_{S^2} \mu(r\omega) \, d\omega
      = \frac{1}{4\pi} \int_{S^2} \nabla\mu(r\omega) \cdot \omega \, d\omega
      = \langle \partial_r\mu \rangle_r .
\end{equation}
With $4\pi r^2 \langle \rho \rangle_r = \partial_r m_t$ from the first part, the drift is the first term of~\eqref{eq:closure-error}.
\end{proofpart}
\begin{proofpart}[The noise]
Adding and subtracting $\sqrt{2D_0\langle \rho \rangle_r}$ under the integral, and using that $\langle \rho \rangle_r = \partial_r m_t/(4\pi r^2)$ is constant upon the sphere,
\begin{equation}
    \oint_{\partial B_r} \sqrt{2D_0\rho} \; d\zeta^\delta_t \cdot \hat{\mathbf{n}} \, dS
      = \sqrt{\frac{2D_0 \, \partial_r m_t}{4\pi r^2}} \oint_{\partial B_r} d\zeta^\delta_t \cdot \hat{\mathbf{n}} \, dS + d\mathcal{E}_{\mathcal{B}}(r)_t ,
\end{equation}
which is the noise term of~\eqref{eq:closure-error}.
\end{proofpart}
\begin{proofpart}[The drift remainder]
The Cauchy--Schwarz inequality for the mean over the sphere, and the definitions~\eqref{eq:asymmetries}, give
\begin{equation}
    |\mathcal{E}_{\mathcal{A}}(r)| \le 4\pi r^2 D_0\beta \, \langle (\rho - \langle \rho \rangle_r)^2 \rangle_r^{1/2} \, \bigl\langle (\partial_r\mu - \langle \partial_r\mu \rangle_r)^2 \bigr\rangle_r^{1/2}
      = 4\pi r^2 D_0\beta \, \langle \rho \rangle_r \, a_\rho \, |\langle \partial_r\mu \rangle_r| \, a_\mu ,
\end{equation}
which is the bound of the statement upon the drift remainder, by the second part.
\end{proofpart}
\begin{proofpart}[The noise remainder]
Write $b \coloneqq \sqrt{2D_0\rho} - \sqrt{2D_0\langle \rho \rangle_r}$ upon $\partial B_r$, at time $s$. The remainder $\mathcal{E}_{\mathcal{B}}(r)$ is a stochastic integral in time, and its quadratic variation is the time integral of the variance of its integrand. The components of the noise field have the covariation $g_\delta(|\mathbf{x}-\mathbf{y}|)\,ds$ at two points, by Lemma~\ref{lemma:mollified-noise}, so $[\mathcal{E}_{\mathcal{B}}(r)]_t = \int_0^t q_s(r) \, ds$ with
\begin{equation}
    q_s(r) = \oint_{\partial B_r} \oint_{\partial B_r} b(\mathbf{x}) \, b(\mathbf{y}) \, g_\delta(|\mathbf{x}-\mathbf{y}|) \, \hat{\mathbf{n}}(\mathbf{x}) \cdot \hat{\mathbf{n}}(\mathbf{y}) \, dS_{\mathbf{x}} \, dS_{\mathbf{y}} .
\end{equation}
Using $|\hat{\mathbf{n}}(\mathbf{x}) \cdot \hat{\mathbf{n}}(\mathbf{y})| \le 1$, the inequality $|b(\mathbf{x})b(\mathbf{y})| \le \frac12 b(\mathbf{x})^2 + \frac12 b(\mathbf{y})^2$, and the symmetry of the kernel in $\mathbf{x}$ and $\mathbf{y}$,
\begin{equation}
    q_s(r)
      \le \oint_{\partial B_r} b(\mathbf{x})^2 \Bigl( \oint_{\partial B_r} |g_\delta(|\mathbf{x}-\mathbf{y}|)| \, dS_{\mathbf{y}} \Bigr) dS_{\mathbf{x}} .
\end{equation}
For $\mathbf{x}$ upon the sphere, take spherical coordinates about $\mathbf{x}/r$, so that $|\mathbf{x}-\mathbf{y}|^2 = 2r^2(1-\cos\gamma)$ and $u \, du = r^2 \sin\gamma \, d\gamma$ for $u = |\mathbf{x}-\mathbf{y}|$. The inner integral does not depend upon $\mathbf{x}$,
\begin{equation}
    \oint_{\partial B_r} |g_\delta(|\mathbf{x}-\mathbf{y}|)| \, dS_{\mathbf{y}}
      = 2\pi r^2 \int_0^\pi |g_\delta(u)| \, \sin\gamma \, d\gamma
      = 2\pi \int_0^{2r} |g_\delta(u)| \, u \, du
      = \Gamma_\delta(r) \, V_\delta(r) .
\end{equation}
Finally $(\sqrt{\rho} - \sqrt{\langle \rho \rangle_r})^2 = (\rho - \langle \rho \rangle_r)^2/(\sqrt{\rho} + \sqrt{\langle \rho \rangle_r})^2 \le (\rho - \langle \rho \rangle_r)^2/\langle \rho \rangle_r$, so that
\begin{equation}
    \oint_{\partial B_r} b^2 \, dS \le \frac{2D_0}{\langle \rho \rangle_r} \oint_{\partial B_r} (\rho - \langle \rho \rangle_r)^2 \, dS
      = 2D_0 \, 4\pi r^2 \, \langle \rho \rangle_r \, a_\rho^2 .
\end{equation}
Combining the three displays gives the bound of the statement upon $q_s(r)$. The density of the quadratic variation of the retained noise term is the square $2D_0\,\partial_r m_t/(4\pi r^2) = 2D_0\langle \rho \rangle_r$ of its coefficient times that of $\oint_{\partial B_r} d\zeta^\delta_t \cdot \hat{\mathbf{n}} \, dS$, which is $4\pi r^2 V_\delta(r)$ by Definition~\ref{def:radial-noise-norm-ref} and Lemma~\ref{lemma:radial-covariance}, and the product is $2D_0\langle \rho \rangle_r \, 4\pi r^2 V_\delta(r)$.
\end{proofpart}
\begin{proofpart}[The symmetric case]
If $\rho_t$ is spherically symmetric, it is constant upon every sphere, so $\rho_t = \langle \rho_t \rangle_r$ there and $a_\rho = 0$. Its chemical potential is then spherically symmetric, by Lemma~\ref{lemma:radial-chemical-potential}, so $\partial_r\mu$ is constant upon every sphere, $a_\mu = 0$, and $\partial_r\langle \mu \rangle_r = \partial_r\mu$. The integrands of $\mathcal{E}_{\mathcal{A}}$ and $\mathcal{E}_{\mathcal{B}}$ vanish.
\end{proofpart}
\end{proof}
In mesoscopic nucleation theory, nucleation is the growth or the decay of a localised excess of mass over the background density, and this excess is taken to be spherically symmetric, as the domain and the background density are. The remainders of Proposition~\ref{prop:shell-integration} are then neglected. We state this closure as an assumption.
\begin{assumption}[Spherical-mean closure of the mass dynamics]
\label{ass:spherical-closure}
The density field remains close to spherically symmetric over the time interval of interest. That is, the asymmetries $a_\rho$ and $a_\mu$ of Definition~\ref{def:spherical-asymmetry} are small for $t < \tau_\alpha$, with the stopping time of Assumption~\ref{ass:strict-positivity}, and we close the equation of Proposition~\ref{prop:shell-integration} by replacing the density upon each sphere by its mean over that sphere. This neglects the remainders $\mathcal{E}_{\mathcal{A}}$ and $\mathcal{E}_{\mathcal{B}}$, together with the error of evaluating the chemical potential at the spherically symmetric density with the same cumulative mass.
\end{assumption}
\noindent{}Proposition~\ref{prop:shell-integration} does not establish that the asymmetries remain small as the dynamics evolves. In that sense, Assumption~\ref{ass:spherical-closure} is the premise that the theory rests on. Let us now derive the evolution equation for $m_t$ by applying the integral operator directly to the three-dimensional SPDE, and see in detail how the geometric scaling factors cancel exactly.
\begin{proposition}[Dynamics of cumulative mass]
\label{prop:mass-dynamics}
Under Assumption~\ref{ass:spherical-closure}, let $\rho_t \in \mathcal{H}_{\text{rad}}$ be a spherically symmetric density field satisfying the It\^{o} stochastic partial differential equation~\eqref{eq:full_spde_ito}, which in the present notation reads
\begin{equation}
    d\rho_t = \nabla \cdot \left( D_0 \rho \nabla (\beta \mu) \right) dt + \varepsilon \nabla \cdot \left( \sqrt{2D_0 \rho} \; d\zeta^\delta_t \right).
\end{equation}
The cumulative mass profile $m_t(r) = \int_0^r \rho_t(s) 4\pi s^2 ds$ evolves according to the scalar stochastic differential equation:
\begin{equation}
    dm_t(r) = \bigl[\mathcal{A}_m(m_t) - \varepsilon^2 \mathfrak{S}_m^\delta(r)\bigr] \, dt + \varepsilon \mathcal{B}_m(m_t) \circ d\xi^\delta_t,
\end{equation}
where the drift and diffusion operators are:
\begin{align}
    \mathcal{A}_m(m) &= - D_0 \beta \frac{\partial m}{\partial r} \frac{\delta \mathcal{F}_m}{\delta m}, \label{eq:mass-drift-explicit} \\
    \mathcal{B}_m(m) &= \sqrt{2 D_0 \frac{\partial m}{\partial r}}, \label{eq:mass-diffusion-explicit}
\end{align}
and $\xi^\delta_t(r)$ is the scalar radial noise field.
\end{proposition}
\begin{proof}
Because $\mathcal{T}$ is a linear spatial integral operator, its second Fr\'{e}chet derivative vanishes identically, so that the It\^{o} formula for the change of variable $\rho \mapsto \mathcal{T}\rho$ carries no second-order term and reduces to $dm_t = \mathcal{T}[d\rho_t]$. No correction is generated by the change of variable itself, and every correction which appears below is generated instead by a change of stochastic convention, for which the relevant object is the first derivative of the diffusion operator rather than the second derivative of the map. We apply $\mathcal{T}$ to the drift and diffusion terms separately.
\begin{proofpart}[The drift term]
Applying the mass operator to the deterministic divergence term yields a volume integral over the ball $B_r$:
\begin{equation}
    \mathcal{A}_m(r) = \int_{B_r} \nabla \cdot \left( D_0 \rho(\mathbf{x}) \nabla (\beta \mu(\mathbf{x})) \right) d\mathbf{x}.
\end{equation}
Applying the Divergence Theorem transforms this volume integral into a surface integral over the sphere $\partial B_r$:
\begin{equation}
    \mathcal{A}_m(r) = \oint_{\partial B_r} \left( D_0 \rho(\mathbf{s}) \nabla (\beta \mu(\mathbf{s})) \right) \cdot \hat{\mathbf{n}} \, dS.
\end{equation}
Because the density $\rho$ and the chemical potential $\mu$ are spherically symmetric, the gradient $\nabla \mu$ points strictly in the radial direction $\hat{\mathbf{n}}$. Therefore, the inner product simplifies to $\nabla \mu \cdot \hat{\mathbf{n}} = \frac{\partial \mu}{\partial r}$. Since $\rho$ and $\frac{\partial \mu}{\partial r}$ depend only on the radial coordinate $r$, they are constant over the surface of integration and can be factored out:
\begin{equation}
    \mathcal{A}_m(r) = D_0 \beta \frac{\partial \mu(r)}{\partial r} \oint_{\partial B_r} \rho(\mathbf{s}) \, dS.
\end{equation}
The remaining surface integral equates to the radial derivative of the cumulative mass, $\frac{\partial m}{\partial r} = 4\pi r^2 \rho(r) = \oint_{\partial B_r} \rho(\mathbf{s}) \, dS$. Substituting this definition yields:
\begin{equation}
    \mathcal{A}_m(r) = D_0 \beta \frac{\partial m}{\partial r} \frac{\partial \mu}{\partial r}.
\end{equation}
Applying the duality relation $\frac{\partial \mu}{\partial r} = - \frac{\delta \mathcal{F}_m}{\delta m}$ derived in Lemma \ref{lemma:mu-duality} gives the drift of the mass process in the It\^{o} frame, since $\mathcal{T}$ is linear and therefore commutes with the It\^{o} integral:
\begin{equation}
    \mathcal{A}_m(r) = - D_0 \beta \frac{\partial m}{\partial r} \frac{\delta \mathcal{F}_m}{\delta m}.
\end{equation}
The projection of Section~\ref{sec:projection} requires the Stratonovich form, and converting subtracts the correction of Lemma~\ref{lemma:ito-strat-correction} below, which is the term carried by the drift in the statement.
\end{proofpart}
\begin{proofpart}[The diffusion term]
Applying the mass operator to the stochastic flux divergence:
\begin{equation}
    \varepsilon \mathcal{B}_m(m_t) d\xi^\delta_t = \varepsilon \int_{B_r} \nabla \cdot \left( \sqrt{2D_0 \rho(\mathbf{x})} d\zeta^\delta_t(\mathbf{x}) \right) d\mathbf{x}.
\end{equation}
Applying the Divergence Theorem:
\begin{equation}
    \varepsilon \mathcal{B}_m(m_t) d\xi^\delta_t = \varepsilon \oint_{\partial B_r} \sqrt{2D_0 \rho(\mathbf{s})} d\zeta^\delta_t(\mathbf{s}) \cdot \hat{\mathbf{n}} \, dS.
\end{equation}
Since $\rho$ is spherically symmetric, the scalar term $\sqrt{2D_0 \rho(r)}$ is constant over $\partial B_r$ and can be extracted from the integral. Using the inverse density relation (Remark \ref{rem:inverse-density}), we express $\rho$ strictly in terms of $m$:
\begin{equation}
    \varepsilon \mathcal{B}_m(m_t) d\xi^\delta_t = \varepsilon \sqrt{2D_0 \left( \frac{1}{4\pi r^2} \frac{\partial m}{\partial r} \right)} \oint_{\partial B_r} d\zeta^\delta_t(\mathbf{s}) \cdot \hat{\mathbf{n}} \, dS.
\end{equation}
The remaining surface integral is precisely the integrated radial noise $d\tilde{W}^\delta(r, t)$. Since the variance of a surface integral of a fluctuating field scales with the area $4\pi r^2$ of the sphere, its standard deviation scales as $\sqrt{4\pi r^2}$. We define the normalised scalar radial noise $\xi^\delta(r,t)$ via the square-root-of-area scaling $d\tilde{W}^\delta(r, t) = \sqrt{4\pi r^2} \, d\xi^\delta(r, t)$. Substituting this relation yields:
\begin{align}
    \varepsilon \mathcal{B}_m(m_t) d\xi^\delta_t &= \varepsilon \frac{\sqrt{2D_0 \frac{\partial m}{\partial r}}}{\sqrt{4\pi r^2}} \left( \sqrt{4\pi r^2} \, d\xi^\delta(r, t) \right).
\end{align}
The geometric factors $\sqrt{4\pi r^2}$ cancel exactly in the numerator and denominator, yielding the final operator form:
\begin{equation}
    \varepsilon \mathcal{B}_m(m_t) \, d\xi^\delta_t = \varepsilon \sqrt{2 D_0 \frac{\partial m}{\partial r}} \, d\xi^\delta(r, t).
\end{equation}
This identifies the diffusion coefficient $\mathcal{B}_m$ of the It\^{o} equation, and its Stratonovich form follows from Lemma~\ref{lemma:ito-strat-correction}.
\end{proofpart}
\end{proof}
\begin{lemma}[The explicit It\^{o}--Stratonovich correction]
\label{lemma:ito-strat-correction}
Let the dynamics of the cumulative mass $m_t$ be governed by the diffusion operator $\mathcal{B}_m(m)$ acting as a multiplication operator by $\sigma[m](r) = \sqrt{2D_0 \partial_r m(r)}$ and evaluated under the spatially correlated radial noise field with covariance $C_\delta(r, s)$. The It\^{o}--Stratonovich correction term $\mathfrak{S}_m^\delta(r)$ is generically non-zero and evaluates to:
\begin{equation}
    \mathfrak{S}_m^\delta(r) \coloneqq \frac{1}{2} \operatorname{Tr}_{\mathcal{Q}} \left( D\mathcal{B}_m[m] \cdot \mathcal{B}_m[m] \right) = \frac{D_0}{2\sigma(r)} \left[ \sigma'(r) V_\delta(r) + \frac{1}{2} \sigma(r) V_\delta'(r) \right],
\end{equation}
where $V_\delta(r) \coloneqq C_\delta(r, r)$ is the local variance of the radial noise.
\end{lemma}
\begin{proof}
We proceed via exact evaluation of the trace in the Karhunen--Lo\`{e}ve eigenbasis of the radial covariance operator.
\begin{proofpart}[The Fr\'{e}chet derivative of the diffusion coefficient]
We compute the Fr\'{e}chet derivative $D\sigma[m]: H^1_{\mathrm{loc}}((0,L]) \to L^2_{\mathrm{loc}}((0,L])$.
Consider a perturbation $m \to m + \epsilon \phi$ with $\phi \in H^1_{\mathrm{loc}}((0,L])$. Expanding $\sigma[m + \epsilon\phi](r) = \sqrt{2D_0 \partial_r(m + \epsilon\phi)}$ to first order in $\epsilon$:
\begin{align}
    D\sigma[m](\phi)(r) &= \frac{d}{d\epsilon}\bigg|_{\epsilon=0} \sqrt{2D_0(\partial_r m + \epsilon \partial_r \phi)} = \frac{D_0}{\sigma(r)} \partial_r \phi(r).
\end{align}
The prefactor $D_0/\sigma(r)$ diverges at $r=0$ at rate $r^{-1}$, since $\sigma(r) = \sqrt{2D_0 \cdot 4\pi r^2\rho(r)} = \mathcal{O}(r)$ (Remark~\ref{rem:origin-singularity}). This pointwise singularity does not invalidate the Galerkin evaluation that follows, because the trace is ultimately projected onto the tangent space via the weighted inner product $\langle \cdot, \cdot \rangle_\varpi$, and the integrability of all projected quantities is guaranteed by the compensating vanishing of the tangent basis vectors $\mathbf{e}_k^m(r) = \mathcal{O}(r^3)$ at the origin (Lemma~\ref{lemma:boundary-vanishing} and Remark~\ref{rem:origin-singularity}).
\end{proofpart}
\begin{proofpart}[Galerkin evaluation of the stochastic trace]
Let $\{e_k, \lambda_k\}_{k=1}^\infty$ denote the eigenpairs of the radial covariance operator, so that $C_\delta(r,s) = \sum_{k=1}^\infty \lambda_k e_k(r) e_k(s)$. The Stratonovich correction is defined as the trace
\begin{equation}
    \mathfrak{S}_m^\delta(r) = \frac{1}{2} \sum_{k=1}^\infty \lambda_k \bigl[D\sigma[m]\bigl(\sigma \, e_k\bigr)\bigr](r) \, e_k(r).
\end{equation}
For each $k$, we evaluate the action of the Fr\'{e}chet derivative on the function $\phi_k(r) \coloneqq \sigma(r) e_k(r) \in H^1_{\mathrm{loc}}((0,L])$:
\begin{equation}
    D\sigma[m](\sigma \, e_k)(r) = \frac{D_0}{\sigma(r)} \frac{\partial}{\partial r}\bigl[\sigma(r) e_k(r)\bigr] = \frac{D_0}{\sigma(r)}\bigl[\sigma'(r) e_k(r) + \sigma(r) e_k'(r)\bigr].
\end{equation}
We first retain a finite truncation at order $K$ and collect terms:
\begin{align}
    \mathfrak{S}_{m,K}^\delta(r) &= \frac{D_0}{2\sigma(r)} \sum_{k=1}^K \lambda_k \bigl[\sigma'(r) e_k(r)^2 + \sigma(r) e_k'(r) e_k(r)\bigr] \nonumber\\
    &= \frac{D_0}{2\sigma(r)} \biggl[\sigma'(r) \sum_{k=1}^K \lambda_k e_k(r)^2 + \sigma(r) \sum_{k=1}^K \lambda_k e_k'(r) e_k(r)\biggr].
\end{align}
\end{proofpart}
\begin{proofpart}[Identification with covariance evaluations and passage to the limit]
By the spectral representation of the covariance kernel, the partial sums converge for every $r \in (0,L]$, as we now show:
\begin{equation}
    \sum_{k=1}^K \lambda_k e_k(r)^2 \xrightarrow{K \to \infty} C_\delta(r, r) = V_\delta(r),
\end{equation}
\begin{equation}
    \sum_{k=1}^K \lambda_k e_k'(r) e_k(r) \xrightarrow{K \to \infty} \left.\frac{\partial C_\delta(r, s)}{\partial s}\right|_{s=r} = \frac{1}{2} V_\delta'(r),
\end{equation}
where the second identity follows from Lemma~\ref{lemma:radial-covariance}. Both limits are evaluated through a representation of the radial noise in the noise space, which the proof now constructs. For $h \in \mathcal{U}_0$ let $f_h(r) \coloneqq (4\pi r^2)^{-1/2} \oint_{\partial B_r} (J_\delta h) \cdot \hat{\mathbf{n}} \, dS$ be the normalised flux of $J_\delta h$, so that Definitions~\ref{def:mollifier} and~\ref{def:radial-noise-norm-ref} give $\xi^\delta_t(r) = \sum_k \beta_k(t) f_{h_k}(r)$. Writing the convolution out and exchanging the two integrals by Fubini's theorem,
\begin{equation}
    \begin{aligned}
    f_h(r) &= \langle h, \vartheta_r \rangle_{\mathcal{U}_0} ,
    \\
    \vartheta_r(\mathbf{z}) &\coloneqq \frac{1}{\sqrt{4\pi r^2}} \oint_{\partial B_r} \psi_\delta(\mathbf{x} - \mathbf{z}) \, \hat{\mathbf{n}}(\mathbf{x}) \, dS_{\mathbf{x}}
      = \frac{r}{\sqrt{4\pi}} \int_{S^2} \psi_\delta(r\omega - \mathbf{z}) \, \omega \, d\omega ,
    \end{aligned}
    \label{eq:flux-kernel}
\end{equation}
where $d\omega$ is the surface measure upon the unit sphere. The integrand of the last form is smooth in $(r,\mathbf{z})$ and vanishes for $|\mathbf{z}| > r + \delta$, so the map $r \mapsto \vartheta_r$ is continuously differentiable into $\mathcal{U}_0$, with the derivative $\partial_r\vartheta_r$ obtained by differentiating under the integral sign, and $f_h'(r) = \langle h, \partial_r\vartheta_r \rangle_{\mathcal{U}_0}$. The Cauchy--Schwarz inequality against $d\omega$ gives $\|\vartheta_r\|_{\mathcal{U}_0} \le \sqrt{4\pi}\, r \, \|\psi_\delta\|_{L^2}$, so the map $F : h \mapsto f_h$ is Hilbert--Schmidt into $L^2(0,L)$, with $\sum_k \|F h_k\|^2 = \int_0^L \|\vartheta_r\|^2 \, dr$ by Parseval's identity, and the covariance operator of $\xi^\delta$ is $F F^*$, with eigenpairs $\{\lambda_k, e_k\}$. For $\lambda_k > 0$ set $u_k \coloneqq \lambda_k^{-1/2} F^* e_k$. Then
\begin{equation}
    \langle u_k, u_j \rangle_{\mathcal{U}_0} = (\lambda_k\lambda_j)^{-1/2} \langle F F^* e_k, e_j \rangle = \delta_{kj} ,
    \qquad
    F u_k = \lambda_k^{-1/2} F F^* e_k = \lambda_k^{1/2} e_k ,
\end{equation}
so the $u_k$ are orthonormal, and the continuous representative of $e_k$ satisfies $\lambda_k^{1/2} e_k(r) = \langle u_k, \vartheta_r \rangle$ and $\lambda_k^{1/2} e_k'(r) = \langle u_k, \partial_r\vartheta_r \rangle$. Let $w$ be orthogonal to every $u_k$. Then $\langle e_k, F w \rangle = \langle F^* e_k, w \rangle = 0$ for $\lambda_k > 0$, and also for $\lambda_k = 0$, since then $\|F^* e_k\|^2 = \langle F F^* e_k, e_k \rangle = 0$. The $e_k$ being a complete orthonormal family, $F w = 0$, so the continuous function $f_w$ vanishes, and with it $f_w'$. Hence $\langle w, \vartheta_r \rangle = \langle w, \partial_r\vartheta_r \rangle = 0$ for every $r$, and $\vartheta_r$ and $\partial_r\vartheta_r$ lie in the closed span of the $u_k$. Parseval's identity in that span, polarised for the middle sum, gives the three series absolutely convergent, with
\begin{equation}
    \sum_k \lambda_k e_k(r)^2 = \|\vartheta_r\|^2 ,
    \qquad
    \sum_k \lambda_k e_k'(r) e_k(r) = \langle \vartheta_r, \partial_r\vartheta_r \rangle ,
    \qquad
    \sum_k \lambda_k \bigl[e_k'(r)\bigr]^2 = \|\partial_r\vartheta_r\|^2 .
\end{equation}
Parseval's identity over the basis $\{h_k\}$ gives in the same way $C_\delta(r,s) = \sum_k f_{h_k}(r) f_{h_k}(s) = \langle \vartheta_r, \vartheta_s \rangle$, so that $\|\vartheta_r\|^2 = C_\delta(r,r) = V_\delta(r)$, $\langle \vartheta_r, \partial_r\vartheta_r \rangle = \partial_s C_\delta(r,s)|_{s=r}$, and $\|\partial_r\vartheta_r\|^2 = \partial_r\partial_s C_\delta(r,s)|_{s=r} = -\kappa_\delta(r)$ by Definition~\ref{def:kernel-invariants}. These are the two limits displayed above, and the third sum enters the Cauchy--Schwarz bound $\sum_k \lambda_k |e_k'(r)e_k(r)| \le \bigl(-\kappa_\delta(r)V_\delta(r)\bigr)^{1/2}$ upon the second. In particular the third identity is
\begin{equation}
    -\kappa_\delta(r) = \left.\frac{\partial^2 C_\delta(r,s)}{\partial r\,\partial s}\right|_{s=r}
    = \sum_k \lambda_k \bigl[e_k'(r)\bigr]^2 = \|\partial_r\vartheta_r\|^2 \;\ge\; 0 ,
    \label{eq:kappa-sign}
\end{equation}
so that $\kappa_\delta \le 0$ throughout, which is the reason for the sign in Eq.~\eqref{eq:kappa-convention}. Substituting into the expression for $\mathfrak{S}_{m,K}^\delta$ and passing $K \to \infty$:
\begin{equation}
    \mathfrak{S}_m^\delta(r) = \frac{D_0}{2\sigma(r)} \left[ \sigma'(r) V_\delta(r) + \frac{1}{2} \sigma(r) V_\delta'(r) \right].
\end{equation}
This expression is well defined pointwise for all $r \in (0,L]$. The prefactor $\sigma(r)^{-1}$ grows like $r^{-1}$ at the origin, but $V_\delta$ vanishes there like $r^4$, and Remark~\ref{rem:origin-singularity} shows that the correction itself has no singularity at the origin.
\end{proofpart}
\end{proof}
\subsection{The concentration manifold hypothesis}
\label{subsec:manifold-hypothesis}
Having established the cumulative mass $m(r)$ as the canonical variable for the stochastic dynamics in Section \ref{sec:reduction}, we now formalise the geometry of the dimensional reduction. The reduced configuration space is constructed as the structure induced by the action of the mass operator $\mathcal{T}$ on the physical density manifold of Definition~\ref{def:manifold}, instead of being an arbitrary ansatz. Furthermore, we formalise the spectral gap assumption that justifies the adiabatic elimination of orthogonal modes.
\begin{remark}[Notation]
Recall that $d$ denotes the dimension of the reduced manifold, fixed before Definition~\ref{def:immersion}, while $N$ denotes the total conserved mass (Definition~\ref{def:density-space}).
\end{remark}
\subsubsection{The induced mass manifold}
\label{subsubsec:induced-manifold}
Recall that the physical manifold $\mathcal{M} \subset \mathcal{H}$ is the image of the parameter space $\mathcal{X}$ under the density embedding $\iota$. The cumulative mass operator $\mathcal{T}: \mathcal{H}_{\text{rad}} \to H^1([0, L])$ of Definition~\ref{def:mass-operator} is bounded and injective, but it is not onto. Its inverse carries the factor $(4\pi r^2)^{-1}$, so a profile $m \in H^1([0,L])$ with $m'(0) \neq 0$ does not correspond to any density in $L^2(V)$. What $\mathcal{T}$ provides is an isomorphism onto its image,
\begin{equation}
    \mathcal{M}_{\mathrm{mass}} \coloneqq \mathcal{T}\bigl(\mathcal{H}_{\text{rad}}\bigr) = \Bigl\{ m \in H^1([0,L]) : m(0)=0,\ (4\pi r^2)^{-1}m' \in L^2(V) \Bigr\} ,
    \label{eq:mass-image}
\end{equation}
which is the subspace of mass profiles with a derivative that vanishes fast enough at the centre. Every embedding used in this work takes its values in $\mathcal{M}_{\mathrm{mass}}$, since $M(\cdot;\mathbf{X})$ is generated by a density profile, and $\mathcal{T}$ therefore induces a diffeomorphic manifold structure in the mass representation.
\begin{definition}[The induced mass manifold]
We define the \textit{mass manifold} $\mathcal{M}_m \subset H^1([0, L])$ as the image of the physical manifold $\mathcal{M}$ under the operator $\mathcal{T}$:
\begin{equation}
    \mathcal{M}_m \coloneqq \mathcal{T}(\mathcal{M}) = \{ m \in H^1([0, L]) : m = \mathcal{T}[\rho], \, \rho \in \mathcal{M} \}.
\end{equation}
The corresponding embedding $\iota_m: \mathcal{X} \to \mathcal{M}_m$ is defined by the composition $\iota_m = \mathcal{T} \circ \iota$. We denote the mass profile parameterised by the collective coordinates $\mathbf{X}$ as $M(\cdot; \mathbf{X})$:
\begin{equation}
    M(r; \mathbf{X}) \coloneqq \iota_m(\mathbf{X})(r) = \int_0^r \rho(s; \mathbf{X}) \, 4\pi s^2 \, ds.
\end{equation}
Since the density embedding is a smooth immersion (Definition \ref{def:immersion}), $M(\cdot; \mathbf{X})$ defines a $d$-dimensional differentiable submanifold within the affine space of functions satisfying $m(L)=N$.
\end{definition}
\begin{assumption}[Finite covariance contraction of the embedding curvature]
\label{ass:trace-class-embedding}
We strictly assume that the coordinate embedding $M(\mathbf{X})$ is sufficiently smooth so that its second derivative defines a trace-class operator with respect to the noise covariance $C_\delta$. In particular, the geometric contraction of the manifold's curvature vectors $W_{ij}(r; \mathbf{X}) \coloneqq \frac{\partial^2 M(r; \mathbf{X})}{\partial X^i \partial X^j}$ with the covariance kernel must be finite. This condition guarantees that infinite-dimensional geometric traces, such as the mean curvature of the embedding, do not yield ultraviolet divergences in the analytical expansions.
\end{assumption}
\begin{remark}[Gaussian ansatz and the covariance contraction]
\label{rem:gaussian-ansatz}
Assumption~\ref{ass:trace-class-embedding} is non-vacuous. Consider the one-parameter family of spherically symmetric density profiles given by the Gaussian ansatz
\begin{equation}
    \rho(r; R) = \rho_0 + \frac{A_{\mathrm{exc}}}{(2\pi w^2)^{1/2}} \exp\!\left( -\frac{(r - R)^2}{2w^2} \right), \qquad R \in (0, L),
\end{equation}
where $\rho_0 = \rho_0(R) > 0$ is the background density, fixed by the constraint $M(L;R) = N$ (and hence the family lies in $\mathcal{D}_N$), $A_{\mathrm{exc}} > 0$ is the excess mass amplitude (with a subscript to keep it apart from the moment $A$ of Proposition~\ref{prop:kernel-asymptotics}), and $w > 0$ is the interfacial width. The parameter $R$ is the radius of the concentration shell, and the manifold dimension is $d = 1$.
The cumulative mass profile is $M(r; R) = \int_0^r \rho(s; R) 4\pi s^2 \, ds$, and the second derivative of the embedding with respect to $R$ is
\begin{equation}
    W(r; R) = \frac{\partial^2 M}{\partial R^2}(r; R) = \int_0^r \frac{\partial^2 \rho(s; R)}{\partial R^2} 4\pi s^2 \, ds.
\end{equation}
Since $\partial^2 \rho / \partial R^2$ is a smooth function of $s$ upon $[0,L]$, the curvature vector $W(\cdot; R)$ lies in $C^\infty([0, L])$. The trace-class contraction with the covariance kernel evaluates to
\begin{equation}
    \int_0^L \int_0^L W(r; R) \, C_\delta(r, s) \, W(s; R) \, dr \, ds \le \|W\|_{L^2}^2 \, \|C_\delta\|_{\mathrm{op}} < \infty,
\end{equation}
which is finite because $W$ is square integrable and $C_\delta$ is a bounded operator. This confirms that the Gaussian ansatz satisfies Assumption~\ref{ass:trace-class-embedding} for any $\delta > 0$.
\end{remark}
\subsubsection{Induced tangent structure}
\label{subsubsec:induced-tangent}
The local linearisation of the manifold is described by its tangent bundle, which is the pushforward of the density tangent bundle (Definition~\ref{def:tangent-space}).
\begin{definition}[Mass tangent basis]
Let $\{\mathbf{e}_k\}_{k=1}^d \subset T_{\rho}\mathcal{M}$ be the tangent basis vectors of the density manifold. We define the tangent basis $\{\mathbf{e}_k^m\}_{k=1}^d$ for $\mathcal{M}_m$ via the linear action of $\mathcal{T}$:
\begin{equation}
    \mathbf{e}_k^m(r; \mathbf{X}) \coloneqq \mathcal{T}[\mathbf{e}_k(\cdot; \mathbf{X})](r) = \int_0^r \frac{\partial \rho(s; \mathbf{X})}{\partial X^k} \, 4\pi s^2 \, ds.
\end{equation}
Equivalently, these are the partial derivatives of the mass embedding:
\begin{equation}
    \mathbf{e}_k^m(r; \mathbf{X}) = \frac{\partial M(r; \mathbf{X})}{\partial X^k}.
\end{equation}
\end{definition}
\begin{lemma}[Boundary vanishing]
\label{lemma:boundary-vanishing}
For any configuration $\mathbf{X} \in \mathcal{X}$, the tangent vectors in the mass representation satisfy homogeneous Dirichlet boundary conditions:
\begin{equation}
    \mathbf{e}_k^m(0; \mathbf{X}) = 0 \quad \text{and} \quad \mathbf{e}_k^m(L; \mathbf{X}) = 0.
\end{equation}
Near the origin, moreover, $\mathbf{e}_k^m(r;\mathbf{X}) = \mathcal{O}(r^3)$ and $\partial_r\mathbf{e}_k^m(r;\mathbf{X}) = \mathcal{O}(r^2)$ as $r \to 0^+$. Consequently, the tangent space $T_{\mathbf{X}}\mathcal{M}_m$ is a closed subspace of the Hilbert space $H^1_0([0, L])$.
\end{lemma}
\begin{proof}
We proceed by direct evaluation of the boundary values of the differentiated mass profile.
\begin{proofpart}[At the origin $r=0$]
By definition of the cumulative mass operator (Definition \ref{def:mass-operator}), $M(0; \mathbf{X}) = 0$ for all $\mathbf{X}$, as the volume of a ball of radius zero is null. Differentiating this constant constraint with respect to the parameter $X^k$:
\begin{equation}
    \mathbf{e}_k^m(0; \mathbf{X}) = \frac{\partial}{\partial X^k} M(0; \mathbf{X}) = \frac{\partial}{\partial X^k} (0) = 0.
\end{equation}
\end{proofpart}
\begin{proofpart}[At the boundary $r=L$]
The total mass of the system is a conserved quantity $N$, fixed for all configurations in the density configuration space $\mathcal{D}_N$ (Definition \ref{def:density-space}). Thus, $M(L; \mathbf{X}) = \int_V \rho(\mathbf{r}; \mathbf{X}) d\mathbf{r} = N$. Differentiating this constraint:
\begin{equation}
    \mathbf{e}_k^m(L; \mathbf{X}) = \frac{\partial}{\partial X^k} M(L; \mathbf{X}) = \frac{\partial}{\partial X^k} (N) = 0.
\end{equation}
\end{proofpart}
\begin{proofpart}[The order of vanishing at the origin]
By Definition~\ref{def:mass-operator}, $M(r;\mathbf{X}) = \int_0^r 4\pi s^2\rho(s;\mathbf{X})\,ds$. The parameter derivative $\partial_{X^k}\rho$ is continuous and bounded upon a neighbourhood $(0,r_0]$ of the origin, locally uniformly in $\mathbf{X}$, by the requirement placed upon the embedding, so that differentiation under the integral sign is permitted and gives
\begin{equation}
    \mathbf{e}_k^m(r;\mathbf{X}) = \frac{\partial}{\partial X^k}\int_0^r 4\pi s^2\rho(s;\mathbf{X})\,ds = \int_0^r 4\pi s^2\,\partial_{X^k}\rho(s;\mathbf{X})\,ds .
\end{equation}
Since $\partial_{X^k}\rho(\cdot\,;\mathbf{X})$ is continuous upon $(0,r_0]$, the fundamental theorem of calculus then gives $\partial_r\mathbf{e}_k^m(r;\mathbf{X}) = 4\pi r^2\,\partial_{X^k}\rho(r;\mathbf{X})$. Writing $C_\rho$ for a bound of $|\partial_{X^k}\rho|$ upon $(0,r_0]$, for every $r \in (0,r_0]$,
\begin{equation}
    \bigl|\partial_r\mathbf{e}_k^m(r;\mathbf{X})\bigr| \le 4\pi C_\rho\,r^2 ,
    \qquad
    \bigl|\mathbf{e}_k^m(r;\mathbf{X})\bigr| \le \int_0^r 4\pi s^2 C_\rho\,ds = \frac{4\pi C_\rho}{3}\,r^3 ,
\end{equation}
which are the two orders stated.
\end{proofpart}
Since the basis vectors vanish at the boundaries and are smooth images of $\mathbf{e}_k$, they belong to $H^1_0([0, L])$.
\end{proof}
\subsubsection{Orthogonal decomposition}
\label{subsubsec:orthogonal-decomposition}
In order to define the reduced dynamics, we must decompose the infinite-dimensional state $m$ into a finite-dimensional projection on the manifold and a residual fluctuation. To this end, we proceed constructively, and we distinguish strictly between the affine manifold and its linear tangent bundle.
\begin{definition}[Weighted Hilbert space]
\label{def:weighted-space}
To account for the spatial inhomogeneity of the noise, we introduce a weight function $\varpi(r)$, strictly positive and smooth on the open interval $(0, L]$. The weight selected by Lemma~\ref{lemma:metric-selection} is $\varpi = (4\pi r^2\rho)^{-1}$, which is of order $r^{-2}$ at the centre, so it is not bounded on the closed interval and the weighted space is defined on $(0,L]$. Every pairing used below remains finite, because the tangent fields vanish as $\mathcal{O}(r^3)$ at the origin by Lemma~\ref{lemma:boundary-vanishing}, and Remark~\ref{rem:origin-singularity} treats the point $r = 0$ in detail. The ambient space is accordingly not $L^2(0,L)$ with a new pairing, since not every square-integrable function has finite weighted norm, but the weighted space
\begin{equation}
    \mathcal{H}_\varpi \coloneqq L^2\bigl((0,L),\,\varpi(r)\,dr\bigr) ,
    \label{eq:weighted-space}
\end{equation}
upon which we define the weighted inner product $\langle \cdot, \cdot \rangle_{\varpi}$ and the induced norm $\|\cdot\|_{\varpi}$ by
\begin{equation}
    \langle u, v \rangle_{\varpi} \coloneqq \int_0^L u(r) v(r) \varpi(r) \, dr.
\end{equation}
Being an $L^2$ space over a $\sigma$-finite Borel measure on an interval, $\mathcal{H}_\varpi$ is a separable Hilbert space and therefore a Polish space, which is what the disintegration of Theorem~\ref{thm:disintegration} requires of it.
The weight is a \emph{functional of the state}, $\varpi = \varpi[m]$, not a fixed function of $r$. Lemma~\ref{lemma:metric-selection} below fixes it as $\varpi[m](r) = (\partial_r m(r))^{-1}$, and Theorem~\ref{thm:wasserstein-isometry} confirms the same value from the Wasserstein geometry. The following convention is adopted throughout:
\begin{enumerate}[(i)]
    \item At a manifold point $\mathbf{X}$ we write $\varpi(\cdot\,;\mathbf{X}) \coloneqq \varpi[M(\mathbf{X})]$, and $\mathcal{H}_{\varpi(\cdot;\mathbf{X})}$ is the fibre Hilbert space at that base point. Orthogonality, the projectors $\Pi^{\parallel}_{\mathbf{X}}$ and $\Pi^{\perp}_{\mathbf{X}}$, the normal bundle and the tubular trivialisation are all defined \emph{relative to that base point}, and are held fixed as the fluctuation $\eta$ varies within the fibre.
    \item Two distinct differentiations occur in this work and they obey different conventions. A derivative taken \emph{along} the manifold differentiates the weight together with the field, since $\varpi[m]$ depends upon $m$, and it is that derivative which produces the weight-gradient terms of the geometric anomaly, of the frame correction and of the connection correction. A derivative taken \emph{transversally at a fixed base point} $\mathbf{X}$, which is the one appearing in $D^2\Pi_k$ and in the fibre integration of Section~\ref{sec:kolmogorov-reduction}, holds the chart fixed, so that the tangent fields $\mathbf{e}_k^m(\cdot;\mathbf{X})$ and the fibre metric $\varpi(\cdot;\mathbf{X})$ do not vary. This applies in particular to the second derivative $D^2\Pi_k$ evaluated in Section~\ref{sec:commutativity}, in which the weight is $\varpi(\cdot\,;\mathbf{X})$ and not $\varpi[m]$.
    \item The parameter derivative $\partial_{X^l}\varpi$ means the derivative of $\varpi(\cdot\,;\mathbf{X})$ along the manifold, and is distinct from the spatial derivative $\varpi' = \partial_r \varpi$. The two are related by $\partial_{X^l}\varpi = -\varpi^2\,\partial_r \mathbf{e}^m_l$.
\end{enumerate}
The distinction in item~(iii) is not cosmetic, since the two derivatives carry different indices and different dimensions, and they may not be exchanged.
While the spectral analysis of the Hessian is often performed in the standard $L^2$ topology (corresponding to $\varpi=1$), the stochastic projection requires a weight consistent with the Onsager--Machlup functional. We leave $\varpi$ general for this derivation.
\end{definition}
\begin{definition}[Riemannian metric tensor]
\label{def:metric-tensor}
We define the components of the metric tensor (Gram matrix) $g(\mathbf{X})$ on the tangent space $T_{\mathbf{X}}\mathcal{M}_m$ via the weighted inner product of the basis vectors:
\begin{equation}
    g_{ij}(\mathbf{X}) \coloneqq \langle \mathbf{e}_i^m(\cdot; \mathbf{X}), \mathbf{e}_j^m(\cdot; \mathbf{X}) \rangle_{\varpi} = \int_0^L \frac{\partial M}{\partial X^i} \frac{\partial M}{\partial X^j} \varpi(r) \, dr.
\end{equation}
Since the embedding is an immersion and weights are positive, this matrix remains positive definite and invertible. We denote the components of the inverse matrix by $g^{ij}(\mathbf{X})$.
\end{definition}
\begin{definition}[Tangent projection operators]
\label{def:projectors}
Using the weighted metric tensor, we define the orthogonal projection operators acting on the Hilbert space:
\begin{enumerate}[(i)]
    \item The \emph{parallel projection} $\Pi_{\mathbf{X}}^{\parallel}: \mathcal{H}_\varpi \to T_{\mathbf{X}}\mathcal{M}_m$ is defined by the expansion:
    \begin{equation}
        \Pi_{\mathbf{X}}^{\parallel} u \coloneqq \sum_{i,j=1}^d g^{ij}(\mathbf{X}) \langle u, \mathbf{e}_j^m(\cdot; \mathbf{X}) \rangle_{\varpi} \, \mathbf{e}_i^m(\cdot; \mathbf{X}).
        \label{def:parallel-projector}
    \end{equation}
    \item The \emph{orthogonal projection} $\Pi_{\mathbf{X}}^{\perp}: \mathcal{H}_\varpi \to (T_{\mathbf{X}}\mathcal{M}_m)^{\perp_{\varpi}}$ is the complement:
    \begin{equation}
        \Pi_{\mathbf{X}}^{\perp} u \coloneqq u - \Pi_{\mathbf{X}}^{\parallel} u.
	\label{def:perp-projector}
    \end{equation}
\end{enumerate}
\end{definition}
\begin{lemma}[Projector-Jacobian identity]
\label{lemma:projector-jacobian}
Let $\mathcal{J}_{\iota_m}(\mathbf{X}): \mathbb{R}^d \to \mathcal{H}_\varpi$ be the Jacobian operator of the mass embedding, defined by $\mathcal{J}_{\iota_m}(\mathbf{X}) \mathbf{v} = \sum_{k=1}^d v^k \mathbf{e}_k^m(\cdot; \mathbf{X})$. Let $\mathcal{J}_{\iota_m}^{\dagger}(\mathbf{X})$ denote the adjoint operator with respect to the \textit{weighted} inner product. The parallel projection operator satisfies the factorisation:
\begin{equation}
    \Pi_{\mathbf{X}}^{\parallel} = \mathcal{J}_{\iota_m}(\mathbf{X}) \, [g(\mathbf{X})]^{-1} \, \mathcal{J}_{\iota_m}^{\dagger}(\mathbf{X}).
\end{equation}
\end{lemma}
\begin{proof}
Let $u \in \mathcal{H}_\varpi$ be an arbitrary test function. We evaluate the action of the operator composition on the right-hand side step by step.
\begin{proofpart}[Action of the weighted adjoint Jacobian]
The weighted adjoint operator $\mathcal{J}_{\iota_m}^{\dagger}(\mathbf{X}): \mathcal{H}_\varpi \to \mathbb{R}^d$ is defined by the duality relation $\langle \mathcal{J}\mathbf{v}, u \rangle_{\varpi} = \langle \mathbf{v}, \mathcal{J}^{\dagger}u \rangle_{\mathbb{R}^d}$. This maps the function $u$ to the vector of covariant coefficients:
\begin{equation}
    c_k = [\mathcal{J}_{\iota_m}^{\dagger}(\mathbf{X}) u]_k = \langle u, \mathcal{J}_{\iota_m}(\mathbf{X}) \hat{\mathbf{u}}_k \rangle_{\varpi} = \langle u, \mathbf{e}_k^m(\cdot; \mathbf{X}) \rangle_{\varpi},
\end{equation}
where $\hat{\mathbf{u}}_k$ is the standard basis vector in parameter space.
\end{proofpart}
\begin{proofpart}[Action of the inverse metric]
We apply the inverse metric matrix $[g(\mathbf{X})]^{-1}$ to the coefficient vector $\mathbf{c}$ to obtain the contravariant vector $\mathbf{v} \in \mathbb{R}^d$:
\begin{equation}
    v^i = \sum_{j=1}^d g^{ij}(\mathbf{X}) c_j = \sum_{j=1}^d g^{ij}(\mathbf{X}) \langle u, \mathbf{e}_j^m(\cdot; \mathbf{X}) \rangle_{\varpi}.
\end{equation}
\end{proofpart}
\begin{proofpart}[Action of the Jacobian]
Finally, we apply the Jacobian operator $\mathcal{J}_{\iota_m}(\mathbf{X})$ to the vector $\mathbf{v}$. By linearity:
\begin{equation}
    (\mathcal{J}_{\iota_m}(\mathbf{X}) \mathbf{v})(r) = \sum_{i=1}^d v^i \mathbf{e}_i^m(r; \mathbf{X}).
\end{equation}
Substituting the expression for $v^i$ from the previous step:
\begin{equation}
    (\mathcal{J}_{\iota_m}(\mathbf{X}) \mathbf{v})(r) = \sum_{i=1}^d \mathbf{e}_i^m(r; \mathbf{X}) \sum_{j=1}^d g^{ij}(\mathbf{X}) \langle u, \mathbf{e}_j^m(\cdot; \mathbf{X}) \rangle_{\varpi}  = \Pi_{\mathbf{X}}^{\parallel}u. \quad \text{(Def.~\eqref{def:parallel-projector})}
\end{equation}
\end{proofpart}
\end{proof}
\begin{definition}[Normal bundle]
\label{def:normal-bundle}
We define the weighted normal bundle $\mathcal{N} \subset \mathcal{X} \times \mathcal{H}_\varpi$ as the collection of all pairs $(\mathbf{X}, \eta)$ such that $\eta$ is orthogonal to the tangent space at $\mathbf{X}$, the orthogonality and the ambient space alike being taken in the weighted sense of Definition~\ref{def:weighted-space}:
\begin{equation}
    \mathcal{N} \coloneqq \{ (\mathbf{X}, \eta) \in \mathcal{X} \times \mathcal{H}_\varpi : \Pi_{\mathbf{X}}^{\parallel} \eta = 0 \}.
\end{equation}
\end{definition}
The Inverse Function Theorem gives the decomposition only locally in $\mathbf{X}$, and a uniform radius needs control of the normal injectivity radius.
\begin{assumption}[Uniform curvature bound]
\label{ass:uniform-curvature}
The second fundamental form of the embedding $M: \mathcal{X} \to \mathcal{H}_\varpi$ satisfies
\begin{equation}
    \sup_{\mathbf{X} \in \mathcal{X}} \| \mathbf{I\!I}(\mathbf{X}) \|_{\mathrm{op}} < \infty,
\end{equation}
where $\| \cdot \|_{\mathrm{op}}$ denotes the operator norm from $T_{\mathbf{X}}\mathcal{M}_m \otimes T_{\mathbf{X}}\mathcal{M}_m$ to $\mathcal{N}_{\mathbf{X}}$.
\end{assumption}
\begin{lemma}[Existence of tubular neighbourhood]
\label{lemma:tubular-existence}
Under Assumption~\ref{ass:uniform-curvature}, there exists a radius $\varrho_{\mathrm{tub}} > 0$ and an open neighbourhood $U_{\varrho_{\mathrm{tub}}}$ of the manifold $\mathcal{M}_m$ in $\mathcal{H}_\varpi$ such that for every $m \in U_{\varrho_{\mathrm{tub}}}$, there exists a unique pair $(\mathbf{X}, \eta) \in \mathcal{N}$ with $\|\eta\|_{\varpi} < \varrho_{\mathrm{tub}}$ satisfying $m = M(\cdot; \mathbf{X}) + \eta$.
\end{lemma}
\begin{proof}
We construct the decomposition map and prove that it is a local diffeomorphism. The normal spaces form a family of weighted Hilbert spaces varying with the base point, so that the domain of $\Psi$ is a Hilbert bundle rather than a fixed product of Banach spaces, and the inverse function theorem does not apply to it directly. We therefore fix $\mathbf{X}_0 \in \mathcal{X}$ and choose, as in Definition~\ref{def:local-trivialisation} below, the local trivialisation $\mathcal{R}_{\mathbf{X}}$ of the normal bundle over a neighbourhood of $\mathbf{X}_0$, which identifies every nearby fibre with the single model space $\mathcal{H}_\perp$. In those local coordinates $\Psi$ is a map between fixed Banach spaces and the theorem applies, which is why the conclusion obtained is local in $\mathbf{X}$ and is afterwards made uniform by the compactness of $\overline{\mathcal{X}}$.
\begin{proofpart}[Construction of the map]
Consider the map $\Psi: \mathcal{N} \to \mathcal{H}_\varpi$ defined by addition:
\begin{equation}
    \Psi(\mathbf{X}, \eta) = M(\cdot; \mathbf{X}) + \eta.
\end{equation}
This map is smooth because $M(\mathbf{X})$ is a smooth embedding. We restrict our attention to the ``zero section'' of the bundle, i.e., points where $\eta = 0$, which correspond to the manifold itself.
\end{proofpart}
\begin{proofpart}[Differentiation]
We compute the Fréchet derivative of $\Psi$ at a point on the manifold, $p_0 = (\mathbf{X}, 0)$. Let $(\mathbf{v}, \zeta) \in T_{p_0}\mathcal{N}$ be a tangent vector to the bundle. We evaluate the change in $\Psi$ under a perturbation $\Psi(\mathbf{X} + \mathbf{v}, 0 + \zeta) = M(\cdot; \mathbf{X} + \mathbf{v}) + \zeta$. Using the Taylor expansion of the mass embedding $M(\mathbf{X} + \mathbf{v}) = M(\mathbf{X}) + \mathcal{J}_{\iota_m}(\mathbf{X})\mathbf{v} + o(\|\mathbf{v}\|)$, we have:
\begin{align}
    \Psi(\mathbf{X} + \mathbf{v}, \zeta) - \Psi(\mathbf{X}, 0) &= (M(\mathbf{X}) + \mathcal{J}_{\iota_m}(\mathbf{X})\mathbf{v} + \zeta + o(\|\mathbf{v}\|) )- M(\mathbf{X}) \nonumber \\
    &= \mathcal{J}_{\iota_m}(\mathbf{X})\mathbf{v} + \zeta + o(\|\mathbf{v}\|).
\end{align}
The linear part defines the Fréchet derivative $D\Psi_{(\mathbf{X}, 0)}$:
\begin{equation}
    D\Psi_{(\mathbf{X}, 0)} \cdot (\mathbf{v}, \zeta) = \mathcal{J}_{\iota_m}(\mathbf{X}) \mathbf{v} + \zeta.
\end{equation}
Here, $\mathcal{J}_{\iota_m}(\mathbf{X}) \mathbf{v}$ lies in the tangent space $T_{\mathbf{X}}\mathcal{M}_m$, while $\zeta$ lies in the orthogonal fibre $(T_{\mathbf{X}}\mathcal{M}_m)^{\perp_{\varpi}}$ by the definition of the tangent space to the normal bundle at the zero section.
\end{proofpart}
\begin{proofpart}[Isomorphism property]
We now check if $D\Psi$ is an isomorphism.
For any vector $w \in \mathcal{H}_\varpi$, we can uniquely decompose it using the projection operators (Def. \ref{def:projectors}): $w = w_{\parallel} + w_{\perp}$.
To satisfy $D\Psi \cdot (\mathbf{v}, \zeta) = w$, we must have:
\begin{equation}
    \mathcal{J}_{\iota_m}(\mathbf{X}) \mathbf{v} = w_{\parallel} \quad \text{and} \quad \zeta = w_{\perp}.
\end{equation}
Since $\mathcal{J}_{\iota_m}$ is injective (it is an embedding), $\mathbf{v}$ is uniquely determined by $w_{\parallel}$. Since $\zeta$ is directly $w_{\perp}$, the solution is unique. Thus, $D\Psi_{(\mathbf{X}, 0)}$ is a linear bijection. By the Open Mapping Theorem, it is a topological isomorphism.
\end{proofpart}
\begin{proofpart}[Existence of decomposition]
Since $D\Psi$ is an isomorphism at the zero section (as proven above), the Inverse Function Theorem for Banach spaces implies that $\Psi$ is a local diffeomorphism, and since the closure $\iota_m(\overline{\mathcal{X}})$ is compact with positive reach (Remark~\ref{rem:tubular-scope}), $\Psi$ is injective upon $\{\|\eta\|_\varpi < \varrho_{\mathrm{tub}}\}$ for some $\varrho_{\mathrm{tub}} > 0$. Consequently, there exists an open neighbourhood $U_{\mathcal{N}} \subset \mathcal{N}$ of the zero section and an open neighbourhood $U_{\varrho_{\mathrm{tub}}} \subset \mathcal{H}_\varpi$ of the manifold such that $\Psi: U_{\mathcal{N}} \to U_{\varrho_{\mathrm{tub}}}$ is a bijection.
Thus, let $m \in U_{\varrho_{\mathrm{tub}}}$, and let $(\mathbf{X}^*, \eta^*) \coloneqq \Psi^{-1}(m)$ be the unique pair associated by bijection. By the fundamental property of the inverse map, applying $\Psi$ yields the identity $\Psi(\mathbf{X}^*, \eta^*) = \Psi(\Psi^{-1}(m)) = m$.  Substituting the explicit definition of the map $\Psi(\mathbf{X}, \eta) = M(\cdot; \mathbf{X}) + \eta$, we obtain the orthogonal decomposition equality:
\begin{equation}
    M(\cdot; \mathbf{X}^*) + \eta^* = m.
\end{equation}
Since the pair $(\mathbf{X}^*, \eta^*)$ lies in the domain $\mathcal{N}$ (the normal bundle), the residual $\eta^*$ must satisfy the orthogonality condition $\Pi_{\mathbf{X}^*}^{\parallel} \eta^* = 0$ by definition. The uniqueness of this decomposition is guaranteed by the bijectivity of $\Psi$ on $U_{\varrho_{\mathrm{tub}}}$.
\end{proofpart}
\end{proof}
\begin{remark}[Scope of the tubular neighbourhood]
\label{rem:tubular-scope}
The Inverse Function Theorem for Banach spaces, applied in the proof above, guarantees only the \emph{local} existence of the decomposition, in that for each $\mathbf{X}_0 \in \mathcal{X}$, there exists a neighbourhood $U_{\mathbf{X}_0}$ on which $\Psi$ is invertible. To obtain a \emph{uniform} tubular neighbourhood (a single $\varrho_{\mathrm{tub}} > 0$ valid for all $\mathbf{X} \in \mathcal{X}$), we need a uniform lower bound on the injectivity radius of the normal exponential map.
 Under Assumption~\ref{ass:uniform-curvature}, the principal curvatures are uniformly bounded, the normal exponential map has a positive injectivity radius, and the tubular neighbourhood radius $\varrho_{\mathrm{tub}}$ can be chosen independently of $\mathbf{X}$. The parameter domain $\mathcal{X}$ is open, so we cannot use compactness on $\mathcal{X}$ itself, but we can use it on its closure. Since $\overline{\mathcal{X}}$ is compact and the embedding extends smoothly to it, the image $\iota_m(\overline{\mathcal{X}})$ is a compact embedded submanifold of $\mathcal{H}_\varpi$ (the norms $\|\cdot\|_{\varpi(\cdot;\mathbf{X})}$, $\mathbf{X} \in \overline{\mathcal{X}}$, are uniformly equivalent, since the density is bounded above and below on $\overline{\mathcal{X}}$), and a compact embedded submanifold has positive reach. A uniform lower bound on the normal injectivity radius therefore exists, and $\varrho_{\mathrm{tub}}$ may be chosen independently of $\mathbf{X}$. Note that the curvature bound alone would not suffice for this conclusion, since a bound on the second fundamental form controls the local geometry but does not prevent distant portions of an embedded manifold from approaching one another. It is the compactness of the closure, together with injectivity of the embedding, that supplies the uniform radius. The assumption becomes non-trivial only for families of manifolds with dimension $d \to \infty$ or with unbounded extrinsic curvature.
\end{remark}
\begin{lemma}[Implicit definition of collective coordinates]
\label{lemma:orthogonality-condition}
The unique coordinate $\mathbf{X}$ identified by Lemma~\ref{lemma:tubular-existence} satisfies the weighted orthogonality condition
\begin{equation}
    \langle m - M(\cdot; \mathbf{X}), \mathbf{e}_k^m(\cdot; \mathbf{X}) \rangle_{\varpi(\cdot;\mathbf{X})} = 0, \quad k=1,\dots,d,
\end{equation}
in which the weight is taken at the base point itself, as Definition~\ref{def:weighted-space}~(i) prescribes. Equivalently, $\mathbf{Y} = \mathbf{X}$ is a stationary point of the weighted distance functional $\mathbf{Y}\mapsto\|m-M(\cdot;\mathbf{Y})\|_{\varpi(\cdot;\mathbf{X})}^2$, the weight being held at $\mathbf{X}$.
\end{lemma}
\begin{proof}
By Lemma \ref{lemma:tubular-existence}, $m = M(\mathbf{X}) + \eta$ with $\eta$ in the fibre of the normal bundle. By Definition \ref{def:normal-bundle}, this means $\Pi_{\mathbf{X}}^{\parallel} \eta = 0$. Using the explicit form of the projector (Def. \ref{def:projectors}), vanishing parallel projection implies $\langle \eta, \mathbf{e}_k^m \rangle_{\varpi} = 0$ for all basis vectors, which is the stated condition.
The condition is exactly the stationarity condition of the weighted distance functional with the weight held at $\mathbf{X}$. Differentiating $\mathbf{Y} \mapsto \|m - M(\cdot;\mathbf{Y})\|_{\varpi(\cdot;\mathbf{X})}^2$ gives
\begin{equation}
    \frac{\partial}{\partial Y^k}\bigl\|m - M(\cdot;\mathbf{Y})\bigr\|_{\varpi(\cdot;\mathbf{X})}^2
    = -2\bigl\langle m - M(\cdot;\mathbf{Y}), \mathbf{e}_k^m(\cdot;\mathbf{Y})\bigr\rangle_{\varpi(\cdot;\mathbf{X})} ,
\end{equation}
which at $\mathbf{Y} = \mathbf{X}$ vanishes for every $k$ precisely when the orthogonality condition holds. By Lemma \ref{lemma:tubular-existence} the pair $(\mathbf{X},\eta)$ is unique within the tubular neighbourhood, so $\mathbf{X}$ is the only point of that neighbourhood which satisfies the condition.
\end{proof}
%%%
\subsubsection{The stiffness hypothesis}
\label{subsubsec:stiffness}
The orthogonal decomposition constructed in Section \ref{subsubsec:orthogonal-decomposition} provides the static coordinate framework. However, the validity of the reduction to the manifold $\mathcal{M}_m$ relies on a dynamical property, namely that the manifold act as a stable attractor for the stochastic trajectories. To formalise it, we examine the linear stability of fluctuations in the weighted norm.
\begin{definition}[Linearised evolution operator]
\label{def:linearised-operator}
We define the \textit{linearised evolution operator} $\mathcal{L}_{\mathbf{X}}$ at a manifold point $M(\mathbf{X})$ as the Fréchet derivative of the drift:
\begin{equation}
    \mathcal{L}_{\mathbf{X}} \coloneqq D\mathcal{A}_m[M(\mathbf{X})].
\end{equation}
(Here, $D$ denotes the Fr\'{e}chet derivative in the Banach space.)
In other words, $\mathcal{L}_{\mathbf{X}}$ describes the first-order expansion of the dynamics in the weighted topology:
\begin{equation}
    \mathcal{A}_m(M+\eta) = \mathcal{A}_m(M) + \mathcal{L}_{\mathbf{X}} \eta + o(\|\eta\|_{\varpi}).
\end{equation}
We further define the \textit{normal stability operator} $\mathcal{L}_{\perp}$, also written $\mathcal{L}_\perp(\mathbf{X})$ when the base point must be displayed, as the restriction of this linearised dynamics to the weighted orthogonal subspace:
\begin{equation}
    \mathcal{L}_{\perp} \eta \coloneqq \Pi_{\mathbf{X}}^{\perp} \left( \mathcal{L}_{\mathbf{X}} \eta \right) = \Pi_{\mathbf{X}}^{\perp}\, D\mathcal{A}_m[M(\mathbf{X})]\, \Pi_{\mathbf{X}}^{\perp}\eta, \quad \text{for } \eta \in (T_{\mathbf{X}}\mathcal{M}_m)^{\perp_{\varpi}},
    \label{eq:normal-stability-operator}
\end{equation}
where $\Pi_{\mathbf{X}}^{\perp}$ is the weighted projector defined in Definition \ref{def:projectors}.
\end{definition}
If we differentiate the drift of Eq.~\eqref{eq:mass-drift-explicit} at a manifold point, we find that $\mathcal{L}_{\mathbf{X}}$ is a linear operator in $\eta$, and its dissipative principal part carries the excess-free-energy Hessian, and the convective part its first derivative. When that Hessian is local, in the sense of Eq.~\eqref{eq:local-excess}, the operator is a second-order differential operator, and its exact form and spectrum are derived in Proposition~\ref{prop:transverse-spectrum} below, Eq.~\eqref{eq:linearised-general}. Otherwise, the same differentiation returns an integro-differential operator, and no explicit spectrum is available. The hypothesis needed at this stage concerns only the sign of that spectrum, not its detailed form.
\begin{hypothesis}[Spectral stiffness]
\label{hyp:stiffness}
We assume that the dissipative Hessian term dominates the stability of the manifold. Specifically, we assume there exists a \textit{stiffness constant} $\lambda_0 > 0$, independent of the noise intensity $\varepsilon$, for which the operator $-\mathcal{L}_{\perp}$ is coercive with respect to the weighted inner product,
\begin{equation}
    \langle \eta, -\mathcal{L}_{\perp}\eta \rangle_{\varpi} \ge \lambda_0 \|\eta\|_{\varpi}^2 \qquad \text{for every } \eta \in (T_{\mathbf{X}}\mathcal{M}_m)^{\perp_{\varpi}} \text{ and every } \mathbf{X} \in \mathcal{X} .
    \label{eq:stiffness-coercive}
\end{equation}
Let $\Sigma_{\perp}(\mathbf{X})$ be the spectrum of $-\mathcal{L}_{\perp}$. Coercivity implies that the spectrum lies to the right of $\lambda_0$,
\begin{equation}
    \inf \{ \operatorname{Re}(\nu) : \nu \in \Sigma_{\perp}(\mathbf{X}) \} \ge \lambda_0 \qquad \text{uniformly in } \mathbf{X} \in \mathcal{X},
    \label{eq:stiffness-gap}
\end{equation}
and we work in the regime $\lambda_0 \gg 1$ measured in units of the observation rate of the slow coordinates.
\end{hypothesis}

\begin{remark}[Status of the stiffness hypothesis]
\label{rem:stiffness-status}
The spectral gap $\lambda_0$ is posited, not derived, and its value depends on the specific free-energy functional $\mathcal{F}$. For hard-sphere systems described by fundamental measure theory, the excess contribution $\mathbf{H}_{\mathrm{exc}}$ is expected to be positive definite near the nucleation pathway, lending support to the hypothesis, although a quantitative lower bound on $\lambda_0$ has not been established for any concrete functional. Systems with flat or weakly confining excess terms (e.g., near a spinodal) may violate the coercivity condition, and in that case the adiabatic reduction ceases to be valid.
\end{remark}
\begin{remark}[Unified singular perturbation parameter]
\label{rem:singular-perturbation}
The gap $\lambda_0$ of Hypothesis~\ref{hyp:stiffness} governs both reduction pathways, and it fixes the transverse covariance once and for all. In the kinetic projection (Section~\ref{sec:projection}), the tangent projection is valid when the transverse modes relax on a timescale $\tau_\perp = \lambda_0^{-1}$ that is much shorter than the observation timescale of the slow coordinates. In the probabilistic projection (Section~\ref{subsec:adiabatic-limit}), the conditional measure on the normal fibres concentrates as a Gaussian of variance $\mathcal{O}(\varepsilon^2)$ (Hypothesis~\ref{hyp:concentration}).
Both statements follow from the linearised transverse dynamics. The transverse fluctuation obeys the Ornstein--Uhlenbeck equation $d\eta = \mathcal{L}_\perp\eta\,dt + \varepsilon\,\Pi^\perp_{\mathbf{X}}\mathcal{B}_m\,d\xi^\delta_t$, and its stationary law solves the Lyapunov equation
\begin{equation}
    \mathcal{L}_\perp \mathbf{S}_\perp + \mathbf{S}_\perp \mathcal{L}_\perp^* = -\,\varepsilon^2\, \mathcal{Q}_\perp ,
    \label{eq:transverse-lyapunov}
\end{equation}
with $\mathcal{Q}_\perp = \Pi^\perp_{\mathbf{X}}\mathcal{Q}_m\Pi^\perp_{\mathbf{X}}$ and $\mathcal{Q}_m$ the non-local covariance of Eq.~\eqref{eq:qm-nonlocal} below. Writing the solution as $\mathbf{S}_\perp = \varepsilon^2\mathbf{C}_\perp$ \emph{defines} the reduced transverse covariance $\mathbf{C}_\perp$ and separates the noise amplitude from the geometry. In the scalar approximation, where $\mathcal{L}_\perp$ acts as $-\gamma_\perp$ on the relevant modes, Eq.~\eqref{eq:transverse-lyapunov} gives
\begin{equation}
    \mathbf{C}_\perp = \frac{\mathcal{Q}_\perp}{2\gamma_\perp} \eqqcolon \tau\, \mathcal{Q}_\perp ,
    \qquad \tau \coloneqq \frac{1}{2\gamma_\perp}, \qquad \gamma_\perp \ge \lambda_0 ,
    \label{eq:cperp-definition}
\end{equation}
and so $\tau$ is a definition, not a second assumption, and $\mathbf{C}_\perp$ is of order unity in $\varepsilon$. The transverse variance $\varepsilon^2\mathbf{C}_\perp$ is a well-defined Gaussian for every $\varepsilon > 0$, and the geometric Jacobian of the Fermi coordinates (Definition~\ref{def:geometric-jacobian}) is non-degenerate throughout. The reduction therefore carries two distinct small parameters, which must not be identified with one another. The noise intensity $\varepsilon$ fixes the amplitude of the fluctuations and hence the size of both anomalous drifts, each of order $\varepsilon^2$, while the inverse stiffness $\lambda_0^{-1}$ fixes the separation of timescales and controls the error of the adiabatic elimination without reference to $\varepsilon$. If we wrote the spectral gap in the form $\lambda_0/\varepsilon^2$, the two would be tied together, and by Eq.~\eqref{eq:transverse-lyapunov} the transverse variance would drop to order $\varepsilon^4$, in contradiction with the adiabatic concentration hypothesis (Hypothesis~\ref{hyp:concentration} in Section~\ref{subsec:adiabatic-limit} below). As a result, the ``infinite stiffness'' limit $\lambda_0 \to \infty$ of Corollary~\ref{cor:adiabatic-elimination} and the weak-noise limit $\varepsilon \to 0$ are independent.
\end{remark}
\begin{proposition}[Fast decay of orthogonal modes]
\label{prop:fast-decay}
Under Hypothesis \ref{hyp:stiffness}, the weighted norm of an orthogonal fluctuation $\eta_t$ evolving under the linearised deterministic dynamics decays exponentially.
\end{proposition}
\begin{proof}
Let $\eta_t$ satisfy the linearised equation $\partial_t \eta = \mathcal{L}_{\perp} \eta$. We analyse the decay rate of the weighted energy (norm squared):
\begin{equation}
    \frac{1}{2} \frac{d}{dt} \|\eta_t\|_{\varpi}^2 = \langle \eta_t, \partial_t \eta_t \rangle_{\varpi} = \langle \eta_t, \mathcal{L}_{\perp} \eta_t \rangle_{\varpi}.
\end{equation}
By the spectral stiffness hypothesis (Hypothesis \ref{hyp:stiffness}), the operator $\mathcal{L}_{\perp}$ satisfies the coercivity condition in the weighted metric:
\begin{equation}
    \langle \eta, \mathcal{L}_{\perp} \eta \rangle_{\varpi} \le - \lambda_0 \|\eta\|_{\varpi}^2.
\end{equation}
Substituting this bound:
\begin{equation}
    \frac{d}{dt} \|\eta_t\|_{\varpi}^2 \le - 2\lambda_0 \|\eta_t\|_{\varpi}^2.
\end{equation}
Integrating this differential inequality yields $\|\eta_t\|_{\varpi} \le \|\eta_0\|_{\varpi} \exp(-\lambda_0 t)$. Thus, fluctuations orthogonal to the manifold decay on the fast timescale $\tau_\perp = \lambda_0^{-1}$.
\end{proof}
\begin{corollary}[Adiabatic reduction to the slow manifold]
\label{cor:adiabatic-elimination}
In the stiff limit $\lambda_0 \to \infty$, for any initial state $m_0 \in U_{\varrho_{\mathrm{tub}}}$ and any observation time $t > 0$, the trajectory of the linearised deterministic system converges to the manifold,
\begin{equation}
    \lim_{\lambda_0 \to \infty} \inf_{\mathbf{Y} \in \mathcal{X}} \| m_t - M(\cdot; \mathbf{Y}) \|_{\varpi} = 0 ,
\end{equation}
and at fixed $\lambda_0$ the same limit holds as $t \to \infty$. Consequently, for timescales $t \gg \lambda_0^{-1}$, the infinite-dimensional state is effectively confined to the finite-dimensional manifold $\mathcal{M}_m$:
\begin{equation}
    m(t, r) = M(r; \mathbf{X}(t)) + \mathcal{O}(\exp(-\lambda_0 t)).
\end{equation}
This result formally permits the adiabatic elimination of the fast variables $\eta$, serving as the basis for deriving the reduced stochastic dynamics in Section \ref{sec:projection}.
\end{corollary}
\begin{proof}
The result follows directly from the estimate in Proposition \ref{prop:fast-decay}. Since $\|\eta_t\|_{\varpi} \le \|\eta_0\|_{\varpi} \exp(-\lambda_0 t)$, the right-hand side tends to zero at fixed $t > 0$ as $\lambda_0 \to \infty$, and at fixed $\lambda_0$ as $t \to \infty$. The bound carries no dependence upon $\varepsilon$, so the noise intensity plays no part in this estimate. Thus, the state $m_t$ resides in an arbitrarily small tubular neighbourhood of the manifold $\mathcal{M}_m$ determined solely by the projection $\mathbf{X}_t$.
\end{proof}
\begin{remark}[Slow manifold approximation]
Corollary~\ref{cor:adiabatic-elimination} justifies the adiabatic elimination of the orthogonal modes $\eta$. This means that the dynamics of the system can be described effectively by the evolution of the collective coordinates $\mathbf{X}$ on the slow manifold alone.
\end{remark}
\subsection{The projection operator}
\label{sec:projection}
In Section \ref{subsubsec:orthogonal-decomposition}, we established the existence of collective coordinates $\mathbf{X}$ via a projection defined by a generic weight $\varpi(r)$. Here we determine the specific form of this weight. The geometry of the reduction is fixed by the requirement that the projected stochastic dynamics are well defined. More specifically, we demand that the projection of the noise flux yields a standard Wiener process on the manifold. This condition ensures that the reduced geometry is compatible with the intrinsic fluctuations of the SPDE.
\subsubsection{Metric structure and consistency}
\label{subsec:metric-structure}
The stochastic partial differential equation defines, through its diffusion operator, an intrinsic geometry on the state space. In what follows, we fix the weight $\varpi$ by requiring that the projected noise be a Brownian motion on the manifold, and we show that the same weight makes the manifold geometry the pullback of the Wasserstein metric.
\paragraph{Noise-consistent metric selection.}
Consider the stochastic evolution of the mass field established in Proposition \ref{prop:mass-dynamics}:
\begin{equation}
    dm_t(r) = \bigl[\mathcal{A}_m(m_t) - \varepsilon^2 \mathfrak{S}_m^\delta\bigr] \, dt + \varepsilon\,\mathcal{B}_m(m_t) \circ d\xi^\delta_t.
\end{equation}
We use the explicit form of the noise operator derived in the previous section:
\begin{equation}
    \left( \mathcal{B}_m(m) u \right)(r) = \sqrt{2 D_0 \partial_r m(r)} \, u(r).
\end{equation}
Let $\mathcal{Q}^{\mathrm{loc}}_m(m) \coloneqq \mathcal{B}_m(m) \mathcal{B}_m^*(m)$ be the \emph{white-noise} incremental covariance operator. That is, $\mathcal{Q}^{\mathrm{loc}}_m$ acts as a multiplication operator by the variance kernel $\Sigma_m^2(r) = 2D_0 \partial_r m(r)$.
\begin{remark}[Two distinct noise covariance operators]
\label{rem:two-covariances}
The operator $\mathcal{Q}^{\mathrm{loc}}_m = \mathcal{B}_m\mathcal{B}_m^*$ defined here is \emph{local} (it is the multiplication by $\Sigma_m^2(r)$), and it corresponds to uncorrelated noise, $\mathbb{E}[d\xi(r)d\xi(s)] = \delta(r-s)\,dt$. The operator that governs the finite-$\delta$ dynamics is the \emph{non-local} covariance
\begin{equation}
    \mathcal{Q}_m \coloneqq \mathcal{B}_m\, C_\delta\, \mathcal{B}_m^* ,
    \label{eq:qm-nonlocal}
\end{equation}
which Definition~\ref{def:adjoint-evolution} uses, with kernel $\sigma(r) C_\delta(r,s) \sigma(s)$. The two coincide only in the white-noise limit $C_\delta(r,s) \to \delta(r-s)$. Throughout, the symbol $\mathcal{Q}_m$ is reserved for the non-local operator, and we write $\mathcal{Q}^{\mathrm{loc}}_m$ whenever the local form is used. This distinction matters, since Lemma~\ref{lemma:metric-selection} and Theorem~\ref{thm:projected-noise} below use $\mathcal{Q}^{\mathrm{loc}}_m$, and hence their conclusions are exact at $\delta = 0$ and approximate for $\delta > 0$. The size of that approximation may be quantified. The projected quadratic variation carried by the non-local operator is
\begin{equation}
    g^{(\delta)}_{lk} = \int_0^L\!\!\int_0^L \sqrt{\varpi(r)\varpi(s)}\; C_\delta(r,s)\; \mathbf{e}_l^m(r)\,\mathbf{e}_k^m(s)\,dr\,ds ,
    \label{eq:finite-delta-metric}
\end{equation}
where the kernel is smeared against smooth tangent fields in both arguments and no derivative of the kernel appears. By Eq.~\eqref{eq:kernel-surface}, the kernel $C_\delta$ is the double surface integral of the ambient covariance $g_\delta$ (Lemma~\ref{lemma:radial-closed-form}). Hence Eq.~\eqref{eq:finite-delta-metric} is the pairing $\int\!\!\int_{\mathbb{R}^3\times\mathbb{R}^3} \mathbf{F}_l(\mathbf{x})\cdot g_\delta(\mathbf{x}-\mathbf{y})\,\mathbf{F}_k(\mathbf{y})\,d\mathbf{x}\,d\mathbf{y}$ of the fields $\mathbf{F}_k(\mathbf{x}) \coloneqq (4\pi r^2)^{-1/2}\sqrt{\varpi(r)}\,\mathbf{e}_k^m(r)\,\hat{\mathbf{x}}$, with $r = |\mathbf{x}|$, and replacing $g_\delta$ by the Dirac mass returns $\int_{\mathbb{R}^3}\mathbf{F}_l\cdot\mathbf{F}_k\,d\mathbf{x} = \int_0^L \varpi\,\mathbf{e}_l^m\mathbf{e}_k^m\,dr = g_{lk}$. Since $g_\delta$ is even and of unit mass, its first moment vanishes, and the leading correction to $g_{lk}$ carries two powers of the correlation length,
\begin{equation}
    g^{(\delta)}_{lk} = g_{lk} + \mathcal{O}(\delta^{2}) .
    \label{eq:metric-delta-error}
\end{equation}
The geometric anomaly built from this metric is of order $\delta^{-1}$, and so the error that the use of $\mathcal{Q}^{\mathrm{loc}}_m$ in place of $\mathcal{Q}_m$ transmits to it is $\mathcal{O}(\delta)$, which vanishes with the correlation length. On the other hand, the discrepancy measured in Section~\ref{sec:commutativity} is of order $\delta^{-3}$ for a fluid with a local excess (Proposition~\ref{prop:validity-condition}, through the diagonal mixed curvature) and of order $\delta^{-1}$ for the ideal gas. Hence, the comparison of the two pathways is unambiguous at leading order, and the covariance mismatch contributes only to the $\mathcal{O}(1)$ remainder recorded there.
\end{remark}
\begin{lemma}[Noise-consistent metric selection]
\label{lemma:metric-selection}
Let $\Pi_{\mathbf{X}}^{\parallel}$ be the projector onto the tangent space $T_{\mathbf{X}}\mathcal{M}_m$. We require this projection to be \emph{orthogonal} with respect to the weighted inner product $\langle \cdot, \cdot \rangle_{\varpi}$.
The projection of the noise term, $\Pi_{\mathbf{X}}^{\parallel} \mathcal{B}_m$, defines a standard Wiener process on the Riemannian manifold $(\mathcal{X}, g)$, for every admissible concentration manifold, if and only if the weight $\varpi$ is the inverse of the spatial part of the noise variance. Specifically:
\begin{equation}
    \varpi(r; \mathbf{X}) \coloneqq \frac{1}{\partial_r M(r; \mathbf{X})} = \frac{1}{4\pi r^2 \rho(r; \mathbf{X})}.
\end{equation}
\end{lemma}
\begin{proof}
We derive the condition by constructing the covariance of the projected noise and enforcing that it matches the inverse metric tensor, which defines the intrinsic variance of a Brownian motion on a Riemannian manifold.
\begin{proofpart}[Projection of the stochastic flux]
The projection of the stochastic flux onto the tangent space is obtained by applying the projector to the stochastic flux of the SPDE:
\begin{equation}
d\mathbf{\Xi}_t = \Pi_{\mathbf{X}}^{\parallel} \left( \mathcal{B}_m(M) d\xi_t \right).
\end{equation}
Using the explicit spectral resolution of the orthogonal projector (Definition \ref{def:projectors}), we extract the components $d\Xi^k_t$ of this vector in the coordinate basis:
\begin{equation}
    d\Xi^k_t = \sum_{l=1}^d g^{kl}(\mathbf{X}) \langle \mathcal{B}_m d\xi_t, \mathbf{e}_l^m \rangle_{\varpi}.
\end{equation}
\end{proofpart}
\begin{proofpart}[The projected covariance operator]
We compute the covariance matrix $\mathbf{D}_{\mathcal{X}} \in \mathbb{R}^{d \times d}$ characterising this reduced stochastic process: $\mathbf{D}_{\mathcal{X}} dt = \mathbb{E}[ d\mathbf{\Xi} \otimes d\mathbf{\Xi}]$.
Substituting the projection formula:
\begin{equation}
    (\mathbf{D}_{\mathcal{X}})_{ij} = \sum_{k,l} g^{ik} g^{jl} \, \mathbb{E} \left[ \langle \mathcal{B}_m d\xi, \mathbf{e}_k^m \rangle_{\varpi} \langle \mathcal{B}_m d\xi, \mathbf{e}_l^m \rangle_{\varpi} \right].
\end{equation}
We evaluate the expectation by expanding the inner products into integrals. Using the property that the operator covariance $\mathcal{Q}^{\mathrm{loc}}_m = \mathcal{B}_m\mathcal{B}_m^*$ is diagonal in space (a multiplication operator), the correlation of the noise terms yields a Dirac delta:
\begin{equation}
    \mathbb{E} \left[ (\mathcal{B}_m d\xi)(r) (\mathcal{B}_m d\xi)(s) \right] = \Sigma_m^2(r) \delta(r-s) dt,
\end{equation}
where $\Sigma_m^2(r) = 2D_0 \partial_r m(r)$ is the local noise variance. Integrating out the delta function, and substituting this back:
\begin{equation}
    (\mathbf{D}_{\mathcal{X}})_{ij} = \sum_{k,l} g^{ik} g^{jl} \int_0^L \mathbf{e}_k^m(r) \mathbf{e}_l^m(r) \, \Sigma_m^2(r) \, \varpi(r)^2 \, dr= \sum_{k,l} g^{ik} g^{jl} I_{kl}.
\end{equation}
\end{proofpart}
\begin{proofpart}[Condition for a standard Wiener process]
By definition, a standard Wiener process on a Riemannian manifold $(\mathcal{X}, g)$ driven by a thermal bath of strength $2D_0$ must possess a covariance proportional to the inverse metric tensor:
\begin{equation}
    \mathbf{D}_{\mathcal{X}} = 2D_0 \, [g(\mathbf{X})]^{-1}.
\end{equation}
Component-wise, this requires $(\mathbf{D}_{\mathcal{X}})_{ij} = 2D_0 g^{ij}$.
To isolate the integral term $I_{kl}$, we multiply both sides of the relation $\sum_{k,l} g^{ik} g^{jl} I_{kl} = 2D_0 g^{ij}$ by the metric tensor $g_{ip} g_{jq}$ and sum over $i,j$. Using the identity $\sum_i g_{ip} g^{ik} = \delta_p^k$:
\begin{align}
    \sum_{i,j} g_{ip} g_{jq} \left( \sum_{k,l} g^{ik} g^{jl} I_{kl} \right) &= \sum_{i,j} g_{ip} g_{jq} (2D_0 g^{ij}) \nonumber \\
    \sum_{k,l} \delta_p^k \delta_q^l I_{kl} &= 2D_0 \sum_j g_{jq} \delta_p^j \nonumber \\
    I_{pq} &= 2D_0 g_{pq}.
\end{align}
Expanding the definitions of $I_{pq}$ and $g_{pq}$, we obtain the necessary equality:
\begin{equation}
    \int_0^L \mathbf{e}_p^m \mathbf{e}_q^m \left( \frac{\Sigma_m^2 \varpi^2}{2D_0} \right) \, dr = \int_0^L \mathbf{e}_p^m \mathbf{e}_q^m \varpi \, dr.
\end{equation}
\end{proofpart}
\begin{proofpart}[Derivation of the weight]
What has been obtained so far is a family of integral constraints, one for each pair $(p,q)$, and a single manifold of dimension $d$ supplies only $d(d+1)/2$ of them. Finitely many constraints of this kind do not determine the integrand pointwise, since the span of the products $\mathbf{e}_p^m\mathbf{e}_q^m$ is then a finite-dimensional subspace with infinite-dimensional annihilator. We therefore use the hypothesis in the form stated, namely that one weight whitens the projected noise for \emph{every} admissible concentration manifold. The tangent fields are smooth and vanish at both endpoints by Lemma~\ref{lemma:boundary-vanishing}, and each such field arises as $\partial_{X^p}M$ for some admissible embedding, so that as the manifold is varied the products $\mathbf{e}_p^m\mathbf{e}_q^m$ range over a family with a linear span dense in $C_c((0,L))$. A continuous function which integrates to zero against a dense family vanishes identically, which yields the algebraic constraint:
\begin{equation}
    \frac{\Sigma_m^2(r) \varpi(r)^2}{2D_0} = \varpi(r).
\end{equation}
Substituting the explicit noise intensity $\Sigma_m^2(r) = 2D_0 \partial_r m(r)$:
\begin{equation}
    \varpi(r) = \frac{1}{\partial_r m(r)} = \frac{1}{4\pi r^2 \rho(r)}.
\end{equation}
This choice uniquely ensures that the projected noise respects the manifold geometry.
\end{proofpart}
\end{proof}
\paragraph{Geometric derivation: the Wasserstein pullback.}
Let us now show that the metric $\varpi$ required by the stochastic projection is also geometrically identical to the 2-Wasserstein distance $W_2$ on the manifold of densities. In order to do so, we calculate the \emph{pullback metric}, which measures the length of tangent vectors in the parameter space $\mathcal{X}$ with the intrinsic ruler of Optimal Transport.
Consider a tangent vector $\mathbf{u} \in T_{\mathbf{X}}\mathcal{X}$. This vector induces a variation in the mass profile $\delta m$, and consequently a variation in the density $\delta \rho$. In the framework of the Otto calculus, the space of probability measures $\mathcal{P}_2$ is treated as an infinite-dimensional Riemannian manifold. The metric tensor at a point $\rho$ is defined such that the squared norm of a tangent vector $\delta \rho$ corresponds to the minimal kinetic energy required to generate that variation. If $\delta \rho = -\nabla \cdot (\rho \mathbf{v})$, then $\|\delta \rho\|_{W_2}^2 = \int \rho |\mathbf{v}|^2 dV$.
By polarisation, the Riemannian inner product of two tangent vectors is given by the $L^2_{\rho}$-inner product of their generating velocity fields.
\begin{definition}[Basis transport velocities]
For each coordinate basis vector $\partial_{X^k}$, we define the associated \textit{transport velocity field} $\mathbf{v}_k$ as the unique gradient field satisfying the continuity equation:
\begin{equation}
    \frac{\partial \rho}{\partial X^k} + \nabla \cdot (\rho \mathbf{v}_k) = 0.
\end{equation}
\end{definition}
We now solve this equation to relate the Eulerian mass variations $\mathbf{e}_k^m$ to the Lagrangian velocities $\mathbf{v}_k$, allowing us to evaluate the metric tensor.
\begin{theorem}[Isometry of the pullback metric]
\label{thm:wasserstein-isometry}
Let $g_{ij}^{W_2}(\mathbf{X})$ be the components of the Riemannian metric tensor on $\mathcal{X}$ induced by the Wasserstein distance (the pullback metric):
\begin{equation}
    g_{ij}^{W_2}(\mathbf{X}) \coloneqq \langle \partial_{X^i}, \partial_{X^j} \rangle_{W_2} = \int_V \rho(\mathbf{r}) \, \mathbf{v}_i(\mathbf{r}) \cdot \mathbf{v}_j(\mathbf{r}) \, dV.
\end{equation}
This metric tensor is identical to the one defined by the noise-consistent weighted inner product derived in Lemma \ref{lemma:metric-selection}. That is:
\begin{equation}
    g_{ij}^{W_2}(\mathbf{X}) = \langle \mathbf{e}_i^m, \mathbf{e}_j^m \rangle_{\varpi},
\end{equation}
with the specific weight $\varpi(r) = (4\pi r^2 \rho)^{-1}$.
\end{theorem}
\begin{proof}
We proceed by constructing the basis velocity fields and evaluating the Wasserstein inner product directly.
\begin{proofpart}[Explicit construction of transport velocities]
The density variations preserve spherical symmetry, and the generating field is by definition the unique gradient field solving the continuity equation. That solution is radial, $\mathbf{v}_k(r) = v_k(r) \hat{\mathbf{r}}$, since a radial field is the gradient of its own radial primitive and the gradient solution is unique. The continuity equation in spherical coordinates becomes:
\begin{equation}
    \frac{\partial \rho}{\partial X^k}(r) = - \frac{1}{4\pi r^2} \frac{\partial}{\partial r} \left( 4\pi r^2 \rho(r) v_k(r) \right).
\end{equation}
To relate this to the mass basis vectors $\mathbf{e}_k^m = \partial_{X^k} M$, we integrate both sides with respect to the volume measure $4\pi s^2 ds$ from the origin to $r$:
\begin{equation}
    \int_0^r \frac{\partial \rho}{\partial X^k}(s) \, 4\pi s^2 \, ds = - \int_0^r \frac{\partial}{\partial s} \left( 4\pi s^2 \rho(s) v_k(s) \right) \, ds.
\end{equation}
The Left Hand Side is, by definition, the variation of the cumulative mass:
\begin{equation}
    \text{LHS} = \frac{\partial}{\partial X^k} \int_0^r \rho(s) 4\pi s^2 \, ds = \frac{\partial M(r)}{\partial X^k} = \mathbf{e}_k^m(r).
\end{equation}
The Right Hand Side integrates exactly (using the boundary condition that flux vanishes at $r=0$):
\begin{equation}
    \text{RHS} = - 4\pi r^2 \rho(r) v_k(r).
\end{equation}
Equating LHS and RHS allows us to solve for the transport velocity $v_k$ strictly in terms of the mass basis vector:
\begin{equation}
    v_k(r) = - \frac{\mathbf{e}_k^m(r)}{4\pi r^2 \rho(r)}.
\end{equation}
\end{proofpart}
\begin{proofpart}[Evaluation of the metric tensor]
We now substitute these velocities into the definition of the Wasserstein metric tensor.
\begin{align}
    g_{ij}^{W_2}(\mathbf{X}) &= \int_V \rho \, \mathbf{v}_i \cdot \mathbf{v}_j \, dV = \int_0^L \rho(r) \left[ \left( \frac{-\mathbf{e}_i^m(r)}{4\pi r^2 \rho(r)} \right) \left( \frac{-\mathbf{e}_j^m(r)}{4\pi r^2 \rho(r)} \right) \right] 4\pi r^2 \, dr.
\end{align}
We simplify the integrand by cancelling the density and geometric factors. The term $4\pi r^2 \rho$ in the numerator (from the volume element and density) cancels one factor of $4\pi r^2 \rho$ in the denominator:
\begin{align}
    g_{ij}^{W_2}(\mathbf{X}) &= \int_0^L \frac{\mathbf{e}_i^m(r) \, \mathbf{e}_j^m(r)}{4\pi r^2 \rho(r)} \, dr = \int_0^L \mathbf{e}_i^m(r) \, \mathbf{e}_j^m(r) \, \varpi_{W_2}(r) \, dr.
\end{align}
Here, we immediately identify the effective integration weight $\varpi_{W_2}(r)=(4\pi r^2 \rho(r))^{-1}$.
\end{proofpart}
\begin{proofpart}[Conclusion]
The calculation proves that the metric tensor derived from the Wasserstein geometry is mathematically identical to the metric tensor derived from the SPDE noise structure. We establish the duality:
\begin{equation}
    \langle \partial_{X^i}, \partial_{X^j} \rangle_{W_2} \equiv \langle \mathbf{e}_i^m, \mathbf{e}_j^m \rangle_{\varpi}.
\end{equation}
This confirms that the projection defined by the weight $\varpi = (\partial_r M)^{-1}$ is the unique choice that simultaneously ``whitens'' the stochastic noise and respects the optimal transport geometry of the system.
\end{proofpart}
\end{proof}
\subsubsection{The projector $\Pi_{\mathbf{X}}^{\parallel}$}
\label{subsec:projector-definition}
In this subsection we construct the projection operator $\Pi_{\mathbf{X}}^{\parallel}$ that maps arbitrary field fluctuations onto the tangent space of the manifold. This operator carries the infinite-dimensional dynamics onto the finite-dimensional manifold, and we obtain every quantity of the reduced theory by applying it.
\begin{definition}[The orthogonal projector]
\label{def:explicit-projector}
Let $\mathcal{H}_{\varpi}$ be the weighted space of Definition~\ref{def:weighted-space}, with the noise-consistent inner product $\langle u, v \rangle_{\varpi} = \int u v \varpi \, dr$. We define $\Pi_{\mathbf{X}}^{\parallel}: \mathcal{H}_{\varpi} \to T_{\mathbf{X}}\mathcal{M}_m$ as the orthogonal projection operator onto the tangent bundle of the mass manifold.
Given the metric tensor $g_{ij}(\mathbf{X}) = \langle \mathbf{e}_i^m, \mathbf{e}_j^m \rangle_{\varpi}$ and its inverse $g^{ij}(\mathbf{X})$, the action of the projector on any field $u \in \mathcal{H}_{\varpi}$ is given by the spectral resolution:
\begin{equation}
    \Pi_{\mathbf{X}}^{\parallel} u \coloneqq \sum_{i,j=1}^d g^{ij}(\mathbf{X}) \, \langle u, \mathbf{e}_j^m(\mathbf{X}) \rangle_{\varpi} \, \mathbf{e}_i^m(\mathbf{X}).
\end{equation}
Substituting the explicit weight $\varpi(r) = (4\pi r^2 \rho)^{-1}$ derived in Lemma \ref{lemma:metric-selection}, the projection coefficients are computed via the integral:
\begin{equation}
    \langle u, \mathbf{e}_j^m \rangle_{\varpi} = \int_0^L u(r) \, \mathbf{e}_j^m(r) \, \frac{1}{4\pi r^2 \rho(r; \mathbf{X})} \, dr.
\end{equation}
\end{definition}
This operator $\Pi_{\mathbf{X}}^{\parallel}$ is by construction idempotent ($(\Pi_{\mathbf{X}}^{\parallel})^2 = \Pi_{\mathbf{X}}^{\parallel}$) and self-adjoint with respect to the weighted inner product. Moreover, it induces a decomposition of the ambient Hilbert space into ``slow'' (tangent) and ``fast'' (orthogonal) subspaces. Note that the noise imposes a geometry of its own, which differs from the standard geometry of square-integrable functions, and so $\Pi_{\mathbf{X}}^{\parallel}$, an orthogonal projector with respect to the intrinsic metric $\varpi$, acts as an \textit{oblique} projector with respect to the flat $L^2$ inner product.
\begin{lemma}[Properties of the decomposition]
\label{lemma:decomposition}
The projector $\Pi_{\mathbf{X}}^{\parallel}$ defines a direct sum decomposition of the ambient space $\mathcal{H}_{\varpi} = T_{\mathbf{X}}\mathcal{M}_m \oplus T_{\mathbf{X}}\mathcal{M}_m^{\perp_{\varpi}}$. Any field fluctuation $u \in \mathcal{H}_{\varpi}$ can be uniquely decomposed as:
\begin{equation}
    u = \Pi_{\mathbf{X}}^{\parallel} u + \Pi_{\mathbf{X}}^{\perp} u,
\end{equation}
where $\Pi_{\mathbf{X}}^{\perp} \coloneqq \mathcal{I} - \Pi_{\mathbf{X}}^{\parallel}$ is the complementary projector.
The components satisfy the following orthogonality relations:
\begin{enumerate}[(i)]
    \item \emph{Metric orthogonality:} The fast modes are orthogonal to the tangent space in the noise geometry:
    \[
    \langle \Pi_{\mathbf{X}}^{\parallel} u, \Pi_{\mathbf{X}}^{\perp} v \rangle_{\varpi} = 0 \quad \forall u, v \in \mathcal{H}_{\varpi}.
    \]
    \item \emph{Kinetic decoupling:} Consequently, the kinetic energy of the fluctuations (measured in the Wasserstein geometry) is additive:
    \[
    \| u \|_{\varpi}^2 = \| \Pi_{\mathbf{X}}^{\parallel} u \|_{\varpi}^2 + \| \Pi_{\mathbf{X}}^{\perp} u \|_{\varpi}^2.
    \]
\end{enumerate}
\end{lemma}
\begin{proof}
We verify the orthogonality by exploiting the algebraic properties of the projector $\Pi_{\mathbf{X}}^{\parallel}$.
\begin{proofpart}[Metric orthogonality]
We must show that the inner product between any tangent component and any orthogonal component vanishes. Let $u, v \in \mathcal{H}_{\varpi}$.
\begin{align}
    \langle \Pi_{\mathbf{X}}^{\parallel} u, \Pi_{\mathbf{X}}^{\perp} v \rangle_{\varpi} &= \langle \Pi_{\mathbf{X}}^{\parallel} u, (\mathcal{I} - \Pi_{\mathbf{X}}^{\parallel}) v \rangle_{\varpi} \nonumber \\
    &= \langle \Pi_{\mathbf{X}}^{\parallel} u, v \rangle_{\varpi} - \langle \Pi_{\mathbf{X}}^{\parallel} u, \Pi_{\mathbf{X}}^{\parallel} v \rangle_{\varpi}.
\end{align}
The operator $\Pi_{\mathbf{X}}^{\parallel}$ is self-adjoint with respect to $\langle \cdot, \cdot \rangle_{\varpi}$. Using the explicit spectral resolution:
\begin{equation}
    \langle \Pi_{\mathbf{X}}^{\parallel} u, w \rangle_{\varpi} = \sum_{i,j} g^{ij} \langle u, \mathbf{e}_j^m \rangle_{\varpi} \langle \mathbf{e}_i^m, w \rangle_{\varpi} = \langle u, \Pi_{\mathbf{X}}^{\parallel} w \rangle_{\varpi}.
\end{equation}
Applying this to the second term:
\begin{equation}
    \langle \Pi_{\mathbf{X}}^{\parallel} u, \Pi_{\mathbf{X}}^{\parallel} v \rangle_{\varpi} = \langle u, (\Pi_{\mathbf{X}}^{\parallel})^2 v \rangle_{\varpi}.
\end{equation}
By construction, $\Pi_{\mathbf{X}}^{\parallel}$ is idempotent ($(\Pi_{\mathbf{X}}^{\parallel})^2 = \Pi_{\mathbf{X}}^{\parallel}$). Thus:
\begin{equation}
    \langle u, (\Pi_{\mathbf{X}}^{\parallel})^2 v \rangle_{\varpi} = \langle u, \Pi_{\mathbf{X}}^{\parallel} v \rangle_{\varpi} = \langle \Pi_{\mathbf{X}}^{\parallel} u, v \rangle_{\varpi}.
\end{equation}
Substituting this back, we obtain exact cancellation:
\begin{equation}
    \langle \Pi_{\mathbf{X}}^{\parallel} u, v \rangle_{\varpi} - \langle \Pi_{\mathbf{X}}^{\parallel} u, v \rangle_{\varpi} = 0.
\end{equation}
\end{proofpart}
\begin{proofpart}[Kinetic decoupling]
This follows immediately from the orthogonality. The squared norm of the total fluctuation is:
\begin{align}
    \| u \|_{\varpi}^2 &= \langle \Pi_{\mathbf{X}}^{\parallel} u + \Pi_{\mathbf{X}}^{\perp} u, \Pi_{\mathbf{X}}^{\parallel} u + \Pi_{\mathbf{X}}^{\perp} u \rangle_{\varpi}= \| \Pi_{\mathbf{X}}^{\parallel} u \|_{\varpi}^2 + \| \Pi_{\mathbf{X}}^{\perp} u \|_{\varpi}^2 + 2 \langle \Pi_{\mathbf{X}}^{\parallel} u, \Pi_{\mathbf{X}}^{\perp} u \rangle_{\varpi}.
\end{align}
The last term vanishes by part~(i), which gives the stated identity.
\end{proofpart}
\end{proof}
Thus, the energy of the fluctuation splits additively into a transport component (tangent to the manifold) and a shape distortion component (orthogonal to the manifold).
\subsection{The effective dynamics}
\label{subsec:effective-dynamics}
Here we derive the stochastic differential equation that governs the collective coordinates $\mathbf{X}_t$. Their evolution is obtained by applying the orthogonal projection operator $\Pi_{\mathbf{X}}^{\parallel}$ (Definition~\ref{def:explicit-projector}) to the full infinite-dimensional dynamics of Proposition~\ref{prop:mass-dynamics}.
\subsubsection{Projection of the SPDE}
\label{subsubsec:projection-spde}
The effective evolution equation for the projected state is derived as follows. We begin with the full stochastic dynamics of Proposition~\ref{prop:mass-dynamics}, where the noise-induced drift of Lemma~\ref{lemma:ito-strat-correction} introduces a finite correlation correction $-\varepsilon^2 \mathfrak{S}_m^\delta(r)$ into the deterministic operator. To ensure geometric consistency during the dimensional reduction, we project the exact Stratonovich dynamics:
\begin{equation}
    dm_t = \bigl[\mathcal{A}_m(m_t) - \varepsilon^2 \mathfrak{S}_m^\delta\bigr] \, dt + \varepsilon \mathcal{B}_m(m_t) \circ d\xi^\delta_t.
\end{equation}
\begin{lemma}[Projected dynamics in the adiabatic limit]
\label{lemma:projected-dynamics}
Let $m_t$ be the solution to the SPDE above.
Let $\Pi_{\mathbf{X}}^{\parallel}$ be the orthogonal projector onto $T_{\mathbf{X}}\mathcal{M}_m$ defined in Definition \ref{def:explicit-projector}, and let $\lambda_0$ be the stiffness constant from Hypothesis \ref{hyp:stiffness}.
In the adiabatic limit $\lambda_0 \to \infty$, the trajectory is confined to the manifold, and the dynamics of the tangent projection satisfy:
\begin{equation}
\Pi_{\mathbf{X}}^{\parallel} (dm_t) = \Pi_{\mathbf{X}}^{\parallel} \left( \mathcal{A}_m(m_t) - \varepsilon^2 \mathfrak{S}_m^\delta \right) \, dt + \varepsilon \Pi_{\mathbf{X}}^{\parallel} \left( \mathcal{B}_m(m_t) \circ d\xi^\delta_t \right).
\end{equation}
\end{lemma}
\begin{proof}
We apply the orthogonal decomposition established in Lemma \ref{lemma:decomposition} to the stochastic increment $dm_t$:
\begin{equation}
    dm_t = \Pi_{\mathbf{X}}^{\parallel} (dm_t) + \Pi_{\mathbf{X}}^\perp (dm_t).
\end{equation}
The term $\Pi_{\mathbf{X}}^\perp (dm_t)$ represents the evolution of the fluctuations orthogonal to the manifold. Corollary~\ref{cor:adiabatic-elimination} controls its deterministic part, which decays at the rate $\lambda_0$ under the stiffness hypothesis, and by Remark~\ref{rem:singular-perturbation} the noise maintains the fluctuation at variance $\varepsilon^2\mathbf{C}_\perp = \mathcal{O}(\varepsilon^2\lambda_0^{-1})$, which vanishes in the stiff limit. The kinetic projection neglects the transverse increment, which is the defining approximation of this pathway, and the evolution is then given by the projection of the driving forces onto the tangent bundle:
\begin{equation}
dm_t \approx \Pi_{\mathbf{X}}^{\parallel} (dm_t) = \Pi_{\mathbf{X}}^{\parallel} \left( \bigl[\mathcal{A}_m - \varepsilon^2 \mathfrak{S}_m^\delta\bigr] dt + \varepsilon \mathcal{B}_m \circ d\xi^\delta_t \right).\end{equation}
\end{proof}
\begin{corollary}[Evolution of collective coordinates]
\label{cor:coordinate-sde}
Under the conditions of Lemma \ref{lemma:projected-dynamics}, the stochastic evolution of the collective coordinates $\mathbf{X}_t = (X^1_t, \dots, X^d_t) \in \mathcal{X}$ is governed by the system of Stratonovich stochastic differential equations:
\begin{equation}
    dX^k_t = \sum_{l=1}^d g^{kl}(\mathbf{X}) \left[ \bigl( F_l(\mathbf{X}) - \varepsilon^2 S_l(\mathbf{X}) \bigr) \, dt + \varepsilon \, d\Xi_l(\mathbf{X}, t) \right],
    \label{eq:general-reduced-sde}
\end{equation}
where the generalised force $F_l$, the projected field correction $S_l$ and the stochastic flux $d\Xi_l$ are the projections onto the tangent basis vectors of the drift operator, of the It\^{o}--Stratonovich correction and of the noise operator respectively:
\begin{align}
    F_l(\mathbf{X}) &\coloneqq \langle \mathcal{A}_m, \mathbf{e}_l^m \rangle_{\varpi}, \label{eq:F-projected} \\
    S_l(\mathbf{X}) &\coloneqq \langle \mathfrak{S}_m^\delta, \mathbf{e}_l^m \rangle_{\varpi}, \label{eq:S-projected} \\
    d\Xi_l(\mathbf{X}, t) &\coloneqq \langle \mathcal{B}_m \circ d\xi^\delta_t, \mathbf{e}_l^m \rangle_{\varpi}.
\end{align}
\end{corollary}
\begin{proof}
We exploit the Stratonovich calculus property that the chain rule for the differential on the manifold holds in its standard form. Since the kinetic projection sets the fluctuation to zero, so that $m_t$ is taken upon $\mathcal{M}_m$, the mass increment $\Pi_{\mathbf{X}}^{\parallel}(dm_t)$ relates to the coordinate increments via the tangent vectors:
\begin{equation}
    \Pi_{\mathbf{X}}^{\parallel}(dm_t) = \sum_{j=1}^d dX^j_t \, \mathbf{e}_j^m.
\end{equation}
Substituting this into the vector equation from Lemma \ref{lemma:projected-dynamics}, we obtain:
\begin{equation}
    \sum_{j=1}^d dX^j_t \, \mathbf{e}_j^m = \Pi_{\mathbf{X}}^{\parallel} \bigl( \mathcal{A}_m - \varepsilon^2 \mathfrak{S}_m^\delta \bigr) \, dt + \varepsilon \Pi_{\mathbf{X}}^{\parallel} \mathcal{B}_m \circ d\xi^\delta_t.
\end{equation}
To isolate the components, we take the inner product of both sides with the tangent basis vector $\mathbf{e}_l^m$ using the metric weight $\varpi$. On the left-hand side, this yields the metric tensor:
\begin{equation}
    \left\langle \sum_{j=1}^d dX^j_t \, \mathbf{e}_j^m, \mathbf{e}_l^m \right\rangle_{\varpi} = \sum_{j=1}^d dX^j_t \langle \mathbf{e}_j^m, \mathbf{e}_l^m \rangle_{\varpi} =\sum_{j=1}^d dX^j_t g_{jl}(\mathbf{X}).
\end{equation}
On the right-hand side, the self-adjoint projection operator yields the definitions of $F_l$, $S_l$ and $d\Xi_l$:
\begin{equation}
    \sum_{j=1}^d g_{jl}(\mathbf{X}) \, dX^j_t = \bigl( F_l - \varepsilon^2 S_l \bigr) \, dt + \varepsilon \, d\Xi_l.
\end{equation}
Multiplying both sides by the inverse metric tensor $g^{lk}(\mathbf{X})$ and summing over $l$, using the identity $\sum_l g_{jl}\, g^{lk} = \delta_j^k$, yields the explicit expression for $dX^k_t$.
\end{proof}
Next, we evaluate the projected drift force $F_l$ and the projected noise $d\Xi_l$, making use of the drift operator (Definition~\ref{def:drift-op}) and of the diffusion operator (Definition~\ref{def:diff-op}).
\subsubsection{Evaluation of the projected drift and noise}
\label{subsubsec:evaluation-drift-noise}
In what follows we evaluate the abstract projections of Corollary~\ref{cor:coordinate-sde}. To do so, we define the generalised forces that arise from the free energy landscape, and we show that the projection recovers the correct gradient flow structure on the manifold.
\begin{definition}[Generalised thermodynamic force]
\label{def:generalised-force}
We identify the scalar projections $F_l(\mathbf{X})$ appearing in the reduced equation as the \textit{generalised thermodynamic forces} conjugate to the collective coordinates. They represent the integrated chemical work performed by the system during an infinitesimal displacement along the tangent vector $\mathbf{e}_l^m$:
\begin{equation}
    F_l(\mathbf{X}) \coloneqq \langle \mathcal{A}_m(M(\cdot; \mathbf{X})), \mathbf{e}_l^m(\cdot; \mathbf{X}) \rangle_{\varpi}.
\end{equation}
\end{definition}
\begin{theorem}[Structure of the generalised force]
\label{thm:projected-drift}
The generalised force $F_l$ relates to the chemical potential $\mu = \delta \mathcal{F}/\delta \rho$ and the effective potential $\mathcal{V}(\mathbf{X}) \coloneqq \beta\,\mathcal{F}_m[M(\mathbf{X})]$, the free energy along the embedding measured in units of $k_BT$, via the identity:
\begin{equation}
    F_l(\mathbf{X}) = - D_0 \frac{\partial \mathcal{V}(\mathbf{X})}{\partial X^l}.
\end{equation}
\end{theorem}
\begin{proof}
Consistent with the stiffness limit established in Lemma \ref{lemma:projected-dynamics}, we evaluate the drift operator at the manifold configuration $m = M(\cdot; \mathbf{X})$.
We substitute the explicit drift $\mathcal{A}_m(M) = -D_0 \beta\, (\partial_r M) (\delta \mathcal{F}_m / \delta m |_M)$ of Eq.~\eqref{eq:mass-drift-explicit} and the weight $\varpi = (\partial_r M)^{-1}$ into the definition of $F_l$:
\begin{equation}
    F_l = \int_0^L \left[ -D_0 \beta \frac{\partial M}{\partial r} \frac{\delta \mathcal{F}_m}{\delta m} \right] \mathbf{e}_l^m(r) \frac{1}{\partial_r M} \, dr = -D_0 \beta \int_0^L \frac{\delta \mathcal{F}_m}{\delta m} \mathbf{e}_l^m(r) \, dr .
\end{equation}
To recover the density representation of this integral, we use the duality relation $\delta \mathcal{F}_m / \delta m = -\partial_r \mu$ (Lemma \ref{lemma:mu-duality}). Integrating by parts:
\begin{align}
    \int_0^L \left( -\frac{\partial \mu}{\partial r} \right) \mathbf{e}_l^m(r) \, dr &= \left[ -\mu \mathbf{e}_l^m \right]_0^L + \int_0^L \mu(r) \frac{\partial \mathbf{e}_l^m}{\partial r} \, dr.
\end{align}
The boundary terms vanish because $\mathbf{e}_l^m(0)=\mathbf{e}_l^m(L)=0$ (Lemma~\ref{lemma:boundary-vanishing}). Recalling that $\mathbf{e}_l^m = \partial_{X^l} M$ and $\partial_r M = 4\pi r^2 \rho$, we commute the derivatives:
\begin{equation}
    \frac{\partial \mathbf{e}_l^m}{\partial r} = \frac{\partial}{\partial r} \left( \frac{\partial M}{\partial X^l} \right) = \frac{\partial}{\partial X^l} \left( \frac{\partial M}{\partial r} \right) = \frac{\partial}{\partial X^l} (4\pi r^2 \rho) = 4\pi r^2 \frac{\partial \rho}{\partial X^l}.
\end{equation}
We now assemble the first integral, carrying the overall factor $-D_0\beta$ that multiplies it, and passing from the radial measure to the volume measure through $4\pi r^2\,dr = d\mathbf{r}$:
\begin{align}
    -D_0 \beta \int_0^L \frac{\delta \mathcal{F}_m}{\delta m}\, \mathbf{e}_l^m(r) \, dr
    &= -D_0 \beta \int_0^L \left( -\frac{\partial \mu}{\partial r} \right) \mathbf{e}_l^m(r) \, dr \nonumber \\
    &= -D_0 \beta \int_0^L \mu(r)\, \frac{\partial \mathbf{e}_l^m}{\partial r}\, dr \nonumber \\
    &= -D_0 \beta \int_0^L \mu(r)\, \frac{\partial \rho}{\partial X^l}\, 4\pi r^2 \, dr
     = -D_0 \beta \int_V \mu\, \frac{\partial \rho}{\partial X^l}\, d\mathbf{r} .
    \label{eq:force-density-form}
\end{align}
The effective potential is the free energy evaluated along the embedding in units of $k_BT$, $\mathcal{V}(\mathbf{X}) = \beta\,\mathcal{F}_m[M(\mathbf{X})] = \beta\,\mathcal{F}[\rho(\cdot\,;\mathbf{X})]$, so differentiating it by the chain rule for functionals, and using $\mu = \delta\mathcal{F}/\delta\rho$, returns $\beta$ times the integral just obtained:
\begin{equation}
    \frac{\partial \mathcal{V}}{\partial X^l}
    = \beta \int_V \frac{\delta \mathcal{F}}{\delta \rho}(\mathbf{r})\, \frac{\partial \rho(\mathbf{r};\mathbf{X})}{\partial X^l}\, d\mathbf{r}
    = \beta \int_V \mu\, \frac{\partial \rho}{\partial X^l}\, d\mathbf{r} .
    \label{eq:V-chain-rule}
\end{equation}
Combining the two previous displays gives the generalised force in the form stated:
\begin{equation}
    F_l = -D_0\, \frac{\partial \mathcal{V}}{\partial X^l} .
\end{equation}
\end{proof}
\begin{remark}[Geometric nature of the thermodynamic force]
\label{rem:geometric-force}
One has to distinguish between the thermodynamic force components $F_l$, which are covariant derivatives of a scalar potential ($F_l = -D_0\partial_l \mathcal{V}$), and the resulting drift vector components. The force $\mathbb{F} = -D_0d\mathcal{V}$ is a \emph{covector field} (or 1-form) in the cotangent space $T^*_X\mathcal{M}_m$, while physical motion is a vector in the tangent space $T_X\mathcal{M}_m$, and the inverse metric tensor $g^{kl}$ maps the one to the other, by raising the index, $v^k = \sum_l g^{kl} F_l$. In differential geometry this map is the \textit{musical isomorphism} ($\sharp : T^*\mathcal{M}_m \to T\mathcal{M}_m$). Acting in this way, the inverse metric is a state-dependent \emph{mobility tensor}, which quantifies the local ``stiffness'' of the manifold and determines how easily the system responds to the applied thermodynamic gradient in different directions.
\end{remark}
\begin{theorem}[Quadratic variation of the projected flux]
\label{thm:projected-noise}
The projected stochastic fluxes $d\Xi_l(\mathbf{X}, t)$ constitute a state-dependent multiplicative noise. Their diffusion tensor (quadratic variation per unit time) is determined strictly by the manifold metric:
\begin{equation}
    \left[ d\Xi_l, d\Xi_k \right]_t = \bigl(2 D_0 \, g_{lk}(\mathbf{X}) + \mathcal{O}(\delta^{2})\bigr)\, dt .
    \label{eq:projected-noise-qv}
\end{equation}
The identity is exact in the white-noise limit. For $\delta > 0$ the exact quadratic variation is
\begin{equation}
    \left[ d\Xi_l, d\Xi_k \right]_t = 2D_0 \int_0^L\!\!\int_0^L \sqrt{\varpi(r)\varpi(s)}\; C_\delta(r,s)\, \mathbf{e}^m_l(r)\, \mathbf{e}^m_k(s)\, dr\,ds\; dt ,
    \label{eq:projected-noise-exact}
\end{equation}
which is $2D_0\,g^{(\delta)}_{lk}$ in the notation of Eq.~\eqref{eq:finite-delta-metric} and reduces to $2D_0 g_{lk}$ when $C_\delta(r,s) \to \delta(r-s)$. The departure from that limit at fixed $\delta > 0$ is $\mathcal{O}(\delta^{2})$ by Eq.~\eqref{eq:metric-delta-error}, so that the metric used below differs from the exact finite-$\delta$ one by a quantity vanishing with the correlation length.
This quadratic variation defines the diffusion geometry of the reduced stochastic process, serving as the basis for the Gaussian substitution in the Itô representation (see Section \ref{subsec:noise-equivalence}).
\end{theorem}
\begin{proof}
\begin{proofpart}[The finite-$\delta$ covariance]
We compute the covariance of the projected noise increments defined in Corollary \ref{cor:coordinate-sde}. In the limit of infinite stiffness ($\lambda_0 \to \infty$), we evaluate the operators at the manifold configuration $m = M(\cdot; \mathbf{X})$. The covariance integral is given by the expectation of the product of the projections:
\begin{equation}
    C_{lk} \, dt = \int_0^L \!\! \int_0^L \mathbb{E} \left[ (\mathcal{B}_m d\xi^\delta)(r) \cdot (\mathcal{B}_m d\xi^\delta)(s) \right] \, \mathbf{e}_l^m(r) \mathbf{e}_k^m(s) \varpi(r) \varpi(s) \, dr \, ds.
\end{equation}
Using the explicit diffusion operator $\mathcal{B}_m u = \sqrt{2D_0 \partial_r M} \, u$ (Eq.~\ref{eq:mass-diffusion-explicit}) and the spatial covariance $\mathbb{E}[d\xi^\delta(r) d\xi^\delta(s)] = C_\delta(r,s)\,dt$ of Lemma~\ref{lemma:radial-covariance}, the kernel becomes
\begin{equation}
    \mathbb{E} \left[ (\mathcal{B}_m d\xi^\delta)(r) \cdot (\mathcal{B}_m d\xi^\delta)(s) \right] = 2D_0 \sqrt{\partial_r M(r)\, \partial_r M(s)}\; C_\delta(r,s) \, dt .
\end{equation}
Substituting this into the double integral, and writing the square root through the weight $\varpi = (\partial_r M)^{-1}$ of Lemma~\ref{lemma:metric-selection} as $\sqrt{\partial_rM(r)\,\partial_rM(s)} = 1/\sqrt{\varpi(r)\varpi(s)}$, the two remaining factors of $\varpi$ combine with it to leave a single power:
\begin{align}
    C_{lk}
    &= 2D_0 \int_0^L\!\!\int_0^L \frac{C_\delta(r,s)}{\sqrt{\varpi(r)\varpi(s)}}\, \mathbf{e}_l^m(r)\mathbf{e}_k^m(s)\, \varpi(r)\varpi(s) \, dr\, ds \nonumber \\
    &= 2D_0 \int_0^L\!\!\int_0^L \sqrt{\varpi(r)\varpi(s)}\; C_\delta(r,s)\, \mathbf{e}_l^m(r)\, \mathbf{e}_k^m(s)\, dr\, ds ,
    \label{eq:projected-noise-exact-proof}
\end{align}
which is~\eqref{eq:projected-noise-exact}.
\end{proofpart}
\begin{proofpart}[The white-noise limit]
As $\delta \to 0$ the covariance kernel converges to $C_\delta(r,s) \to \delta(r-s)$ in the sense of distributions, and the integration over $s$ in~\eqref{eq:projected-noise-exact-proof} collapses onto $s = r$:
\begin{equation}
    C_{lk} = 2D_0 \int_0^L \varpi(r)\, \mathbf{e}_l^m(r)\, \mathbf{e}_k^m(r) \, dr = 2D_0\, g_{lk}(\mathbf{X}) ,
\end{equation}
the last equality being the definition of the metric tensor, Definition~\ref{def:metric-tensor}. This is~\eqref{eq:projected-noise-qv}, and the deviation from it at finite $\delta$ is $\mathcal{O}(\delta^{2})$ by Eq.~\eqref{eq:metric-delta-error}, the ambient kernel $g_\delta$, even and of unit mass, having vanishing first moment.
\end{proofpart}
\end{proof}
\subsubsection{Equivalence of the projected noise}
\label{subsec:noise-equivalence}
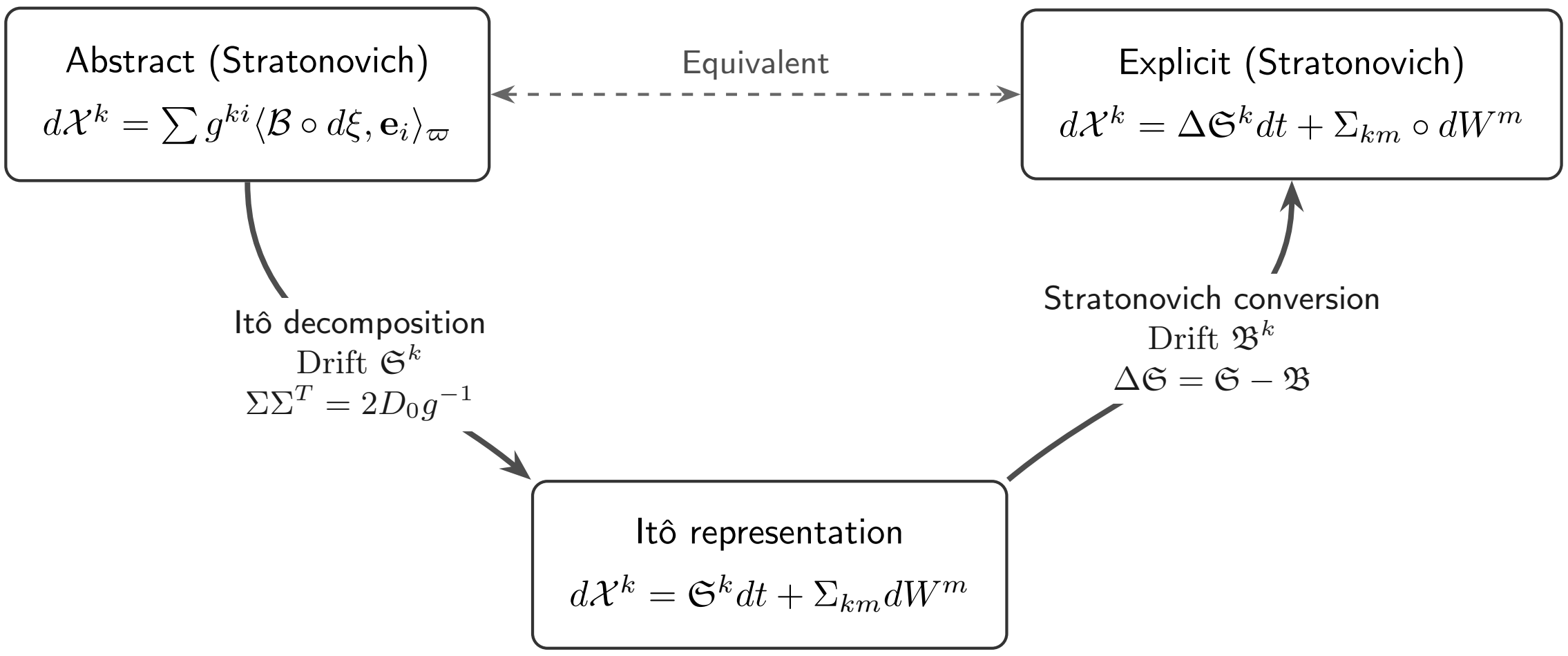
\begin{figure}[ht]
    \centering
\begin{tikzpicture}[
    % Global styles
    box/.style={
        rectangle,
        draw=black!80,
        thick,
        rounded corners=4pt,
        fill=white,
        align=center,
        inner sep=10pt,
        minimum width=2pt,
        minimum height=5pt,
        font=\small
    },
    arrow/.style={
        ->,
        >={Stealth[length=3mm]},
        line width=1.5pt,
        color=black!70
    },
    label_text/.style={
        font=\footnotesize, % No italics
        align=center,
        fill=white,
        text=black!90,
        inner sep=3pt
    }
]
    % Nodes arranged in a V-shape
    % 1. Start (Top Left)
    \node[box] (start) at (-5, 0) {
        \textsf{Abstract (Stratonovich)} \\[0.5em]
        $d\mathcal{X}^k = \sum g^{ki} \langle \mathcal{B} \circ d\xi, \mathbf{e}_i \rangle_{\varpi}$
    };
    % 3. End (Top Right)
    \node[box] (end) at (5, 0) {
        \textsf{Explicit (Stratonovich)} \\[0.5em]
        $d\mathcal{X}^k = \Delta \mathfrak{S}^k dt + \Sigma_{km} \circ dW^m$
    };
    % 2. Intermediate (Bottom Center)
    \node[box] (ito) at (0, -4.5) {
        \textsf{It\^{o} representation} \\[0.5em]
        $d\mathcal{X}^k = \mathfrak{S}^k dt + \Sigma_{km} dW^m$
    };
    % Curved Arrows
    % Path 1: Downward (Start -> Ito)
    \draw[arrow] (start.south) to [out=270, in=140] node[label_text, pos=0.55] {
        \textsf{It\^{o} decomposition} \\
        Drift $\mathfrak{S}^k$ \\
        $\Sigma \Sigma^T = 2D_0 g^{-1}$
    } (ito.north west);
    % Path 2: Upward (Ito -> End)
    \draw[arrow] (ito.north east) to [out=40, in=270] node[label_text, pos=0.55] {
        \textsf{Stratonovich conversion} \\
        Drift $\mathfrak{B}^k$ \\
        $\Delta \mathfrak{S} = \mathfrak{S} - \mathfrak{B}$
    } (end.south);
    % Equivalence Line
    \draw[dashed, <->, >=Stealth, thick, color=black!60] (start.east) -- node[above, font=\footnotesize\sffamily, color=black!70] {Equivalent} (end.west);
\end{tikzpicture}
    \caption{The conversion cycle used to identify the anomalous force. The projected increment is written first in the implicit Stratonovich form which the projection returns, then in It\^{o} form, and finally reassembled as a standard multiplicative noise in the explicit Stratonovich form. The difference between the two deterministic drifts met on the way is the anomalous force $\Delta\mathfrak{S}^k$.}
    \label{fig:conversion-cycle}
\end{figure}

As we have shown in Corollary~\ref{cor:coordinate-sde}, the stochastic evolution of the coordinate $X^k$ is driven by the contravariant projected noise term. With the notation $d\Xi_i = \langle \mathcal{B}_m \circ d\xi^\delta_t, \mathbf{e}_i^m \rangle_{\varpi}$ for the covariant flux of Theorem~\ref{thm:projected-noise}, the coordinate increment is given by:
\begin{equation}
    d\mathcal{X}^k \coloneqq \sum_{i=1}^d g^{ki}(\mathbf{X}) \, d\Xi_i(\mathbf{X}, t).
\end{equation}
Let us now prove that this complex functional noise is statistically equivalent to a standard multiplicative noise plus a deterministic drift correction. To determine this correction, we follow the conversion cycle of Figure~\ref{fig:conversion-cycle}.

\begin{theorem}[Equivalent fluctuation dynamics]
\label{thm:equivalent-fluctuations}
The projected Stratonovich noise term $d\mathcal{X}^k$ is equivalent in law to the sum of a deterministic anomalous force $\Delta \mathfrak{S}^k$ and a standard geometric noise process:
\begin{equation}
    d\mathcal{X}^k = \Delta \mathfrak{S}^k(\mathbf{X}) \, dt + \sum_{m=1}^d \Sigma^{km}(\mathbf{X}) \circ dW_m(t).
\end{equation}
The diffusion amplitude $\mathbf{\Sigma}(\mathbf{X})$ is fixed by the condition $\mathbf{\Sigma}\mathbf{\Sigma}^T = 2D_0 \mathbf{g}^{-1}$ only up to an orthogonal factor on the right, since $\mathbf{\Sigma}\mathbf{O}$ satisfies the same condition for every orthogonal $\mathbf{O}(\mathbf{X})$. Such a change of factor rotates the driving Brownian motion and alters $\Delta\mathfrak{S}^k$ by exactly the difference of the two conversion drifts, so that the law of $\mathcal{X}$ is unchanged. We take throughout the symmetric positive-definite root, which the condition determines uniquely, and the expression below is written for that choice. The anomalous force $\Delta\mathfrak{S}^k$ is given by:
\begin{align}
    \Delta\mathfrak{S}^k(\mathbf{X}) &= \frac{1}{2}\sum_{l,m} \Sigma^{km} \frac{\partial \Sigma^{lm}}{\partial X^l} - \frac{D_0}{2} \sum_{l,j,m} g^{kl} \frac{\partial g^{jm}}{\partial X^l} g_{mj} \nonumber\\
    &\quad + \frac{D_0}{2} \sum_{i,j,l} \left( g^{ki} g^{lj} - g^{kl} g^{ij} \right) \int_0^L \frac{\partial \varpi}{\partial X^l} \mathbf{e}_i^m \mathbf{e}_j^m \, dr,
\end{align}
with the weight variation $\frac{\partial \varpi}{\partial X^l} = -\frac{1}{4\pi r^2 \rho^2} \frac{\partial \rho}{\partial X^l}$.
\end{theorem}
\begin{proof}
The proof determines the covariance, computes the It\^{o} drift of the projection, and removes the It\^{o} drift of the target geometric form.
\begin{proofpart}[Covariance equivalence]
We first determine the required diffusion tensor for the coordinate process. From the definition $d\mathcal{X}^k = \sum_i g^{ki} d\Xi_i$, the quadratic variation is:
\begin{equation}
    [d\mathcal{X}^k, d\mathcal{X}^l]_t = \sum_{p,q} g^{kp} g^{lq} [d\Xi_p, d\Xi_q]_t.
\end{equation}
Substituting the result from Theorem \ref{thm:projected-noise}, $[d\Xi_p, d\Xi_q]_t = 2D_0 g_{pq} dt$, we obtain:
\begin{equation}
    [d\mathcal{X}^k, d\mathcal{X}^l]_t = 2D_0 \sum_{p,q} g^{kp} g^{lq} g_{pq} \, dt = 2D_0 \sum_{p} g^{kp} \delta^l_p \, dt = 2D_0 g^{kl} \, dt.
\end{equation}
For the target process $\sum_m \Sigma^{km} dW_m$ to match this quadratic variation, the amplitude matrix must satisfy the condition $(\mathbf{\Sigma}\mathbf{\Sigma}^T)^{kl} = \sum_m \Sigma^{km} \Sigma^{lm} = 2D_0 g^{kl}$. This condition is necessary and does not by itself fix $\mathbf{\Sigma}$, the symmetric positive-definite root being the choice made in the statement.
\end{proofpart}
\begin{proofpart}[The It\^{o} drift of the integral noise ($\mathfrak{S}^k$)]
Let $\eta_i = \langle \mathcal{B}_m d\xi, \mathbf{e}_i^m \rangle_{\varpi}$ be the It\^{o} part of the projected flux. The total noise term is $\zeta^k = \sum_i g^{ki} \eta_i$. The conversion to Stratonovich introduces the drift correction:
\begin{equation}
    \mathfrak{S}^k = \frac{1}{2} \sum_{l=1}^d \mathbb{E} \left[ \frac{\partial \zeta^k}{\partial X^l} \zeta^l \right].
\end{equation}
We expand the derivative of $\zeta^k$ by the product rule, substitute the result into the expectation, and evaluate the first of the two resulting correlations with $\mathbb{E}[\eta_i \zeta^l] = 2D_0\,\delta_i^l$, which is the metric definition itself:
\begin{align}
    \frac{\partial \zeta^k}{\partial X^l}
    &= \frac{\partial g^{ki}}{\partial X^l}\,\eta_i + g^{ki}\,\frac{\partial \eta_i}{\partial X^l} , \nonumber \\
    \mathfrak{S}^k
    &= \frac{1}{2} \sum_{l,i} \frac{\partial g^{ki}}{\partial X^l}\,\mathbb{E}\bigl[\eta_i\, \zeta^l\bigr]
     + \frac{1}{2} \sum_{l,i} g^{ki}\,\mathbb{E}\!\left[ \frac{\partial \eta_i}{\partial X^l}\, \zeta^l \right] \nonumber \\
    &= D_0 \sum_{l} \frac{\partial g^{kl}}{\partial X^l}
     + \frac{1}{2} \sum_{l,i,j} g^{ki} g^{lj}\, \mathbb{E}\!\left[ \frac{\partial \eta_i}{\partial X^l}\, \eta_j \right]
    \eqqcolon \mathfrak{S}^k_A + \mathfrak{S}^k_B ,
\end{align}
where the second line has been re-expressed by writing $\zeta^l = \sum_j g^{lj}\eta_j$.
\end{proofpart}
\begin{proofpart}[Explicit expansion of $\mathfrak{S}_B^k$]
We evaluate the correlation in $\mathfrak{S}^k_B$. The derivative of the noise projection is:
\begin{equation}
    \frac{\partial \eta_i}{\partial X^l} = \frac{\partial}{\partial X^l}\!\int_0^L \!\sigma\,\varpi\,\mathbf{e}_i^m \, d\xi \, dr = \sqrt{2D_0}\!\int_0^L\!\left[ \varpi^{1/2} \frac{\partial \mathbf{e}_i^m}{\partial X^l} + \tfrac{1}{2}\varpi^{-1/2}\frac{\partial \varpi}{\partial X^l} \mathbf{e}_i^m \right] d\xi \, dr .
\end{equation}
The diffusion amplitude is differentiated together with the tangent field and the weight, all three depending upon $\mathbf{X}$ through $M(\mathbf{X})$. Writing $\sigma\varpi = \sqrt{2D_0\,\varpi}$ before differentiating makes this explicit and supplies the factor $\tfrac12$ carried by the weight-gradient term. Correlating with $\eta_j = \sqrt{2D_0}\int \varpi^{1/2}\mathbf{e}_j^m\,d\xi\,dr$:
\begin{equation}
    \mathbb{E} \left[ \frac{\partial \eta_i}{\partial X^l} \eta_j \right] = 2D_0\int_0^L\left[\varpi\,\frac{\partial^2 M}{\partial X^l \partial X^i}\,\mathbf{e}_j^m + \tfrac{1}{2}\frac{\partial\varpi}{\partial X^l}\mathbf{e}_i^m\mathbf{e}_j^m\right]dr .
\end{equation}
Dividing through by the factor $2D_0$ carried by the noise amplitude:
\begin{equation}
    \frac{1}{2D_0} \mathbb{E} \left[ \frac{\partial \eta_i}{\partial X^l} \eta_j \right] = \int_0^L \varpi \frac{\partial^2 M}{\partial X^l \partial X^i} \mathbf{e}_j^m \, dr + \frac{1}{2}\int_0^L \frac{\partial \varpi}{\partial X^l} \mathbf{e}_i^m \mathbf{e}_j^m \, dr,
\end{equation}
To handle the second derivative of the mass, we differentiate the metric definition $g_{ij} = \int \varpi \mathbf{e}_i^m \mathbf{e}_j^m$:
\begin{equation}
    \frac{\partial g_{ij}}{\partial X^l} = \int \frac{\partial \varpi}{\partial X^l} \mathbf{e}_i^m \mathbf{e}_j^m + \int \varpi \left( \frac{\partial^2 M}{\partial X^l \partial X^i} \mathbf{e}_j^m + \mathbf{e}_i^m \frac{\partial^2 M}{\partial X^l \partial X^j} \right).
\end{equation}
We abbreviate the two integrals that occur, writing $M_{li} \coloneqq \partial^2 M/\partial X^l \partial X^i$ for the second derivative of the embedding:
\begin{equation}
    P_{lij} \coloneqq \int_0^L \varpi\, M_{li}\, \mathbf{e}^m_j \,dr ,
    \qquad
    T_{lij} \coloneqq \int_0^L \frac{\partial \varpi}{\partial X^l}\, \mathbf{e}^m_i\, \mathbf{e}^m_j \,dr ,
    \label{eq:PT-definitions}
\end{equation}
so that the previous display reads $\partial_l g_{ij} = T_{lij} + P_{lij} + P_{lji}$. Second derivatives commute, so $M_{li} = M_{il}$ and $P$ is symmetric under exchange of its first two indices. We form the Christoffel combination and expand each of its three metric derivatives by this rule:
\begin{align}
    \mathsf{S}_{lij}
    &\coloneqq \partial_l g_{ij} + \partial_i g_{lj} - \partial_j g_{il} \nonumber \\
    &= \bigl( T_{lij} + P_{lij} + P_{lji} \bigr)
     + \bigl( T_{ilj} + P_{ilj} + P_{ijl} \bigr)
     - \bigl( T_{jil} + P_{jil} + P_{jli} \bigr) .
    \label{eq:S-expanded}
\end{align}
The six $P$ terms reduce in pairs. Applying the symmetry in the first two indices in the form $P_{ilj} = P_{lij}$, $P_{ijl} = P_{jil}$ and $P_{jli} = P_{lji}$, two cancellations occur and one term is doubled:
\begin{align}
    P_{lij} + P_{lji} + P_{ilj} + P_{ijl} - P_{jil} - P_{jli}
    &= P_{lij} + P_{lji} + P_{lij} + P_{jil} - P_{jil} - P_{lji} \nonumber \\
    &= 2\,P_{lij} .
    \label{eq:P-collapse}
\end{align}
Consequently $\mathsf{S}_{lij} = T_{lij} + T_{ilj} - T_{jil} + 2P_{lij}$, and solving this relation for $P_{lij}$ isolates the target integral:
\begin{equation}
    \int \varpi \frac{\partial^2 M}{\partial X^l \partial X^i} \mathbf{e}_j^m
    = \frac{1}{2}\bigl(\partial_l g_{ij} + \partial_i g_{lj} - \partial_j g_{il}\bigr)
    - \frac{1}{2}\bigl(T_{lij} + T_{ilj} - T_{jil}\bigr).
    \label{eq:christoffel-isolation}
\end{equation}
All three weight-gradient terms are retained, since each arises from one of the three metric derivatives in the Christoffel combination, and none may be dropped.
Substituting this back into $\mathfrak{S}^k_B$ separates the geometric metric derivatives from the weight gradients:
\begin{align}
    \mathfrak{S}_B^k &= D_0 \sum_{l,i,j} g^{ki} g^{lj} \left[ \frac{1}{2} \mathsf{S}_{lij} \right] \nonumber \\
    & \quad + \frac{D_0}{2} \sum_{l,i,j} g^{ki} g^{lj} \left( T_{lij} - T_{lij} - T_{ilj} + T_{jil} \right).
\end{align}
\end{proofpart}
\begin{proofpart}[The explicit It\^{o} drift ($\mathfrak{S}^k$)]
We simplify the first line of $\mathfrak{S}_B$, namely $D_0 \sum_{l,i,j} g^{ki} g^{lj}\,\tfrac{1}{2} \mathsf{S}_{lij}$, by contracting the three metric derivatives of $\mathsf{S}_{lij}$ one at a time. The first and the third cancel against each other, which renaming the summation indices $j \leftrightarrow l$ in the third contraction and using the symmetry $g^{jl} = g^{lj}$ shows,
\begin{align}
    -\sum_{i,l,j} g^{ki} g^{lj}\, \partial_j g_{il}
    &= -\sum_{i,j,l} g^{ki} g^{jl}\, \partial_l g_{ij}
     = -\sum_{i,l,j} g^{ki} g^{lj}\, \partial_l g_{ij} ,
    \label{eq:S-first-third}
\end{align}
which is the exact negative of the first contraction. Only the second survives. To evaluate it we differentiate the identity $\sum_j g_{ij} g^{jn} = \delta^n_i$ with respect to $X^l$, obtaining
\begin{equation}
    \partial_l g_{ij} = -\sum_{m,n} g_{im}\, g_{jn}\, \partial_l g^{mn} ,
    \label{eq:inverse-metric-derivative}
\end{equation}
and substituting this into the second contraction, where one metric pair collapses through $\sum_l g^{lj} g_{lm} = \delta^j_m$:
\begin{align}
    \sum_{i,l,j} g^{ki} g^{lj}\, \partial_i g_{lj}
    &= -\sum_{i,l,j,m,n} g^{ki}\, \bigl( g^{lj} g_{lm} \bigr)\, g_{jn}\, \partial_i g^{mn} \nonumber \\
    &= -\sum_{i,j,n} g^{ki}\, g_{jn}\, \partial_i g^{jn}
     = -\sum_{l,j,m} g^{kl}\, g_{mj}\, \partial_l g^{jm} ,
    \label{eq:S-second}
\end{align}
the last equality being a renaming of the summation indices $i \to l$ and $n \to m$. This contribution involves the metric alone and carries no weight gradient.
We can now state the total drift $\mathfrak{S}^k = \mathfrak{S}_A^k + \mathfrak{S}^k_B$ that would appear if we remained in the It\^{o} interpretation:
\begin{align}
    \mathfrak{S}^k &= D_0 \sum_{l} \frac{\partial g^{kl}}{\partial X^l} - \frac{D_0}{2} \sum_{l,j,m} g^{kl} (\partial_l g^{jm}) g_{mj} \nonumber \\
    & \quad + \frac{D_0}{2} \sum_{l,i,j} \left( g^{ki} g^{lj} - g^{kl} g^{ij} \right) \int_0^L \frac{\partial \varpi}{\partial X^l} \mathbf{e}_i^m \mathbf{e}_j^m \, dr.
\label{eq:explicit-ito-drift}
\end{align}
\end{proofpart}
\begin{proofpart}[Returning to Stratonovich ($\mathfrak{S}^k \to \Delta \mathfrak{S}^k$)]
To write the final equation in Stratonovich form with the standard multiplicative noise $\sum \Sigma^{km} \circ dW_m$, we must subtract the It\^{o}-to-Stratonovich drift generated by this specific matrix choice. The standard correction is:
\begin{equation}
\label{eq:ito-stratonovich-correction}
    \mathfrak{B}^k = \frac{1}{2} \sum_{l,m} \Sigma^{lm} \frac{\partial \Sigma^{km}}{\partial X^l}.
\end{equation}
The final anomalous force is the difference $\Delta\mathfrak{S}^k = \mathfrak{S}^k - \mathfrak{B}^{k}$.
We rewrite the first term of $\mathfrak{S}^k$, namely $\mathfrak{S}_A^k$, using the relation $2D_0 g^{kl} = \sum_m \Sigma^{km} \Sigma^{lm}$:
\begin{equation}
    \mathfrak{S}_A^k = \frac{1}{2} \sum_{l,m} \left( \frac{\partial \Sigma^{km}}{\partial X^l} \Sigma^{lm} + \Sigma^{km} \frac{\partial \Sigma^{lm}}{\partial X^l} \right) = \mathfrak{B}^k + \frac{1}{2} \sum_{l,m} \Sigma^{km} \frac{\partial \Sigma^{lm}}{\partial X^l}.
\end{equation}
The term $\mathfrak{B}^k$ cancels out exactly, leaving only the second part:
\begin{equation}
\mathfrak{S}_A^k - \mathfrak{B}^k = \frac{1}{2} \sum_{l,m} \Sigma^{km} \frac{\partial \Sigma^{lm}}{\partial X^l}.
\end{equation}
\end{proofpart}
\begin{proofpart}[The integral part]
Finally, we assemble the integral remainder from $\mathfrak{S}_B$, the bracket of which carries the three weight-gradient terms
\begin{equation}
    T_{lij} - T_{lij} - T_{ilj} + T_{jil} = T_{jil} - T_{ilj} .
    \label{eq:T-bracket}
\end{equation}
Two of these are brought to the index pattern of the first by renaming the summation indices. In the third term we set $l' = j$ and $j' = l$, and use the symmetry $g^{j'l'} = g^{l'j'}$:
\begin{align}
    \sum_{i,l,j} g^{ki} g^{lj}\, T_{jil}
    &= \sum_{i,l',j'} g^{ki} g^{j'l'}\, T_{l' i j'}
     = \sum_{i,l',j'} g^{ki} g^{l'j'}\, T_{l' i j'}
     = \sum_{i,l,j} g^{ki} g^{lj}\, T_{lij} .
    \label{eq:T-third}
\end{align}
In the second term the derivative index of $T$ is $i$, and it is contracted with $g^{ki}$, so that setting $l' = i$, $i' = l$ and $j' = j$ carries the derivative index into the first slot:
\begin{align}
    \sum_{i,l,j} g^{ki} g^{lj}\, T_{ilj}
    &= \sum_{l',i',j'} g^{k l'} g^{i' j'}\, T_{l' i' j'}
     = \sum_{i,l,j} g^{kl} g^{ij}\, T_{lij} .
    \label{eq:T-second}
\end{align}
The three contractions therefore collect onto the single tensor
\begin{align}
    \sum_{i,l,j} g^{ki} g^{lj} \bigl( T_{jil} - T_{ilj} \bigr)
    &= \sum_{i,l,j} \bigl( g^{ki} g^{lj} - g^{kl} g^{ij} \bigr) T_{lij} .
    \label{eq:anomaly-tensor-collect}
\end{align}
Collecting,
\begin{equation}
    \mathcal{I}^k = \frac{D_0}{2} \sum_{i,j,l} \bigl(g^{ki} g^{lj} - g^{kl} g^{ij}\bigr) \int_0^L \frac{\partial \varpi}{\partial X^l} \mathbf{e}_i^m \mathbf{e}_j^m \, dr
    = \frac{D_0}{2} \sum_{i,j,l} \bigl(g^{ki} g^{lj} - g^{kl} g^{ij}\bigr) T_{lij}.
    \label{eq:anomaly-tensor}
\end{equation}
The relative weight of the two contractions is fixed by the weight-gradient terms of Eq.~\eqref{eq:christoffel-isolation} together with the factor $\tfrac12$ generated by differentiating the diffusion amplitude, and the antisymmetric combination obtained here agrees with the anomalous force of Lutsko~\cite{lutsko2012mean}.
It remains to express the weight gradient through the tangent basis. We differentiate $\varpi = (4\pi r^2\rho)^{-1}$ with respect to $X^l$, then use $\varpi^{-1} = \partial_r M = 4\pi r^2\rho$ to recognise the derivative of the embedding, and finally $\partial_{X^l}(\partial_r M) = \partial_r \mathbf{e}^m_l$:
\begin{align}
    \frac{\partial \varpi}{\partial X^l}
    &= -\frac{1}{(4\pi r^2\rho)^2}\,\frac{\partial (4\pi r^2 \rho)}{\partial X^l}
     = -\frac{1}{4\pi r^2 \rho^2}\, \frac{\partial \rho}{\partial X^l} \nonumber \\
    &= -\varpi^2\, \frac{\partial}{\partial X^l}\bigl( \partial_r M \bigr)
     = -\varpi^2\, \partial_r \mathbf{e}^m_l .
    \label{eq:weight-gradient}
\end{align}
Substituting this into~\eqref{eq:anomaly-tensor} gives the expression stated in the theorem.
\end{proofpart}
\end{proof}
\begin{remark}[Physical interpretation of the anomalous force]
\label{disc:anomalous-force}
The deterministic drift $\Delta \mathfrak{S}^k$ is the correction required if the reduced dynamics are to carry the equilibrium statistics of the parent field theory, up to the noise-induced shift of Proposition~\ref{prop:invariant-measures}, and it originates in the failure of the stochastic integration to commute with the state-dependent projection $\Pi_{\mathbf{X}}^{\parallel}$. It has two terms with distinct geometric origins, and we keep them apart, since only the second of them is an anomaly in the sense used throughout this work.
The terms in $\partial_l g^{jm}$ contract the Christoffel symbols ($g^{lm}\Gamma_{lm}^k$) and turn the generator into the Laplace--Beltrami operator of $g$, and hence the reduced diffusion is the Brownian motion of the metric $g$. This contribution comes from the curvature of the coordinate grid itself. It can be read as a diffusive centrifugal force, in the sense that inertial forces arise in classical mechanics whenever generalised coordinates describe a curved path, since the drift compensates for that curvature and the coordinate noise carries no drift of its own.
What about the terms in $\partial\varpi$? They are of a different character, and together with the tangential projection of the field-level conversion term they make up what Theorem~\ref{thm:geo-operator} isolates as the geometric anomaly $\mathfrak{s}_a$. They come from the inhomogeneity of the inner product weight along the manifold, not from the rotation of the tangent frame within the ambient space, as Remark~\ref{rem:flat-manifold-test} establishes by exhibiting a manifold with an identically vanishing second fundamental form, on which $\mathfrak{s}_a$ does not vanish in general. Without this contribution, the projected noise would implicitly sample a flat Euclidean space, and the reduced dynamics would relax to the na\"{\i}ve density $\exp(-\mathcal{V}/\varepsilon^2)$ instead of the shifted stationary measure $\exp(-\mathcal{V}/\varepsilon^2 - \Phi_\xi/D_0)$ of Proposition~\ref{prop:invariant-measures} (understood as a density against the Riemannian volume measure $\mathfrak{m}_g$).
\end{remark}
\begin{corollary}[The geometric Langevin equations]
\label{cor:geometric-langevin}
Combining the thermodynamic driving forces (Theorem \ref{thm:projected-drift}) with the equivalent fluctuation dynamics (Theorem \ref{thm:equivalent-fluctuations}), the evolution of the collective coordinates is governed by the following equivalent stochastic differential equations.
\paragraph{Stratonovich representation.}
This form preserves the standard chain rule and transforms covariantly (like a standard vector field) under coordinate changes on the manifold $\mathcal{M}_m$.
\begin{equation}
    dX^k_t = \sum_{l=1}^d g^{kl}(\mathbf{X}) \bigl( F_l(\mathbf{X}) - S_l(\mathbf{X}) \bigr) \, dt + \Delta \mathfrak{S}^k(\mathbf{X}) \, dt + \sum_{m=1}^d \Sigma^{km}(\mathbf{X}) \circ dW_m(t).
\end{equation}
Here, $F_l = -D_0\partial \mathcal{V} / \partial X^l$ is the thermodynamic force (Theorem \ref{thm:projected-drift}), $S_l$ is the projected field correction of Eq.~\eqref{eq:S-projected}, and $\Delta \mathfrak{S}^k$ is the anomalous force (Theorem \ref{thm:equivalent-fluctuations}). The correction $S_l$ is written here without the factor $\varepsilon^2$ which Corollary~\ref{cor:coordinate-sde} carries, in the same convention as $\Delta\mathfrak{S}^k$ and $\Sigma^{km}$ beside it, the factor being restored to all of them together by Definition~\ref{def:weak-noise}.
\paragraph{It\^{o} representation.}
This form separates the martingale term from the total drift, making it suitable for deriving the Fokker--Planck equation.
\begin{equation}
    dX^k_t = \sum_{l=1}^d g^{kl}(\mathbf{X}) \bigl( F_l(\mathbf{X}) - S_l(\mathbf{X}) \bigr) \, dt + \mathfrak{S}^k(\mathbf{X})dt + \sum_{m=1}^d \Sigma^{km}(\mathbf{X}) \, dW_m(t).
\end{equation}
The Itô drift $\mathfrak{S}^k$ (Eq.~\ref{eq:explicit-ito-drift}) equals the sum of the anomalous force $\Delta \mathfrak{S}^k$ and the standard Stratonovich-to-It\^{o} correction (Eq.~\ref{eq:ito-stratonovich-correction}) $\mathfrak{B}^k(\mathbf{X}) = \frac{1}{2}  \sum_{l,m} \Sigma^{lm}(\mathbf{X}) \frac{\partial \Sigma^{km}(\mathbf{X})}{\partial X^l}.$
\end{corollary}
\begin{proof}
The Stratonovich form follows directly by summing the projected deterministic drift of Corollary~\ref{cor:coordinate-sde}, the two parts of which are evaluated in Theorem \ref{thm:projected-drift} and in Eq.~\eqref{eq:S-projected}, and the effective stochastic dynamics (Theorem \ref{thm:equivalent-fluctuations}). The It\^{o} form is obtained by applying the standard conversion formula to the diffusion term $\Sigma^{km} \circ dW_m$, adding the resulting correction $\mathfrak{B}^k$ to the anomalous force $\Delta \mathfrak{S}^k$ to yield the total It\^{o} drift $\mathfrak{S}^k$.
\end{proof}
\begin{remark}[Intrinsic metric-dual representation]
The dynamics can be formulated as the response of the manifold to a total driving semimartingale in the cotangent bundle. Let $\sharp: T^*\mathcal{M}_m \to T\mathcal{M}_m$ (index raising) and $\flat: T\mathcal{M}_m \to T^*\mathcal{M}_m$ (index lowering) denote the musical isomorphisms induced by the metric $g$. The evolution equation is the metric dual of the sum of the total deterministic and stochastic force 1-forms:
\begin{equation}
    dX_t = \left( \mathbb{F}_{\text{tot}} \, dt + \mathbb{W}_t \right)^{\sharp},
\end{equation}
with the total deterministic force 1-form $\mathbb{F}_{\text{tot}} = \mathbb{F} - S + (\Delta \mathfrak{S})^\flat$ (where $\mathbb{F} = - D_0 \, d\mathcal{V}$ and $S$ is the 1-form with components $S_l$), and the stochastic force 1-form $\mathbb{W}_t = \sum_{m=1}^d \omega_m \circ dW_m(t)$ (where $\omega_m$ is the one-form dual to the $m$-th column of $\boldsymbol{\Sigma}$), and both are elements of $T^*\mathcal{M}_m$. In this formulation, the inverse metric $g^{-1}$ acts as the universal mobility tensor, which converts the abstract total impulse $(\mathbb{F}_{\text{tot}} \, dt + \mathbb{W}_t)$ in the cotangent bundle into physical displacement in the tangent bundle.
\end{remark}
\subsection{The induced geometric evolution}
\label{subsec:induced-evolution}
We now characterise the evolution of the effective process on the manifold $\mathcal{M}_m$. Instead of postulating a conservation law, we derive the infinitesimal generator associated with the geometric Langevin equations (Corollary~\ref{cor:geometric-langevin}), and we show that it is equivalent to the intrinsic Laplace--Beltrami operator.
\begin{definition}[Induced generator]
Let $\phi \in C^\infty_c(\mathcal{M}_m)$ be a smooth test function. The infinitesimal generator $\mathcal{L}: C^\infty_c(\mathcal{M}_m) \to C^\infty(\mathcal{M}_m)$ of the process $\mathbf{X}_t$ is defined by the It\^{o} formula applied to the SDE derived in Corollary \ref{cor:geometric-langevin}:
\begin{equation}
    \mathcal{L}\phi \coloneqq \sum_{k=1}^d \mathcal{A}^k(\mathbf{X}) \frac{\partial \phi}{\partial X^k} + \frac{1}{2} \sum_{k,l=1}^d \mathcal{D}^{kl}(\mathbf{X}) \frac{\partial^2 \phi}{\partial X^k \partial X^l},
\end{equation}
where the diffusion tensor is $\mathcal{D}^{kl} = 2D_0 g^{kl}$ and the total It\^{o} drift vector is $\mathcal{A}^k = \sum_l g^{kl} \bigl( F_l - S_l \bigr) + \mathfrak{S}^k$.
\end{definition}
The following theorem establishes that the explicit form of the drift $\mathfrak{S}^k$ derived in the previous section ensures the operator is geometrically intrinsic.
\begin{theorem}[Geometric operator equivalence]
\label{thm:geo-operator}
The generator $\mathcal{L}$ is strictly equivalent to the twisted Laplace--Beltrami operator on $(\mathcal{M}_m, g)$:
\begin{equation}
    \mathcal{L}\phi = D_0 \Delta_g \phi - D_0\langle \nabla \mathcal{V}, \nabla \phi \rangle_g + \langle \mathfrak{s}_a, \nabla \phi \rangle_g.
\end{equation}
Here, $\langle u, v \rangle_g \equiv \sum_{i,j} g_{ij}u^i v^j$ denotes the Riemannian metric inner product, $\Delta_g$ is the Laplace--Beltrami operator, and $\mathfrak{s}_a$ is the geometric anomaly vector field, defined as the whole of the reduced drift carried neither by the Laplace--Beltrami operator nor by the gradient of $\mathcal{V}$. It receives two contributions of the same order, the weighted projection of the noise measure gradient and the tangential projection of the field-level It\^{o}--Stratonovich correction of Lemma~\ref{lemma:ito-strat-correction}, and is given in closed form by
\begin{equation}
    \mathfrak{s}_a^k = \frac{D_0}{2} \sum_{i,j,l} \left( g^{ki} g^{lj} - g^{kl} g^{ij} \right) \int_0^L \frac{\partial \varpi}{\partial X^l} \mathbf{e}_i^m(r)\, \mathbf{e}_j^m(r) \, dr
    \;-\; \sum_{l} g^{kl} S_l .
    \label{eq:geo-anomaly-explicit}
\end{equation}
The second contribution is present because the kinetic pathway projects the Stratonovich equation and converts to It\^{o} afterwards, so that it carries the field-level conversion term through the projection, whereas the probabilistic pathway converts at the level of the field and never encounters it.
\end{theorem}

\begin{proof}
We proceed by expanding the Laplace--Beltrami operator in local coordinates and comparing it term-by-term with the generator $\mathcal{L}$.
\begin{proofpart}[Expansion of the Laplace--Beltrami operator]
By definition, the action of $\Delta_g$ on a scalar $\phi$ is given by the divergence of the gradient:
\begin{equation}
    \Delta_g \phi = \frac{1}{\sqrt{|g|}} \sum_{k,l} \frac{\partial}{\partial X^k} \left( \sqrt{|g|} g^{kl} \frac{\partial \phi}{\partial X^l} \right),
\end{equation}
where $|g| \equiv \det(g_{ij})$. Applying the product rule to the term inside the divergence:
\begin{align}
    \Delta_g \phi &= \frac{1}{\sqrt{|g|}} \sum_{k,l} \left[ \frac{\partial}{\partial X^k}(\sqrt{|g|} g^{kl}) \frac{\partial \phi}{\partial X^l} + \sqrt{|g|} g^{kl} \frac{\partial^2 \phi}{\partial X^k \partial X^l} \right] \nonumber \\
    &= \sum_{k,l} g^{kl} \frac{\partial^2 \phi}{\partial X^k \partial X^l} + \sum_{l} \left[ \sum_k \frac{\partial g^{kl}}{\partial X^k} + \sum_k g^{kl} \frac{\partial \ln \sqrt{|g|}}{\partial X^k} \right]\frac{\partial \phi}{\partial X^l}\\
    &=\sum_{k,l} g^{kl} \frac{\partial^2 \phi}{\partial X^k \partial X^l} + \sum_{l}\mathcal{T}^l \frac{\partial \phi}{\partial X^l}
\end{align}
Multiplying by $D_0$, we identify the diffusion term $D_0 \sum g^{kl} \partial^2_{kl}\phi$, which matches the second-order term in $\mathcal{L}$ exactly. The proof thus reduces to showing that the drift coefficient of $\mathcal{L}$ matches $D_0 \mathcal{T}^l$ plus the potential and anomaly terms.
\end{proofpart}
\begin{proofpart}[Analysis of the generator drift]
The total drift in $\mathcal{L}$ is $\mathcal{A}^k = \sum_l g^{kl} \bigl( F_l - S_l \bigr) + \mathfrak{S}^k$, in which $F_l$ is the generalised force of Theorem \ref{thm:projected-drift} and $S_l$ is the projected field correction of Eq.~\eqref{eq:S-projected}, carried into the reduced equation by the Stratonovich drift. We substitute the explicit expression for the It\^{o} drift $\mathfrak{S}^k$ derived in Eq. \eqref{eq:explicit-ito-drift} of Theorem \ref{thm:equivalent-fluctuations}:
\begin{equation}
    \mathfrak{S}^k =D_0 \sum_{l} \frac{\partial g^{kl}}{\partial X^l}- \frac{D_0}{2} \sum_{l,j,m} g^{kl} \frac{\partial g^{jm}}{\partial X^l} g_{mj} + \mathcal{I}^k,
\end{equation}
where $\mathcal{I}^k$ is the first contribution to $\mathfrak{s}_a^k$ in Eq.~\eqref{eq:geo-anomaly-explicit}.
To relate the second term to the metric determinant, we use the Jacobi formula for the derivative of the inverse matrix. For any invertible matrix $M(x)$, $\partial_x \ln (\det M) = \text{Tr}(M^{-1} \partial_x M)$. Applied to the inverse metric $g^{-1}$ (with components $g^{jm}$):
\begin{equation}
    \sum_{j,m} g_{mj} \frac{\partial g^{jm}}{\partial X^l} = \text{Tr}\left( g \frac{\partial g^{-1}}{\partial X^l} \right) = \frac{\partial \ln (\det g^{-1})}{\partial X^l}.
\end{equation}
Since $\det g^{-1} = (\det g)^{-1} = |g|^{-1}$, we have:
\begin{equation}
    \frac{\partial \ln |g|^{-1}}{\partial X^l} = - \frac{\partial \ln |g|}{\partial X^l} = - 2 \frac{\partial \ln \sqrt{|g|}}{\partial X^l}.
\end{equation}
Substituting this identity back into the second term of $\mathfrak{S}^k$:
\begin{equation}
   -\frac{D_0}{2} \sum_{l,j,m} g^{kl} \frac{\partial g^{jm}}{\partial X^l} g_{mj}  = -\frac{D_0}{2} \sum_l g^{kl} \left( - 2 \frac{\partial \ln \sqrt{|g|}}{\partial X^l} \right) = D_0 \sum_l g^{kl} \frac{\partial \ln \sqrt{|g|}}{\partial X^l}.
\end{equation}
\end{proofpart}
\begin{proofpart}[Synthesis]
The total It\^{o} drift, taken together with the projected field correction which the Stratonovich drift carries, is therefore:
\begin{equation}
    \mathfrak{S}^k - \sum_l g^{kl} S_l = D_0 \sum_{l} \frac{\partial g^{kl}}{\partial X^l} + D_0 \sum_l g^{kl} \frac{\partial \ln \sqrt{|g|}}{\partial X^l} + \mathfrak{s}_a^k,
\end{equation}
an identity between the unscaled objects, the factor $\varepsilon^2$ being supplied to $\mathfrak{S}^k$, to $S_l$ and hence to $\mathfrak{s}_a^k$ alike by Definition~\ref{def:weak-noise}, so that the two contributions to $\mathfrak{s}_a^k$ are $\mathcal{I}^k$ and $-\sum_l g^{kl}S_l$.
Comparing this to the Laplace--Beltrami geometric drift $D_0 \mathcal{T}^k$, we observe an exact match for the first two terms. The generator action is thus:
\begin{align}
    \mathcal{L}\phi &= \frac{1}{2}\sum_{k,l} \bigl(2 D_0 g^{kl}\bigr) \partial^2_{kl} \phi + \sum_k \Bigl( \sum_l g^{kl} F_l + D_0 \mathcal{T}^k + \mathfrak{s}_{a}^k \Bigr) \partial_k \phi \nonumber \\
    &= D_0 \Bigl( \sum_{k,l} g^{kl} \partial^2_{kl} \phi + \sum_k \mathcal{T}^k \partial_k \phi \Bigr) + \sum_{k,l} g^{kl} F_l \partial_k \phi + \sum_k \mathfrak{s}_a^k \partial_k \phi\nonumber\\
    &= D_0\Delta_g \phi + \sum_{k,l} g^{kl} F_l \partial_k \phi + \sum_k \mathfrak{s}_a^k\partial_k \phi
\end{align}
\end{proofpart}
\noindent{}We identify the remaining terms geometrically. First,
the term $\sum g^{kl} F_l \partial_k \phi$ is the action of the vector field $-D_0\nabla \mathcal{V}$, which in metric notation is $-D_0\langle \nabla \mathcal{V}, \nabla \phi \rangle_g$. The projected field correction $-\sum_l g^{kl}S_l$ is a drift of the same order as the anomaly, and is carried within $\mathfrak{s}_a$ accordingly. Second, the term $\sum \mathfrak{s}_a^k \partial_k \phi$ is the direct contraction of the vector field $\mathfrak{s}_a$ with the differential $d\phi$. Since the gradient vector is defined as $(\nabla \phi)^k = \sum g^{km} \partial_m \phi$, we have the identity $\mathfrak{s}_a^k \partial_k \phi = \langle \mathfrak{s}_a, \nabla \phi \rangle_g$.
Thus, the operator takes the invariant form:
\begin{equation}
    \mathcal{L}\phi = D_0 \Delta_g \phi - D_0\langle \nabla \mathcal{V}, \nabla \phi \rangle_g + \langle \mathfrak{s}_{a}, \nabla \phi \rangle_g.
\end{equation}
\end{proof}
The covariance kernel enters the anomaly through the projected field correction $S_l$, and Section~\ref{sec:commutativity} carries that dependence through the invariants $V_\delta$ and $V_\delta'$. The quadratic variation is independent of the kernel, since Theorem~\ref{thm:projected-noise} evaluates it in the white-noise limit.
\begin{remark}[The anomaly is not curvature-induced: a flat-manifold test]
\label{rem:flat-manifold-test}
Since $\mathfrak{s}_a$ arises from the geometry of the reduction, it is natural to attribute it to the extrinsic curvature of $\mathcal{M}_m$. We test this attribution on a manifold of vanishing second fundamental form, where we find that it is incorrect and isolate the mechanism that is actually responsible. Let us take an affine family
\begin{equation}
    M(r;\mathbf{X}) = M_0(r) + \sum_{k=1}^d X^k \phi_k(r),
\end{equation}
where the tangent vectors $\mathbf{e}^m_k = \phi_k$ are independent of $\mathbf{X}$, the second fundamental form vanishes identically, and $\mathcal{M}_m$ is a flat affine subspace. Nevertheless, the weight still varies along the manifold, since by Definition~\ref{def:weighted-space},
\begin{equation}
    \frac{\partial\varpi}{\partial X^l} = -\varpi^2\,\partial_r\phi_l \;\neq\; 0 ,
\end{equation}
and hence the contraction $T_{lij}$ of Eq.~\eqref{eq:PT-definitions} does not vanish, and so the geometric anomaly, built from $T_{lij}$ (through Eq.~\eqref{eq:anomaly-tensor}) and from the projected field correction $S_l$, does not vanish in general. In the case $d = 1$, the two contractions of $T_{lij}$ in Eq.~\eqref{eq:anomaly-tensor} cancel, and the anomaly is then carried by $S_l$ alone, where the weight enters through the noise amplitude $\sigma = (2D_0/\varpi)^{1/2}$. The anomaly therefore survives, in general, on a manifold of identically vanishing extrinsic curvature.

This means that the anomaly is not determined by the second fundamental form of $\mathcal{M}_m$. It is produced instead by the state dependence of the noise amplitude, a distinct mechanism, which enters through $\varpi = (\partial_r M)^{-1}$ (both through its variation in $r$ and through its variation along the manifold). For this reason the second fundamental form does not appear in $\mathfrak{s}_a$, and we call the geometric anomaly noise-induced, not curvature-dependent.
\end{remark}

\subsubsection{The weak-noise limit and invariant measures}
\label{subsubsec:thermodynamic-limit}
The generator $\mathcal{L}$ derived in Theorem \ref{thm:geo-operator} describes the effective evolution of observables on the manifold. To recover the macroscopic covariant description, we must pass to the dual picture (the Fokker--Planck equation) and quantify the magnitude of the geometric anomalies relative to the thermodynamic forces. To this end, we formally introduce the \emph{weak-noise limit}, by scaling the magnitude of the projected fluctuations. In this way, and consistently with the Freidlin--Wentzell theory of large deviations, we can treat the geometric anomalies as higher-order perturbations.
\begin{definition}[Weak-noise scaling convention]
\label{def:weak-noise}
Let $\varepsilon > 0$ be a dimensionless parameter characterising the noise intensity (physically $\varepsilon = 1$, and $\varepsilon^2 = N^{-1}$ when the density is measured in units of $N$). We write the stochastic increment of Corollary \ref{cor:coordinate-sde} as:
\begin{equation}
    d\mathcal{X}^k_\varepsilon = \varepsilon\,\sum_{i=1}^d g^{ki}(\mathbf{X}) \,d\Xi_i(\mathbf{X}, t).
\end{equation}
This scaling implies the following hierarchy for the coefficients:
\begin{enumerate}[(i)]
    \item \emph{Diffusion tensor:} The quadratic variation scales quadratically with $\varepsilon$:
    \begin{equation}
        \mathcal{D}^{kl}_\varepsilon = \varepsilon^2 \mathcal{D}^{kl} = 2 \varepsilon^2 D_0 g^{kl}.
    \end{equation}
    \item \emph{Anomalous drift:} The It\^{o} drift $\mathfrak{S}^k$, its part $\mathfrak{B}^k$ and the geometric anomaly $\mathfrak{s}_a^k$ arise from the contraction of the noise coefficients (e.g., $\Sigma \partial \Sigma$). Consequently, they scale quadratically with the noise amplitude:
    \begin{equation}
        \mathfrak{S}^k_\varepsilon = \varepsilon^2 \mathfrak{S}^k \sim \mathcal{O}(\varepsilon^2).
    \end{equation}
    The projected field correction $S_l$ scales in the same way, $S_{l,\varepsilon} = \varepsilon^2 S_l$, being quadratic in the noise amplitude through Lemma~\ref{lemma:ito-strat-correction}, which is why Corollary~\ref{cor:geometric-langevin} and the results after it carry both without the factor.
    \item \emph{Thermodynamic force:} The effective potential $\mathcal{V}$ is a static property of the free energy landscape, independent of the dynamical noise scaling:
    \begin{equation}
        F_k \sim \mathcal{O}(1).
    \end{equation}
\end{enumerate}
\end{definition}
Let us define the functional space for the probability density. Let $\mathfrak{m}_g$ be the Riemannian volume measure on $\mathcal{M}_m$, denoted locally by $d\mathfrak{m}_g = \sqrt{g} \, d\mathbf{X}$ (where $g \equiv \det(g_{ij})$). We denote by $L^2(\mathfrak{m}_g)$ the Hilbert space of square-integrable functions with respect to this measure, equipped with the inner product $\langle u, v \rangle_{L^2(\mathfrak{m}_g)} = \int_{\mathcal{M}_m} u\,v\, d\mathfrak{m}_g.$
\begin{theorem}[The general Fokker--Planck equation]
\label{thm:general-fpe}
Let $\rho(\mathbf{X}, t)$ be the probability density of the process with respect to the measure $\mathfrak{m}_g$, with no flux through the boundary of $\mathcal{M}_m$. The exact evolution of $\rho$ is governed by the $L^2(\mathfrak{m}_g)$-adjoint of the scaled generator $\mathcal{L}_\varepsilon$. Under Definition~\ref{def:weak-noise}, the density evolves according to:
\begin{equation}
    \frac{\partial \rho}{\partial t} = \nabla_g \cdot \left( \varepsilon^2 D_0 \nabla \rho + D_0\rho \nabla \mathcal{V} - \varepsilon^2 \rho \mathfrak{s}_a \right),
\end{equation}
where $\nabla_g \cdot$ denotes the Riemannian divergence operator and $\nabla$ the gradient vector field.
\end{theorem}
\begin{proof}
We derive the scaled adjoint operator $\mathcal{L}_\varepsilon^*$. By definition, $\mathcal{L}_\varepsilon^*$ satisfies the duality relation $\langle \psi, \mathcal{L}_\varepsilon \phi \rangle_{L^2(\mathfrak{m}_g)} = \langle \mathcal{L}_\varepsilon^*\psi, \phi \rangle_{L^2(\mathfrak{m}_g)}$ for any test functions $\phi, \psi$.
Substituting the scaled coefficients into the geometric form of Theorem \ref{thm:geo-operator}:
\begin{equation}
    \mathcal{L}_\varepsilon \phi = \varepsilon^2 D_0 \Delta_g \phi - D_0\langle \nabla \mathcal{V}, \nabla \phi \rangle_g + \varepsilon^2 \langle \mathfrak{s}_a, \nabla \phi \rangle_g.
\end{equation}
Using the divergence theorem (Green's identity on the manifold) to integrate by parts with respect to $d\mathfrak{m}_g$:
\begin{proofpart}[Diffusion term]
The Laplace--Beltrami operator is self-adjoint with respect to $\mathfrak{m}_g$:
    \begin{equation}
        \int \psi (\varepsilon^2 D_0 \Delta_g \phi) \, d\mathfrak{m}_g = \int \phi (\varepsilon^2 D_0 \Delta_g \psi) \, d\mathfrak{m}_g.
    \end{equation}
\end{proofpart}
\begin{proofpart}[Drift terms]
For the thermodynamic drift ($-D_0\nabla \mathcal{V}$) and the anomalous drift ($\varepsilon^2 \mathfrak{s}_a$), we examine the general property of a vector field $W$:
\begin{align}
    \int_{\mathcal{M}_m} \psi \langle W, \nabla \phi \rangle_g \, d\mathfrak{m}_g
    &= \int_{\mathcal{M}_m} \left[ \nabla_g \cdot (\psi \phi W) - \phi \nabla_g \cdot (\psi W) \right] d\mathfrak{m}_g \nonumber \\
    &= \oint_{\partial \mathcal{M}_m} \psi \phi \langle W, \hat{n} \rangle_g \, dS_g - \int_{\mathcal{M}_m} \phi \nabla_g \cdot (\psi W) \, d\mathfrak{m}_g.
\end{align}
By the no-flux condition of the statement, the boundary integral vanishes. Letting $W = -D_0\nabla \mathcal{V} + \varepsilon^2 \mathfrak{s}_a$ exactly yields the negative divergence of the flux:
    \begin{equation}
        \int \psi \langle -D_0\nabla \mathcal{V} + \varepsilon^2 \mathfrak{s}_a, \nabla \phi \rangle_g \, d\mathfrak{m}_g = - \int \phi \, \nabla_g \cdot \left( \psi [-D_0\nabla \mathcal{V} + \varepsilon^2 \mathfrak{s}_a] \right) \, d\mathfrak{m}_g.
    \end{equation}
\end{proofpart}
\noindent{}Replacing the test function $\psi$ with the density $\rho$, the exact evolution equation is:
\begin{equation}
    \frac{\partial \rho}{\partial t} = \mathcal{L}_\varepsilon^* \rho = \nabla_g \cdot \left( \varepsilon^2 D_0 \nabla \rho + D_0\rho \nabla \mathcal{V} - \varepsilon^2 \rho \mathfrak{s}_a \right).
\end{equation}
\end{proof}
\begin{corollary}[The covariant Fokker--Planck equation in the adiabatic limit]
\label{cor:covariant-fpe}
Discarding the anomalous flux, which is of order $\varepsilon^2$ relative to the thermodynamic one, and retaining the diffusive flux, which is of the same order but carries the highest derivative, the equation of Theorem~\ref{thm:general-fpe} reduces to the covariant Fokker--Planck equation, which retains the large-deviation exponent $\mathcal{V}/\varepsilon^2$ of the stationary density but, by Proposition~\ref{prop:invariant-measures}~(iii) below, misses a factor of order unity in it:
\begin{equation}
    \frac{\partial \rho}{\partial t} = \nabla_g \cdot \left( \varepsilon^2 D_0 \nabla \rho + D_0\rho \nabla \mathcal{V} \right).
\end{equation}
\end{corollary}
\begin{proof}
We analyse the order of magnitude of the flux terms in Theorem \ref{thm:general-fpe} under the limit $\varepsilon \ll 1$.
The total drift is the sum of the $\mathcal{O}(1)$ thermodynamic force and the $\mathcal{O}(\varepsilon^2)$ anomalous correction. The ratio of the anomalous flux to the thermodynamic flux vanishes asymptotically:
\begin{equation}
    \frac{\| \varepsilon^2 \rho \mathfrak{s}_a \|}{\| D_0\rho \nabla \mathcal{V} \|} \sim \varepsilon^2 \to 0.
\end{equation}
Consequently, the drift is $\mathbf{F}_0 = -D_0\nabla \mathcal{V}$ at leading order.
However, the diffusion term $\varepsilon^2 D_0 \nabla \rho$, despite being $\mathcal{O}(\varepsilon^2)$, represents the singular perturbation (the highest spatial derivative) that defines the parabolic nature of the equation. Dropping it would reduce the system to the deterministic Liouville equation $\partial_t \rho = \nabla_g \cdot (D_0\rho \nabla \mathcal{V})$, which cannot satisfy the boundary conditions or describe the Gaussian broadening of the distribution. We therefore retain the diffusion term while neglecting the drift correction, a truncation which is not uniform, since the neglected drift changes the stationary density by the factor of order unity of Proposition~\ref{prop:invariant-measures}~(iii).
\end{proof}
\begin{proposition}[Invariant measures and detailed balance]
\label{prop:invariant-measures}
Let $d\mu^{\varepsilon}_{\text{inv}}$ be a stationary measure of the dynamics for a given $\varepsilon$, and $\rho^{\varepsilon}_{\text{inv}} = d\mu^{\varepsilon}_{\text{inv}}/d\mathfrak{m}_g$ its Radon--Nikodym derivative. The following assertions hold:
\begin{enumerate}[(i)]
    \item {General Case (Broken Detailed Balance):} For finite $\varepsilon$, if the 1-form $\mathfrak{s}_a^\flat$ is not closed ($d\mathfrak{s}_a^\flat \neq 0$, equivalently $\nabla_{[i}\mathfrak{s}_{a,j]} \neq 0$), the measure satisfies a generalised stationarity condition driven by a solenoidal flux $\mathbf{J}_s$:
    \begin{equation}
        \varepsilon^2 D_0 \nabla \ln \rho^{\varepsilon}_{\text{inv}} + D_0\nabla \mathcal{V} - \varepsilon^2 \mathfrak{s}_a = \frac{\mathbf{J}_s}{\rho^{\varepsilon}_{\text{inv}}}.
    \end{equation}
    \item {Conservative Anomaly (Noise-Induced Potential):} If the anomaly is conservative such that $\mathfrak{s}_a = -\nabla \Phi_\xi$, and if the exponent below is integrable against the Riemannian volume, so that the right-hand side normalises to a probability measure, then the measure below is stationary and satisfies detailed balance. The equilibrium measure $d\mu^{\varepsilon}_{\text{eq}}$ takes the Gibbs form shifted by the noise potential:
    \begin{equation}
        d\mu^{\varepsilon}_{\text{eq}} \propto \exp\left( -\frac{\mathcal{V}}{\varepsilon^2} - \frac{\Phi_\xi}{D_0} \right) d\mathfrak{m}_g.
    \end{equation}
    \item {Weak-noise concentration:} The noise-induced potential carries no power of $\varepsilon$ and is therefore not a small correction to the density. What it does not affect is the large-deviation exponent, which is governed by $\mathcal{V}/\varepsilon^2$ alone, so that
    \begin{equation}
        d\mu^{\varepsilon}_{\text{eq}} \;\propto\; \exp\left( -\frac{\mathcal{V}}{\varepsilon^2} \right) \exp\left( -\frac{\Phi_\xi}{D_0} \right) d\mathfrak{m}_g .
    \end{equation}
    The factor $e^{-\Phi_\xi/D_0}$ is of order unity and does not tend to one, so the family does not converge to a finite-temperature Gibbs density. As $\varepsilon \to 0$ it concentrates, in the sense of Laplace, upon the minimisers of $\mathcal{V}$, the noise-induced potential leaving that rate function untouched while modifying the prefactor and the relative weights carried by degenerate minimisers.
\end{enumerate}
\end{proposition}
\begin{proof}
The proof of (i) and (iii) follows from the flux vanishing condition derived in Theorem \ref{thm:general-fpe}, with a flux that retains the anomaly that Corollary \ref{cor:covariant-fpe} discards. For (ii), assume $\mathfrak{s}_a = -\nabla \Phi_\xi$. The stationarity condition $\nabla_g \cdot \mathbf{J} = 0$ is satisfied by the detailed-balance solution $\mathbf{J}=0$, which gives
\begin{equation}
    \varepsilon^2 D_0 \nabla \ln \rho + D_0\nabla \mathcal{V} + \varepsilon^2 \nabla\Phi_\xi = 0 ,
    \qquad\text{hence}\qquad
    \ln \rho = -\frac{\mathcal{V}}{\varepsilon^2} - \frac{\Phi_\xi}{D_0} + C .
\end{equation}
We exponentiate this relation and absorb the integration constant into the normalisation, which gives the shifted stationary density and therefore the measure asserted in part~(ii):
\begin{align}
    \rho &= e^{C}\exp\!\left( -\frac{\mathcal{V}}{\varepsilon^2} - \frac{\Phi_\xi}{D_0} \right) ,
    \nonumber \\
    d\mu^{\varepsilon}_{\text{eq}} &\propto \exp\!\left( -\frac{\mathcal{V}}{\varepsilon^2} - \frac{\Phi_\xi}{D_0} \right) d\mathfrak{m}_g .
\end{align}
The two exponents carry different powers of $\varepsilon$, the thermodynamic one scaling as $\varepsilon^{-2}$ and the noise-induced one being of order unity, so that their ratio vanishes in the weak-noise limit:
\begin{equation}
    \frac{\Phi_\xi/D_0}{\mathcal{V}/\varepsilon^{2}}
    = \frac{\varepsilon^{2}\,\Phi_\xi}{D_0\,\mathcal{V}}
    \xrightarrow[\;\varepsilon \to 0\;]{} 0 .
\end{equation}
The shift is therefore subdominant in the large-deviation exponent, as part~(iii) states, whereas the factor $e^{-\Phi_\xi/D_0}$ itself carries no power of $\varepsilon$ and remains of order unity, so that the stationary family does not reduce to the canonical measure. The exponent is dimensionless, because $[\Phi_\xi] = [D_0] = L^2T^{-1}$ and $\mathcal{V} = \beta\,\mathcal{F}_m[M]$ is measured in units of $k_BT$ by Theorem~\ref{thm:projected-drift}, the factor $D_0$ being carried by the thermodynamic force $F_l = -D_0\partial_l\mathcal{V}$.
\end{proof}
%%%%%%%%%%
\section{Dimensional reduction of the forward Kolmogorov equation}
\label{sec:kolmogorov-reduction}
Having reduced the stochastic flow by projection in Section~\ref{sec:reduction}, we now carry out the second reduction, which proceeds along the other edge of the diagram of Figure~\ref{fig:reduction_diagram} and reaches the manifold without passing through the trajectories at all. Here we project the functional forward Kolmogorov equation associated with the full infinite-dimensional mass SPDE onto $\mathcal{M}_m$ by integrating out the orthogonal degrees of freedom. It does not settle, however, whether both routes reach the same generator, and Section~\ref{sec:commutativity} takes up this question. We show that, under a measure-theoretic analogue of the stiffness hypothesis of Section~\ref{subsubsec:stiffness} (namely the adiabatic concentration of Hypothesis~\ref{hyp:concentration}), this probabilistic reduction yields an effective generator of the same form as the one derived in Section~\ref{subsec:induced-evolution}, and that it differs from that one only in the anomalous drift it carries. This drift has three sources, and the curvature of the drift field is one of them, through the second variation $D^2\mathcal{A}_m$ averaged against the transverse fluctuations. This contribution is zero for the ideal gas, since the drift is then linear in the mass. The Fermi coordinates also depend on the fluctuation, and the projection of the drift inherits that dependence together with the displacement of the mean of the fibre law by the mean curvature of the manifold. Finally, as the weight depends on the base point, the connection induced by the projection differs from the Levi-Civita connection of $g$. Note that the Jacobian $J_\Psi$ of the Fermi coordinates cancels between the conditional profile and the reference measure, and the curvature of the tube acts instead through that mean.
\subsection{The functional forward Kolmogorov equation}
\label{subsec:infinite-kolmogorov}
Let us begin with the evolution equation for the probability law of the cumulative mass field $m_t$. Let $\mu_t$ be the Borel probability measure on the Hilbert space $\mathcal{H}_{\varpi_0}$ that represents the law of the process $m_t$ at time $t$. Here $\varpi_0 \coloneqq (4\pi r^2)^{-1}$, and the norm of this space is equivalent to every weighted norm of Definition~\ref{def:weighted-space} for densities bounded above and below. The dynamics of $m_t$ are governed by the SPDE of Proposition~\ref{prop:mass-dynamics} (under the spherical-mean closure of Assumption~\ref{ass:spherical-closure}).
\begin{remark}[Starting from the It\^{o} representation]
\label{rem:ito-start}
In Section~\ref{sec:reduction} we formulated the dynamics in the Stratonovich sense, in order to preserve the chain rule during the geometric projection. However, in order to construct the Kolmogorov generator, we have to go back to the It\^{o} representation. As we have shown in Lemma~\ref{lemma:ito-strat-correction}, the conversion adds the Stratonovich correction $+\varepsilon^2 \mathfrak{S}_m^\delta$. Since the Stratonovich drift of Proposition~\ref{prop:mass-dynamics} is $\mathcal{A}_m - \varepsilon^2 \mathfrak{S}_m^\delta$, the conversion cancels this term exactly and returns the bare thermodynamic drift $\mathcal{A}_m$ for the Kolmogorov generator. The corresponding It\^{o} SPDE is given by:
\begin{equation}
    dm_t = \mathcal{A}_m(m_t) \, dt + \varepsilon \mathcal{B}_m(m_t) \, d\xi^\delta_t.
\end{equation}
The noise intensity $\varepsilon$ is that of Definition~\ref{def:weak-noise}.
\end{remark}
\begin{lemma}[The mass covariance is trace class]
\label{lemma:mass-covariance-trace}
Let $m$ be a state with $\partial_r m > 0$ upon $(0,L]$ and weight $\varpi = \varpi[m]$, and let $\mathcal{Q}_m = \mathcal{B}_m C_\delta \mathcal{B}_m^*$ of Eq.~\eqref{eq:qm-nonlocal} be the covariance operator upon $\mathcal{H}_\varpi$ of the noise increment $\mathcal{B}_m\,d\xi^\delta_t$, so that
\begin{equation}
    \langle \mathcal{Q}_m u, v \rangle_\varpi \, dt = \mathbb{E}\bigl[\langle \mathcal{B}_m d\xi^\delta_t, u \rangle_\varpi \, \langle \mathcal{B}_m d\xi^\delta_t, v \rangle_\varpi\bigr] , \qquad u, v \in \mathcal{H}_\varpi .
    \label{eq:mass-covariance-form}
\end{equation}
\begin{enumerate}[(i)]
    \item The operator $\mathcal{Q}_m$ is non-negative, self-adjoint and trace class upon $\mathcal{H}_\varpi$, with
    \begin{equation}
        \operatorname{Tr}_\varpi(\mathcal{Q}_m) = 2D_0 \int_0^L V_\delta(r)\,dr .
        \label{eq:mass-covariance-trace}
    \end{equation}
    \item For every bounded linear operator $A$ from $\mathcal{H}_\varpi$ into a separable Hilbert space $\mathcal{K}$, the operator $A \mathcal{Q}_m A^*$ is non-negative, self-adjoint and trace class upon $\mathcal{K}$.
\end{enumerate}
\end{lemma}
\begin{proof}
Throughout, $\{f_j\}$ is a complete orthonormal basis of $\mathcal{H}_\varpi$ and $\sigma = \sqrt{2D_0\,\partial_r m}$ is the amplitude of $\mathcal{B}_m$.
\begin{proofpart}[Non-negativity and self-adjointness]
The right-hand side of~\eqref{eq:mass-covariance-form} is symmetric in $u$ and $v$, so that $\langle \mathcal{Q}_m u, v \rangle_\varpi = \langle u, \mathcal{Q}_m v \rangle_\varpi$, and setting $v = u$ gives
\begin{equation}
    \langle \mathcal{Q}_m u, u \rangle_\varpi \, dt = \mathbb{E}\bigl[\langle \mathcal{B}_m d\xi^\delta_t, u \rangle_\varpi^2\bigr] \ge 0 .
\end{equation}
\end{proofpart}
\begin{proofpart}[The trace series]
Setting $u = v = f_j$ in~\eqref{eq:mass-covariance-form}, summing over $j$, exchanging the sum with the expectation by monotone convergence, and applying Parseval's identity,
\begin{equation}
    \sum_j \langle \mathcal{Q}_m f_j, f_j \rangle_\varpi \, dt
    = \mathbb{E}\Bigl[\sum_j \langle \mathcal{B}_m d\xi^\delta_t, f_j \rangle_\varpi^2\Bigr]
    = \mathbb{E}\bigl[\|\mathcal{B}_m d\xi^\delta_t\|_\varpi^2\bigr] .
\end{equation}
Writing the norm as an integral, exchanging the expectation with it by Fubini's theorem, and using $\mathbb{E}[d\xi^\delta(r,t)^2] = C_\delta(r,r)\,dt$ from Lemma~\ref{lemma:radial-covariance}, then $\sigma^2\varpi = 2D_0$ and $V_\delta(r) = C_\delta(r,r)$ from Definition~\ref{def:kernel-invariants},
\begin{equation}
    \begin{split}
    \mathbb{E}\bigl[\|\mathcal{B}_m d\xi^\delta_t\|_\varpi^2\bigr]
    &= \int_0^L \sigma(r)^2\,\mathbb{E}\bigl[d\xi^\delta(r,t)^2\bigr]\,\varpi(r)\,dr
    \\
    &= \int_0^L \sigma(r)^2\,\varpi(r)\,C_\delta(r,r)\,dr\;dt
    = 2D_0 \int_0^L V_\delta(r)\,dr\;dt .
    \end{split}
\end{equation}
\end{proofpart}
\begin{proofpart}[Finiteness]
By Eq.~\eqref{eq:V-exact}, $V_\delta(r) = 2\pi\int_0^{2r} g_\delta(u)\,u\,du - \pi r^{-2}\int_0^{2r} g_\delta(u)\,u^3\,du$. For $0 \le u \le 2r$ we have $u^3/r^2 \le 4u$, so that
\begin{equation}
    |V_\delta(r)| \le 2\pi\int_0^{2r} |g_\delta(u)|\,u\,du + 4\pi\int_0^{2r} |g_\delta(u)|\,u\,du \le 6\pi\int_0^\infty |g_\delta(u)|\,u\,du ,
\end{equation}
which is finite, $g_\delta$ being continuous with compact support by Definition~\ref{def:mollifier}. The series of the previous part is therefore finite, and Definition~\ref{def:hilbert-schmidt}, with the non-negativity and self-adjointness of the first part, makes $\mathcal{Q}_m$ trace class with the trace~\eqref{eq:mass-covariance-trace}.
\end{proofpart}
\begin{proofpart}[The compressed covariance]
Let $\{h_k\}$ be a complete orthonormal basis of $\mathcal{K}$. The operator $A\mathcal{Q}_m A^*$ is self-adjoint, and $\langle A\mathcal{Q}_m A^* h, h \rangle_{\mathcal{K}} = \langle \mathcal{Q}_m A^* h, A^* h \rangle_\varpi \ge 0$ by the first part. Setting $u = v = A^* h_k$ in~\eqref{eq:mass-covariance-form}, summing over $k$, and applying Parseval's identity in $\mathcal{K}$,
\begin{equation}
    \sum_k \langle A\mathcal{Q}_m A^* h_k, h_k \rangle_{\mathcal{K}} \, dt
    = \mathbb{E}\Bigl[\sum_k \langle A\mathcal{B}_m d\xi^\delta_t, h_k \rangle_{\mathcal{K}}^2\Bigr]
    = \mathbb{E}\bigl[\|A\mathcal{B}_m d\xi^\delta_t\|_{\mathcal{K}}^2\bigr]
    \le \|A\|^2\,\operatorname{Tr}_\varpi(\mathcal{Q}_m)\,dt ,
\end{equation}
which is finite by the previous part, and Definition~\ref{def:hilbert-schmidt} makes $A\mathcal{Q}_m A^*$ trace class.
\end{proofpart}
\end{proof}
\begin{definition}[Functional adjoint evolution]
\label{def:adjoint-evolution}
The time evolution of the measure $\mu_t$ is defined by duality with the observables of the system, consistent with the definition of the weak solution given in Definition \ref{def:ffpe}. Let $\Phi \in C_b^2(\mathcal{H}_{\varpi})$ be a cylindrical test functional. The evolution is given by:
\begin{equation}
    \frac{d}{dt} \int_{\mathcal{H}_{\varpi}} \Phi(m) \, d\mu_t(m) = \int_{\mathcal{D}_\alpha} \mathcal{L}_\varepsilon \Phi(m) \, d\mu_t(m),
\end{equation}
in which $\mathcal{D}_\alpha$ is the set of Definition~\ref{def:ffpe}, read in the mass representation through $\rho = (4\pi r^2)^{-1}\partial_r m$. The operator $\mathcal{L}_\varepsilon$ is the infinitesimal generator for the It\^{o} process:
\begin{equation}
    \mathcal{L}_\varepsilon \Phi(m) = \langle \mathcal{A}_m(m), D\Phi(m) \rangle_{\varpi} + \frac{\varepsilon^2}{2} \operatorname{Tr}_{\varpi} \left( \mathcal{Q}_m(m) D^2\Phi(m) \right).
\end{equation}
Here, $\mathcal{Q}_m = \mathcal{B}_m C_\delta \mathcal{B}_m^*$ is the diffusion covariance operator incorporating the finite-$\delta$ kernel.
\end{definition}
\begin{remark}[The choice of topology]
\label{rem:weighted-topology}
In Section~\ref{sec:preliminaries} we introduced the general $L^2$ topology, but in this section we work strictly in the weighted Hilbert space $\mathcal{H}_{\varpi}$ of Section~\ref{subsubsec:orthogonal-decomposition}. This change is needed, since the projector $\Pi_{\mathbf{X}}^{\parallel}$ (Definition~\ref{def:explicit-projector}) is orthogonal with respect to the weighted inner product $\langle \cdot, \cdot \rangle_\varpi$, and Lemma~\ref{lemma:metric-selection} shows that with this weight the projected noise is a standard Wiener process on $(\mathcal{X}, g)$. If we used the standard $L^2$ inner product here, off-diagonal cross-terms would appear in the projected diffusion tensor, and the geometric interpretation would become more complicated.
\end{remark}
\subsection{Tubular coordinates and the Fermi map}
\label{subsec:tubular-coordinates}
To separate the slow dynamics on the manifold from the fast transverse fluctuations, we must parameterise the Hilbert space in a neighbourhood of the manifold. To this end, we use the orthogonal decomposition of Section~\ref{subsubsec:orthogonal-decomposition} to construct a curvilinear coordinate system.
\subsubsection{The geometry of the tubular neighbourhood}
Let us recall from Section~\ref{subsubsec:induced-manifold} that the mass manifold $\mathcal{M}_m$ is the image of the embedding $M: \mathcal{X} \to \mathcal{H}_{\varpi}$. The tangent space $T_{\mathbf{X}}\mathcal{M}_m$ is then spanned by the basis vectors $\mathbf{e}_k^m(\mathbf{X}) = \partial_{X^k} M(\mathbf{X})$.
\begin{definition}[The Fermi map]
\label{def:fermi-map}
Let $\mathcal{N}$ be the weighted normal bundle of Definition~\ref{def:normal-bundle}. We define the \textit{Fermi map} $\Psi: \mathcal{N} \to \mathcal{H}_{\varpi}$ by the affine addition:
\begin{equation}
    \Psi(\mathbf{X}, \eta) \coloneqq M(\cdot; \mathbf{X}) + \eta,
\end{equation}
where $\mathbf{X} \in \mathcal{X}$ are the collective coordinates and $\eta \in \mathcal{H}_{\varpi}$ satisfies the orthogonality condition $\langle \eta, \mathbf{e}_k^m(\mathbf{X}) \rangle_{\varpi} = 0$ for all $k=1, \dots, d$.
\end{definition}
The Fermi map, in turn, gives the coordinate transformation that we need for the reduction. By Lemma~\ref{lemma:tubular-existence}, $\Psi$ is a diffeomorphism from a neighbourhood of the zero section of $\mathcal{N}$ onto the tubular neighbourhood $U_{\varrho_{\mathrm{tub}}}$ of the manifold, and by Lemma~\ref{lemma:orthogonality-condition} the coordinates $(\mathbf{X}, \eta) = \Psi^{-1}(m)$ of every $m \in U_{\varrho_{\mathrm{tub}}}$ are those fixed by the weighted orthogonality condition. For the Kolmogorov reduction we need the explicit structure of the Jacobian of $\Psi$, in order to transform the volume measure and the differential operators.
\begin{lemma}[The Fermi Jacobian block structure]
\label{lemma:fermi-jacobian}
Let $p_0 = (\mathbf{X}, 0)$ be a point on the zero section of the normal bundle $\mathcal{N}$. Let $(\mathbf{v}, \zeta) \in T_{p_0}\mathcal{N}$ be a tangent vector, where $\mathbf{v} \in T_{\mathbf{X}}\mathcal{X} \cong \mathbb{R}^d$ and $\zeta \in \mathcal{N}_{\mathbf{X}} = (T_{\mathbf{X}}\mathcal{M}_m)^{\perp}$.
The Fr\'{e}chet derivative $D\Psi_{(\mathbf{X}, 0)}$ is the linear operator defined by:
\begin{equation}
    D\Psi_{(\mathbf{X}, 0)} (\mathbf{v}, \zeta) = \mathcal{J}_{\iota_m}(\mathbf{X}) \mathbf{v} + \zeta.
\end{equation}
The operator $\mathcal{J}_{\iota_m}$ is the same pushforward that factorises the tangential projector in Lemma~\ref{lemma:projector-jacobian}. Consequently, at the manifold interface, the Jacobian admits a block-diagonal decomposition with respect to the splitting $\mathcal{H}_\varpi = T_{\mathbf{X}}\mathcal{M}_m \oplus (T_{\mathbf{X}}\mathcal{M}_m)^{\perp}$:
\begin{equation}
    D\Psi_{(\mathbf{X}, 0)} = \mathcal{J}_{\iota_m}(\mathbf{X}) \oplus \mathcal{I}_{\perp},
\end{equation}
where $\mathcal{J}_{\iota_m}$ maps parameter space to the tangent space, and $\mathcal{I}_{\perp}$ is the identity on the normal fibre.
\end{lemma}
\begin{proof}
\begin{proofpart}[The Fr\'{e}chet derivative]
We proceed by the definition of the Fr\'{e}chet derivative. We must show that the linear operator $L(\mathbf{v}, \zeta) \coloneqq \mathcal{J}_{\iota_m}(\mathbf{X})\mathbf{v} + \zeta$ satisfies the limit condition:
\begin{equation}
    \lim_{\|(\mathbf{v}, \zeta)\| \to 0} \frac{\| \Psi(\mathbf{X}+\mathbf{v}, \zeta) - \Psi(\mathbf{X}, 0) - L(\mathbf{v}, \zeta) \|_{\varpi}}{\|(\mathbf{v}, \zeta)\|} = 0.
\end{equation}
Substituting the definition $\Psi(\mathbf{X}, \eta) = M(\mathbf{X}) + \eta$:
\begin{align}
    \Psi(\mathbf{X}+\mathbf{v}, \zeta) - \Psi(\mathbf{X}, 0) &= M(\mathbf{X}+\mathbf{v}) + \zeta - M(\mathbf{X}).
\end{align}
Since the embedding $M$ is a smooth immersion (Definition \ref{def:immersion}), it is Fr\'{e}chet differentiable. Its Taylor expansion is:
\begin{equation}
    M(\mathbf{X}+\mathbf{v}) - M(\mathbf{X}) = D M(\mathbf{X}) \mathbf{v} + R(\mathbf{v}),
\end{equation}
where the remainder satisfies $\|R(\mathbf{v})\|_{\varpi} = o(\|\mathbf{v}\|)$. Recall from Lemma~\ref{lemma:projector-jacobian} that the derivative of the embedding is precisely the Jacobian $\mathcal{J}_{\iota_m}(\mathbf{X}) \mathbf{v}$.
Substituting this into the numerator of the difference quotient:
\begin{align}
    \| \Psi(\mathbf{X}+\mathbf{v}, \zeta) - \Psi(\mathbf{X}, 0) - L(\mathbf{v}, \zeta) \|_{\varpi} &= \| \mathcal{J}_{\iota_m}(\mathbf{X}) \mathbf{v} + R(\mathbf{v}) + \zeta - (\mathcal{J}_{\iota_m}(\mathbf{X})\mathbf{v} + \zeta) \|_{\varpi} \nonumber \\
    &= \| R(\mathbf{v}) \|_{\varpi}.
\end{align}
Dividing by the norm of the increment $\|(\mathbf{v}, \zeta)\| \ge \|\mathbf{v}\|$:
\begin{equation}
    \lim \frac{\| R(\mathbf{v}) \|_{\varpi}}{\|(\mathbf{v}, \zeta)\|} \le \lim \frac{\| R(\mathbf{v}) \|_{\varpi}}{\|\mathbf{v}\|} = 0.
\end{equation}
Thus, the derivative is $L$.
\end{proofpart}
\begin{proofpart}[The block structure]
To prove the block structure, we verify the ranges. The term $\mathcal{J}_{\iota_m}(\mathbf{X}) \mathbf{v} = \sum v^k \mathbf{e}_k^m$ lies in the tangent space $T_{\mathbf{X}}\mathcal{M}_m$. The term $\zeta$ is a tangent vector to the fibre $\mathcal{N}_{\mathbf{X}}$, so it lies in $(T_{\mathbf{X}}\mathcal{M}_m)^{\perp}$. Thus, the operator maps the direct sum of the input spaces to the direct sum of the orthogonal subspaces of $\mathcal{H}_\varpi$ without mixing.
\end{proofpart}
\end{proof}
\begin{remark}[Necessity of the block structure]
\label{rem:jacobian-necessity}
Note that the Jacobian (evaluated at $\eta=0$) is block-diagonal, and this structure is essential for the reduction. To leading order, variations in the order parameters ($\mathbf{X}$) and variations in the fluctuation ($\eta$) are then orthogonal, and we can factorise the integration measure in the vicinity of the manifold, and we separate, in this way, the calculation of the marginal distribution (on $\mathcal{M}_m$) from the integration over the normal modes (on the fibres).
\end{remark}
\subsubsection{Disintegration of the probability measure}
\label{subsubsec:measure-disintegration}
The Fermi map $\Psi$ decomposes the state space $\mathcal{H}_{\varpi}$, and hence, for the reduction, it gives the coordinate framework. Two technical difficulties arise, however, when the probability measure $\mu_t$ is decomposed:
\begin{enumerate}[(i)]
    \item \emph{Topological:} since $\mathcal{H}_{\varpi}$ is infinite-dimensional, it carries no non-zero, locally finite, translation-invariant Borel measure.     \item \emph{Geometric:} the subspace of orthogonal fluctuations $\mathcal{N}_{\mathbf{X}}$ rotates as the order parameter $\mathbf{X}$ evolves, and so the normal bundle $\mathcal{N}$ is ``twisted'' and the space is not globally a simple product.
\end{enumerate}
To solve the second problem, we first ``untwist'' the geometry locally. In this procedure, known as \textit{local trivialisation}, we identify the varying normal fibres with a single, fixed Hilbert space (the model fibre), and this fixed space then carries the Gaussian reference measure.
\begin{definition}[Local trivialisation and the model fibre]
\label{def:local-trivialisation}
We fix an abstract separable Hilbert space $\mathcal{H}_{\perp}$, termed the \textit{model normal fibre}. For $\mathbf{X}_0 \in \mathcal{X}$, a \emph{local trivialisation} of the normal bundle $\mathcal{N}$ of Definition~\ref{def:normal-bundle} over an open neighbourhood $U \subset \mathcal{X}$ of $\mathbf{X}_0$ is a smooth map
\begin{equation}
    \chi: \mathcal{N}|_{U} \to U \times \mathcal{H}_{\perp}, \quad (\mathbf{X}, \eta) \mapsto (\mathbf{X}, \mathcal{R}_{\mathbf{X}}\eta),
\end{equation}
in which each $\mathcal{R}_{\mathbf{X}}: \mathcal{N}_{\mathbf{X}} \to \mathcal{H}_{\perp}$ is a unitary isomorphism, the fibre $\mathcal{N}_{\mathbf{X}} = (T_{\mathbf{X}}\mathcal{M}_m)^{\perp_{\varpi}}$ carrying the weighted inner product at the base point $\mathbf{X}$.
\end{definition}
Such a trivialisation exists about every point. The tangent spaces form a subbundle of rank $d$ of the product $\mathcal{X} \times \mathcal{H}_\varpi$, and the normal fibres are their orthogonal complements in the inner products $\langle\cdot,\cdot\rangle_{\varpi(\cdot;\mathbf{X})}$, which depend smoothly upon $\mathbf{X}$. Theorem~3.1 of Chapter~VII of Lang~\cite{lang1999fundamentals}, together with the construction of orthogonal complements which follows it, shows that such a complement is a Hilbert bundle, which is to say that it admits local trivialisations by unitary isomorphisms. The fibres are infinite-dimensional and separable, and hence each of them is unitarily isomorphic to $\mathcal{H}_\perp$. In the remainder of this section we work within such a trivialisation, and we treat the state space locally as the product $U \times \mathcal{H}_{\perp}$.
\begin{example}[Trivialisation of radial fluctuations]
\label{ex:trivialisation}
To clarify the necessity of the trivialisation, consider a particle constrained to a circle $\mathcal{M} = S^1$ in $\mathbb{R}^2$, parameterised by angle $\theta$:
\begin{enumerate}[(i)]
    \item \emph{Twisted frame:} The normal space $\mathcal{N}_\theta$ is the radial line spanned by the vector $\mathbf{n}(\theta) = (\cos \theta, \sin \theta)^T$. A normal fluctuation is $\eta = \xi \mathbf{n}(\theta)$. The direction of this fluctuation depends on the state $\theta$. Consequently, the covariance matrix of the noise in the ambient coordinates is $Q(\theta) = \sigma^2 \mathbf{n}(\theta)\mathbf{n}(\theta)^T$, which varies along the manifold.
    \item \emph{Trivialisation:} We choose the model fibre $\mathcal{H}_{\perp} \cong \mathbb{R}$. The map $\mathcal{R}_{\theta}$ is the rotation that aligns $\mathbf{n}(\theta)$ with the fixed basis vector $\mathbf{e}_1$.
    \item \emph{Result:} In the trivialised frame, the fluctuation is simply the scalar amplitude $\xi \in \mathbb{R}$. The transformed covariance is then given by $\mathcal{R}_{\theta} Q(\theta) \mathcal{R}_{\theta}^* = \sigma^2$.
\end{enumerate}
In the trivialised frame, the noise appears as a \textit{constant} Gaussian measure on the model fibre, even though it rotates geometrically in the ambient space.
\end{example}
\begin{remark}[Geometric interpretation]
From a physical point of view, the model fibre $\mathcal{H}_{\perp}$ represents the internal reference frame of the noise. While the physical orientation of the orthogonal fluctuations changes as the system moves along the curved manifold $\mathcal{M}_m$, the trivialisation map $\mathcal{R}_{\mathbf{X}}$ acts as a \textit{moving frame} that rotates these fluctuations into a standard orientation. We can therefore define a single Gaussian reference measure for the noise, taken at a reference point $\mathbf{X}_0$, while the variation of the projected covariance along the manifold is controlled by Assumption~\ref{ass:fibre-homogeneity}.
\end{remark}
With the product structure established, we can now define the reference measure required to construct the probability densities.
\begin{definition}[The reference measure $\lambda$]
\label{def:reference-measure}
We define the reference measure $\lambda$ on the trivialised product space $U \times \mathcal{H}_{\perp}$ as the tensor product $\lambda \coloneqq \lambda_{\mathcal{X}} \otimes \lambda_{\perp}$, where:
\begin{enumerate}[(i)]
    \item $\lambda_{\mathcal{X}}$ is the Riemannian volume measure on the parameter space, induced by the metric $g(\mathbf{X})$ defined in Definition \ref{def:metric-tensor}:
    \begin{equation}
        d\lambda_{\mathcal{X}}(\mathbf{X}) \coloneqq d\mathfrak{m}_g(\mathbf{X}) = \sqrt{\det g(\mathbf{X})} \, d\mathbf{X}.
    \end{equation}
    \item $\lambda_{\perp}$ is the centred Gaussian measure on the model fibre $\mathcal{H}_{\perp}$ with covariance $s\,\mathcal{Q}_{\perp}(\mathbf{X}_0)$, for a fixed reference point $\mathbf{X}_0 \in U$ and a fixed scale $s > 0$, where, for every $\mathbf{X} \in U$,
    \begin{equation}
        \mathcal{Q}_{\perp}(\mathbf{X}) \coloneqq \mathcal{R}_{\mathbf{X}} \Pi_{\mathbf{X}}^{\perp} \mathcal{Q}_m \Pi_{\mathbf{X}}^{\perp} \mathcal{R}_{\mathbf{X}}^*
    \end{equation}
    is the covariance of the noise projected onto the normal fibre and carried into the model fibre, $\mathcal{Q}_m$ being evaluated at the manifold point $M(\mathbf{X})$:
    \begin{equation}
        \lambda_{\perp} \coloneqq \mathrm{N}\bigl(0, s\,\mathcal{Q}_{\perp}(\mathbf{X}_0)\bigr).
    \end{equation}
\end{enumerate}
\end{definition}
\noindent{}The measure $\lambda_\perp$ exists by Proposition~2.18 of Da~Prato and Zabczyk~\cite{daprato2014stochastic}, since $\mathcal{Q}_\perp(\mathbf{X}_0)$ is trace class by Lemma~\ref{lemma:mass-covariance-trace}\,(ii).
\begin{assumption}[Homogeneity of the fibre covariance]
\label{ass:fibre-homogeneity}
For every $\mathbf{X} \in U$, the covariances of Definition~\ref{def:reference-measure} satisfy $\mathcal{Q}_\perp(\mathbf{X})^{1/2}(\mathcal{H}_\perp) = \mathcal{Q}_\perp(\mathbf{X}_0)^{1/2}(\mathcal{H}_\perp) \eqqcolon \mathcal{H}_{\perp,0}$, and the operator
\begin{equation}
    \bigl(\mathcal{Q}_\perp(\mathbf{X}_0)^{-1/2}\mathcal{Q}_\perp(\mathbf{X})^{1/2}\bigr)\bigl(\mathcal{Q}_\perp(\mathbf{X}_0)^{-1/2}\mathcal{Q}_\perp(\mathbf{X})^{1/2}\bigr)^* - I
    \label{eq:fibre-homogeneity}
\end{equation}
is Hilbert--Schmidt upon $\mathcal{H}_{\perp,0}$.
\end{assumption}
\begin{remark}[The condition of the Feldman--H\'ajek theorem]
\label{rem:fibre-homogeneity}
Assumption~\ref{ass:fibre-homogeneity} is the condition of the Feldman--H\'ajek theorem (Theorem~2.25 of Da~Prato and Zabczyk~\cite{daprato2014stochastic}) for the centred Gaussian measures with covariances $\mathcal{Q}_\perp(\mathbf{X})$ and $\mathcal{Q}_\perp(\mathbf{X}_0)$ to be equivalent, the square roots being those of that theorem. Note that the covariance varies along the manifold through the rotation of the normal fibre and, at finite $\delta$, through the dependence of the weight upon the base point, and the assumption requires that variation to be of Hilbert--Schmidt size relative to $\mathcal{Q}_\perp(\mathbf{X}_0)$.
\end{remark}
\begin{assumption}[Absolute continuity of the law]
\label{ass:absolute-continuity}
The law $\mu_t$, read in the trivialised chart, is absolutely continuous with respect to the reference measure $\lambda$ of Definition~\ref{def:reference-measure}, $\mu_t \ll \lambda$, for every $t$ in the interval considered.
\end{assumption}
\begin{remark}[Absolute continuity and densities]
\label{rem:absolute-continuity}
Assumption~\ref{ass:absolute-continuity} is a hypothesis. Girsanov's theorem does not give it, and every statement obtained from the disintegration below is conditional on it. Indeed, Girsanov's theorem concerns the equivalence of path measures under a change of drift at fixed diffusion coefficient, and it does not deliver absolute continuity of the time-$t$ law with respect to a product reference measure upon an infinite-dimensional state space. The hypothesis also fixes $s$, since by the Feldman--H\'ajek theorem two centred Gaussian measures upon an infinite-dimensional space are either equivalent or singular, and $Q$ and $cQ$ give singular Gaussian measures unless $c = 1$, the operator of condition~\eqref{eq:fibre-homogeneity} being then $(c-1)I$, which is Hilbert--Schmidt only for $c = 1$. The conditional law of Hypothesis~\ref{hyp:concentration} has covariance $\varepsilon^2\mathbf{C}_\perp(\mathbf{X})$, which Hypothesis~\ref{hyp:kinematic-covariance} fixes as $\varepsilon^2\tau\,\mathcal{Q}_\perp(\mathbf{X})$, and so Assumption~\ref{ass:absolute-continuity} can hold only where $\varepsilon^2\tau(\mathbf{X}) = s$ for every $\mathbf{X} \in U$, that is only upon charts where the transverse relaxation time is constant, and then only under Assumption~\ref{ass:fibre-homogeneity} and only if the mean displacement $\hat\eta(\mathbf{X})$ of Hypothesis~\ref{hyp:concentration} lies in $\mathcal{H}_{\perp,0}$, which is condition~(ii) of the Feldman--H\'ajek theorem, Theorem~2.25 of~\cite{daprato2014stochastic}. The Radon--Nikodym theorem then gives a density $\varrho_t: U \times \mathcal{H}_{\perp} \to [0,\infty)$, unique up to $\lambda$-null sets, such that
\begin{equation}
    d\mu_t(\mathbf{X}, \eta) = \varrho_t(\mathbf{X}, \eta) \, d\lambda(\mathbf{X}, \eta).
\end{equation}
In this framework, the density $\varrho_t$ represents the re-weighting of the Gaussian noise background by the nonlinear drift forces.
\end{remark}
We now introduce the machinery to separate the marginal and conditional statistics of this density.
\begin{definition}[Disintegration of measures]
\label{def:disintegration}
Let $\mathcal{E}$ be a Polish space and $\mu$ a probability measure on $\mathcal{E}$. Let $\pi: \mathcal{E} \to \mathcal{S}$ be a measurable map to a Polish space $\mathcal{S}$. A \textit{disintegration} of $\mu$ with respect to $\pi$ is a pair $(\bar{\mu}, \{\nu^{\mathbf{X}}\}_{\mathbf{X} \in \mathcal{S}})$ where:
\begin{enumerate}[(i)]
    \item $\bar{\mu} = \mu \circ \pi^{-1}$ is the pushforward measure on $\mathcal{S}$ (the marginal).
    \item $\{\nu^{\mathbf{X}}\}_{\mathbf{X}}$ is a family of probability measures on $\mathcal{E}$ such that for every bounded measurable $f$:
    \begin{equation}
        \int_{\mathcal{E}} f(u) \, d\mu(u) = \int_{\mathcal{S}} \left( \int_{\pi^{-1}(\mathbf{X})} f(u) \, d\nu^{\mathbf{X}}(u) \right) d\bar{\mu}(\mathbf{X}).
    \end{equation}
\end{enumerate}
\end{definition}
Let us now prove that, for measures admitting a density with respect to a product reference (Definition~\ref{def:reference-measure}), the disintegration takes a specific functional form.
\begin{theorem}[Disintegration theorem for product-like spaces]
\label{thm:disintegration}
Let $\mathcal{E} \cong \mathcal{S} \times \mathcal{Y}$ be a product of Polish spaces equipped with the reference product measure $\lambda = \lambda_{\mathcal{S}} \otimes \lambda_{\mathcal{Y}}$. Let $\mu$ be a probability measure on $\mathcal{E}$ such that $\mu \ll \lambda$, with Radon--Nikodym derivative $\varrho(x, y) = \frac{d\mu}{d\lambda}(x,y)$.
Then a disintegration, unique up to $\bar{\mu}$-null sets, exists in which the marginal measure $\bar{\mu}$ on $\mathcal{S}$ and the conditional measures $\{\nu^{x}\}_{x \in \mathcal{S}}$ on $\mathcal{Y}$ are defined by the density relations:
\begin{enumerate}[(i)]
    \item The marginal density $\bar{\varrho}(x)$ is the integral of the full density over the fibre with respect to the fibre reference:
    \begin{equation}
        \bar{\varrho}(x) = \int_{\mathcal{Y}} \varrho(x, y) \, d\lambda_{\mathcal{Y}}(y).
    \end{equation}
    Then $d\bar{\mu}(x) = \bar{\varrho}(x) \, d\lambda_{\mathcal{S}}(x)$.
    \item The conditional measure $\nu^{x}$ on the fibre is given by the normalised profile:
    \begin{equation}
        d\nu^{x}(y) = \frac{\varrho(x, y)}{\bar{\varrho}(x)} \, d\lambda_{\mathcal{Y}}(y).
    \end{equation}
\end{enumerate}
This decomposition holds for $\bar{\mu}$-almost every $x$.
\end{theorem}
\begin{proof}
We construct the disintegration using the properties of the Radon--Nikodym density.
\begin{proofpart}[Marginal consistency]
Let $A \subset \mathcal{S}$ be a measurable set. By definition of the pushforward and the density:
\begin{align}
    \bar{\mu}(A) &= \mu(A \times \mathcal{Y}) = \int_{A \times \mathcal{Y}} \varrho(x, y) \, d\lambda(x, y) \nonumber \\
    &= \int_A \left( \int_{\mathcal{Y}} \varrho(x, y) \, d\lambda_{\mathcal{Y}}(y) \right) d\lambda_{\mathcal{S}}(x) \quad (\text{Fubini's Theorem}).
\end{align}
Defining $\bar{\varrho}(x) \coloneqq \int_{\mathcal{Y}} \varrho(x, y) \, d\lambda_{\mathcal{Y}}(y)$, we see that $\bar{\mu}(A) = \int_A \bar{\varrho}(x) \, d\lambda_{\mathcal{S}}(x)$. Thus, $\bar{\varrho}(x)$ is the density of the marginal measure with respect to the base reference $\lambda_{\mathcal{S}}$.
\end{proofpart}
\begin{proofpart}[Disintegration formula]
Let $f(x, y)$ be a bounded measurable function. We integrate against the full measure:
\begin{align}
    \int_{\mathcal{E}} f(x, y) \, d\mu &= \int_{\mathcal{S}} \int_{\mathcal{Y}} f(x, y) \varrho(x, y) \, d\lambda_{\mathcal{Y}}(y) \, d\lambda_{\mathcal{S}}(x).
\end{align}
We multiply and divide the integrand by the marginal density $\bar{\varrho}(x)$ (noting that $\bar{\varrho}(x) > 0$ on the support of $\bar{\mu}$):
\begin{align}
    \int_{\mathcal{E}} f \, d\mu &= \int_{\mathcal{S}} \bar{\varrho}(x) \left[ \int_{\mathcal{Y}} f(x, y) \frac{\varrho(x, y)}{\bar{\varrho}(x)} \, d\lambda_{\mathcal{Y}}(y) \right] d\lambda_{\mathcal{S}}(x).
\end{align}
The term in brackets is exactly the expectation of $f$ with respect to the measure $d\nu^{x}(y) \coloneqq \frac{\varrho(x,y)}{\bar{\varrho}(x)} d\lambda_{\mathcal{Y}}(y)$. The outer integral is with respect to $d\bar{\mu}(x) = \bar{\varrho}(x) d\lambda_{\mathcal{S}}(x)$. Thus, we recover the disintegration formula:
\begin{equation}
    \int_{\mathcal{E}} f \, d\mu = \int_{\mathcal{S}} \left( \int_{\mathcal{Y}} f(x, y) \, d\nu^{x}(y) \right) d\bar{\mu}(x).
\end{equation}
\end{proofpart}
\end{proof}
We now finalise the construction by applying this theorem to the trivialised tubular neighbourhood.
\begin{proposition}[Application to the tubular neighbourhood]
\label{prop:measure-decomposition}
Let $\Psi$ be the Fermi map and $\chi$ the local trivialisation defined in Definition \ref{def:local-trivialisation}. We apply Theorem \ref{thm:disintegration} to the measure $\mu_t$ in this chart. The disintegration of $\mu_t$ is given by:
\begin{equation}
    d\mu_t(\mathbf{X}, \eta) = \bar{\varrho}_t(\mathbf{X}) \, d\mathfrak{m}_g(\mathbf{X}) \otimes d\nu_t^{\mathbf{X}}(\eta),
\end{equation}
where $\bar{\varrho}_t$ is the marginal density on the manifold and $\nu_t^{\mathbf{X}}$ is the conditional probability measure on the model fibre $\mathcal{H}_{\perp}$.
\end{proposition}
\begin{proof}
We proceed by direct substitution of the reference measure structure into the Radon--Nikodym definition.
\begin{proofpart}[Representation of the full measure]
By Assumption~\ref{ass:absolute-continuity} and Remark~\ref{rem:absolute-continuity}, the measure $\mu_t$ is defined by its density $\varrho_t$ against the reference $\lambda$:
\begin{equation}
    d\mu_t(\mathbf{X}, \eta) = \varrho_t(\mathbf{X}, \eta) \, d\lambda(\mathbf{X}, \eta).
\end{equation}
Substitute the explicit tensor definition of the reference measure $\lambda = \lambda_{\mathcal{X}} \otimes \lambda_{\perp}$ from Definition \ref{def:reference-measure}:
\begin{equation}
    d\mu_t(\mathbf{X}, \eta) = \varrho_t(\mathbf{X}, \eta) \, \left( d\mathfrak{m}_g(\mathbf{X}) \otimes d\lambda_{\perp}(\eta) \right).
\end{equation}
\end{proofpart}
\begin{proofpart}[Construction of the marginal]
We define the marginal density $\bar{\varrho}_t(\mathbf{X})$ as the projection of the full density onto the fibre reference measure:
\begin{equation}
    \bar{\varrho}_t(\mathbf{X}) \coloneqq \int_{\mathcal{H}_{\perp}} \varrho_t(\mathbf{X}, \eta) \, d\lambda_{\perp}(\eta).
\end{equation}
This allows us to rewrite the full density identically as:
\begin{equation}
    \varrho_t(\mathbf{X}, \eta) = \bar{\varrho}_t(\mathbf{X}) \cdot \frac{\varrho_t(\mathbf{X}, \eta)}{\bar{\varrho}_t(\mathbf{X})}.
\end{equation}
\end{proofpart}
\begin{proofpart}[Identification of the conditional measure]
Substitute this factorisation back into the measure expression:
\begin{align}
    d\mu_t(\mathbf{X}, \eta) &= \left[ \bar{\varrho}_t(\mathbf{X}) \frac{\varrho_t(\mathbf{X}, \eta)}{\bar{\varrho}_t(\mathbf{X})} \right] \, d\mathfrak{m}_g(\mathbf{X}) \otimes d\lambda_{\perp}(\eta) \nonumber \\
    &= \left( \bar{\varrho}_t(\mathbf{X}) \, d\mathfrak{m}_g(\mathbf{X}) \right) \otimes \left( \frac{\varrho_t(\mathbf{X}, \eta)}{\bar{\varrho}_t(\mathbf{X})} \, d\lambda_{\perp}(\eta) \right).
\end{align}
We identify the second factor as the conditional measure $\nu_t^{\mathbf{X}}$ on the fibre:
\begin{equation}
    d\nu_t^{\mathbf{X}}(\eta) \coloneqq \frac{\varrho_t(\mathbf{X}, \eta)}{\bar{\varrho}_t(\mathbf{X})} \, d\lambda_{\perp}(\eta).
\end{equation}
By construction, this measure is normalised for every $\mathbf{X}$:
\begin{equation}
    \int_{\mathcal{H}_{\perp}} d\nu_t^{\mathbf{X}}(\eta) = \int_{\mathcal{H}_{\perp}} \frac{\varrho_t(\mathbf{X}, \eta)}{\bar{\varrho}_t(\mathbf{X})} \, d\lambda_{\perp}(\eta) = \frac{1}{\bar{\varrho}_t(\mathbf{X})} \int_{\mathcal{H}_{\perp}} \varrho_t(\mathbf{X}, \eta) \, d\lambda_{\perp}(\eta) = 1,
\end{equation}
where the last equality uses the definition $\bar{\varrho}_t(\mathbf{X}) = \int_{\mathcal{H}_{\perp}} \varrho_t(\mathbf{X}, \eta) \, d\lambda_{\perp}(\eta)$. Thus, the disintegration is exact:
\begin{equation}
    d\mu_t = \bar{\varrho}_t(\mathbf{X}) \, d\mathfrak{m}_g(\mathbf{X}) \otimes d\nu_t^{\mathbf{X}}(\eta).
\end{equation}
\end{proofpart}
\end{proof}
%%%%%%%%%%%
\subsection{The geometric measure decomposition}
\label{subsec:measure-decomposition}
In Proposition~\ref{prop:measure-decomposition} we established the disintegration of the physical measure $\mu_t$ with respect to the reference measure $\lambda = \lambda_{\mathcal{X}} \otimes \lambda_{\perp}$ in the trivialised coordinate chart. To use this decomposition in the functional forward Kolmogorov equation (Definition~\ref{def:adjoint-evolution}), we need the relationship between the Fermi coordinates and the reference measure $\lambda$. The geometric Jacobian of the Fermi map relates the Fermi coordinates to the reference measure, through the \emph{Jacobian determinant} of the coordinate transformation, although it cancels from every fibre average (Remark~\ref{rem:volume-no-drift}). The Fermi map $\Psi: U \times \mathcal{H}_{\perp} \to \mathcal{H}_{\varpi}$ parametrises the physical state $m = \Psi(\mathbf{X}, \eta)$, and the linear properties of its derivative $D\Psi$ govern the distortion of the phase space volume. To determine the density $\varrho_t$, we therefore compute the Fr\'{e}chet derivative of the Fermi map first, and then its determinant.
\subsubsection{The shape operator and curvature}
In the trivialised chart of Definition~\ref{def:local-trivialisation}, the fibre coordinate $\eta$ belongs to the fixed model fibre $\mathcal{H}_\perp$, not to the normal space at $\mathbf{X}$, and the Fermi map of Definition~\ref{def:fermi-map} reads $\Psi(\mathbf{X}, \eta) = M(\mathbf{X}) + \mathcal{R}_{\mathbf{X}}^{-1} \eta$. This is the form that we use from here on. To compute its derivative with respect to the base point $\mathbf{X}$, we analyse the variation of the inverse trivialisation map $\mathcal{R}_{\mathbf{X}}^{-1}$. This operator maps a fixed vector from the model fibre $\mathcal{H}_{\perp}$ to a rotating vector in the physical normal bundle $(T_{\mathbf{X}}\mathcal{M}_m)^{\perp_{\varpi}}$. The shape operator encodes the geometric dependence of the normal fibres on the tangent directions.
\begin{definition}[The shape operator]
\label{def:shape-operator}
Let $\mathbf{n} \in (T_{\mathbf{X}}\mathcal{M}_m)^{\perp_{\varpi}}$ be a normal vector at $\mathbf{X}$, and let $\tilde{\mathbf{n}}$ be a smooth local extension of it, defined on an open neighbourhood of $\mathbf{X}$ in $\mathcal{M}_m$, with $\tilde{\mathbf{n}}(\mathbf{X}) = \mathbf{n}$ and $\tilde{\mathbf{n}}(\mathbf{Y}) \in (T_{\mathbf{Y}}\mathcal{M}_m)^{\perp_{\varpi}}$. The \textit{shape operator}, or Weingarten map, $A_{\mathbf{n}}: T_{\mathbf{X}}\mathcal{M}_m \to T_{\mathbf{X}}\mathcal{M}_m$ is the linear operator defined by the tangential projection of the directional derivative of this normal field. For any tangent vector $\mathbf{v} \in T_{\mathbf{X}}\mathcal{M}_m$:
\begin{equation}
    A_{\mathbf{n}} \mathbf{v} \coloneqq - \Pi_{\mathbf{X}}^{\parallel} \left( D_{\mathbf{v}} \tilde{\mathbf{n}} \right),
\end{equation}
where $D_{\mathbf{v}}$ denotes the Fr\'{e}chet derivative in the ambient space $\mathcal{H}_{\varpi}$, that is, the ordinary derivative of $\tilde{\mathbf{n}}(\cdot\,;\mathbf{X})$ with respect to $\mathbf{X}$, holding the weight $\varpi$ fixed.
\end{definition}
\noindent{}This is the weighted analogue of the shape operator of Definition~\ref{def:second-fundamental-form}. Note that, since the weight $\varpi(\cdot\,;\mathbf{X})$ entering the inner product is itself a function of $\mathbf{X}$, the coordinate components of $A_{\mathbf{n}}$ do not reduce to a pairing of $\mathbf{n}$ against the coordinate derivatives $\partial_{X^j}\mathbf{e}_l^m$ alone, and Lemma~\ref{lemma:weingarten-weighted} below establishes their expression together with the correction that this dependence introduces. That expression involves $\mathbf{n}$ only through its value at $\mathbf{X}$, and hence $A_{\mathbf{n}}$ does not depend on the choice of extension.
\begin{lemma}[Weingarten duality under a state-dependent weight]
\label{lemma:weingarten-weighted}
Let $\mathbf{n} \in (T_{\mathbf{X}}\mathcal{M}_m)^{\perp_{\varpi}}$ be a normal vector at $\mathbf{X}$, and let $\{\mathbf{e}_k^m\}_{k=1}^d$ be the tangent basis of Definition~\ref{def:metric-tensor}. Then, for all $k, l = 1, \dots, d$,
\begin{equation}
    \langle A_{\mathbf{n}}\mathbf{e}_l^m, \mathbf{e}_k^m \rangle_\varpi = \langle \mathbf{I\!I}(\mathbf{e}_l^m,\mathbf{e}_k^m), \mathbf{n}\rangle_\varpi - \int_0^L \mathbf{n}\,\mathbf{e}_k^m\,\varpi^2\,\partial_r\mathbf{e}_l^m \, dr .
    \label{eq:An-II-corrected}
\end{equation}
The classical Weingarten duality~\eqref{eq:classical-weingarten-duality} is recovered when the ambient weight is independent of the base point, since $\partial_{X^l}\varpi = 0$ then removes the integral correction.
\end{lemma}
\begin{proof}
We differentiate the orthogonality relation that defines the normal bundle, as in the classical derivation of the Weingarten equation given by Lee~\cite{lee2018riemannian}, Proposition~8.4, and we identify every term that the differentiation produces.
\begin{proofpart}[The full product rule]
By Definition~\ref{def:normal-bundle}, the normal field obeys $\langle \mathbf{n}(\mathbf{X}), \mathbf{e}_k^m(\mathbf{X})\rangle_{\varpi(\cdot;\mathbf{X})} = 0$ at every base point $\mathbf{X}$ and for every $k$. We differentiate this identity with respect to $X^l$. Written as an explicit integral, the integrand carries \emph{three} $\mathbf{X}$-dependent factors, namely $\mathbf{n}$, $\mathbf{e}_k^m$, and the weight $\varpi$ itself, so the ordinary product rule produces three terms:
\begin{equation}
    0 = \frac{\partial}{\partial X^l}\int_0^L \mathbf{n}\,\mathbf{e}_k^m\,\varpi\,dr
    = \langle \partial_l\mathbf{n}, \mathbf{e}_k^m\rangle_\varpi + \langle \mathbf{n}, \partial_l\mathbf{e}_k^m\rangle_\varpi + \int_0^L \mathbf{n}\,\mathbf{e}_k^m\,(\partial_l\varpi)\,dr .
    \label{eq:full-product-rule}
\end{equation}
The first two terms are those retained by a naive product rule, in which the weight is treated as inert. The third is generated by the $\mathbf{X}$-dependence of $\varpi$, and it does not vanish in general.
\end{proofpart}
\begin{proofpart}[Identification of the three terms]
We treat each term of~\eqref{eq:full-product-rule} in turn. For the first, Definition~\ref{def:shape-operator} gives $A_{\mathbf{n}}\mathbf{e}_l^m = -\Pi^{\parallel}_{\mathbf{X}}(\partial_l\mathbf{n})$, and the projector may be dropped when pairing against $\mathbf{e}_k^m$, because $\mathbf{e}_k^m$ is itself tangential:
\begin{equation}
    \langle A_{\mathbf{n}}\mathbf{e}_l^m,\mathbf{e}_k^m\rangle_\varpi = -\langle\partial_l\mathbf{n},\mathbf{e}_k^m\rangle_\varpi .
    \label{eq:term-one-id}
\end{equation}
For the second, Definition~\ref{def:second-fundamental-form} gives $\mathbf{I\!I}(\mathbf{e}_l^m,\mathbf{e}_k^m) = \Pi^{\perp}_{\mathbf{X}}(\partial_l\mathbf{e}_k^m)$, and here the projector may be dropped because pairing against the normal vector $\mathbf{n}$ already annihilates the tangential part:
\begin{equation}
    \langle\mathbf{n},\partial_l\mathbf{e}_k^m\rangle_\varpi = \langle \mathbf{I\!I}(\mathbf{e}_l^m,\mathbf{e}_k^m),\mathbf{n}\rangle_\varpi .
    \label{eq:term-two-id}
\end{equation}
For the third, the weight-gradient identity $\partial_l\varpi = -\varpi^2\,\partial_r\mathbf{e}_l^m$ established in Definition~\ref{def:weighted-space} gives
\begin{equation}
    \int_0^L \mathbf{n}\,\mathbf{e}_k^m\,(\partial_l\varpi)\,dr = -\int_0^L \mathbf{n}\,\mathbf{e}_k^m\,\varpi^2\,\partial_r\mathbf{e}_l^m\,dr .
    \label{eq:term-three-id}
\end{equation}
\end{proofpart}
\begin{proofpart}[Rearrangement]
Substituting the three identifications~\eqref{eq:term-one-id},~\eqref{eq:term-two-id} and~\eqref{eq:term-three-id} into~\eqref{eq:full-product-rule}:
\begin{equation}
    0 = -\langle A_{\mathbf{n}}\mathbf{e}_l^m,\mathbf{e}_k^m\rangle_\varpi + \langle \mathbf{I\!I}(\mathbf{e}_l^m,\mathbf{e}_k^m),\mathbf{n}\rangle_\varpi - \int_0^L \mathbf{n}\,\mathbf{e}_k^m\,\varpi^2\,\partial_r\mathbf{e}_l^m\,dr .
\end{equation}
Solving this for $\langle A_{\mathbf{n}}\mathbf{e}_l^m,\mathbf{e}_k^m\rangle_\varpi$ gives~\eqref{eq:An-II-corrected}. If the weight does not depend on the base point, then $\partial_l\varpi = 0$, the third term of~\eqref{eq:full-product-rule} is absent from the outset, and the same rearrangement returns the classical duality~\eqref{eq:classical-weingarten-duality}.
\end{proofpart}
\end{proof}
\begin{lemma}[Derivative of the trivialisation]
\label{lemma:trivialisation-derivative}
Consider the inverse trivialisation, which carries a fixed vector $\eta \in \mathcal{H}_{\perp}$ to the normal vector $\mathcal{R}_{\mathbf{X}}^{-1}\eta$ at $\mathbf{X}$. Its partial Fr\'{e}chet derivative with respect to $\mathbf{X}$ satisfies the projection identity:
\begin{equation}
    \Pi_{\mathbf{X}}^{\parallel} \left( D_{\mathbf{X}} [\mathcal{R}_{\mathbf{X}}^{-1} \eta] \cdot \mathbf{v} \right) = - A_{\mathcal{R}_{\mathbf{X}}^{-1}\eta} \mathbf{v}.
\end{equation}
\end{lemma}
\begin{proof}
Let $\mathbf{n}(\mathbf{X}) \coloneqq \mathcal{R}_{\mathbf{X}}^{-1} \eta$ for the fixed $\eta \in \mathcal{H}_\perp$. Since $\mathcal{R}_{\mathbf{Y}}^{-1}$ maps $\mathcal{H}_\perp$ into $(T_{\mathbf{Y}}\mathcal{M}_m)^{\perp_\varpi}$ at every base point $\mathbf{Y}$, the assignment $\mathbf{Y} \mapsto \mathbf{n}(\mathbf{Y})$ is itself a smooth section of the normal bundle defined on the chart $U$, and in particular it is a valid smooth extension of its own value $\mathbf{n} = \mathbf{n}(\mathbf{X})$ to a neighbourhood of $\mathbf{X}$, in the sense required by Definition~\ref{def:shape-operator}. Applying that definition directly, with $\tilde{\mathbf{n}} = \mathbf{n}(\cdot)$ itself as the extension, gives
\begin{equation}
    A_{\mathbf{n}} \mathbf{v} \coloneqq -\Pi_{\mathbf{X}}^{\parallel}\bigl(D_{\mathbf{v}}\mathbf{n}\bigr) = -\Pi_{\mathbf{X}}^{\parallel}\bigl(D_{\mathbf{X}}[\mathcal{R}_{\mathbf{X}}^{-1}\eta]\cdot\mathbf{v}\bigr) ,
\end{equation}
which rearranges to the stated identity. No further computation is required, the identity being the defining formula for $A_{\mathbf{n}}$ evaluated at $\mathbf{n} = \mathcal{R}_{\mathbf{X}}^{-1}\eta$.
\end{proof}
\subsubsection{The geometric Jacobian}
We now evaluate the volume distortion factor $J_{\Psi}$ arising from the change of variables. We first provide the definition of this object.
\begin{definition}[Geometric Jacobian]
\label{def:geometric-jacobian}
Let $D\Psi_{(\mathbf{X}, \eta)}: T_{\mathbf{X}}\mathcal{X} \times \mathcal{H}_{\perp} \to \mathcal{H}_{\varpi}$ be the Fr\'{e}chet derivative of the Fermi map, the domain carrying the inner product of $g(\mathbf{X})$ on $T_{\mathbf{X}}\mathcal{X}$ and that of $\mathcal{H}_{\perp}$ on the fibre, and the codomain the weighted inner product $\langle \cdot, \cdot \rangle_{\varpi}$. The \textit{geometric Jacobian} is
\begin{equation}
    J_{\Psi}(\mathbf{X}, \eta) \coloneqq \sqrt{ \det \left( [D\Psi_{(\mathbf{X}, \eta)}]^* D\Psi_{(\mathbf{X}, \eta)} \right) },
\end{equation}
the adjoint being taken with respect to these inner products, and the determinant being the Fredholm determinant of Theorem~3.7 of Simon~\cite{simon2005trace}.
\end{definition}
\noindent{}The determinant exists, since $[D\Psi]^* D\Psi$ differs from the identity by an operator of finite rank, as the block form of $D\Psi$ established in the proof of Theorem~\ref{thm:jacobian-expansion} shows.
\begin{theorem}[Expansion of the geometric Jacobian]
\label{thm:jacobian-expansion}
Let $J_{\Psi}(\mathbf{X}, \eta)$ be the Geometric Jacobian defined above. For small fluctuations $\|\eta\|_{\mathcal{H}_{\perp}} \ll 1$, the Jacobian admits the asymptotic expansion:
\begin{equation}
    J_{\Psi}(\mathbf{X}, \eta) = 1 - \operatorname{Tr}(A_{\mathcal{R}_{\mathbf{X}}^{-1}\eta}) + \mathcal{O}(\|\eta\|_{\mathcal{H}_{\perp}}^2).
\end{equation}
Here, $\operatorname{Tr}(A_{\mathbf{n}})$ is the mean curvature of the manifold in the direction $\mathbf{n}$ less the weight-gradient correction of Eq.~\eqref{eq:An-II-corrected}.
\end{theorem}
\begin{proof}
We compute the determinant by constructing the matrix representation of the operator $D\Psi_{(\mathbf{X}, \eta)}$ in the orthogonal splitting $\mathcal{H}_{\varpi} = T_{\mathbf{X}}\mathcal{M}_m \oplus (T_{\mathbf{X}}\mathcal{M}_m)^{\perp_{\varpi}}$, and its block-diagonal form on the zero section is established in Lemma~\ref{lemma:fermi-jacobian}.
\begin{proofpart}[Computation of the Fr\'{e}chet derivative]
The Fermi map is given by the sum $\Psi(\mathbf{X}, \eta) = M(\mathbf{X}) + \mathcal{R}_{\mathbf{X}}^{-1}\eta$. We evaluate its derivative acting on a generic tangent vector $(\mathbf{v}, \zeta) \in T_{\mathbf{X}}\mathcal{X} \times \mathcal{H}_{\perp}$.
Using the linearity of differentiation:
\begin{equation}
    D\Psi_{(\mathbf{X}, \eta)}(\mathbf{v}, \zeta) = D_{\mathbf{X}}[M(\mathbf{X})]\mathbf{v} + D_{\mathbf{X}}[\mathcal{R}_{\mathbf{X}}^{-1}\eta]\mathbf{v} + D_{\eta}[\mathcal{R}_{\mathbf{X}}^{-1}\eta]\zeta.
\end{equation}
We identify each term. The derivative of the embedding, $D_{\mathbf{X}}[M(\mathbf{X})]\mathbf{v} = \mathcal{J}_{\iota_m}(\mathbf{X})\mathbf{v}$, is the pushforward of the parameter velocity, which lies in $T_{\mathbf{X}}\mathcal{M}_m$. The derivative in the fibre variable, $D_{\eta}[\mathcal{R}_{\mathbf{X}}^{-1}\eta]\zeta = \mathcal{R}_{\mathbf{X}}^{-1}\zeta$, is the operator itself applied to the perturbation $\zeta$, since $\mathcal{R}_{\mathbf{X}}^{-1}$ is linear in $\eta$, and its image lies in the normal bundle $(T_{\mathbf{X}}\mathcal{M}_m)^{\perp_{\varpi}}$, since $\zeta \in \mathcal{H}_{\perp}$. The derivative along the manifold, $D_{\mathbf{X}}[\mathcal{R}_{\mathbf{X}}^{-1}\eta]\mathbf{v}$, represents the variation of the normal vector $\mathbf{n} = \mathcal{R}_{\mathbf{X}}^{-1}\eta$ due to transport along the manifold, and it generally has components in both the tangent and the normal spaces.
Thus, the total derivative is:
\begin{equation}
    D\Psi_{(\mathbf{X}, \eta)}(\mathbf{v}, \zeta) = \mathcal{J}_{\iota_m}(\mathbf{X})\mathbf{v} + D_{\mathbf{X}}[\mathcal{R}_{\mathbf{X}}^{-1}\eta]\mathbf{v} + \mathcal{R}_{\mathbf{X}}^{-1}\zeta.
\end{equation}
\end{proofpart}
\begin{proofpart}[Projection onto the orthogonal basis]
We now project this vector result onto the splitting $T_{\mathbf{X}}\mathcal{M}_m \oplus (T_{\mathbf{X}}\mathcal{M}_m)^{\perp_{\varpi}}$ to form the operator matrix blocks. The upper-left block maps the tangent space into itself, and
\begin{equation}
    \Pi_{\mathbf{X}}^{\parallel} \left( \mathcal{J}_{\iota_m}(\mathbf{X})\mathbf{v} \right) = \mathcal{J}_{\iota_m}(\mathbf{X})\mathbf{v}.
\end{equation}
Since $\mathcal{J}_{\iota_m}$ is the identity isomorphism from the parameter tangent space to the manifold tangent space, this contribution is the identity operator $\mathcal{I}_{\parallel}$. Lemma~\ref{lemma:trivialisation-derivative} supplies the remaining tangential contribution, $\Pi_{\mathbf{X}}^{\parallel} D_{\mathbf{X}}[\mathbf{n}]\mathbf{v} = -A_{\mathbf{n}}\mathbf{v}$, so that the upper-left block is $\mathcal{I}_{\parallel} - A_{\mathbf{n}}$.

The upper-right block maps the normal space into the tangent space. Since $\mathcal{R}_{\mathbf{X}}^{-1}$ maps $\mathcal{H}_{\perp}$ isometrically to the normal space $(T_{\mathbf{X}}\mathcal{M}_m)^{\perp_{\varpi}}$, and $\Pi_{\mathbf{X}}^{\parallel}$ projects onto the tangent space, the image is orthogonal to the projection range, so that
\begin{equation}
    \Pi_{\mathbf{X}}^{\parallel} \left( \mathcal{R}_{\mathbf{X}}^{-1}\zeta \right) = 0
\end{equation}
and the upper-right block is exactly zero.

The lower-right block maps the normal space into itself. Since $\mathcal{R}_{\mathbf{X}}^{-1}\zeta$ is already normal, the projector $\Pi_{\mathbf{X}}^{\perp}$ acts upon it as the identity, and the map is the unitary isomorphism between fibres. Relative to the reference measure $\lambda_{\perp}$ (which corresponds to the basis in $\mathcal{H}_{\perp}$), the lower-right block is the identity $\mathcal{I}_{\perp}$.

The lower-left block maps the tangent space into the normal space. It is $\Pi_{\mathbf{X}}^{\perp} \left( D_{\mathbf{X}}[\mathcal{R}_{\mathbf{X}}^{-1}\eta]\mathbf{v} \right)$, which represents the covariant derivative of the normal vector projected back onto the normal bundle (the normal connection $\nabla^{\perp}$), and we denote it by $\mathcal{D}^{\perp}_{\mathbf{n}}$.
\end{proofpart}
\begin{proofpart}[Determinant calculation]
Combining the projection results from the previous step, the linear operator $D\Psi_{(\mathbf{X}, \eta)}$ admits the exact matrix representation
\begin{equation}
    D\Psi_{(\mathbf{X}, \eta)} = \begin{pmatrix} \mathcal{I}_{\parallel} - A_{\mathbf{n}} & 0 \\ \mathcal{D}^{\perp}_{\mathbf{n}} & \mathcal{I}_{\perp} \end{pmatrix}
\end{equation}
with respect to the orthogonal decomposition $\mathcal{H}_{\varpi} = T_{\mathbf{X}}\mathcal{M}_m \oplus \mathcal{N}_{\mathbf{X}}$. Writing $T \coloneqq \mathcal{I}_{\parallel} - A_{\mathbf{n}}$ and $D \coloneqq \mathcal{D}^{\perp}_{\mathbf{n}}$, the operator of Definition~\ref{def:geometric-jacobian} is
\begin{equation}
    [D\Psi_{(\mathbf{X}, \eta)}]^* D\Psi_{(\mathbf{X}, \eta)} = \begin{pmatrix} T^*T + D^*D & D^* \\ D & \mathcal{I}_{\perp} \end{pmatrix} ,
\end{equation}
and eliminating the lower-right block by its Schur complement gives
\begin{equation}
    \det\bigl([D\Psi]^* D\Psi\bigr) = \det(\mathcal{I}_{\perp})\,\det\bigl(T^*T + D^*D - D^*D\bigr) = \det(T^*T) = |\det T|^2 ,
\end{equation}
the Fredholm determinant of the identity being one and the remaining determinant that of a $d \times d$ matrix. The off-diagonal term $D$ therefore does not enter, and the whole Jacobian reduces to
\begin{equation}
    J_{\Psi}(\mathbf{X}, \eta) = \bigl| {\det}_{\mathbb{R}^d}\bigl( I_d - A_{\mathbf{n}} \bigr) \bigr| .
    \label{eq:fermi-jacobian-finite}
\end{equation}
For small fluctuations $\|\mathbf{n}\|_{\varpi} \ll 1$, the linear operator $A_{\mathbf{n}}$ is small in norm, so that $\det(I_d - A_{\mathbf{n}})$ is close to one and in particular positive, which allows the absolute value to be dropped. Since the tangent space $T_{\mathbf{X}}\mathcal{M}_m \cong \mathbb{R}^d$ is finite-dimensional, $A_{\mathbf{n}}$ is a $d \times d$ matrix and its trace is a finite sum, and the trace-class Assumption~\ref{ass:trace-class-embedding} is relevant only for the second fundamental form viewed as an operator on the infinite-dimensional normal bundle, which enters the higher-order corrections.
We expand the determinant using the functional calculus identity $\ln \det(B) = \operatorname{Tr}(\ln B)$. Setting $B = \mathcal{I}_{\parallel} - A_{\mathbf{n}}$:
\begin{align}
    \ln J_{\Psi} &= \operatorname{Tr} \left( \ln (\mathcal{I}_{\parallel} - A_{\mathbf{n}}) \right) \nonumber \\
    &= \operatorname{Tr} \left( - A_{\mathbf{n}} - \frac{1}{2} A_{\mathbf{n}}^2 - \mathcal{O}(\|A_{\mathbf{n}}\|^3) \right).
\end{align}
Exponentiating the result and linearising the exponential function gives the asymptotic form:
\begin{equation}
    J_{\Psi}(\mathbf{X}, \eta) = 1 - \operatorname{Tr}(A_{\mathbf{n}}) + \mathcal{O}(\|\mathbf{n}\|_{\varpi}^2).
\end{equation}
Substituting the definition $\mathbf{n} = \mathcal{R}_{\mathbf{X}}^{-1}\eta$ and using the isometry property $\|\mathbf{n}\|_{\varpi} = \|\eta\|_{\mathcal{H}_{\perp}}$, we obtain the stated expansion.
\end{proofpart}
\end{proof}
\subsubsection{The factorised density ansatz}
\label{subsubsec:density-factorisation}
We substitute the geometric Jacobian into the disintegration formula to obtain the explicit form of the density $\varrho_t$.
\begin{definition}[Intrinsic conditional profile]
\label{def:intrinsic-profile}
Let $\nu_t^{\mathbf{X}}$ be the conditional measure on the fibre defined in Proposition \ref{prop:measure-decomposition}. We define the \textit{intrinsic conditional profile} $\Gamma_t(\mathbf{X}, \eta)$ as the density of this conditional measure with respect to the Gaussian reference $\lambda_{\perp}$, rescaled by the geometric volume factor $J_{\Psi}^{-1}$:
\begin{equation}
    \Gamma_t(\mathbf{X}, \eta) \coloneqq \frac{d\nu_t^{\mathbf{X}}}{d\lambda_{\perp}}(\eta) \cdot J_{\Psi}(\mathbf{X}, \eta).
\end{equation}
\end{definition}
\noindent{}This profile captures the purely statistical deviation of the fluctuations from the Gaussian reference, separating out the deterministic volume contraction effects encoded in $J_{\Psi}$.
\begin{proposition}[Geometric Gibbs factorisation]
\label{prop:gibbs-factorisation}
Let $\mu_t$ be a solution of the functional forward Kolmogorov equation of Definition~\ref{def:adjoint-evolution}, and assume Assumption~\ref{ass:absolute-continuity}. Its Radon--Nikodym density $\varrho_t(\mathbf{X}, \eta)$ with respect to the reference measure $\lambda$ factorises as:
\begin{equation}
    \varrho_t(\mathbf{X}, \eta) = \bar{\varrho}_t(\mathbf{X}) \, \Gamma_t(\mathbf{X}, \eta) \, J_{\Psi}(\mathbf{X}, \eta)^{-1},
\end{equation}
where:
\begin{enumerate}[(i)]
    \item $\bar{\varrho}_t(\mathbf{X})$ is the marginal density on the manifold (Proposition \ref{prop:measure-decomposition}).
    \item $J_{\Psi}^{-1}$ is the geometric correction factor derived in Theorem \ref{thm:jacobian-expansion}.
    \item $\Gamma_t(\mathbf{X}, \eta)$ is the intrinsic conditional profile (Definition \ref{def:intrinsic-profile}), which satisfies the normalisation condition:
    \begin{equation}
        \int_{\mathcal{H}_{\perp}} \Gamma_t(\mathbf{X}, \eta) \, J_{\Psi}(\mathbf{X}, \eta)^{-1} \, d\lambda_{\perp}(\eta) = 1.
    \end{equation}
\end{enumerate}
\end{proposition}
\begin{proof}
Theorem~\ref{thm:disintegration}, together with Definition~\ref{def:disintegration}, writes the full measure as the marginal on the base carrying the conditional on the fibre, $d\mu_t = d\bar{\mu}_t \otimes d\nu_t^{\mathbf{X}}$. We express both sides as densities against the product reference measure $\lambda = \lambda_{\mathcal{X}} \otimes \lambda_{\perp}$ and identify the factors of $d\lambda_{\mathcal{X}}\,d\lambda_{\perp}$:
\begin{equation}
    \varrho_t \, d\lambda_{\mathcal{X}}\, d\lambda_{\perp}
    = \bigl( \bar{\varrho}_t \, d\lambda_{\mathcal{X}} \bigr) \otimes \left( \frac{d\nu_t^{\mathbf{X}}}{d\lambda_{\perp}}\, d\lambda_{\perp} \right) .
\end{equation}
Cancelling the common factor $d\lambda_{\mathcal{X}}\,d\lambda_{\perp}$ from the two sides, therefore,
\begin{equation}
    \varrho_t
    = \bar{\varrho}_t \, \frac{d\nu_t^{\mathbf{X}}}{d\lambda_{\perp}}
     = \bar{\varrho}_t \, \Gamma_t \, J_{\Psi}^{-1} ,
    \label{eq:gibbs-density-identification}
\end{equation}
the final equality being the definition of the intrinsic profile, $\Gamma_t \coloneqq (d\nu^{\mathbf{X}}_t/d\lambda_\perp)\, J_\Psi$, given in Definition~\ref{def:intrinsic-profile}. The normalisation stated in the proposition is the assertion that the conditional measure is a probability measure on the fibre, and the Jacobian cancels between the profile and its inverse:
\begin{align}
    1 = \int_{\mathcal{H}_\perp} d\nu_t^{\mathbf{X}}
    &= \int_{\mathcal{H}_\perp} \frac{d\nu_t^{\mathbf{X}}}{d\lambda_{\perp}}\, d\lambda_{\perp}
     = \int_{\mathcal{H}_\perp} \Gamma_t \, J_{\Psi}^{-1}\, d\lambda_{\perp} .
    \label{eq:gibbs-normalisation}
\end{align}
The same cancellation identifies the fibre average with the conditional expectation, since for every bounded measurable $f$:
\begin{align}
    \mathbb{E}^\perp[f]
    = \int_{\mathcal{H}_\perp} f \, \Gamma_t \, J_\Psi^{-1}\, d\lambda_\perp
    = \int_{\mathcal{H}_\perp} f \, \frac{d\nu^{\mathbf{X}}_t}{d\lambda_\perp}\, d\lambda_\perp
    = \int_{\mathcal{H}_\perp} f \, d\nu^{\mathbf{X}}_t .
    \label{eq:gibbs-conditional-expectation}
\end{align}
The factor $J_\Psi$ is carried separately only in order that the reference measure $\lambda_\perp$ not depend upon the base point, and~\eqref{eq:gibbs-conditional-expectation} shows that it cancels from every fibre average.
\end{proof}

\begin{remark}[The volume factor and the mean of the fibre law]
\label{rem:volume-no-drift}
The factor $J_{\Psi}^{-1}$ is the density that the curvature of the tube contributes to the conditional measure against the flat reference measure $\lambda_\perp$, and it would enter the dynamics as an entropic potential $\ln J_{\Psi}$ if the fibre average were taken against $\lambda_\perp$. By Eq.~\eqref{eq:gibbs-conditional-expectation}, however, the average is taken against $\nu_t^{\mathbf{X}}$, and so $J_{\Psi}$ cancels from it. Still, the curvature of the tube has an effect, since it acts on $\nu_t^{\mathbf{X}}$ itself and displaces its mean by $\varepsilon^2\hat\eta$ under Hypothesis~\ref{hyp:concentration}. In other words, the last two terms of Eq.~\eqref{eq:frame-correction} are the drift that this displacement produces.
\end{remark}
\subsection{The probabilistic projection}
\label{subsec:probabilistic-projection}
Having decomposed the probability density $\varrho_t$ into a marginal component on the manifold ($\bar{\varrho}_t$) and a conditional profile on the fibres ($\Gamma_tJ_\Psi^{-1}$), we now derive the effective evolution equation for the marginal density $\bar{\varrho}_t$. For this we project the functional forward Kolmogorov equation (Definition~\ref{def:adjoint-evolution}) onto the tangent space of the manifold, in two steps.
In the \emph{operator decomposition} we express the generator $\mathcal{L}_\varepsilon$ in the Fermi coordinates $(\mathbf{X}, \eta)$. Here we evaluate the differential operators (the gradient and the Laplacian) in the curved geometry of the tubular neighbourhood, while tracking the metric perturbations induced by the curvature. In the \emph{fibre averaging} we then integrate the evolution equation over the normal fibres with respect to the reference measure $\lambda_{\perp}$. The two steps together generate an effective drift correction of their own, and we set this quantity against the one produced by the stochastic projection when we compare the two pathways in Section~\ref{sec:commutativity}.
\subsubsection{The generator in Fermi coordinates}
To express the generator $\mathcal{L}_\varepsilon \Phi = \langle \mathcal{A}_m, D\Phi \rangle + \frac{\varepsilon^2}{2} \operatorname{Tr}(\mathcal{Q}_m D^2\Phi)$ in the local chart, we must compute the components of the gradient and the Hessian relative to the basis $\{\partial_{X^i}, \partial_{\eta}\}$. For this we need the exact metric tensor of the tubular neighbourhood.
Before we derive the metric, we establish the decomposition of the derivatives of the normal frame, which ensures that the terms appearing in the metric expansion are well-defined geometric objects.
\begin{lemma}[Orthogonal decomposition of derivatives]
\label{lemma:weingarten-decomposition}
Let $\mathbf{X} \in \mathcal{X}$ be a coordinate vector parametrising the point $M(\mathbf{X}) \in \mathcal{M}_m$. Let $\mathbf{u} \in T_{\mathbf{X}}\mathcal{M}_m$ be a tangent vector and $\mathbf{n} \in \Gamma(\mathcal{N})$ be a smooth section of the normal bundle (a normal field).
The Fr\'{e}chet derivative of the normal field $D_{\mathbf{u}}\mathbf{n}$ in the ambient space $\mathcal{H}_{\varpi}$ admits a unique orthogonal decomposition:
\begin{equation}
    D_{\mathbf{u}}\mathbf{n} = \Pi_{\mathbf{X}}^{\parallel}(D_{\mathbf{u}}\mathbf{n}) + \Pi_{\mathbf{X}}^{\perp}(D_{\mathbf{u}}\mathbf{n}).
\end{equation}
Furthermore, the tangential component is $-A_{\mathbf{n}}\mathbf{u}$, with $A_{\mathbf{n}}$ the shape operator of Definition~\ref{def:shape-operator}, and Lemma~\ref{lemma:weingarten-weighted} expresses the components of $A_{\mathbf{n}}$ through the second fundamental form and the dependence of the weight upon the base point.
\end{lemma}
\begin{proof}
We proceed in two steps, establishing the existence of the decomposition via orthogonal projection and then characterising the tangential component via the geometry of the embedding.
\begin{proofpart}[Existence and uniqueness of the split]
Since $\mathcal{H}_{\varpi}$ is a Hilbert space and the tangent space $T_{\mathbf{X}}\mathcal{M}_m$ is a closed linear subspace (isomorphic to $\mathbb{R}^d$), the orthogonal projection theorem guarantees that the space admits the direct sum decomposition $\mathcal{H}_{\varpi} = T_{\mathbf{X}}\mathcal{M}_m \oplus \mathcal{N}_{\mathbf{X}}$.
Consequently, the identity operator decomposes as $\mathcal{I} = \Pi_{\mathbf{X}}^{\parallel} + \Pi_{\mathbf{X}}^{\perp}$. Applying this operator to the vector $D_{\mathbf{u}}\mathbf{n}$ (which exists in $\mathcal{H}_{\varpi}$ because $\mathbf{n}$ is a smooth field), we obtain the identity:
\begin{equation}
    D_{\mathbf{u}}\mathbf{n} = \Pi_{\mathbf{X}}^{\parallel}(D_{\mathbf{u}}\mathbf{n}) + \Pi_{\mathbf{X}}^{\perp}(D_{\mathbf{u}}\mathbf{n}).
\end{equation}
This establishes that the decomposition into parallel and normal components exists and is unique.
\end{proofpart}
\begin{proofpart}[Characterisation of the tangential component]
Definition~\ref{def:shape-operator} defines the shape operator by exactly this tangential projection, $A_{\mathbf{n}}\mathbf{u} \coloneqq -\Pi_{\mathbf{X}}^{\parallel}(D_{\mathbf{u}}\tilde{\mathbf{n}})$, for any smooth local extension $\tilde{\mathbf{n}}$ of $\mathbf{n}$ with $\tilde{\mathbf{n}}(\mathbf{X}) = \mathbf{n}$. Since $\mathbf{n}$ is itself already given as a smooth section of the normal bundle, it serves as its own extension along any path through $\mathbf{X}$, and comparing with the decomposition of the previous part identifies the tangential component directly:
\begin{equation}
    \Pi_{\mathbf{X}}^{\parallel}(D_{\mathbf{u}}\mathbf{n}) = -A_{\mathbf{n}}\mathbf{u} .
\end{equation}
This component is non-zero precisely when $A_{\mathbf{n}} \neq 0$, and Lemma~\ref{lemma:weingarten-weighted} relates $A_{\mathbf{n}}$ to the second fundamental form, together with the correction which the dependence of the weight upon the base point introduces.
\end{proofpart}
\end{proof}
\begin{remark}[Consistency with trivialisation derivative]
This general decomposition explains the result obtained in Lemma~\ref{lemma:trivialisation-derivative}. In that lemma, we computed the derivative of the specific normal field $\mathbf{n}(\mathbf{X}) = \mathcal{R}_{\mathbf{X}}^{-1} \eta$, and the result $\Pi_{\mathbf{X}}^{\parallel} D_{\mathbf{X}}[\mathcal{R}_{\mathbf{X}}^{-1}\eta]\mathbf{v} = -A_{\mathbf{n}}\mathbf{v}$ is exactly the tangential component identified here, applied to the trivialised normal vector.
\end{remark}
\begin{corollary}[The Weingarten decomposition]
\label{corollary:weingarten-decomposition}
Using Definition \ref{def:shape-operator} and Lemma \ref{lemma:weingarten-decomposition}, the total Fr\'{e}chet derivative of a normal field $\mathbf{n} \in \Gamma(\mathcal{N})$ along a tangent vector $\mathbf{u}$ decomposes into:
\begin{equation}
    D_{\mathbf{u}} \mathbf{n} = -A_{\mathbf{n}}\mathbf{u} + \nabla^{\perp}_{\mathbf{u}} \mathbf{n}.
\end{equation}
\end{corollary}
Based on Lemma \ref{lemma:weingarten-decomposition}, we now define the normal connection, which names the second of these two projection components.
\begin{definition}[The normal connection]
\label{def:geometry-objects}
Using the decomposition established in Lemma~\ref{lemma:weingarten-decomposition}, we define the \emph{normal connection} $\nabla^{\perp}: T_{\mathbf{X}}\mathcal{M}_m \times \Gamma(\mathcal{N}) \to \Gamma(\mathcal{N})$ as the normal projection of the derivative of a normal field,
\begin{equation}
    \nabla^{\perp}_{\mathbf{u}} \mathbf{n} \coloneqq \Pi_{\mathbf{X}}^{\perp} (D_{\mathbf{u}} \mathbf{n}) .
\end{equation}
\end{definition}
\noindent{}Hence the two components of the decomposition are the shape operator (Definition~\ref{def:shape-operator}) and the normal connection, and the second fundamental form entering them is that of Definition~\ref{def:second-fundamental-form}, taken with the weighted inner product at the base point. Next, we define the metric that the embedding into the ambient Hilbert space induces on the product space $U \times \mathcal{H}_{\perp}$.
\begin{definition}[The induced metric tensor on the bundle]
\label{def:induced-metric-bundle}
The \textit{induced metric tensor} $\mathfrak{G}$ is the pullback of the ambient inner product $\langle \cdot, \cdot \rangle_{\varpi}$ via the Fermi map $\Psi: U \times \mathcal{H}_{\perp} \to \mathcal{H}_{\varpi}$.
It is a map that assigns to every point $(\mathbf{X}, \eta)$ in the coordinate chart a positive-definite operator on the tangent space of the product manifold, identifying $T_{(\mathbf{X}, \eta)}(U \times \mathcal{H}_{\perp}) \cong T_{\mathbf{X}}\mathcal{M}_m \oplus \mathcal{H}_{\perp}$. That is, for any two vectors $\mathbf{V}, \mathbf{W}$ in the product tangent space:
\begin{equation}
    \mathfrak{G}_{(\mathbf{X}, \eta)}(\mathbf{V}, \mathbf{W}) = \langle D\Psi_{(\mathbf{X}, \eta)}[\mathbf{V}], D\Psi_{(\mathbf{X}, \eta)}[\mathbf{W}] \rangle_{\varpi}.
\end{equation}
We represent $\mathfrak{G}$ as a block operator acting on the orthogonal decomposition $T_{\mathbf{X}}\mathcal{M}_m \oplus \mathcal{H}_{\perp}$, capturing the geometry of the full tubular neighbourhood.
\end{definition}
\begin{remark}[Coordinate representation and indices]
In order to perform explicit calculations, we introduce a coordinate basis, and we denote indices for the full product space by Greek letters $\mu, \nu \in \mathcal{I}$, where the index set $\mathcal{I} = \{1, \dots, d\} \sqcup \mathbb{N}$ is the disjoint union of:
\begin{enumerate}[(i)]
    \item Tangential indices $i, j \in \{1, \dots, d\}$ corresponding to the basis $\mathbf{e}_i^m = \partial_{X^i} M$.
    \item Fibre indices $\alpha, \beta \in \mathbb{N}$ corresponding to the basis of $\mathcal{H}_{\perp}$.
\end{enumerate}
In this notation, the metric components are $\mathfrak{G}_{\mu \nu} = \langle \partial_\mu \Psi, \partial_\nu \Psi \rangle_{\varpi}$.
When restricted to the zero section ($\eta = 0$), the tangential block $\mathfrak{G}_{ij}$ coincides exactly with the Riemannian metric $g_{ij}$ on $\mathcal{M}_m$.
\end{remark}
\begin{lemma}[The inverse metric tensor]
\label{lemma:inverse-metric}
Let $\mathfrak{G}(\mathbf{X}, \eta)$ be the induced metric tensor defined above. For small fluctuations $\|\eta\|_{\mathcal{H}_{\perp}} \ll 1$, its inverse $\mathfrak{G}^{-1}(\mathbf{X}, \eta)$ admits the following block expansion:
\begin{equation}
    \mathfrak{G}^{-1}(\mathbf{X}, \eta) = \begin{pmatrix} g^{ij} + \sum_{k,l} g^{ik}\, \Delta_{kl}[\mathbf{n}] \, g^{lj} & \mathcal{O}(\|\eta\|) \\ \mathcal{O}(\|\eta\|) & I + \mathcal{O}(\|\eta\|^2) \end{pmatrix} + \mathcal{O}(\|\eta\|^2) ,
\end{equation}
where, writing $\mathbf{n} = \mathcal{R}_{\mathbf{X}}^{-1}\eta$,
\begin{equation}
    \Delta_{kl}[\mathbf{n}] \coloneqq 2\langle \mathbf{I\!I}(\mathbf{e}_k^m,\mathbf{e}_l^m), \mathbf{n} \rangle_\varpi - \int_0^L \mathbf{n}\,\varpi^2 \bigl( \mathbf{e}_l^m\,\partial_r\mathbf{e}_k^m + \mathbf{e}_k^m\,\partial_r\mathbf{e}_l^m \bigr) \, dr ,
    \label{eq:metric-correction-Delta}
\end{equation}
the second term being the correction generated by the dependence of the weight $\varpi(\cdot\,;\mathbf{X})$ on the base point, established in the proof below.
\end{lemma}
\begin{proof}
We construct the metric components $\mathfrak{G}_{\mu\nu} = \langle \partial_\mu \Psi, \partial_\nu \Psi \rangle_{\varpi}$ by expanding the derivatives of the Fermi map $\Psi(\mathbf{X}, \eta) = M(\mathbf{X}) + \mathbf{n}(\mathbf{X})$, where $\mathbf{n} = \mathcal{R}_{\mathbf{X}}^{-1}\eta$.
\begin{proofpart}[Expansion of the tangent basis vectors]
Recall $\mathbf{e}_i(\mathbf{X}) = \partial_{X^i} M(\mathbf{X})$ are the coordinate basis vectors of the tangent space $T_{\mathbf{X}}\mathcal{M}_m$. Differentiating the map $\Psi$ with respect to the coordinate $X^i$, and applying the Weingarten decomposition (Corollary \ref{corollary:weingarten-decomposition}) to the term $\partial_{X^i}\mathbf{n}\equiv D_{\mathbf{e}_i} \mathbf{n}$:
\begin{align}
    \partial_{X^i} \Psi &= \partial_{X^i} M + \partial_{X^i} \mathbf{n} \nonumber \\
    &= \mathbf{e}_i + \left( -A_{\mathbf{n}}\mathbf{e}_i + \nabla^{\perp}_{\mathbf{e}_i} \mathbf{n} \right) \nonumber \\
    &= (\mathcal{I} - A_{\mathbf{n}}) \mathbf{e}_i + \nabla^{\perp}_{\mathbf{e}_i} \mathbf{n}.
\end{align}
The first term lies in the tangent space $T_{\mathbf{X}}\mathcal{M}_m$, whereas the second lies in the normal fibre $\mathcal{N}_{\mathbf{X}}$.
\end{proofpart}
\begin{proofpart}[Computation of the tangential metric block $\mathfrak{G}_{ij}$]
We form the inner product of two tangent basis vectors. Using the orthogonality of the tangent and normal subspaces (Lemma \ref{lemma:weingarten-decomposition}), the cross terms vanish:
\begin{align}
    \mathfrak{G}_{ij} &= \langle \partial_{X^i} \Psi, \partial_{X^j} \Psi \rangle_{\varpi} \nonumber \\
    &= \langle (\mathcal{I} - A_{\mathbf{n}}) \mathbf{e}_i + \nabla^{\perp}_{\mathbf{e}_i} \mathbf{n}, (\mathcal{I} - A_{\mathbf{n}}) \mathbf{e}_j + \nabla^{\perp}_{\mathbf{e}_j} \mathbf{n} \rangle_{\varpi} \nonumber \\
    &= \langle (\mathcal{I} - A_{\mathbf{n}}) \mathbf{e}_i, (\mathcal{I} - A_{\mathbf{n}}) \mathbf{e}_j \rangle_{\varpi} + \langle \nabla^{\perp}_{\mathbf{e}_i} \mathbf{n}, \nabla^{\perp}_{\mathbf{e}_j} \mathbf{n} \rangle_{\varpi}.
\end{align}
We expand the first term bilinearly. Since $A_{\mathbf{n}}$ is linear in $\mathbf{n}$, and $\mathbf{n} = \mathcal{R}_{\mathbf{X}}^{-1}\eta$ is linear in $\eta$, every term carrying two factors of $A_{\mathbf{n}}$ is $\mathcal{O}(\|\eta\|^2)$. The connection term $\langle \nabla^{\perp}_{\mathbf{e}_i}\mathbf{n}, \nabla^{\perp}_{\mathbf{e}_j}\mathbf{n} \rangle_{\varpi}$ is quadratic for the same reason, the covariant derivative $\nabla^{\perp}_{\mathbf{e}_i}$ being linear in $\mathbf{n}$ and hence in $\eta$. Thus, to linear order:
\begin{align}
    \mathfrak{G}_{ij} &= \langle \mathbf{e}_i, \mathbf{e}_j \rangle_{\varpi} - \langle A_{\mathbf{n}}\mathbf{e}_i, \mathbf{e}_j \rangle_{\varpi} - \langle \mathbf{e}_i, A_{\mathbf{n}}\mathbf{e}_j \rangle_{\varpi} + \mathcal{O}(\|\eta\|^2) \nonumber \\
    &= g_{ij} - \langle A_{\mathbf{n}}\mathbf{e}_i, \mathbf{e}_j \rangle_{\varpi} - \langle A_{\mathbf{n}}\mathbf{e}_j, \mathbf{e}_i \rangle_{\varpi} + \mathcal{O}(\|\eta\|^2) ,
    \label{eq:Gij-presymmetrised}
\end{align}
the second line using $\langle \mathbf{e}_i, A_{\mathbf{n}}\mathbf{e}_j\rangle_\varpi = \langle A_{\mathbf{n}}\mathbf{e}_j,\mathbf{e}_i\rangle_\varpi$, which holds by symmetry of the inner product alone and does not require $A_{\mathbf{n}}$ to be self-adjoint. The distinction matters here, because $A_{\mathbf{n}}$ is in fact \emph{not} self-adjoint with respect to $\langle\cdot,\cdot\rangle_\varpi$, since the weight $\varpi(\cdot\,;\mathbf{X})$ entering that pairing depends on the base point $\mathbf{X}$. The two terms in~\eqref{eq:Gij-presymmetrised} must therefore be kept separately, and may not be combined into $2\langle A_{\mathbf{n}}\mathbf{e}_i,\mathbf{e}_j\rangle_\varpi$. We now relate each of them to the second fundamental form. We do this with Lemma~\ref{lemma:weingarten-weighted}, which supplies the weighted Weingarten duality~\eqref{eq:An-II-corrected} together with the integral correction generated by the base-point dependence of $\varpi$. Substituting~\eqref{eq:An-II-corrected}, once with $(l,k)=(i,j)$ and once with $(l,k)=(j,i)$, into~\eqref{eq:Gij-presymmetrised}, and using the symmetry $\mathbf{I\!I}(\mathbf{e}_i^m,\mathbf{e}_j^m) = \mathbf{I\!I}(\mathbf{e}_j^m,\mathbf{e}_i^m)$, which holds because $\partial_i\mathbf{e}_j^m = \partial_j\mathbf{e}_i^m$ (mixed partials of $M$ commute):
\begin{equation}
    \mathfrak{G}_{ij} = g_{ij} - 2\langle \mathbf{I\!I}_{ij}, \mathbf{n} \rangle_{\varpi} + \int_0^L \mathbf{n}\,\varpi^2 \bigl( \mathbf{e}_j^m\,\partial_r\mathbf{e}_i^m + \mathbf{e}_i^m\,\partial_r\mathbf{e}_j^m \bigr) \, dr + \mathcal{O}(\|\eta\|^2) .
    \label{eq:Gij-corrected}
\end{equation}
\end{proofpart}
\begin{proofpart}[Computation of the cross metric block $\mathfrak{G}_{i \alpha}$]
With $\{u_\alpha\}$ an orthonormal basis of $\mathcal{H}_\perp$, the fibre basis vectors are the images of that basis under the unitary map, $\partial_{\eta^\alpha} \Psi = \mathcal{R}_{\mathbf{X}}^{-1} u_\alpha$. We compute the inner product with the tangent vectors:
\begin{align}
    \mathfrak{G}_{i \alpha} &= \langle \partial_{X^i} \Psi, \partial_{\eta^\alpha} \Psi \rangle_{\varpi} \nonumber \\
    &= \langle (\mathcal{I} - A_{\mathbf{n}}) \mathbf{e}_i + \nabla^{\perp}_{\mathbf{e}_i} \mathbf{n}, \mathcal{R}_{\mathbf{X}}^{-1} u_\alpha \rangle_{\varpi} \nonumber \\
    &= \langle (\mathcal{I} - A_{\mathbf{n}}) \mathbf{e}_i, \mathcal{R}_{\mathbf{X}}^{-1} u_\alpha \rangle_{\varpi} + \langle \nabla^{\perp}_{\mathbf{e}_i} \mathbf{n}, \mathcal{R}_{\mathbf{X}}^{-1} u_\alpha \rangle_{\varpi}.
\end{align}
The first inner product vanishes by orthogonality of $\operatorname{Im}(\mathcal{I} - A_{\mathbf{n}})$ and $\operatorname{Im}(\mathcal{R}_{\mathbf{X}}^{-1})$. The surviving term involves the connection $\nabla^{\perp}$ acting on $\mathbf{n}$. Since $\mathbf{n}$ is linear in $\eta$, this term is of order $\mathcal{O}(\|\eta\|)$.
\end{proofpart}
\begin{proofpart}[Block inversion]
We assemble the full metric matrix $\mathfrak{G}$ and compute its inverse. Let $\mathfrak{G}_\parallel \coloneqq \bigl(g_{ij} - \Delta_{ij}[\mathbf{n}]\bigr)$, with $\Delta_{ij}[\mathbf{n}]$ the right-hand correction of~\eqref{eq:Gij-corrected} exactly as defined in~\eqref{eq:metric-correction-Delta}, be the top-left block, and let $\mathfrak{G}_\times$ denote the cross-term block of order $\mathcal{O}(\|\eta\|)$. The normal block is the identity $\mathcal{I}_{\perp}$ (plus $\mathcal{O}(\eta^2)$ corrections).
\begin{equation}
    \mathfrak{G} = \begin{pmatrix} \mathfrak{G}_\parallel & \mathfrak{G}_\times^* \\ \mathfrak{G}_\times & \mathcal{I}_{\perp} \end{pmatrix} + \mathcal{O}(\|\eta\|^2).
\end{equation}
We apply the block inversion formula. The top-left block of the inverse is given by the inverse of the Schur complement $\mathfrak{G}_\parallel - \mathfrak{G}_\times^* \mathcal{I}_{\perp}^{-1} \mathfrak{G}_\times$. Since $\mathfrak{G}_\times \sim \mathcal{O}(\|\eta\|)$, the product $\mathfrak{G}_\times^* \mathfrak{G}_\times \sim \mathcal{O}(\|\eta\|^2)$, so the Schur complement equals $\mathfrak{G}_\parallel$ to linear order.
We compute $\mathfrak{G}_\parallel^{-1}$ using the Neumann series $(g - \Delta)^{-1} \approx g^{-1} + g^{-1} \Delta g^{-1}$:
\begin{align}
    (\mathfrak{G}^{-1})_{ij} &\approx \bigl(g_{ij} - \Delta_{ij}[\mathbf{n}]\bigr)^{-1} \nonumber \\
    &= g^{ij} + \sum_{k,l} g^{ik}\, \Delta_{kl}[\mathbf{n}]\, g^{lj} + \mathcal{O}(\|\eta\|^2).
\end{align}
This confirms the form of the tangential block stated in the lemma. The off-diagonal terms of the inverse are $-\mathfrak{G}_\parallel^{-1} \mathfrak{G}_\times^*$, which remain $\mathcal{O}(\|\eta\|)$. We retain these in the error term notation.
\end{proofpart}
\end{proof}
\begin{definition}[The connection correction]
\label{def:connection-correction}
For $i, j, k = 1, \dots, d$, the \emph{projected connection coefficients} at $\mathbf{X}$ are
\begin{equation}
    \tilde{\Gamma}^k_{ij}(\mathbf{X}) \coloneqq \sum_{l=1}^d g^{kl}(\mathbf{X})\,\bigl\langle \partial_{X^i}\mathbf{e}_j^m, \mathbf{e}_l^m \bigr\rangle_\varpi ,
    \label{eq:projected-connection}
\end{equation}
the inner product being that of the base point, and the \emph{connection correction} is the vector field $\mathfrak{c}$ upon $\mathcal{X}$ with components
\begin{equation}
    \mathfrak{c}^k(\mathbf{X}) \coloneqq D_0 \sum_{i,j=1}^d g^{ij}(\mathbf{X})\bigl(\Gamma^k_{ij}(\mathbf{X}) - \tilde{\Gamma}^k_{ij}(\mathbf{X})\bigr) ,
    \label{eq:connection-correction}
\end{equation}
with $\Gamma^k_{ij} = \tfrac12\sum_{l} g^{kl}(\partial_{X^i} g_{jl} + \partial_{X^j} g_{il} - \partial_{X^l} g_{ij})$ the Christoffel symbols of $g$.
\end{definition}
\begin{proposition}[Decomposition of the diffusion generator]
\label{prop:diffusion-decomposition}
Let $\Phi \in C_b^2(\mathcal{H}_{\varpi})$ be a test functional expressed in the Fermi coordinates $(\mathbf{X}, \eta)$, and consider the diffusion component of the generator computed with the local covariance $\mathcal{Q}^{\mathrm{loc}}_m$ of Remark~\ref{rem:two-covariances},
\begin{equation}
    \mathcal{L}_2^\varepsilon \Phi \coloneqq \frac{\varepsilon^2}{2} \operatorname{Tr}_{\mathcal{H}_\varpi}(\mathcal{Q}^{\mathrm{loc}}_m D^2\Phi).
\end{equation}
In the tubular neighbourhood of the manifold this operator admits the local decomposition
\begin{align}
    \mathcal{L}_2^\varepsilon \Phi &= \varepsilon^2 D_0 \Delta_g \Phi + \varepsilon^2 \sum_{k=1}^d \mathfrak{c}^k \frac{\partial \Phi}{\partial X^k} \nonumber \\
    &\quad + \varepsilon^2 D_0 \operatorname{Tr}_{\perp}(D_{\eta}^2 \Phi) - \varepsilon^2 D_0 \langle \mathbf{H}, \nabla_\eta \Phi \rangle_{\perp} + \mathcal{R}_{\kappa} \Phi ,
    \label{eq:diffusion-decomposition}
\end{align}
in which $\Delta_g$ is the Laplace--Beltrami operator on $\mathcal{M}_m$, $\mathfrak{c}$ is the connection correction of Definition~\ref{def:connection-correction}, $\mathbf{H} \coloneqq \sum_{i,j} g^{ij}\,\mathcal{R}_{\mathbf{X}}\mathbf{I\!I}(\mathbf{e}_i^m, \mathbf{e}_j^m)$ is the mean curvature vector carried into the model fibre, and $\mathcal{R}_{\kappa}$ vanishes at $\eta = 0$, its coefficients depending smoothly upon $\eta$, with second-order part in the tangential directions equal, to linear order in $\eta$, to
\begin{equation}
    \varepsilon^2 D_0 \sum_{i,j,k,l=1}^d g^{ik}\,\bigl(\Delta_{kl}[\mathbf{n}] + E_{kl}[\mathbf{n}]\bigr)\,g^{lj} \, \frac{\partial^2 \Phi}{\partial X^i \partial X^j} ,
\end{equation}
where $\Delta_{kl}[\mathbf{n}]$ is the linear metric correction of Eq.~\eqref{eq:metric-correction-Delta} and $E_{kl}[\mathbf{n}] \coloneqq \int_0^L \mathbf{e}_k^m\mathbf{e}_l^m\varpi^2\,\partial_r\mathbf{n}\,dr$ carries the dependence of the noise amplitude upon the state.
\end{proposition}
\begin{proof}
To perform the decomposition, we calculate the trace in the natural coordinate basis of the tubular neighbourhood, $\{\partial_{\mu}\}_{\mu \in \mathcal{I}}$, where indices $\mu, \nu$ run over both tangential directions ($i,j$) and fibre directions ($\alpha, \beta$). Writing $(D^2\Phi)_{\mu\nu} \coloneqq D^2\Phi(\partial_\mu\Psi, \partial_\nu\Psi)$ for the components of the second derivative upon the coordinate vectors, the diffusion generator is the contraction of the noise tensor components $\mathcal{Q}^{\mu\nu}$ with them:
\begin{equation}
    \mathcal{L}_2^\varepsilon \Phi = \frac{\varepsilon^2}{2} \sum_{\mu, \nu \in \mathcal{I}} \mathcal{Q}^{\mu\nu} (D^2 \Phi)_{\mu\nu}.
\end{equation}
\begin{proofpart}[Identification of the noise tensor]
We first establish the explicit form of the noise tensor in these coordinates. Recall from Lemma \ref{lemma:metric-selection} that the weight function $\varpi$ was specifically chosen to ``whiten'' the noise on the tangent space.
We first compute the \emph{covariant} components, obtained by pairing the noise increment with the coordinate basis vectors $\partial_\mu\Psi$ in the weighted inner product. Each of the two pairings contributes one factor of $\varpi$, so the covariance carries $\varpi^2$:
\begin{equation}
    \mathcal{Q}_{\mu\nu} \coloneqq \frac{1}{dt}\,\mathbb{E}\bigl[\langle \mathcal{B}_m d\xi, \partial_\mu\Psi\rangle_\varpi \,\langle \mathcal{B}_m d\xi, \partial_\nu\Psi\rangle_\varpi\bigr]
    = \int \partial_\mu \Psi\, \partial_\nu \Psi \,\Sigma_m^2\, \varpi^2 \, dr .
\end{equation}
The operator $\mathcal{Q}^{\mathrm{loc}}_m$ is multiplication by $\Sigma_m^2 = 2D_0\,\partial_r m$, and at $m = M(\mathbf{X}) + \mathbf{n}$, with $\varpi = \varpi(\cdot\,;\mathbf{X})$ the weight of the base point, $\Sigma_m^2 = (2D_0/\varpi)(1 + \varpi\,\partial_r\mathbf{n})$. One factor of $\varpi$ therefore cancels in the first term, which is the induced metric tensor $\mathfrak{G}_{\mu\nu}$ of Definition~\ref{def:induced-metric-bundle}:
\begin{equation}
    \mathcal{Q}_{\mu\nu} = 2D_0 \int \partial_\mu \Psi\, \partial_\nu \Psi \, \varpi \, dr + 2D_0 \int \partial_\mu \Psi\, \partial_\nu \Psi \, \varpi^2\,\partial_r\mathbf{n} \, dr = 2D_0 \, \mathfrak{G}_{\mu\nu}(\mathbf{X}, \eta) + \mathcal{O}(\|\eta\|).
    \label{eq:noise-tensor-covariant}
\end{equation}
The coordinate diffusion tensor is the contravariant object, obtained by raising both indices with $\mathfrak{G}^{-1}$, exactly as $\mathcal{D}^{kl} = 2D_0 g^{kl}$ was obtained in Section~\ref{sec:reduction}:
\begin{equation}
    \mathcal{Q}^{\mu\nu} = \sum_{\mu', \nu' \in \mathcal{I}} \mathfrak{G}^{\mu\mu'}\,\mathfrak{G}^{\nu\nu'}\,\mathcal{Q}_{\mu'\nu'} = 2D_0 \, \mathfrak{G}^{\mu\nu}(\mathbf{X}, \eta) + \mathcal{O}(\|\eta\|).
    \label{eq:noise-tensor-contravariant}
\end{equation}
Thus, at $\eta = 0$ the noise tensor is exactly proportional to the inverse metric of the tubular neighbourhood. Both factors of $\varpi$ in the integrand of~\eqref{eq:noise-tensor-covariant} are essential, since with a single factor the integral would not reproduce $\mathfrak{G}_{\mu\nu}$. We now substitute the asymptotic expansion of $\mathfrak{G}^{\mu\nu}$ derived in Lemma~\ref{lemma:inverse-metric}. The tangential block is $\mathcal{Q}^{ij} = 2D_0 \bigl(g^{ij}(\mathbf{X}) + \sum_{k,l} g^{ik}(\Delta_{kl}[\mathbf{n}] + E_{kl}[\mathbf{n}]) g^{lj} + \mathcal{O}(\|\eta\|^2)\bigr)$, with $E_{kl}[\mathbf{n}] = \int_0^L \mathbf{e}_k^m\mathbf{e}_l^m\varpi^2\,\partial_r\mathbf{n}\,dr$ from the second term of~\eqref{eq:noise-tensor-covariant}, the normal block is $\mathcal{Q}^{\alpha\beta} = 2D_0 (\delta^{\alpha\beta} + \mathcal{O}(\|\eta\|))$, and the cross terms $\mathcal{Q}^{i\alpha} = \mathcal{O}(\|\eta\|)$ vanish at $\eta = 0$.
\end{proofpart}
\begin{proofpart}[The second derivatives in the Fermi coordinates]
We write $F(\mathbf{X}, \eta) \coloneqq \Phi(\Psi(\mathbf{X}, \eta))$. The second derivative along vector fields is $D^2\Phi(\mathbf{u},\mathbf{v}) = \mathbf{u}(\mathbf{v}\Phi) - (D_{\mathbf{u}}\mathbf{v})\Phi$, so that upon the coordinate vectors
\begin{equation}
    (D^2\Phi)_{\mu\nu} = \partial_\mu\partial_\nu F - D\Phi(\partial_\mu\partial_\nu\Psi) .
\end{equation}
At $\eta = 0$ we have $\partial_{X^i}\partial_{X^j}\Psi = \partial_{X^i}\mathbf{e}_j^m$, which the projectors of Definition~\ref{def:explicit-projector} split, by Eq.~\eqref{eq:projected-connection} and Definition~\ref{def:second-fundamental-form}, as
\begin{equation}
    \partial_{X^i}\mathbf{e}_j^m = \sum_{k=1}^d \tilde{\Gamma}^k_{ij}\,\mathbf{e}_k^m + \mathbf{I\!I}(\mathbf{e}_i^m, \mathbf{e}_j^m) .
\end{equation}
Since $D\Phi(\mathbf{e}_k^m) = \partial_{X^k}F$ and $D\Phi\bigl(\mathbf{I\!I}(\mathbf{e}_i^m, \mathbf{e}_j^m)\bigr) = \langle \nabla_\eta F, \mathcal{R}_{\mathbf{X}}\mathbf{I\!I}(\mathbf{e}_i^m, \mathbf{e}_j^m) \rangle_\perp$, the tangential block at $\eta = 0$ is
\begin{equation}
    (D^2\Phi)_{ij} = \frac{\partial^2 F}{\partial X^i \partial X^j} - \sum_{k=1}^d \tilde{\Gamma}^k_{ij}\,\frac{\partial F}{\partial X^k} - \bigl\langle \nabla_\eta F, \mathcal{R}_{\mathbf{X}}\mathbf{I\!I}(\mathbf{e}_i^m, \mathbf{e}_j^m) \bigr\rangle_\perp .
\end{equation}
The fibre coordinate enters $\Psi$ affinely, so that $\partial_{\eta^\alpha}\partial_{\eta^\beta}\Psi = 0$ and the normal block is $(D^2\Phi)_{\alpha\beta} = \partial_{\eta^\alpha}\partial_{\eta^\beta} F$.
\end{proofpart}
\begin{proofpart}[The tangential diffusion]
At $\eta = 0$ the tangential block of the noise tensor is $2D_0 g^{ij}$, and the contraction with the tangential block of the previous part gives
\begin{equation}
    \frac{\varepsilon^2}{2}\sum_{i,j=1}^d 2D_0\, g^{ij} (D^2\Phi)_{ij}
    = \varepsilon^2 D_0 \sum_{i,j=1}^d g^{ij}\Bigl(\frac{\partial^2 F}{\partial X^i \partial X^j} - \sum_{k=1}^d \tilde{\Gamma}^k_{ij}\frac{\partial F}{\partial X^k}\Bigr) - \varepsilon^2 D_0 \langle \mathbf{H}, \nabla_\eta F \rangle_\perp .
    \label{eq:tangential-diffusion-projected}
\end{equation}
\end{proofpart}
\begin{proofpart}[The Laplace--Beltrami operator and the connection correction]
The proof of Theorem~\ref{thm:geo-operator} writes the Laplace--Beltrami operator as $\Delta_g F = \sum_{k,l} g^{kl}\partial^2_{X^kX^l}F + \sum_l \mathcal{T}^l \partial_{X^l}F$, with $\mathcal{T}^l = \sum_k \partial_k g^{kl} + \sum_k g^{kl}\partial_k \ln\sqrt{|g|}$. Using $\partial_k g^{kl} = -\sum_{a,b} g^{ka}g^{lb}\partial_k g_{ab}$ and $\partial_k \ln\sqrt{|g|} = \tfrac12\sum_{a,b} g^{ab}\partial_k g_{ab}$,
\begin{equation}
    \mathcal{T}^l = -\sum_{k,a,b} g^{ka}g^{lb}\,\partial_k g_{ab} + \frac12\sum_{k,a,b} g^{kl}g^{ab}\,\partial_k g_{ab} = -\sum_{i,j=1}^d g^{ij}\,\Gamma^l_{ij} ,
\end{equation}
the last equality following from the formula for the Christoffel symbols in Definition~\ref{def:connection-correction} after the renamings $(k,a,b) \to (i,j,m)$ in the first sum and $(k,a,b) \to (m,i,j)$ in the second. Hence
\begin{equation}
    \sum_{i,j=1}^d g^{ij}\Bigl(\frac{\partial^2 F}{\partial X^i \partial X^j} - \sum_{k=1}^d \tilde{\Gamma}^k_{ij}\frac{\partial F}{\partial X^k}\Bigr) = \Delta_g F + \frac{1}{D_0}\sum_{k=1}^d \mathfrak{c}^k \frac{\partial F}{\partial X^k} ,
\end{equation}
by Definition~\ref{def:connection-correction}. The weight-gradient form of $\mathfrak{c}$ follows from Eq.~\eqref{eq:christoffel-isolation}, which, with its indices renamed $(l,i,j) \to (i,j,l)$ and multiplied by $g^{kl}$, reads $\tilde{\Gamma}^k_{ij} = \Gamma^k_{ij} - \tfrac12\sum_l g^{kl}(T_{ijl} + T_{jil} - T_{lji})$. Substituting, and using the symmetry of $T_{lij}$ in its last two indices,
\begin{equation}
    \mathfrak{c}^k = \frac{D_0}{2}\sum_{i,j,l=1}^d g^{ij}g^{kl}\bigl(T_{ijl} + T_{jil} - T_{lji}\bigr) = \frac{D_0}{2}\sum_{i,j,l=1}^d \bigl(2g^{ki}g^{lj} - g^{kl}g^{ij}\bigr)T_{lij} .
    \label{eq:connection-correction-T}
\end{equation}
\end{proofpart}
\begin{proofpart}[The mean curvature term]
The last term of~\eqref{eq:tangential-diffusion-projected} is $-\varepsilon^2 D_0 \langle \mathbf{H}, \nabla_\eta F \rangle_\perp$, with $\mathbf{H}$ as in the statement, and it is exact at $\eta = 0$. For a test functional of the base point alone it vanishes, since $\nabla_\eta F = 0$.
\end{proofpart}
\begin{proofpart}[The normal diffusion]
At $\eta = 0$ the normal block of the noise tensor is $2D_0\delta^{\alpha\beta}$, and its contraction with the normal block of the second part gives $\varepsilon^2 D_0 \sum_\alpha \partial^2_{\eta^\alpha\eta^\alpha}F = \varepsilon^2 D_0 \operatorname{Tr}_\perp(D^2_\eta F)$.
\end{proofpart}
\begin{proofpart}[The remainder]
For $\eta \neq 0$, the coefficients $\mathcal{Q}^{\mu\nu}(\mathbf{X},\eta)$ and the vectors $\partial_\mu\partial_\nu\Psi$ depend smoothly upon $\eta$. Subtracting their values at $\eta = 0$, which the previous parts have evaluated, leaves $\mathcal{R}_\kappa$, which therefore vanishes at $\eta = 0$. By Lemma~\ref{lemma:inverse-metric} and the identification of the noise tensor, $\mathcal{Q}^{ij} - 2D_0 g^{ij} = 2D_0\sum_{k,l} g^{ik}(\Delta_{kl}[\mathbf{n}] + E_{kl}[\mathbf{n}])g^{lj} + \mathcal{O}(\|\eta\|^2)$, which multiplies the tangential second derivatives and gives the stated second-order part.
\end{proofpart}
Summing the contributions of the third to the seventh parts gives~\eqref{eq:diffusion-decomposition}.
\end{proof}
\noindent{}The kinetic generator of Theorem~\ref{thm:commutativity} uses the same local covariance, so the two generators are compared on equal terms, and the finite-$\delta$ covariance is treated in Section~\ref{sec:commutativity}.
\subsubsection{The averaged evolution equation}
We now derive the exact evolution equation for the marginal density by integrating out the fibres.
\begin{theorem}[The marginal equation]
\label{thm:marginal-equation}
Let the generator of Definition~\ref{def:adjoint-evolution} be taken with the local covariance of Proposition~\ref{prop:diffusion-decomposition}, and assume that the weak equation of Definition~\ref{def:ffpe} extends to the functionals $\phi\circ\mathbf{X}$ with $\phi$ smooth and compactly supported in the chart, which the cylindrical class of that definition does not contain, and that the part of $\mu_t$ outside $\mathcal{D}_\alpha$, of mass $\mathbb{P}(\tau_\alpha \le t)$, is negligible at the order retained, so that the right-hand side of that equation may be integrated over the whole tubular neighbourhood. The marginal density $\bar{\varrho}_t(\mathbf{X})$ then satisfies, in the weak sense, the effective transport equation on the manifold $\mathcal{M}_m$
\begin{equation}
    \frac{\partial \bar{\varrho}_t}{\partial t} = \nabla_g \cdot \Bigl( - \bigl[\mathbb{E}^{\perp}[\overline{\mathbf{F}}] + \varepsilon^2 \mathfrak{c}\bigr] \bar{\varrho}_t + \varepsilon^2 D_0 \nabla_g \bar{\varrho}_t \Bigr) + \bar{\mathcal{R}}_t ,
    \label{eq:marginal-equation}
\end{equation}
where $\mathbb{E}^{\perp}$ is the expectation with respect to the conditional measure $\nu_t^{\mathbf{X}}$ of Proposition~\ref{prop:measure-decomposition}, $\mathfrak{c}$ is the connection correction of Definition~\ref{def:connection-correction}, and the \emph{effective drift force} is
\begin{equation}
    \overline{\mathbf{F}}^k(\mathbf{X}, \eta) = \sum_{l=1}^d \bigl[(g - \mathsf{K}[\mathbf{n}])^{-1}\bigr]_{kl} \bigl\langle \mathcal{A}_m(\Psi(\mathbf{X}, \eta)), \mathbf{e}_l^m \bigr\rangle_{\varpi} , \qquad \mathbf{n} = \mathcal{R}_{\mathbf{X}}^{-1}\eta ,
    \label{eq:effective-drift-force}
\end{equation}
the matrix $g - \mathsf{K}[\mathbf{n}]$ having entry $g_{lk} - \mathsf{K}_{lk}[\mathbf{n}]$ in row $l$ and column $k$, with
\begin{equation}
    \mathsf{K}_{lk}[\mathbf{n}] \coloneqq \bigl\langle \mathbf{n}, \partial_{X^k}\mathbf{e}_l^m \bigr\rangle_\varpi + \int_0^L \mathbf{n}\, \mathbf{e}_l^m\, \frac{\partial\varpi}{\partial X^k}\, dr .
    \label{eq:frame-matrix}
\end{equation}
The averaged remainder $\bar{\mathcal{R}}_t$ is defined by $\int_{\mathcal{X}} \phi\, \bar{\mathcal{R}}_t \, d\mathfrak{m}_g = \int_{\mathcal{X}} \mathbb{E}^{\perp}[\mathcal{R}_\kappa \phi]\, \bar{\varrho}_t\, d\mathfrak{m}_g$ for every test function $\phi$ of the base point alone, with $\mathcal{R}_\kappa$ the remainder of Proposition~\ref{prop:diffusion-decomposition}.
\end{theorem}
\begin{proof}
We derive this result by integrating the functional forward Kolmogorov equation $\langle \partial_t \mu_t, \Phi \rangle = \langle \mu_t, \mathcal{L}_\varepsilon \Phi \rangle$ against the test functional $\Phi \coloneqq \phi \circ \mathbf{X}$, with $\phi$ a function of the base point alone and $\mathbf{X}(m)$ the base point of the Fermi coordinates of $m$, so that $\Phi(\Psi(\mathbf{X},\eta)) = \phi(\mathbf{X})$, $\nabla_{\eta} \Phi = 0$ and $D_{\eta}^2 \Phi = 0$.
\begin{proofpart}[The weak formulation with density factorisation]
We substitute the measure disintegration $d\mu_t = \bar{\varrho}_t d\mathfrak{m}_g \otimes d\nu_t^{\mathbf{X}}$ into the left-hand side:
\begin{equation}
    \text{LHS} = \frac{d}{dt} \int_{\mathcal{X}} \phi(\mathbf{X}) \bar{\varrho}_t(\mathbf{X}) \, d\mathfrak{m}_g(\mathbf{X}).
\end{equation}
For the right-hand side, we use the density factorisation established in Proposition \ref{prop:gibbs-factorisation}, $\varrho_t = \bar{\varrho}_t \Gamma_t J_{\Psi}^{-1}$, integrated against the reference measure $\lambda = \mathfrak{m}_g \otimes \lambda_{\perp}$:
\begin{equation}
    \text{RHS} = \int_{\mathcal{X}} \left[ \int_{\mathcal{H}_{\perp}} (\mathcal{L}_\varepsilon \Phi)(\mathbf{X}, \eta) \, \Gamma_t(\mathbf{X}, \eta) J_{\Psi}^{-1}(\mathbf{X}, \eta) \, d\lambda_{\perp}(\eta) \right] \bar{\varrho}_t(\mathbf{X}) \, d\mathfrak{m}_g(\mathbf{X}).
\end{equation}
\end{proofpart}
\begin{proofpart}[The drift component]
The drift term of the generator is the derivative of $\Phi$ along the drift, $\langle \mathcal{A}_m, D\Phi \rangle_\varpi = D\Phi(m)[\mathcal{A}_m(m)]$, and the chain rule gives $D\Phi(m)[w] = \sum_{k} \partial_{X^k}\phi\; D\mathbf{X}^k(m)[w]$. The base point is fixed by the orthogonality condition of Lemma~\ref{lemma:orthogonality-condition}, $\langle m - M(\mathbf{X}(m)), \mathbf{e}_l^m(\mathbf{X}(m)) \rangle_{\varpi(\cdot;\mathbf{X}(m))} = 0$ for $l = 1, \dots, d$. Differentiating it at $m = \Psi(\mathbf{X},\eta)$ along $w$, the dependence upon $m$ contributes $\langle w, \mathbf{e}_l^m \rangle_\varpi$, and the dependence upon the base point contributes, since $m - M(\mathbf{X}) = \mathbf{n}$,
\begin{equation}
    \begin{split}
    \frac{\partial}{\partial X^k}\bigl\langle m - M(\mathbf{X}), \mathbf{e}_l^m(\mathbf{X}) \bigr\rangle_{\varpi(\cdot;\mathbf{X})}
    &= -\bigl\langle \mathbf{e}_k^m, \mathbf{e}_l^m \bigr\rangle_\varpi + \bigl\langle \mathbf{n}, \partial_{X^k}\mathbf{e}_l^m \bigr\rangle_\varpi + \int_0^L \mathbf{n}\,\mathbf{e}_l^m\,\frac{\partial\varpi}{\partial X^k}\,dr
    \\
    &= -g_{lk} + \mathsf{K}_{lk}[\mathbf{n}] .
    \end{split}
\end{equation}
Hence $\sum_{k} \bigl(g_{lk} - \mathsf{K}_{lk}[\mathbf{n}]\bigr)\,D\mathbf{X}^k(m)[w] = \langle w, \mathbf{e}_l^m \rangle_\varpi$, and solving,
\begin{equation}
    D\mathbf{X}^k(m)[w] = \sum_{l=1}^d \bigl[(g - \mathsf{K}[\mathbf{n}])^{-1}\bigr]_{kl} \bigl\langle w, \mathbf{e}_l^m \bigr\rangle_\varpi ,
\end{equation}
the matrix being invertible within the tubular neighbourhood, where Lemma~\ref{lemma:tubular-existence} gives the base point as a smooth function of $m$. The drift term is therefore $\sum_{k} \overline{\mathbf{F}}^k\,\partial_{X^k}\phi$, with $\overline{\mathbf{F}}$ as in~\eqref{eq:effective-drift-force}. At $\eta = 0$ the matrix $\mathsf{K}$ vanishes, and $\overline{\mathbf{F}}^k$ reduces to $\sum_{l} g^{kl}\langle \mathcal{A}_m(M), \mathbf{e}_l^m \rangle_\varpi$.
\end{proofpart}
\begin{proofpart}[The diffusion component]
Proposition~\ref{prop:diffusion-decomposition} applies to $\Phi$, for which $\nabla_\eta\Phi = 0$ and $D^2_\eta\Phi = 0$, so that its mean curvature and normal terms vanish and $\mathcal{L}_2^\varepsilon\Phi = \varepsilon^2 D_0\Delta_g\phi + \varepsilon^2\sum_{k} \mathfrak{c}^k\,\partial_{X^k}\phi + \mathcal{R}_\kappa\Phi$. Combining the two components,
\begin{equation}
    \mathcal{L}_\varepsilon\Phi = \sum_{k=1}^d \bigl(\overline{\mathbf{F}}^k + \varepsilon^2\mathfrak{c}^k\bigr)\frac{\partial \phi}{\partial X^k} + \varepsilon^2 D_0 \Delta_g\phi + \mathcal{R}_\kappa\Phi .
\end{equation}
\end{proofpart}
\begin{proofpart}[Fibre integration]
We substitute the expanded generator into the right-hand side of the first part. By Eq.~\eqref{eq:gibbs-conditional-expectation}, the fibre integral against $\Gamma_t J_{\Psi}^{-1}\,d\lambda_{\perp}$ is the expectation $\mathbb{E}^{\perp}$, and with the normalisation $\mathbb{E}^{\perp}[1] = 1$ the terms which do not depend upon $\eta$ leave it:
\begin{equation}
    \text{RHS} = \int_{\mathcal{X}} \Bigl( \sum_{k=1}^d \bigl(\mathbb{E}^{\perp}[\overline{\mathbf{F}}^k] + \varepsilon^2\mathfrak{c}^k\bigr)\frac{\partial \phi}{\partial X^k} + \varepsilon^2 D_0\Delta_g\phi + \mathbb{E}^{\perp}[\mathcal{R}_\kappa\Phi] \Bigr) \bar{\varrho}_t\, d\mathfrak{m}_g .
\end{equation}
\end{proofpart}
\begin{proofpart}[Integration by parts on the manifold]
We now isolate the equation for $\bar{\varrho}_t$ by moving derivatives from $\phi$ to the density using Green's identity on $(\mathcal{M}_m, g)$. For the drift term,
\begin{equation}
    \int_{\mathcal{X}} \Bigl(\sum_{k=1}^d \bigl(\mathbb{E}^{\perp}[\overline{\mathbf{F}}^k] + \varepsilon^2\mathfrak{c}^k\bigr)\frac{\partial \phi}{\partial X^k}\Bigr) \bar{\varrho}_t \, d\mathfrak{m}_g = - \int_{\mathcal{X}} \phi \, \nabla_g \cdot \Bigl( \bigl[\mathbb{E}^{\perp}[\overline{\mathbf{F}}] + \varepsilon^2\mathfrak{c}\bigr] \bar{\varrho}_t \Bigr) \, d\mathfrak{m}_g.
\end{equation}
For the diffusion term, since $\Delta_g$ is self-adjoint with respect to $\mathfrak{m}_g$,
\begin{equation}
    \int_{\mathcal{X}} (\Delta_g \phi) (\varepsilon^2 D_0 \bar{\varrho}_t) \, d\mathfrak{m}_g = \int_{\mathcal{X}} \phi \, \Delta_g (\varepsilon^2 D_0 \bar{\varrho}_t) \, d\mathfrak{m}_g.
\end{equation}
For the remainder, the definition of $\bar{\mathcal{R}}_t$ in the statement gives $\int_{\mathcal{X}} \mathbb{E}^{\perp}[\mathcal{R}_\kappa\Phi]\, \bar{\varrho}_t \, d\mathfrak{m}_g = \int_{\mathcal{X}} \phi\, \bar{\mathcal{R}}_t \, d\mathfrak{m}_g$.
Equating the LHS $\frac{d}{dt} \int \phi \bar{\varrho}_t d\mathfrak{m}_g$ with the processed RHS for arbitrary $\phi$, and writing $\varepsilon^2 D_0 \Delta_g \bar{\varrho}_t = \nabla_g\cdot(\varepsilon^2 D_0 \nabla_g \bar{\varrho}_t)$, we recover~\eqref{eq:marginal-equation}.
\end{proofpart}
\end{proof}
\subsection{The adiabatic limit and concentration of measure}
\label{subsec:adiabatic-limit}
The marginal equation of Theorem~\ref{thm:marginal-equation} describes the transport of the marginal density $\bar{\varrho}_t$ subject to the effective force $\mathbb{E}^{\perp}[\overline{\mathbf{F}}]$, and this force depends on the unclosed conditional measure $\nu_t^{\mathbf{X}}$. We close it by the adiabatic separation of time-scales, assuming that the fast variables $\eta$ equilibrate instantaneously to the local geometry of the fibre, and we perform the asymptotic expansion of the effective force in the limit $\varepsilon \to 0$. In this way we derive the explicit form of the corrections of order $\varepsilon^2$ produced by the fibre average (the statistical anomaly and the frame correction), and we keep them in full, since the comparison of Section~\ref{sec:commutativity} requires precisely these terms.
\subsubsection{The concentration hypothesis}
The conditional measure $\nu_t^{\mathbf{X}}$ is taken to be governed by the linear stability of the manifold (Hypothesis~\ref{hyp:stiffness}). In the limit of weak noise the transverse fluctuations become Gaussian, and the mean curvature of the manifold displaces their mean from it at order $\varepsilon^2$.
\begin{hypothesis}[Adiabatic concentration]
\label{hyp:concentration}
Let $\nu_\varepsilon^{\mathbf{X}}$ be the conditional probability measure on the fibre $\mathcal{H}_{\perp}$ at configuration $\mathbf{X}$ in the limit of timescale separation. Under the stiffness gap $\lambda_0 \gg 1$ (Hypothesis~\ref{hyp:stiffness} and Remark~\ref{rem:singular-perturbation}), the fast Ornstein--Uhlenbeck dynamics on the normal fibre equilibrate on the timescale $\tau_\perp = \lambda_0^{-1}$. We assume that the deterministic drift is tangent to the manifold, $\Pi^{\perp}_{\mathbf{X}}\mathcal{A}_m(M(\mathbf{X})) = 0$ for every $\mathbf{X} \in \mathcal{X}$, so that the fibre law is centred in the absence of noise, and that:
\begin{enumerate}[(i)]
    \item \emph{Gaussian limit:} The measure converges to the Gaussian distribution $\mathrm{N}(\varepsilon^2\hat\eta(\mathbf{X}), \varepsilon^2 \mathbf{C}_{\perp}(\mathbf{X}))$, where $\mathbf{C}_{\perp}(\mathbf{X})$ is the transverse covariance operator on $\mathcal{H}_{\perp}$ and the \emph{mean displacement} $\hat\eta(\mathbf{X}) \in \mathcal{H}_{\perp}$ is the solution of
    \begin{equation}
        \mathcal{L}_\perp\,\hat\eta = D_0\,\mathbf{H} ,
        \label{eq:mean-displacement}
    \end{equation}
    with $\mathbf{H}$ the mean curvature vector of Proposition~\ref{prop:diffusion-decomposition} and $\mathcal{L}_\perp$ read in the model fibre through the trivialisation $\mathcal{R}_{\mathbf{X}}$.
    \item \emph{Moment hierarchy:} The moments satisfy:
    \begin{align}
        \mathbb{E}^{\perp}[\eta] &= \varepsilon^2\hat\eta(\mathbf{X}) + \mathcal{O}(\varepsilon^4), \\
        \mathbb{E}^{\perp}[\eta \otimes \eta] &= \varepsilon^2 \mathbf{C}_{\perp}(\mathbf{X}) + \mathcal{O}(\varepsilon^4), \\
        \mathbb{E}^{\perp}[\eta \otimes \eta \otimes \eta] &= \mathcal{O}(\varepsilon^4), \\
        \mathbb{E}^{\perp}[\|\eta\|^k] &= \mathcal{O}(\varepsilon^k) \quad \text{for } k \ge 3.
    \end{align}
\end{enumerate}
\end{hypothesis}
\begin{remark}[The mean displacement]
\label{rem:mean-displacement}
We now check that the displacement~\eqref{eq:mean-displacement} is the stationary mean of the transverse dynamics. By Proposition~\ref{prop:diffusion-decomposition}, the generator carries the drift $-\varepsilon^2 D_0\langle \mathbf{H}, \nabla_\eta\Phi\rangle_\perp$ in the fibre, and hence the Ornstein--Uhlenbeck equation of Remark~\ref{rem:singular-perturbation} acquires a constant forcing,
\begin{equation}
    d\eta = \bigl[\mathcal{L}_\perp\eta - \varepsilon^2 D_0\,\mathbf{H}\bigr]\,dt + \varepsilon\,\Pi^\perp_{\mathbf{X}}\mathcal{B}_m\,d\xi^\delta_t ,
\end{equation}
and in the stationary state this gives $\mathcal{L}_\perp\mathbb{E}^\perp[\eta] = \varepsilon^2 D_0\,\mathbf{H}$. When the drift is not tangent to the manifold, the same argument gives a displacement $-\mathcal{L}_\perp^{-1}\Pi^\perp_{\mathbf{X}}\mathcal{A}_m(M)$ of order unity in $\varepsilon$, and we exclude it through the tangency condition of Hypothesis~\ref{hyp:concentration} (the kinetic projection neglects it as well). Under the scalar closure of Hypothesis~\ref{hyp:kinematic-covariance} we therefore have $\hat\eta = -2\tau D_0\,\mathbf{H}$, since $\mathcal{L}_\perp$ then acts as $-1/(2\tau)$. This displacement is of first order in the relaxation time, as $\mathbf{C}_\perp$ is, and we do not keep the corrections of second order in $\tau$ produced by the curvature of the drift and by the slow motion of the base point.
\end{remark}
\begin{remark}[Formal status of the adiabatic limit]
\label{rem:adiabatic-formal-status}
Note that the concentration hypothesis above is a \emph{formal asymptotic assumption}, not a theorem. Converting it into a rigorous singular perturbation result would require: (i)~establishing tightness of the family of conditional measures $\{\nu_\varepsilon^{\mathbf{X}}\}_{\varepsilon > 0}$ in a suitable Wasserstein or total-variation topology on $\mathcal{H}_\perp$, (ii)~proving weak convergence $\nu_\varepsilon^{\mathbf{X}} \Rightarrow \mathcal{N}(\varepsilon^2\hat\eta, \varepsilon^2 \mathbf{C}_\perp)$ via a Khasminskii-type averaging theorem for infinite-dimensional diffusions, and (iii)~controlling the remainder $\mathcal{O}(\varepsilon^k)$ uniformly in $\mathbf{X}$. Analogous results are, however, classical, both for finite-dimensional systems~\cite{freidlin2012random} and at the abstract level of perturbed operator semigroups~\cite{kurtz1973semigroup}, and extensions to semilinear SPDEs have been developed by Cerrai and Freidlin~\cite{cerrai2009averaging}. The full convergence proof for the present non-linear, multiplicative-noise, infinite-dimensional setting lies beyond the scope of this work. We therefore read the expansion derived below as a formal matched-asymptotic calculation. Moreover, the comparison of Theorem~\ref{thm:commutativity} cannot be offered as a~posteriori evidence for the validity of the expansion, because that comparison is itself derived \emph{from} the expansion, and so the two are not independent.
\end{remark}
\subsubsection{Derivation of the effective drift}
We now calculate the explicit form of the effective force $\mathbb{E}^{\perp}[\overline{\mathbf{F}}]$, with $\overline{\mathbf{F}}$ the effective drift force of Eq.~\eqref{eq:effective-drift-force}.
\begin{lemma}[Expansion of the effective drift]
\label{lemma:mechanical-drift-explicit}
Under Hypothesis~\ref{hyp:concentration}, the effective drift force of Theorem~\ref{thm:marginal-equation} satisfies
\begin{equation}
    \mathbb{E}^{\perp}\bigl[\overline{\mathbf{F}}^k\bigr] = F_0^k(\mathbf{X}) + \varepsilon^2 \bigl( \Upsilon^k(\mathbf{X}) + \mathfrak{f}^k(\mathbf{X}) \bigr) + \mathcal{O}(\varepsilon^4) .
    \label{eq:effective-drift-expansion}
\end{equation}
Here $F_0^k = \sum_l g^{kl}F_l$, with $F_l(\mathbf{X}) = -D_0\,\partial_{X^l}\mathcal{V}$ the standard force. The \emph{statistical anomaly} is $\Upsilon^k = \sum_l g^{kl}\Upsilon_l$, with
\begin{equation}
    \Upsilon_l(\mathbf{X}) = \frac{1}{2} \operatorname{Tr}_{\perp} \Bigl( \mathbf{C}_{\perp}(\mathbf{X}) \bigl[ \Pi_l^{\parallel} D^2 \mathcal{A}_m[M(\mathbf{X})] \bigr] \Bigr) ,
    \label{eq:statistical-anomaly}
\end{equation}
where $\Pi_l^{\parallel}$ pairs the ambient vector with $\mathbf{e}_l^m$ in the weighted inner product at the base point. The \emph{frame correction} is
\begin{equation}
    \begin{split}
    \mathfrak{f}^k(\mathbf{X}) &= \sum_{i,j,l=1}^d g^{ki}g^{jl} \Bigl( \langle \mathbf{C}_\perp \kappa_{ij}, p_l \rangle_\perp + \sum_{a,b=1}^d g^{ab} \langle \mathbf{C}_\perp \kappa_{ia}, \kappa_{bj} \rangle_\perp F_l \Bigr)
    \\
    &\quad + \sum_{l=1}^d g^{kl}\langle p_l, \hat\eta \rangle_\perp + \sum_{i,j,l=1}^d g^{ki}\langle \kappa_{ij}, \hat\eta \rangle_\perp\, g^{jl}F_l ,
    \end{split}
    \label{eq:frame-correction}
\end{equation}
with $\kappa_{ij}, p_l \in \mathcal{H}_\perp$ the vectors representing the linear functionals $\eta \mapsto \mathsf{K}_{ij}[\mathcal{R}_{\mathbf{X}}^{-1}\eta]$ of Eq.~\eqref{eq:frame-matrix} and $\eta \mapsto \langle D\mathcal{A}_m[M](\mathcal{R}_{\mathbf{X}}^{-1}\eta), \mathbf{e}_l^m \rangle_\varpi$, the last two terms being the drift which the mean displacement $\hat\eta$ of Hypothesis~\ref{hyp:concentration} produces.
\end{lemma}
\begin{proof}
We expand the two factors of Eq.~\eqref{eq:effective-drift-force} about $\eta = 0$, with $\mathbf{n} = \mathcal{R}_{\mathbf{X}}^{-1}\eta$, multiply the expansions, and take the expectation order by order in $\eta$.
\begin{proofpart}[The frame matrix]
The matrix $\mathsf{K}[\mathbf{n}]$ of Eq.~\eqref{eq:frame-matrix} is linear in $\mathbf{n}$, so that $\mathsf{K}[\mathbf{n}] = \mathcal{O}(\|\eta\|)$. Writing $g - \mathsf{K} = g(I - g^{-1}\mathsf{K})$ and summing the Neumann series of $(I - g^{-1}\mathsf{K})^{-1}$,
\begin{equation}
    \begin{split}
    (g - \mathsf{K}[\mathbf{n}])^{-1} &= (I - g^{-1}\mathsf{K}[\mathbf{n}])^{-1} g^{-1}
    \\
    &= g^{-1} + g^{-1}\mathsf{K}[\mathbf{n}]\,g^{-1} + g^{-1}\mathsf{K}[\mathbf{n}]\,g^{-1}\mathsf{K}[\mathbf{n}]\,g^{-1} + \mathcal{O}(\|\eta\|^3) ,
    \end{split}
\end{equation}
and in components
\begin{equation}
    \bigl[(g - \mathsf{K}[\mathbf{n}])^{-1}\bigr]_{kl} = g^{kl} + \sum_{i,j=1}^d g^{ki}\mathsf{K}_{ij}[\mathbf{n}]\,g^{jl} + \sum_{i,j,a,b=1}^d g^{ki}\mathsf{K}_{ia}[\mathbf{n}]\,g^{ab}\mathsf{K}_{bj}[\mathbf{n}]\,g^{jl} + \mathcal{O}(\|\eta\|^3) .
\end{equation}
\end{proofpart}
\begin{proofpart}[The projected drift]
Since $\Psi(\mathbf{X},\eta) = M(\mathbf{X}) + \mathbf{n}$ and $\eta \mapsto \mathbf{n}$ is linear, Taylor's theorem applied to $\mathcal{A}_m$ at $M(\mathbf{X})$ gives
\begin{equation}
    \begin{split}
    \bigl\langle \mathcal{A}_m(\Psi(\mathbf{X},\eta)), \mathbf{e}_l^m \bigr\rangle_\varpi &= F_l(\mathbf{X}) + \bigl\langle D\mathcal{A}_m[M](\mathbf{n}), \mathbf{e}_l^m \bigr\rangle_\varpi
    \\
    &\quad + \frac{1}{2}\bigl\langle D^2\mathcal{A}_m[M](\mathbf{n},\mathbf{n}), \mathbf{e}_l^m \bigr\rangle_\varpi + \mathcal{O}(\|\eta\|^3) ,
    \end{split}
\end{equation}
where the term of order zero is $\langle \mathcal{A}_m(M), \mathbf{e}_l^m \rangle_\varpi = F_l(\mathbf{X})$.
\end{proofpart}
\begin{proofpart}[The product]
Substituting the two expansions into Eq.~\eqref{eq:effective-drift-force} and collecting the terms by their degree in $\eta$, the term of degree zero is $\sum_l g^{kl}F_l = F_0^k$, the terms of degree one are
\begin{equation}
    \sum_{l=1}^d g^{kl}\bigl\langle D\mathcal{A}_m[M](\mathbf{n}), \mathbf{e}_l^m \bigr\rangle_\varpi + \sum_{i,j,l=1}^d g^{ki}\mathsf{K}_{ij}[\mathbf{n}]\,g^{jl}F_l ,
\end{equation}
and the terms of degree two are
\begin{equation}
    \begin{split}
    &\frac{1}{2}\sum_{l=1}^d g^{kl}\bigl\langle D^2\mathcal{A}_m[M](\mathbf{n},\mathbf{n}), \mathbf{e}_l^m \bigr\rangle_\varpi + \sum_{i,j,l=1}^d g^{ki}\mathsf{K}_{ij}[\mathbf{n}]\,g^{jl}\bigl\langle D\mathcal{A}_m[M](\mathbf{n}), \mathbf{e}_l^m \bigr\rangle_\varpi
    \\
    &\qquad + \sum_{i,j,l,a,b=1}^d g^{ki}\mathsf{K}_{ia}[\mathbf{n}]\,g^{ab}\mathsf{K}_{bj}[\mathbf{n}]\,g^{jl}F_l .
    \end{split}
\end{equation}
The remaining terms are of order $\|\eta\|^3$, and they consist of a homogeneous polynomial of degree three in $\eta$ and a remainder of order $\|\eta\|^4$.
\end{proofpart}
\begin{proofpart}[The expectations]
By the definitions of $\kappa_{ij}$ and $p_l$, $\mathsf{K}_{ij}[\mathbf{n}] = \langle \kappa_{ij}, \eta \rangle_\perp$ and $\langle D\mathcal{A}_m[M](\mathbf{n}), \mathbf{e}_l^m \rangle_\varpi = \langle p_l, \eta \rangle_\perp$. By the first moment of Hypothesis~\ref{hyp:concentration}~(ii), the terms of degree one therefore have the mean
\begin{equation}
    \varepsilon^2 \sum_{l=1}^d g^{kl}\langle p_l, \hat\eta \rangle_\perp + \varepsilon^2 \sum_{i,j,l=1}^d g^{ki}\langle \kappa_{ij}, \hat\eta \rangle_\perp\, g^{jl}F_l + \mathcal{O}(\varepsilon^4) ,
    \label{eq:degree-one-mean}
\end{equation}
and the second moment of the same item gives, for any $a, b \in \mathcal{H}_\perp$,
\begin{equation}
    \mathbb{E}^{\perp}\bigl[\langle a, \eta \rangle_\perp \langle b, \eta \rangle_\perp\bigr] = \bigl\langle a \otimes b, \mathbb{E}^{\perp}[\eta \otimes \eta] \bigr\rangle = \varepsilon^2 \langle \mathbf{C}_\perp a, b \rangle_\perp + \mathcal{O}(\varepsilon^4) .
\end{equation}
Applied to the second and the third terms of degree two, this identity gives, at order $\varepsilon^2$,
\begin{equation}
    \varepsilon^2 \sum_{i,j,l=1}^d g^{ki}g^{jl}\langle \mathbf{C}_\perp \kappa_{ij}, p_l \rangle_\perp + \varepsilon^2 \sum_{i,j,l,a,b=1}^d g^{ki}g^{ab}g^{jl}\langle \mathbf{C}_\perp \kappa_{ia}, \kappa_{bj} \rangle_\perp F_l ,
\end{equation}
and adding the mean~\eqref{eq:degree-one-mean} of the terms of degree one gives $\varepsilon^2 \mathfrak{f}^k$, which is Eq.~\eqref{eq:frame-correction}. The first term of degree two is the quadratic form $\eta \mapsto \frac{1}{2}\sum_l g^{kl}[\Pi_l^{\parallel} D^2\mathcal{A}_m[M]](\eta,\eta)$, and the expectation of a quadratic form against a measure with second moment $\varepsilon^2\mathbf{C}_\perp$ is the trace of the form against that covariance, so that
\begin{equation}
    \begin{split}
    \mathbb{E}^{\perp}\Bigl[\frac{1}{2}\sum_{l=1}^d g^{kl}\bigl[\Pi_l^{\parallel} D^2\mathcal{A}_m[M]\bigr](\eta,\eta)\Bigr] &= \frac{\varepsilon^2}{2}\sum_{l=1}^d g^{kl}\operatorname{Tr}_{\perp}\Bigl( \mathbf{C}_{\perp}(\mathbf{X}) \bigl[ \Pi_l^{\parallel} D^2 \mathcal{A}_m[M(\mathbf{X})] \bigr] \Bigr)
    \\
    &= \varepsilon^2 \Upsilon^k ,
    \end{split}
\end{equation}
which is Eq.~\eqref{eq:statistical-anomaly}. The homogeneous polynomial of degree three has a mean of order $\varepsilon^4$ by the third moment of Hypothesis~\ref{hyp:concentration}~(ii), and the remainder has mean of order $\mathbb{E}^{\perp}[\|\eta\|^4] = \mathcal{O}(\varepsilon^4)$ by item~(ii). Adding the three orders gives Eq.~\eqref{eq:effective-drift-expansion}. The explicit evaluation of $\Upsilon$, which requires two further functional derivatives of $\mathcal{A}_m$ and the trace against the projected covariance in kernel form, is deferred to the proof of Theorem~\ref{thm:commutativity} (Step~3), where it is required in full generality and derived there step by step.
\end{proofpart}
\end{proof}
\begin{remark}[The ideal gas carries no statistical anomaly]
\label{rem:ideal-gas-no-anomaly}
In the case of the ideal gas, the chemical potential is given by $\mu = \beta^{-1}\ln\rho$, with $\rho = (4\pi r^2)^{-1}\partial_r m$ (Remark~\ref{rem:inverse-density}). Substituting the drift of Eq.~\eqref{eq:mass-drift-explicit} together with Lemma~\ref{lemma:mu-duality}, which gives $\delta\mathcal{F}_m/\delta m = -\partial_r\mu$,
\begin{equation}
    \mathcal{A}_m(m) = D_0\beta\,(\partial_r m)\,(\partial_r\mu) = D_0\,(\partial_r m)\left[\frac{\partial_r^2 m}{\partial_r m} - \frac{2}{r}\right] = D_0\left[\partial_r^2 m - \frac{2}{r}\,\partial_r m\right] ,
\end{equation}
which is linear in $m$ (it is the diffusion equation $\partial_t\rho = D_0\Delta\rho$ written for the cumulative mass). As a result, the second derivative $D^2\mathcal{A}_m$ vanishes, and so does $\Upsilon$ (Eq.~\eqref{eq:statistical-anomaly}). This means that the statistical anomaly is produced by the excess free energy alone, while the frame correction and the connection correction are present for the ideal gas as well.
\end{remark}
\begin{remark}[Big-$\mathcal{O}$ hierarchy of the effective drift]
\label{rem:big-o-hierarchy}
Combining Lemma~\ref{lemma:mechanical-drift-explicit} with Theorem~\ref{thm:marginal-equation}, the drift of the marginal equation admits the ordering
\begin{equation}
    \mathbb{E}^{\perp}\bigl[\overline{\mathbf{F}}^k\bigr] + \varepsilon^2\mathfrak{c}^k = F_0^k(\mathbf{X}) + \varepsilon^2 \bigl( \Upsilon^k(\mathbf{X}) + \mathfrak{f}^k(\mathbf{X}) + \mathfrak{c}^k(\mathbf{X}) \bigr) + \mathcal{O}(\varepsilon^4) ,
\end{equation}
and the averaged remainder $\bar{\mathcal{R}}_t$ is of order $\varepsilon^4$, since $\mathcal{R}_\kappa$ carries the factor $\varepsilon^2$ and vanishes at $\eta = 0$. Hence its part linear in $\eta$ has a mean of order $\varepsilon^4$ (the mean of $\eta$ is of order $\varepsilon^2$ under Hypothesis~\ref{hyp:concentration}), while the rest is of order $\mathbb{E}^{\perp}[\|\eta\|^2] = \mathcal{O}(\varepsilon^2)$. In the $\mathcal{O}(\varepsilon^2)$ theory, only the thermodynamic force and the three corrections $\Upsilon$, $\mathfrak{f}$ and $\mathfrak{c}$ survive.
\end{remark}
\subsubsection{The reduced Fokker--Planck equation}
Substituting the expansion of Lemma~\ref{lemma:mechanical-drift-explicit} into Theorem~\ref{thm:marginal-equation}, with the order of the averaged remainder recorded in Remark~\ref{rem:big-o-hierarchy}, we obtain the closed equation.
\begin{theorem}[The reduced covariant Fokker--Planck equation]
\label{thm:reduced-fpe}
In the adiabatic limit $\varepsilon \to 0$, retaining terms up to $\mathcal{O}(\varepsilon^2)$, the marginal probability density $\bar{\varrho}_t(\mathbf{X})$ on the manifold $\mathcal{M}_m$ evolves according to:
\begin{equation}
    \frac{\partial \bar{\varrho}_t}{\partial t} = \nabla_g \cdot \left( -\left[ \mathbf{F}_0(\mathbf{X}) + \varepsilon^2 \bigl( \Upsilon(\mathbf{X}) + \mathfrak{f}(\mathbf{X}) + \mathfrak{c}(\mathbf{X}) \bigr) \right] \bar{\varrho}_t + \varepsilon^2 D_0 \nabla_g \bar{\varrho}_t \right).
\end{equation}
The terms are defined as follows:
\begin{enumerate}[(i)]
    \item \emph{Thermodynamic driving force:} $\mathbf{F}_0 = - D_0 \nabla \mathcal{V}$. See Remark~\ref{rem:dropped-strat-term} for the finite-$\delta$ term of Eq.~\eqref{eq:S-projected} that is \emph{not} absorbed into $\mathbf{F}_0$.
    \item \emph{Statistical anomaly:} $\Upsilon(\mathbf{X})$ is the deterministic drift correction arising from the interaction of the transverse noise with the curvature of the drift field, with components $\Upsilon^k = \sum_{l} g^{kl}\Upsilon_l$ and $\Upsilon_l$ the trace of Eq.~\eqref{eq:statistical-anomaly}, in which the second derivative of the projected drift $\langle \mathcal{A}_m, \mathbf{e}_l^m \rangle_{\varpi}$ restricted to the normal fibre is contracted against the transverse covariance, and the inverse metric raises the index that the pairing lowers.
    \item \emph{Frame correction:} $\mathfrak{f}(\mathbf{X})$ is the drift correction of Eq.~\eqref{eq:frame-correction}, arising from the dependence of the Fermi coordinates upon the fluctuation and from the mean displacement of the fibre law.
    \item \emph{Connection correction:} $\mathfrak{c}(\mathbf{X})$ is the drift correction of Definition~\ref{def:connection-correction}, arising from the difference between the projected connection and the Levi-Civita connection.
\end{enumerate}
\end{theorem}
\begin{proof}
Theorem~\ref{thm:marginal-equation} gives Eq.~\eqref{eq:marginal-equation},
\begin{equation}
    \frac{\partial \bar{\varrho}_t}{\partial t} = \nabla_g \cdot \Bigl( - \bigl[\mathbb{E}^{\perp}[\overline{\mathbf{F}}] + \varepsilon^2 \mathfrak{c}\bigr] \bar{\varrho}_t + \varepsilon^2 D_0 \nabla_g \bar{\varrho}_t \Bigr) + \bar{\mathcal{R}}_t ,
    \label{eq:reduced-fpe-start}
\end{equation}
with the fibre-averaged force left unevaluated. Lemma~\ref{lemma:mechanical-drift-explicit} evaluates it, and adding the connection correction, as recorded in Remark~\ref{rem:big-o-hierarchy},
\begin{equation}
    \begin{aligned}
    \mathbb{E}^{\perp}\bigl[\overline{\mathbf{F}}^k\bigr] + \varepsilon^2\mathfrak{c}^k &= F_0^k(\mathbf{X}) + \varepsilon^2\bigl(\Upsilon^k(\mathbf{X}) + \mathfrak{f}^k(\mathbf{X}) + \mathfrak{c}^k(\mathbf{X})\bigr) + \mathcal{O}(\varepsilon^4) ,
    \\
    F_0^k &\coloneqq \sum_l g^{kl}F_l = -D_0\,\nabla^k\mathcal{V} .
    \end{aligned}
    \label{eq:reduced-fpe-force}
\end{equation}
Substituting~\eqref{eq:reduced-fpe-force} into~\eqref{eq:reduced-fpe-start}, with $\bar{\mathcal{R}}_t = \mathcal{O}(\varepsilon^4)$ by the same remark,
\begin{align}
    \frac{\partial \bar{\varrho}_t}{\partial t}
    &= \nabla_g \cdot \Bigl( -\bigl[\mathbf{F}_0(\mathbf{X}) + \varepsilon^2\bigl(\Upsilon(\mathbf{X}) + \mathfrak{f}(\mathbf{X}) + \mathfrak{c}(\mathbf{X})\bigr) + \mathcal{O}(\varepsilon^4)\bigr]\, \bar{\varrho}_t + \varepsilon^2 D_0 \nabla_g \bar{\varrho}_t \Bigr) + \mathcal{O}(\varepsilon^4) \nonumber \\
    &= \nabla_g \cdot \Bigl( -\bigl[\mathbf{F}_0(\mathbf{X}) + \varepsilon^2\bigl(\Upsilon(\mathbf{X}) + \mathfrak{f}(\mathbf{X}) + \mathfrak{c}(\mathbf{X})\bigr)\bigr]\, \bar{\varrho}_t + \varepsilon^2 D_0 \nabla_g \bar{\varrho}_t \Bigr) + \mathcal{O}(\varepsilon^4) ,
\end{align}
and dropping the $\mathcal{O}(\varepsilon^4)$ terms gives the stated equation at the order retained.
\end{proof}
\begin{remark}[Where the finite-$\delta$ drift term belongs]
\label{rem:dropped-strat-term}
Corollary~\ref{cor:coordinate-sde} gives the projected Stratonovich drift as
\begin{equation}
    F_l(\mathbf{X}) - \varepsilon^2 S_l(\mathbf{X}) = -D_0\,\partial_{X^l}\mathcal{V} \;-\; \varepsilon^2\,\bigl\langle \mathfrak{S}_m^\delta, \mathbf{e}^m_l\bigr\rangle_{\varpi},
\end{equation}
while $\mathbf{F}_0$ above retains only the first term. However, the second term is not negligible at the order of the comparison. If we evaluate it with $\sigma = \sqrt{2D_0/\varpi}$, for which $\sigma'/\sigma = -\varpi'/(2\varpi)$, Lemma~\ref{lemma:ito-strat-correction} gives
\begin{equation}
    \bigl\langle \mathfrak{S}^\delta_m, \mathbf{e}^m_l\bigr\rangle_\varpi
    = \frac{D_0}{4}\int_0^L \bigl[-\varpi'\,V_\delta + \varpi\, V_\delta'\bigr]\, \mathbf{e}^m_l \, dr ,
    \label{eq:dropped-term}
\end{equation}
which belongs to the same family of weighted kernel integrals as the geometric anomaly (it is built from $V_\delta$, $V_\delta'$ and the weight together with its derivative) and is of order $\mathcal{O}(\varepsilon^2)$, precisely the order of the comparison. We find this term in $\mathcal{L}_A$ but not in $\mathcal{L}_B$, because the kinetic route projects the Stratonovich equation and converts to It\^{o} afterwards, and so carries the field-level conversion term through the projection, while the probabilistic route converts at the level of the field, where the same term cancels exactly (Remark~\ref{rem:ito-start}). Hence it is carried in the geometric anomaly, which Theorem~\ref{thm:geo-operator} defines as the whole of the reduced drift that lies outside the Laplace--Beltrami operator and the gradient of $\mathcal{V}$. In this way we have $\mathbf{F}_0 = -D_0\nabla\mathcal{V}$ throughout, and Eq.~\eqref{eq:iff-condition} requires no further term. Note that the convention adopted for the parent equation does not affect the criterion. If we wrote~\eqref{eq:full_spde} with the bare Stratonovich drift, this contribution would move from $\mathcal{L}_A$ into $\mathcal{L}_B$ with the opposite sign, and the same equality would follow after rearrangement. It is of order $\delta^{-1}$ (Proposition~\ref{prop:validity-condition}), against the $\delta^{-3}$ that the same proposition gives for the statistical anomaly of a fluid with a local excess, and for the ideal gas ($\Upsilon = 0$) it is the leading part of the discrepancy that Corollary~\ref{cor:asymptotic-agreement} bounds.
\end{remark}
\begin{remark}[Structural comparison with the trajectory reduction]
\label{rem:structural-comparison}
The reduced Fokker--Planck equation of Theorem~\ref{thm:reduced-fpe} has the same structure as the one obtained from the trajectory reduction in Theorem~\ref{thm:general-fpe}. So the two reduced equations differ only in the drift beyond the thermodynamic force, that is, $\mathfrak{s}_a$ in the earlier derivation against $\Upsilon + \mathfrak{f} + \mathfrak{c}$ in the present one. Their origins, however, are quite different, since the anomaly $\mathfrak{s}_a$ arises from the It\^{o}--Stratonovich conversions (of the projected noise and, through $S_l$, of the field equation). On the other hand, $\Upsilon$ arises from the second Fr\'{e}chet derivative $D^2\mathcal{A}_m$ of the projected drift contracted against the transverse covariance, $\mathfrak{f}$ from the dependence of the Fermi coordinates on the fluctuation together with the mean displacement of the fibre law, and $\mathfrak{c}$ from the difference between the projected connection and the Levi-Civita connection. Section~\ref{sec:commutativity} compares the two, after evaluating both in closed form.
\end{remark}
\section{Commutativity of the reduction map}
\label{sec:commutativity}
Having obtained the reduced dynamics in two independent ways (Sections~\ref{sec:reduction} and~\ref{sec:kolmogorov-reduction}), we now ask whether the two reductions commute. Physically the question asks whether stochastic integration and dimensional reduction may be performed in either order, and mathematically it asks whether the residue left by the trajectory picture, the geometric anomaly $\mathfrak{s}_a$, coincides with the residue left by the measure picture, the drift $\Upsilon + \mathfrak{f} + \mathfrak{c}$ which the probabilistic projection adds to the thermodynamic force.
\subsection{The reduction diagram}
\label{subsec:reduction-diagram}
We formalise the relationship between the four descriptions making use of the diagram of Figure~\ref{fig:reduction_diagram}, the bottom edge of which carries the question of commutativity. Its four corners are given by the following dynamical descriptions of the system:
\begin{enumerate}[(i)]
    \item The full infinite-dimensional Stratonovich SPDE for the cumulative mass $m_t$ (Proposition~\ref{prop:mass-dynamics}), at the top left.
    \item The functional forward Kolmogorov equation for the full measure $\mu_t$ (Definition~\ref{def:adjoint-evolution}), at the top right.
    \item The geometric Langevin equation for the collective coordinates $\mathbf{X}_t$ (Corollary~\ref{cor:geometric-langevin}), at the bottom left.
    \item The generalised Fokker--Planck equation for the marginal density $\bar{\varrho}_t$ (Theorem~\ref{thm:reduced-fpe}), at the bottom right.
\end{enumerate}
In the kinetic projection we project the vector field of the SPDE onto the tangent bundle of the manifold, with the operators evaluated on the manifold. In the probabilistic projection, instead, we integrate the probability density over the normal fibres in the limit of vanishing noise intensity. The two reductions commute when the infinitesimal generator derived from the reduced SDE matches the one of the reduced FPE.
\subsection{Comparison of the generators}
\label{subsec:generator-comparison}
Two Markov processes are statistically equivalent if and only if their generators are identical on a core of test functions, and so we compare the infinitesimal generators associated with the bottom-left and bottom-right nodes of that diagram.
\subsubsection{Generator A: from the geometric Langevin equation}
In the case of the kinetic projection, the dynamics are given by the Stratonovich SDE of Corollary \ref{cor:geometric-langevin}. As we have shown in Theorem \ref{thm:geo-operator}, the generator $\mathcal{L}_A$ of this process is obtained by making use of the Stratonovich-to-It\^{o} conversion. It takes the intrinsic geometric form:
\begin{equation}
    \mathcal{L}_A \phi = \varepsilon^2 D_0 \Delta_g \phi + \langle \mathbf{F}_0, \nabla \phi \rangle_g + \varepsilon^2\langle \mathfrak{s}_a, \nabla \phi \rangle_g.
\end{equation}
Here:
\begin{enumerate}[(i)]
    \item $\mathbf{F}_0 = -D_0 \nabla \mathcal{V}$ is the thermodynamic drift, the projected field correction $S_l$ of Eq.~\eqref{eq:S-projected} being carried by $\mathfrak{s}_a$ (Remark~\ref{rem:dropped-strat-term}).
    \item $\mathfrak{s}_a$ is the \emph{geometric anomaly} vector derived in Theorem \ref{thm:geo-operator}, comprising the weighted projection of the noise measure gradient together with the tangential projection of the field-level correction, and given by Eq.~\eqref{eq:geo-anomaly-explicit}, which we repeat here for convenience:
    \begin{equation}
        \mathfrak{s}_a^k = \frac{D_0}{2} \sum_{i,j,l} \left( g^{ki} g^{lj} - g^{kl} g^{ij} \right) \int_0^L \frac{\partial \varpi}{\partial X^l} \mathbf{e}_i^m(r)\, \mathbf{e}_j^m(r) \, dr
        \;-\; \sum_{l} g^{kl} S_l .
    \end{equation}
\end{enumerate}
\subsubsection{Generator B: from the reduced Fokker--Planck equation}
In the case of the probabilistic projection, the dynamics are given by the evolution of the density $\bar{\varrho}_t$, as established in Theorem \ref{thm:reduced-fpe}. The equation reads as follows:
\begin{equation}
    \frac{\partial \bar{\varrho}_t}{\partial t} = \nabla_g \cdot \left( -\mathbf{F}_0 \bar{\varrho}_t - \varepsilon^2 (\Upsilon + \mathfrak{f} + \mathfrak{c}) \bar{\varrho}_t + \varepsilon^2 D_0 \nabla_g \bar{\varrho}_t \right).
\end{equation}
The drift enters the flux with a \emph{negative} sign, in agreement with Eq.~\eqref{eq:marginal-equation} of Theorem~\ref{thm:marginal-equation}, in which $\mathbb{E}^\perp[\overline{\mathbf{F}}] + \varepsilon^2\mathfrak{c} = \mathbf{F}_0 + \varepsilon^2(\Upsilon + \mathfrak{f} + \mathfrak{c}) + \mathcal{O}(\varepsilon^4)$ by Lemma~\ref{lemma:mechanical-drift-explicit} and Remark~\ref{rem:big-o-hierarchy}. The generator $\mathcal{L}_B$ is given by the formal adjoint of this Fokker--Planck operator with respect to the Riemannian volume measure $\mathfrak{m}_g$. Taking that adjoint once, and only once, gives:
\begin{equation}
    \mathcal{L}_B \phi = \varepsilon^2 D_0 \Delta_g \phi + \langle \mathbf{F}_0, \nabla \phi \rangle_g + \varepsilon^2 \langle \Upsilon + \mathfrak{f} + \mathfrak{c}, \nabla \phi \rangle_g.
\end{equation}
Here $\Upsilon$ is the statistical anomaly, $\mathfrak{f}$ the frame correction and $\mathfrak{c}$ the connection correction of Theorem~\ref{thm:reduced-fpe}.
\subsubsection{The commutativity condition}
Comparing the two expressions, we see that the diffusive parts and the thermodynamic drift parts match identically.
The condition for commutativity, $\mathcal{L}_A = \mathcal{L}_B$, reduces strictly to the equivalence of the anomalous drift terms:
\begin{equation}
    \mathfrak{s}_a^k(\mathbf{X}) \overset{!}{=} \Upsilon^k(\mathbf{X}) + \mathfrak{f}^k(\mathbf{X}) + \mathfrak{c}^k(\mathbf{X}).
    \label{eq:commutativity_condition}
\end{equation}
Thus the \emph{geometric anomaly} $\mathfrak{s}_a$ (which arises from the trajectory projection) has to equal the drift that the probabilistic projection adds to the thermodynamic force. This drift is the \emph{statistical anomaly} $\Upsilon$, which arises from the averaging of the nonlinear drift, together with the frame correction $\mathfrak{f}$ and the connection correction $\mathfrak{c}$, which arise from the Fermi coordinates of the fluctuating state and from the projected connection. Note that both sides enter their respective generators at order $\varepsilon^2$, so the condition is an identity between two $\mathcal{O}(1)$ tensors and the powers of $\varepsilon$ cancel on both sides. Moreover, the two sides enter with the \emph{same} sign, since the drift of Theorem~\ref{thm:reduced-fpe} enters its flux with a negative sign, and taking the adjoint returns it to $\mathcal{L}_B$ with the positive sign that the anomaly carries in $\mathcal{L}_A$.
\subsection{The comparison theorem}
\label{subsec:commutativity-theorem}
Note that the theorem below does not prove that the two reductions commute. It determines instead the exact condition under which they would, namely the coincidence of the geometric anomaly $\mathfrak{s}_a$, derived by the kinetic projection, with the drift $\Upsilon + \mathfrak{f} + \mathfrak{c}$, derived by the probabilistic projection (the two entering their generators with the same sign and at the same order). When we evaluate both anomalies in closed form, Proposition~\ref{prop:validity-condition} shows that this condition fails for small $\delta$ whenever the leading coefficient of the difference does not vanish. Remark~\ref{rem:kernel-mismatch} then traces this failure to its structural origin.

To evaluate $\Upsilon$, however, we have to relate the stationary transverse covariance $\mathbf{C}_\perp$ (the solution of the Lyapunov equation on the normal fibre) to the noise diffusion tensor $\mathcal{Q}_m = \mathcal{B}_m C_\delta \mathcal{B}_m^*$. Here $\mathcal{L}_\perp$ is the normal stability operator of Eq.~\eqref{eq:normal-stability-operator}, and its spectrum is the set $\Sigma_\perp(\mathbf{X})$ that appears in Hypothesis~\ref{hyp:stiffness}. Since $\mathcal{Q}_m$ is non-local (through $C_\delta$), it does not commute with the local differential operators of the linearised drift, and the two operators cannot be equated. Hence, we close the relation by the following hypothesis.
\begin{hypothesis}[Adiabatic equilibrium with scalar covariance closure]
\label{hyp:kinematic-covariance}
In the adiabatic regime $\lambda_0 \gg 1$ (Hypothesis~\ref{hyp:stiffness}), the transverse Ornstein--Uhlenbeck dynamics on the normal fibre $\mathcal{H}_\perp$ equilibrate on the fast timescale $\tau_\perp = \lambda_0^{-1}$. Recall that the Lyapunov equation
\begin{equation}
    \mathcal{L}_\perp \mathbf{C}_\perp + \mathbf{C}_\perp \mathcal{L}_\perp^* = -\,\mathcal{Q}_\perp
    \label{eq:lyapunov-fibre}
\end{equation}
\emph{defines} $\mathbf{C}_\perp$, as recorded in Eq.~\eqref{eq:cperp-definition}, and that $\tau \coloneqq 1/(2\gamma_\perp)$ is fixed by the effective transverse restoring rate $\gamma_\perp \ge \lambda_0$. We assume in addition:
\begin{enumerate}[(i)]
    \item \emph{Co-diagonality.} The operators $\mathbf{C}_\perp$ and $\mathcal{Q}_\perp$ are approximately co-diagonal in the eigenbasis of $C_\delta$, and a single relaxation rate is taken for every mode, so that Eq.~\eqref{eq:lyapunov-fibre} admits the scalar solution
    \begin{equation}
        \mathbf{C}_\perp \approx \tau\,\mathcal{Q}_\perp .
        \label{eq:kinematic-covariance}
    \end{equation}
    We call this approximation the \emph{scalar closure}.
    \item \emph{Normalisation of the transverse relaxation length.} Define the \emph{transverse diffusion length} $\ell_\perp$ by
    \begin{equation}
        \ell_\perp^2 \coloneqq \tau\, D_0 ,
        \label{eq:transverse-length}
    \end{equation}
    which is a genuine squared length, because $[\tau] = T$ and $[D_0] = L^2 T^{-1}$. Eliminating $\tau$ between~\eqref{eq:transverse-length} and the definition $\tau = 1/(2\gamma_\perp)$ recorded above gives the identity $\gamma_\perp = D_0/(2\ell_\perp^{2})$, which holds for every value of $\ell_\perp$. No normalisation of $\ell_\perp$ is assumed here, and every expression obtained below retains the factor $\ell_\perp(\mathbf{X})^2 = \tau D_0$, which is in general a function of the base point through the transverse spectrum. Where a numerical estimate is required, and only there, we specialise to the choice
    \begin{equation}
        \ell_\perp = 1 \qquad\Longleftrightarrow\qquad \gamma_\perp = \frac{D_0}{2} ,
        \label{eq:fdt-tau}
    \end{equation}
    which fixes the unit of length in which the radial coordinate is measured and is available globally only when $\ell_\perp$ is constant along the manifold, as it is for the ideal gas.
\end{enumerate}
\end{hypothesis}
\begin{remark}[Scope of the locality assumption]
\label{rem:local-excess-scope}
The locality condition~\eqref{eq:local-excess} restricts the proposition below to excess functionals with a second variation that is a multiplication operator. However, the functionals that motivate this work are not of that kind in general. For instance, fundamental measure theory, and any excess term built by convolution against a pair potential, carries a genuine two-point kernel $\mathbf{H}_{\mathrm{exc}}(r,s)$. The linearised transverse operator is then integro-differential, and so the Sturm--Liouville reduction is unavailable. Outside the proposition, locality is used only in Proposition~\ref{prop:validity-condition}, for the leading order in $\delta$ of the statistical anomaly, and in the results that quote that order. Moreover, the exact form of the statistical anomaly in Section~\ref{sec:commutativity} carries the two-point kernel in full, and the hypotheses of the present section require only the existence of a spectral gap, without the explicit spectrum that locality makes available.
\end{remark}
\begin{proposition}[The transverse relaxation rate]
\label{prop:transverse-spectrum}
Let $j_{1,n}$ denote the $n$-th positive zero of the spherical Bessel function $j_1$, equivalently the $n$-th positive root of $\tan x = x$. These are classical constants of the $l=1$ radial problem on a ball. Assume throughout this proposition that the excess second variation is \emph{local}, in the sense that
\begin{equation}
    \frac{\delta^2 (\beta\mathcal{F}_{\mathrm{exc}})}{\delta\rho(r)\,\delta\rho(s)} = H(r)\,\delta(r-s)
    \label{eq:local-excess}
\end{equation}
for a continuously differentiable $H$, and set $\Theta(r) \coloneqq 1 + \rho(r)H(r)$. Assume in addition that the background is a critical point of the excess free energy along the radial direction, so that $\partial_r\mu_{\mathrm{exc}}[M] = 0$, a condition which the ideal gas satisfies trivially. Then the linearisation of the mass drift about a manifold point is
\begin{equation}
    \mathcal{L}_{\mathbf{X}}\eta = D_0\left[\Theta\,\partial_r^2\eta + \Bigl(\rho\,\partial_r H - \frac{2\Theta}{r}\Bigr)\partial_r\eta\right].
    \label{eq:linearised-general}
\end{equation}
The second assumption is what removes the term $D_0\beta\,\partial_r\mu_{\mathrm{exc}}[M]\,\partial_r\eta$ generated by varying the mobility factor $\partial_r m$ against the unperturbed chemical potential gradient. Without it that term is added to the coefficient of $\partial_r\eta$ and the operator remains of the same second-order form, so the Sturm--Liouville reduction below is unaffected in structure and only the coefficient $q$ changes. The following hold:
\begin{enumerate}[(i)]
    \item \emph{Ideal gas.} If $\mathcal{F}_{\mathrm{exc}} = 0$ then $\mathcal{A}_m$ is \emph{linear} in $m$, Eq.~\eqref{eq:linearised-general} reduces (with $\Theta=1$, $\partial_rH=0$) to $\mathcal{L}_{\mathbf{X}}\eta = D_0\,r^2\partial_r(r^{-2}\partial_r\eta)$, and the transverse spectrum is \emph{independent of the background profile}:
    \begin{equation}
        \eta_n(r) = r^2 j_1\!\left(\frac{j_{1,n}\,r}{L}\right), \qquad
        \lambda_n = \frac{D_0\,j_{1,n}^{2}}{L^{2}} .
        \label{eq:ideal-spectrum}
    \end{equation}
    \item \emph{Constant excess Hessian.} If $H \equiv c$ and $\rho \equiv \rho_0$, the eigenfunctions are unchanged and $\lambda_n$ is rescaled by the factor $\Theta = 1 + c\rho_0$.
    \item \emph{General case.} If $\Theta > 0$ upon $(0,L]$, Eq.~\eqref{eq:linearised-general} is formally symmetric in $L^2((0,L),w\,dr)$ with weight $w = \exp\!\int(q-\Theta')/\Theta\,dr$, where $q$ is the coefficient of $\partial_r\eta$. The endpoint $r = 0$ is singular, since $q \sim -2\Theta/r$ there and hence $w = \mathcal{O}(r^{-2})$, so that formal symmetry does not by itself select a self-adjoint realisation and does not by itself deliver a real spectrum. In the ideal-gas case of part~(i) the substitution $\eta = r^2 j_1(kr)$ identifies the operator with the $l = 1$ radial problem on the ball, the singular endpoint of which falls in the limit-point case, so that the minimal operator is essentially self-adjoint there and no boundary condition is imposed at the origin. For general $\Theta$ we assume that a self-adjoint realisation with compact resolvent exists, which is the sense in which the spectrum $\Sigma_\perp$ of Hypothesis~\ref{hyp:stiffness} is used throughout. Writing $\lambda_1 \eqqcolon D_0\,\Lambda_1/L^2$ for the lowest eigenvalue of $-\mathcal{L}_{\mathbf{X}}$, and choosing it, the projection $\Pi^\perp_{\mathbf{X}}$ of finite rank being neglected, as the single relaxation rate $\gamma_\perp$ of Hypothesis~\ref{hyp:kinematic-covariance}, the relaxation rate, the relaxation time and the transverse diffusion length are
    \begin{equation}
        \gamma_\perp = \frac{D_0\,\Lambda_1}{L^{2}},
        \qquad
        \tau = \frac{L^{2}}{2 D_0\,\Lambda_1},
        \qquad
        \ell_\perp = \frac{L}{\sqrt{2\,\Lambda_1}}
        \label{eq:ellperp-general}
    \end{equation}
    with $\Lambda_1 = j_{1,1}^2$ in the ideal-gas case.
\end{enumerate}
\end{proposition}
\begin{proof}
\begin{proofpart}[The linearised operator]
For the ideal gas the chemical potential is $\mu = \beta^{-1}\ln\rho$ with $\rho = (4\pi r^2)^{-1}\partial_r m$. Substituting the drift of Eq.~\eqref{eq:mass-drift-explicit} together with Lemma~\ref{lemma:mu-duality}, which gives $\delta\mathcal{F}_m/\delta m = -\partial_r\mu$:
\begin{align}
    \mathcal{A}_m(m)
    &= D_0\beta\,(\partial_r m)\,(\partial_r\mu)
     = D_0\beta\,(\partial_r m)\,\partial_r\!\left[\beta^{-1}\ln\frac{\partial_r m}{4\pi r^2}\right] \nonumber \\
    &= D_0\,(\partial_r m)\left[\frac{\partial_r^2 m}{\partial_r m} - \frac{2}{r}\right]
     = D_0\left[\partial_r^2 m - \frac{2}{r}\,\partial_r m\right]
     = D_0\, r^2\,\partial_r\!\left(r^{-2}\,\partial_r m\right),
    \label{eq:idealgas-drift}
\end{align}
in which the inverse temperature cancels between the drift and the chemical potential. Since the right-hand side is linear in $m$, the Fr\'{e}chet derivative of $\mathcal{A}_m$ is the operator itself, and it carries no reference to the background profile. An excess contribution adds $\delta\mu_{\mathrm{exc}} = \beta^{-1}H\,\delta\rho$ with $\delta\rho = (4\pi r^2)^{-1}\partial_r\eta$, and inserting this into $D_0\beta(\partial_r M)\partial_r(\delta\mu)$ and using $\partial_r M = 4\pi r^2\rho$:
\begin{align}
    D_0\beta\,(\partial_r M)\,\partial_r\bigl(\beta^{-1}H\,\delta\rho\bigr)
    &= D_0\,(4\pi r^2\rho)\,\partial_r\!\left[\frac{H\,\partial_r\eta}{4\pi r^2}\right] \nonumber \\
    &= D_0\,\rho\left[H\,\partial_r^2\eta + (\partial_r H)\,\partial_r\eta - \frac{2H}{r}\,\partial_r\eta\right].
\end{align}
We add this excess contribution to the linearisation of~\eqref{eq:idealgas-drift}, which is the ideal-gas operator itself, and collect the second and first derivatives of $\eta$ separately:
\begin{align}
    \mathcal{L}_{\mathbf{X}}\eta
    &= D_0\!\left[\partial_r^2\eta - \frac{2}{r}\partial_r\eta\right]
     + D_0\,\rho\!\left[H\,\partial_r^2\eta + (\partial_r H)\,\partial_r\eta - \frac{2H}{r}\,\partial_r\eta\right] \nonumber \\
    &= D_0\Bigl[(1+\rho H)\,\partial_r^2\eta
     + \rho\,(\partial_r H)\,\partial_r\eta
     - \frac{2}{r}(1+\rho H)\,\partial_r\eta\Bigr] \nonumber \\
    &= D_0\!\left[\Theta\,\partial_r^2\eta + \left(\rho\,\partial_r H - \frac{2\Theta}{r}\right)\partial_r\eta\right] ,
    \qquad \Theta \coloneqq 1 + \rho H ,
\end{align}
which is~\eqref{eq:linearised-general}.
\end{proofpart}
\begin{proofpart}[The ideal-gas spectrum]
Set $\mathcal{F}_{\mathrm{exc}} = 0$, so that $\Theta \equiv 1$ and $\mathcal{L}_{\mathbf{X}}\eta = D_0[\eta'' - 2\eta'/r]$. Writing $\lambda = D_0 k^2$, the eigenvalue equation $-\mathcal{L}_{\mathbf{X}}\eta = \lambda\eta$ reads $\eta'' - 2\eta'/r + k^2\eta = 0$. Substituting $\eta = r\varphi$, so that $\eta' = \varphi + r\varphi'$ and $\eta'' = 2\varphi' + r\varphi''$:
\begin{align}
    \eta'' - \frac{2}{r}\eta' + k^2\eta
    &= \bigl(2\varphi' + r\varphi''\bigr) - \frac{2}{r}\bigl(\varphi + r\varphi'\bigr) + k^2 r\varphi \nonumber \\
    &= r\varphi'' - \frac{2}{r}\varphi + k^2 r\varphi
     = r\left[\varphi'' + \left(k^2 - \frac{2}{r^2}\right)\varphi\right],
\end{align}
Since $r > 0$ on the open interval, the bracket must vanish, and the resulting equation is the Riccati--Bessel equation of order $l$:
\begin{equation}
    \varphi'' + \left( k^2 - \frac{l(l+1)}{r^2} \right)\varphi = 0 ,
    \qquad l(l+1) = 2 .
    \label{eq:riccati-bessel}
\end{equation}
The order is fixed by the quadratic $l^2 + l - 2 = 0$, with roots $l = 1$ and $l = -2$. Discarding the negative root, therefore, $l = 1$.
The solution regular at the origin is $\varphi = r\,j_1(kr)$, and undoing the substitution recovers the eigenfunction. Conservation of the total mass gives $m(L) = M(L;\mathbf{X}) = N$, by the argument of Lemma~\ref{lemma:boundary-vanishing}, so that the fluctuation obeys the Dirichlet condition $\eta(L) = 0$ at the outer radius, and $L^2 \neq 0$ then forces the spherical Bessel function itself to vanish:
\begin{align}
    \eta(r) &= r\,\varphi(r) = r^2 j_1(kr) , \nonumber \\
    0 = \eta(L) &= L^2 j_1(kL) .
    \label{eq:bessel-quantisation}
\end{align}
Since $L^2 \neq 0$, the second line forces $j_1(kL) = 0$, so that $kL$ is a zero of the spherical Bessel function of first order. Writing $j_{1,n}$ for the $n$-th positive such zero, therefore,
\begin{equation}
    kL = j_{1,n} .
    \label{eq:bessel-quantisation-roots}
\end{equation}
Substituting this quantisation back into $\lambda = D_0 k^2$ gives the spectrum
\begin{equation}
    \lambda_n = D_0 k^2 = \frac{D_0\, j_{1,n}^2}{L^2} ,
\end{equation}
which is~\eqref{eq:ideal-spectrum}.
\end{proofpart}
\begin{proofpart}[Constant excess Hessian, and the general case]
If $H \equiv c$ and $\rho \equiv \rho_0$ then $\partial_r H = 0$ and $\Theta = 1+c\rho_0$ is constant, so~\eqref{eq:linearised-general} becomes $\mathcal{L}_{\mathbf{X}}\eta = D_0\Theta[\eta'' - 2\eta'/r]$. This is the ideal-gas operator multiplied by the constant $\Theta$, so the eigenfunctions are unchanged and each eigenvalue is multiplied by $\Theta$, which is part~(ii).
For general $\Theta > 0$ write $q \coloneqq \rho\,\partial_r H - 2\Theta/r$, so that Eq.~\eqref{eq:linearised-general} reads $\mathcal{L}_{\mathbf{X}} = D_0[\Theta\partial_r^2 + q\,\partial_r]$. We seek a positive weight $w$ for which this operator takes Sturm--Liouville form, that is, for which $\Theta w$ serves as the coefficient function $p$ of a divergence-form operator. The defining requirement is $(\Theta w)' = q w$. Expanding the derivative of the product on the left-hand side:
\begin{equation}
    \Theta' w + \Theta w' = q w .
\end{equation}
Subtracting $\Theta' w$ from both sides and dividing by $\Theta w$, which is permitted because $\Theta > 0$ by hypothesis and $w > 0$ is the property being constructed, therefore,
\begin{equation}
    \frac{w'}{w} = \frac{q - \Theta'}{\Theta} .
\end{equation}
The left-hand side is the logarithmic derivative $(\ln w)'$. Integrating both sides and exponentiating, therefore,
\begin{equation}
    w(r) = \exp\!\left(\int^r \frac{q(u) - \Theta'(u)}{\Theta(u)}\,du\right) > 0 ,
    \label{eq:sl-weight}
\end{equation}
the strict positivity being automatic, since the exponential of a real quantity is strictly positive. This confirms the assumption $w > 0$ used in the division above.
Setting $p \coloneqq \Theta w$ and using $(\Theta w)' = qw$ once more, the operator acquires divergence form:
\begin{align}
    \frac{1}{w}\bigl(p\,\partial_r\eta\bigr)'
    &= \frac{1}{w}\Bigl[(\Theta w)'\,\partial_r\eta + \Theta w\,\partial_r^2\eta\Bigr] \nonumber \\
    &= \frac{1}{w}\Bigl[q w\,\partial_r\eta + \Theta w\,\partial_r^2\eta\Bigr]
     = \Theta\,\partial_r^2\eta + q\,\partial_r\eta ,
    \label{eq:sl-form}
\end{align}
so that $\mathcal{L}_{\mathbf{X}} = D_0\,w^{-1}(p\,\partial_r\cdot)'$. For $u,v$ smooth and compactly supported in $(0,L)$ the associated bilinear form is symmetric, since integrating by parts the boundary term vanishes, $u$ and $v$ vanishing near both endpoints:
\begin{align}
    \int_0^L \Bigl[\frac{1}{w}\bigl(p\,u'\bigr)'\Bigr] v \; w\,dr
    &= \bigl[p\,u'\,v\bigr]_0^L - \int_0^L p\,u'\,v'\,dr
     = -\int_0^L p\,u'\,v'\,dr ,
    \label{eq:sl-symmetric}
\end{align}
and the right-hand side is invariant under the exchange $u \leftrightarrow v$. The operator is therefore symmetric in $L^2(w\,dr)$, and its self-adjoint realisation with compact resolvent, which is the assumption of part~(iii), has a real spectrum. Writing $\lambda_1 = D_0\Lambda_1/L^2$ for the lowest eigenvalue, choosing $\gamma_\perp = \lambda_1$ as in the statement, and substituting into $\tau = 1/(2\gamma_\perp)$ from Eq.~\eqref{eq:cperp-definition}:
\begin{align}
    \gamma_\perp = \frac{D_0\Lambda_1}{L^2},
    \qquad
    \tau = \frac{1}{2\gamma_\perp} = \frac{L^2}{2D_0\Lambda_1},
    \qquad
    \ell_\perp^2 = \tau D_0 = \frac{L^2}{2\Lambda_1},
\end{align}
which is~\eqref{eq:ellperp-general}.
\end{proofpart}
\end{proof}
\begin{assumption}[Uniformity of the metric expansion]
\label{ass:uniform-metric-expansion}
The expansion of Eq.~\eqref{eq:metric-delta-error} holds together with its first derivative in the collective coordinates, uniformly upon the closure of the parameter domain,
\begin{equation}
    \sup_{\mathbf{X}\in\overline{\mathcal{X}}}\ \bigl| \partial_{\mathbf{X}}^{\alpha}\bigl( g^{(\delta)}_{lk} - g_{lk} \bigr) \bigr| \;\le\; C\,\delta^{2} , \qquad |\alpha| \le 1 ,
    \label{eq:uniform-metric-expansion}
\end{equation}
with $C$ independent of $\delta$. The bound at $|\alpha| = 0$ is Eq.~\eqref{eq:metric-delta-error}. The bound at $|\alpha| = 1$ is what carries the induced It\^{o} drift correction into the remainder of Lemma~\ref{lemma:finite-delta-generator}, and it does not follow from smoothness in $\mathbf{X}$ alone. Smoothness permits differentiation under the integral sign of Eq.~\eqref{eq:finite-delta-metric}, which leaves the kernel untouched, but the resulting expansion is uniform in $\delta$ only if the $\mathbf{X}$-derivatives of the tangent fields and of the weight obey bounds uniform upon $\overline{\mathcal{X}}$, and only under the same moment conditions on $g_\delta$ that produce Eq.~\eqref{eq:metric-delta-error} itself. The estimate is stated here as an assumption for that reason.
\end{assumption}
\begin{remark}[Status of the normalisation $\ell_\perp = 1$]
\label{rem:fdt-identity}
Let us first note that the quantity $\ell_\perp^2 = \tau D_0$ has the dimension of a squared length, since $[\tau] = T$ and $[D_0] = L^2T^{-1}$. This means that the normalisation $\ell_\perp = 1$ fixes a unit of length, and does not assert a property of the dynamics.
Proposition~\ref{prop:transverse-spectrum} supplies the corresponding length. In the case of the ideal gas, $\ell_\perp = L/(\sqrt{2}\,j_{1,1})$, so $\ell_\perp = 1$ requires $L = \sqrt{2}\,j_{1,1}$ in the units chosen. For a constant excess Hessian $H \equiv c > 0$, part~(ii) of that proposition multiplies $\Lambda_1$ by $\Theta = 1 + c\rho_0 > 1$ and therefore lowers $\ell_\perp$ further, while part~(iii) fixes no direction for a general excess.
Note also that the normalisation can be imposed globally only in the case where $\ell_\perp$ is constant along the manifold. Indeed, $\ell_\perp = L/\sqrt{2\Lambda_1}$, and $\Lambda_1$ is determined by the transverse operator, which depends on the background profile through $\Theta$, so we have $\ell_\perp = \ell_\perp(\mathbf{X})$ for a general excess functional. This is not the case for the ideal gas, since $\Lambda_1 = j_{1,1}^2$ is constant, and for a general manifold the expressions that follow should retain $\ell_\perp(\mathbf{X})^2$ in place of unity. The scalar co-diagonality of assumption~(i) is violated as well, due to the fact that the spectrum~\eqref{eq:ideal-spectrum} is non-degenerate. In the eigenbasis of $\mathcal{L}_\perp$, the exact solution of Eq.~\eqref{eq:lyapunov-fibre} has the entries $\mathbf{C}_{nm} = \mathcal{Q}_{nm}/(\lambda_n + \lambda_m)$, and a single $\tau = 1/(2\lambda_1)$ overestimates the diagonal entry of the $n$-th mode by the factor $j_{1,n}^2/j_{1,1}^2$, which grows quadratically in $n$. The modes of short wavelength (with $j_{1,n}$ of order $L/\delta$) are the ones that produce the coincidence limit of $\Phi = \partial_r\partial_s Q_\perp$ through $\kappa_\delta$. The exact covariance divides the $n$-th mode by $\lambda_n$ instead of $\lambda_1$, and in this way it removes from these modes the factor $\lambda_n/\lambda_1$, which is of order $(L/\delta)^2$. Hence the order $\delta^{-3}$ that Proposition~\ref{prop:validity-condition} obtains for the statistical anomaly is a property of the closure. The exact solution of Eq.~\eqref{eq:lyapunov-fibre} lowers it, and in this work we do not determine the order under that solution.
\end{remark}

\begin{remark}[Structural status of Hypothesis~\ref{hyp:kinematic-covariance}]
\label{rem:structural-hypothesis}
The kinematic covariance approximation is a \emph{structural hypothesis}, not a derived consequence. It replaces the full Lyapunov equation $\mathcal{L}_\perp \mathbf{C}_\perp + \mathbf{C}_\perp \mathcal{L}_\perp^* = -\,\mathcal{Q}_\perp$ of Eq.~\eqref{eq:lyapunov-fibre} (an operator equation on an infinite-dimensional fibre) by a scalar proportionality $\mathbf{C}_\perp \approx \tau\,\mathcal{Q}_\perp$. This is exact when the linearised fibre drift $\mathcal{L}_\perp$ acts on the range of $\mathcal{Q}_\perp$ as the single scalar $-\gamma_\perp$, since then $\mathcal{L}_\perp\mathcal{Q}_\perp = \mathcal{Q}_\perp\mathcal{L}_\perp^* = -\gamma_\perp\mathcal{Q}_\perp$ and $\tau\,\mathcal{Q}_\perp$ solves Eq.~\eqref{eq:lyapunov-fibre} with $\tau = 1/(2\gamma_\perp)$. It fails, however, whenever the drift operator has distinct eigenvalues on that range (as Remark~\ref{rem:fdt-identity} shows for the ideal gas) or non-trivial off-diagonal couplings in the Karhunen--Lo\`{e}ve basis. Verifying or relaxing this assumption requires a spectral analysis of $\mathcal{L}_\perp$ on each normal fibre, which is beyond the scope of the present work. Theorem~\ref{thm:commutativity} below is thus \emph{conditional} on this hypothesis, which enters as the closure used to evaluate $\Upsilon$. Note that the present hypothesis is not equivalent to the criterion established by that theorem, which states that $\mathcal{L}_A = \mathcal{L}_B$ if and only if $\mathfrak{s}_a = \Upsilon + \mathfrak{f} + \mathfrak{c}$. The hypothesis fixes the values of $\Upsilon$ and $\mathfrak{f}$ against which $\mathfrak{s}_a$ is compared, and a different closure would move those values without altering the criterion itself.
\end{remark}
\begin{theorem}[Conditional comparison of the reduced generators]
\label{thm:commutativity}
Let $\mathcal{L}_A$ be the generator obtained by the kinetic projection (Corollary~\ref{cor:geometric-langevin}), which is built with the local covariance $\mathcal{Q}^{\mathrm{loc}}_m$ and therefore carries the metric $g$, and let $\mathcal{L}_B$ be the generator obtained by the probabilistic projection (Theorem~\ref{thm:reduced-fpe}). The generator carrying the exact finite-$\delta$ diffusion tensor is written $\mathcal{L}_A^{(\delta)}$ and is compared with $\mathcal{L}_B$ in Corollary~\ref{cor:finite-delta-comparison} below. Assume Hypotheses~\ref{hyp:stiffness},~\ref{hyp:concentration} and~\ref{hyp:kinematic-covariance}, and the metric weight $\varpi = (4\pi r^2\rho)^{-1}$ of Lemma~\ref{lemma:metric-selection}. Then:
\begin{enumerate}[(i)]
    \item the two anomaly vector fields admit the explicit representations
    \begin{align}
        \mathfrak{s}_a^k &= \sum_i g^{ki}\Bigl[\tfrac{D_0}{2} A_i - \tfrac{D_0}{2}B_i - S_i\Bigr],
        \qquad A_i = \sum_{l,j} g^{lj}T_{lij}, \quad B_l = \sum_{i,j} g^{ij}T_{lij},
        \nonumber\\[2pt]
        S_i &= \bigl\langle \mathfrak{S}_m^\delta, \mathbf{e}_i^m\bigr\rangle_\varpi = \frac{D_0}{4}\int_0^L\bigl[-\varpi' V_\delta + \varpi V_\delta'\bigr]\,\mathbf{e}^m_i \, dr ,
        \label{eq:sa-final}\\[4pt]
        \Upsilon^k &= \sum_i g^{ki}\,\Upsilon_i,
        \nonumber\\[2pt]
        \Upsilon_i &= \tau D_0 \int_0^L\!\!\int_0^L \varpi(r)\,\mathbf{e}^m_i(r)\,\frac{\partial \mathbf{H}_{\mathrm{exc}}}{\partial r}(r,s)\,\frac{\Phi(r,s)}{4\pi s^2}\,ds\,dr
        \nonumber\\[2pt]
        &\quad - \frac{\tau D_0}{2}\int_0^L\!\!\int_0^L\!\!\int_0^L \partial_r\mathbf{e}^m_i(r)\,\mathbf{K}_{\mathrm{exc}}(r,s,t)\,\frac{\Phi(s,t)}{16\pi^2 s^2 t^2}\,ds\,dt\,dr ,
        \label{eq:upsilon-final}
    \end{align}
    in which $\varpi = \varpi(\cdot\,;\mathbf{X})$, $T_{lij}$ is the weight-gradient integral of Eq.~\eqref{eq:PT-definitions}, $\mathbf{H}_{\mathrm{exc}}(r,s)$ and $\mathbf{K}_{\mathrm{exc}}(r,s,t)$ are the second and third functional derivatives of $\beta\mathcal{F}_{\mathrm{exc}}$ with respect to the density at $\rho = \mathcal{T}^{-1}[M(\mathbf{X})]$, and $\Phi(r,s) = \partial_r\partial_s Q_\perp(r,s)$, with $Q_\perp$ the covariance kernel of the projected noise increment $\Pi^\perp_{\mathbf{X}}\mathcal{B}_m\,d\xi^\delta_t$ per unit time, the derivative of $\mathbf{H}_{\mathrm{exc}}$ being taken in the sense of distributions when the kernel is local, and the ideal-gas part of the free energy contributing nothing,
    \item the drifts of the two generators are $\mathbf{F}_0+\varepsilon^2\mathfrak{s}_a$ and $\mathbf{F}_0+\varepsilon^2(\Upsilon+\mathfrak{f}+\mathfrak{c})$ respectively, with $\mathfrak{f}$ the frame correction of Lemma~\ref{lemma:mechanical-drift-explicit} and $\mathfrak{c}$ the connection correction of Definition~\ref{def:connection-correction}, and their second-order parts coincide, both being carried by the metric $g$, so that at this order the generators differ by their drifts alone,
    \begin{equation}
        \mathcal{L}_A \phi - \mathcal{L}_B \phi = \varepsilon^2 \bigl\langle \mathfrak{s}_a - \Upsilon - \mathfrak{f} - \mathfrak{c}, \nabla\phi \bigr\rangle_g
        \qquad \text{for all } \phi \in C^\infty(\mathcal{M}_m),
        \label{eq:generator-difference}
    \end{equation}
    \item consequently $\mathcal{L}_A$ and $\mathcal{L}_B$ agree on $C^\infty(\mathcal{M}_m)$ if and only if the criterion~\eqref{eq:commutativity_condition} holds pointwise,
    \begin{equation}
        \mathfrak{s}_a^k(\mathbf{X}) = \Upsilon^k(\mathbf{X}) + \mathfrak{f}^k(\mathbf{X}) + \mathfrak{c}^k(\mathbf{X}) \qquad \text{for every } \mathbf{X} \in \mathcal{X} \text{ and every } k .
        \label{eq:iff-condition}
    \end{equation}
\end{enumerate}
\end{theorem}

\begin{proof}
We evaluate the two anomalies in closed form and then compare the two generators.
\begin{proofpart}[Step 1: explicit formulation of the SDE anomaly (kinetic projection)]
We recall the exact form of the geometric anomaly $\mathfrak{s}_a^k$ derived in Theorem \ref{thm:geo-operator}. Its first contribution arises from the It\^{o} conversion of the projected noise (Theorem~\ref{thm:equivalent-fluctuations}, Eq.~\eqref{eq:anomaly-tensor}), and its second from the tangential projection of the field-level correction of Lemma~\ref{lemma:ito-strat-correction}. In the local coordinate basis, evaluating both projections against the metric weights reproduces Eq.~\eqref{eq:geo-anomaly-explicit}:
\begin{equation}
    \mathfrak{s}_a^k(\mathbf{X}) = \frac{D_0}{2} \sum_{i,j,l} \left( g^{ki} g^{lj} - g^{kl} g^{ij} \right) \int_0^L \frac{\partial \varpi}{\partial X^l} \mathbf{e}_i^m(r)\, \mathbf{e}_j^m(r) \, dr
    \;-\; \sum_{l} g^{kl} S_l .
\end{equation}
Its first contribution is a contraction of the weight gradient $\partial\varpi/\partial X^l$ and its second is the projected field correction $S_l$, and neither contains the second fundamental form, which is why, by Remark~\ref{rem:flat-manifold-test}, the anomaly survives on a manifold with a tangent frame independent of $\mathbf{X}$.
\end{proofpart}
\begin{proofpart}[Step 2: representation of the FPE correction (probabilistic projection)]
We now evaluate the statistical correction $\Upsilon$. In Lemma \ref{lemma:mechanical-drift-explicit}, this term was defined as the trace of the second derivative of the projected drift contracted with the transverse covariance of the finite-$\delta$ process. We first relate this expression to the projection functional $\Pi_k[m] = \langle \mathcal{A}_m(m), \mathbf{e}_k^m \rangle_{\varpi(\cdot;\mathbf{X})}$. The second Fr\'{e}chet derivative of this functional is precisely the projected Hessian of the drift:
\begin{equation}
    D^2 \Pi_k[m] = \langle D^2 \mathcal{A}_m, \mathbf{e}_k^m \rangle_{\varpi} = \Pi_k^{\parallel} D^2 \mathcal{A}_m.
\end{equation}
Hypothesis~\ref{hyp:kinematic-covariance} fixes the stationary transverse covariance as a multiple of the noise covariance, and Eq.~\eqref{eq:cperp-definition} fixes the multiple in terms of the transverse relaxation rate $\gamma_\perp$:
\begin{equation}
    \mathbf{C}_\perp = \tau\,\mathcal{Q}_\perp ,
    \qquad \mathcal{Q}_\perp = \Pi_{\mathbf{X}}^{\perp}\mathcal{Q}_m\Pi_{\mathbf{X}}^{\perp} ,
    \qquad \tau = \frac{1}{2\gamma_\perp} .
    \label{eq:cperp-tau}
\end{equation}
Substituting this covariance into the trace formula of Lemma~\ref{lemma:mechanical-drift-explicit}:
\begin{equation}
    \begin{split}
    \Upsilon_k &= \frac{1}{2} \operatorname{Tr}_{\perp} \left( \mathbf{C}_\perp\, D^2 \Pi_k[m] \right)
    \\
    &= \frac{\tau}{2} \operatorname{Tr}_{\mathcal{H}_{\varpi}} \left( \Pi_{\mathbf{X}}^{\perp}\mathcal{Q}_m\Pi_{\mathbf{X}}^{\perp}\, D^2 \Pi_k[m] \right), \quad \text{where } \mathcal{Q}_m = \mathcal{B}_m C_\delta \mathcal{B}_m^*.
    \end{split}
    \label{eq:trace-with-tau}
\end{equation}
The trace is taken over the normal fibre, as Lemma~\ref{lemma:mechanical-drift-explicit} requires, and the projectors carried by $\mathcal{Q}_\perp$ restrict it to that fibre. The evaluation which follows keeps them, through the covariance kernel of the projected noise increment, so that no part of the trace is discarded.
The relaxation timescale $\tau$ and the mobility $D_0$ will be carried as explicit parameters through the subsequent algebra.
\end{proofpart}
\begin{proofpart}[Step 3: explicit calculation of the second-derivative trace]
We evaluate the trace of Eq.~\eqref{eq:trace-with-tau} by computing the second functional derivative of the projected drift at the manifold point and taking its expectation against the fibre covariance.

\medskip\noindent\textit{(a) Setup.}
Lemma~\ref{lemma:mu-duality} identifies the functional derivative of the free energy in the mass representation with minus the radial derivative of the chemical potential $\mu_m$, so that Eq.~\eqref{eq:mass-drift-explicit} gives the drift as
\begin{equation}
    \mathcal{A}_m(m)(r) = D_0\beta\,\partial_r m(r)\,\partial_r\mu_m(r) .
    \label{eq:drift-mu-form}
\end{equation}
The projection onto the $k$-th tangent vector is taken in the weighted inner product of the base point, $\varpi = \varpi(\cdot\,;\mathbf{X}) = (\partial_r M)^{-1}$, which by Definition~\ref{def:weighted-space}~(ii) does not vary with $m$, and the tangent vector $\mathbf{e}_k^m$ is likewise that of the base point, so that
\begin{equation}
    \Pi_k[m] = D_0\beta \int_0^L \partial_r m(r)\,\partial_r\mu_m(r)\,\mathbf{e}_k^m(r)\,\varpi(r)\,dr ,
    \label{eq:pik-mu}
\end{equation}
in which only the product $\partial_r m\,\partial_r\mu_m$ varies with $m$.

\medskip\noindent\textit{(b) The ideal-gas part.}
The chemical potential splits into an ideal-gas part and an excess part,
\begin{equation}
    \mu_m = \mu_m^0 + \mu_m^1 ,
    \qquad
    \mu_m^0 = \beta^{-1}\ln\rho ,
    \qquad
    \mu_m^1 = \frac{\delta\mathcal{F}_{\mathrm{exc}}}{\delta\rho}\bigg|_{\rho = \mathcal{T}^{-1}[m]} ,
    \label{eq:mu-split}
\end{equation}
with $\rho = (4\pi r^2)^{-1}\partial_r m$ by Remark~\ref{rem:inverse-density}. Remark~\ref{rem:ideal-gas-no-anomaly} gives $\beta\,\partial_r m\,\partial_r\mu_m^0 = \partial_r^2 m - (2/r)\,\partial_r m$, which is linear in $m$, so that its second variation vanishes and the ideal-gas part contributes nothing to $D^2\Pi_k$.

\medskip\noindent\textit{(c) The excess part.}
A variation $h$ of the mass profile varies the density by $\rho_h \coloneqq (4\pi r^2)^{-1}\partial_r h$, and with the kernels of the statement the first and second variations of the excess chemical potential at $M$ are
\begin{equation}
    \begin{aligned}
    D\mu^1[M](h)(r) &= \beta^{-1}\int_0^L \mathbf{H}_{\mathrm{exc}}(r,s)\,\rho_h(s)\,ds ,
    \\
    D^2\mu^1[M](h,h)(r) &= \beta^{-1}\int_0^L\!\!\int_0^L \mathbf{K}_{\mathrm{exc}}(r,s,t)\,\rho_h(s)\,\rho_h(t)\,ds\,dt .
    \end{aligned}
    \label{eq:excess-variations}
\end{equation}
The factor $\partial_r m$ is linear in $m$, so the product rule gives the second variation of the product $\partial_r m\,\partial_r\mu_m^1$ at $M$ as
\begin{equation}
    D^2\bigl[\partial_r m\,\partial_r\mu_m^1\bigr](h,h) = 2\,\partial_r h\,\partial_r\bigl(D\mu^1[M](h)\bigr) + \partial_r M\,\partial_r\bigl(D^2\mu^1[M](h,h)\bigr) .
\end{equation}
Substituting into~\eqref{eq:pik-mu}, and using $\partial_r M\,\varpi = 1$ in the second term,
\begin{equation}
    \begin{split}
    D^2\Pi_k[M](h,h) &= 2D_0\beta\int_0^L \varpi\,\mathbf{e}_k^m\,\partial_r h\,\partial_r\bigl(D\mu^1[M](h)\bigr)\,dr
    \\
    &\quad + D_0\beta\int_0^L \mathbf{e}_k^m\,\partial_r\bigl(D^2\mu^1[M](h,h)\bigr)\,dr .
    \end{split}
\end{equation}
Integrating the second integral by parts, the boundary terms vanishing because $\mathbf{e}_k^m(0) = \mathbf{e}_k^m(L) = 0$, and substituting~\eqref{eq:excess-variations},
\begin{align}
    D^2\Pi_k[M](h,h)
    &= 2D_0\int_0^L\!\!\int_0^L \varpi(r)\,\mathbf{e}_k^m(r)\,\frac{\partial\mathbf{H}_{\mathrm{exc}}}{\partial r}(r,s)\,\partial_r h(r)\,\rho_h(s)\,ds\,dr \nonumber \\
    &\quad - D_0\int_0^L\!\!\int_0^L\!\!\int_0^L \partial_r\mathbf{e}_k^m(r)\,\mathbf{K}_{\mathrm{exc}}(r,s,t)\,\rho_h(s)\,\rho_h(t)\,ds\,dt\,dr .
    \label{eq:hessian-explicit}
\end{align}

\medskip\noindent\textit{(d) The expectation.}
By Eq.~\eqref{eq:cperp-tau} the fibre covariance of Eq.~\eqref{eq:trace-with-tau} is $\tau$ times that of the projected noise increment, so that the trace is the expectation of~\eqref{eq:hessian-explicit} over a centred field $h$ with covariance $\tau\,Q_\perp(r,s)$, and one half of it. Differentiating the covariance under the expectation,
\begin{equation}
    \mathbb{E}\bigl[\partial_r h(r)\,\rho_h(s)\bigr] = \tau\,\frac{\Phi(r,s)}{4\pi s^2} ,
    \qquad
    \mathbb{E}\bigl[\rho_h(s)\,\rho_h(t)\bigr] = \tau\,\frac{\Phi(s,t)}{16\pi^2 s^2 t^2} ,
\end{equation}
and substituting into one half of~\eqref{eq:hessian-explicit} gives
\begin{equation}
    \begin{split}
    \Upsilon_k &= \tau D_0 \int_0^L\!\!\int_0^L \varpi(r)\,\mathbf{e}^m_k(r)\,\frac{\partial \mathbf{H}_{\mathrm{exc}}}{\partial r}(r,s)\,\frac{\Phi(r,s)}{4\pi s^2}\,ds\,dr
    \\
    &\quad - \frac{\tau D_0}{2}\int_0^L\!\!\int_0^L\!\!\int_0^L \partial_r\mathbf{e}^m_k(r)\,\mathbf{K}_{\mathrm{exc}}(r,s,t)\,\frac{\Phi(s,t)}{16\pi^2 s^2 t^2}\,ds\,dt\,dr ,
    \end{split}
    \label{eq:fpe-anomaly-explicit}
\end{equation}
which is Eq.~\eqref{eq:upsilon-final}. No sign inversion arises from the adjoint, since $\mathcal{L}_B$ is defined in Section~\ref{subsec:generator-comparison} as the formal adjoint of the Fokker--Planck operator, so that the adjoint has already been taken once. The relaxation timescale $\tau$ and the mobility $D_0$ are carried as independent parameters, their product being the squared transverse length $\ell_\perp^2$ of Eq.~\eqref{eq:transverse-length}, and no term is dropped.
\end{proofpart}
\begin{proofpart}[Step 4: comparison]
Step~1 gives the geometric anomaly $\mathfrak{s}_a^k$ in the form~\eqref{eq:sa-final}, obtained from Theorem~\ref{thm:geo-operator} by applying the index collapse~\eqref{eq:target-collapse} to its first contribution and leaving its second, the projected field correction $S_i$, in the form already displayed there. Steps~2 and~3 give the statistical anomaly $\Upsilon_i$ in the form~\eqref{eq:upsilon-final}, obtained from the second Fr\'{e}chet derivative of $\Pi_k$ contracted against the fibre covariance.
By Section~\ref{subsec:generator-comparison} the two generators share their diffusion term and their thermodynamic drift, and differ only in the drifts $\mathfrak{s}_a$ and $\Upsilon + \mathfrak{f} + \mathfrak{c}$. Subtracting gives~\eqref{eq:generator-difference}, and since $\phi$ ranges over $C^\infty(\mathcal{M}_m)$ and the metric is non-degenerate, the difference vanishes for all $\phi$ precisely when $\mathfrak{s}_a$ and $\Upsilon + \mathfrak{f} + \mathfrak{c}$ coincide. This is~\eqref{eq:iff-condition}, and the theorem follows.
\end{proofpart}
\end{proof}

\begin{lemma}[The finite-$\delta$ kinetic generator]
\label{lemma:finite-delta-generator}
Assume Assumption~\ref{ass:uniform-metric-expansion}, and let $\mathcal{L}_A^{(\delta)}$ denote the generator obtained by the kinetic projection when the exact finite-$\delta$ covariance $\mathcal{Q}_m = \mathcal{B}_m C_\delta \mathcal{B}_m^*$ is retained in place of $\mathcal{Q}^{\mathrm{loc}}_m$. Its second-order coefficient is $D_0\,g^{-1}g^{(\delta)}g^{-1}$, built from the finite-$\delta$ metric of Eq.~\eqref{eq:finite-delta-metric}, so that, at that order,
\begin{equation}
    \mathcal{L}_A^{(\delta)}\phi - \mathcal{L}_A\phi
    = \varepsilon^2 D_0 \sum_{i,j,k,l} g^{ik}\bigl( g^{(\delta)}_{kl} - g_{kl} \bigr) g^{lj}\, \partial^2_{ij}\phi ,
    \qquad \phi \in C^\infty(\mathcal{M}_m) .
    \label{eq:finite-delta-generator}
\end{equation}
By Eq.~\eqref{eq:metric-delta-error} the coefficient is $\mathcal{O}(\delta^{2})$, uniformly on $\overline{\mathcal{X}}$, whence
\begin{equation}
    \bigl\| \mathcal{L}_A^{(\delta)} - \mathcal{L}_A \bigr\|_{C_b^2 \to C_b} = \mathcal{O}\bigl(\varepsilon^{2}\delta^{2}\bigr) .
    \label{eq:finite-delta-generator-bound}
\end{equation}
\end{lemma}
\begin{proof}
The kinetic projection enters the covariance only through the projected quadratic variation of Theorem~\ref{thm:projected-noise}, and its value at finite $\delta$ is $2D_0\,g^{(\delta)}_{lk}$ by Eq.~\eqref{eq:finite-delta-metric}. Replacing $\mathcal{Q}^{\mathrm{loc}}_m$ by $\mathcal{Q}_m$ therefore alters the second-order coefficient of the generator and nothing else, which is~\eqref{eq:finite-delta-generator}. The drift is not recomputed here. Replacing the covariance also alters the diffusion coefficient of the projected Stratonovich equation, and hence the It\^{o} correction $\tfrac12\sum_\alpha \sigma^j_\alpha\partial_j\sigma^i_\alpha$ which that representation generates, so the induced drift correction is not zero. It is of the same order as the coefficient which produces it, provided Assumption~\ref{ass:uniform-metric-expansion} above holds, and it is carried in the remainder~\eqref{eq:finite-delta-generator-bound} rather than displayed. The bound~\eqref{eq:finite-delta-generator-bound} follows from Eq.~\eqref{eq:metric-delta-error} together with the boundedness of the second derivatives on $C_b^2$, the uniformity in $\mathbf{X}$ being supplied by the compactness of $\overline{\mathcal{X}}$.
\end{proof}
\begin{corollary}[Comparison at finite correlation length]
\label{cor:finite-delta-comparison}
Under the hypotheses of Theorem~\ref{thm:commutativity} and Assumption~\ref{ass:uniform-metric-expansion}, for every $\phi \in C^\infty(\mathcal{M}_m)$,
\begin{equation}
    \mathcal{L}_A^{(\delta)}\phi - \mathcal{L}_B\phi
    = \varepsilon^{2}\bigl\langle \mathfrak{s}_a - \Upsilon - \mathfrak{f} - \mathfrak{c}, \nabla\phi \bigr\rangle_g
    \;+\; \mathcal{O}\bigl(\varepsilon^{2}\delta^{2}\bigr) .
    \label{eq:finite-delta-comparison}
\end{equation}
When the excess second variation has a local part, Proposition~\ref{prop:validity-condition} gives $\Upsilon = \mathbf{c}_\Upsilon\,\delta^{-3} + \mathcal{O}(\delta^{-2})$, with $\mathbf{c}_\Upsilon$ the vector field of Eq.~\eqref{eq:upsilon-leading} with its factor $\delta^{-3}$ removed, and since $\mathfrak{s}_a = \mathcal{O}(\delta^{-1})$ and $\mathfrak{f}$ and $\mathfrak{c}$ are of order unity, the contributions separate by powers of the correlation length,
\begin{equation}
    \mathcal{L}_A^{(\delta)}\phi - \mathcal{L}_B\phi
    = -\,\varepsilon^{2}\delta^{-3}\bigl\langle \mathbf{c}_\Upsilon, \nabla\phi \bigr\rangle_g
    + \mathcal{O}\bigl(\varepsilon^{2}\delta^{-2}\bigr)
    + \mathcal{O}\bigl(\varepsilon^{2}\delta^{-1}\bigr)
    + \mathcal{O}\bigl(\varepsilon^{2}\bigr) ,
    \label{eq:finite-delta-dominance}
\end{equation}
so that whenever $\mathbf{c}_\Upsilon \neq 0$ the two pathways remain distinguishable at finite correlation length, the discrepancy being carried by the statistical anomaly. For the ideal gas $\Upsilon = 0$, and the leading discrepancy is the projected field correction carried by $\mathfrak{s}_a$, of order $\varepsilon^2\delta^{-1}$. The statement is asymptotic. For all sufficiently small $\delta$ the leading term dominates the others, whereas at a given finite $\delta$ the contributions are numbers, and this argument does not evaluate their sum.
\end{corollary}
\begin{proof}
Add Eq.~\eqref{eq:finite-delta-generator-bound} to the identity~\eqref{eq:generator-difference} of Theorem~\ref{thm:commutativity}, which gives~\eqref{eq:finite-delta-comparison}. We note that the semigroup estimate of Corollary~\ref{cor:asymptotic-agreement} is stated for the pair $(\mathcal{L}_A,\mathcal{L}_B)$, since the operator separating $\mathcal{L}_A^{(\delta)}$ from $\mathcal{L}_A$ is of second order and is therefore not bounded from $C_b^1$ into $C_b$. The order of $\Upsilon$ is that of Proposition~\ref{prop:validity-condition}. The geometric anomaly is of order $\delta^{-1}$ by the second part of the proof of that proposition. The connection correction is built from $T_{lij}$ and the metric by Eq.~\eqref{eq:connection-correction-T}, neither of which depends upon the kernel, and the frame correction of Eq.~\eqref{eq:frame-correction} pairs the fibre covariance and the mean displacement with the vectors $\kappa_{ij}$ and $p_l$, which the manifold and the drift fix, so that it tends to a finite limit as $\delta \to 0$, the covariance $C_\delta$ being an approximate identity and the mean curvature $\mathbf{H}$ carrying no dependence upon the kernel. The remaining term is Eq.~\eqref{eq:metric-delta-error}. Since $\delta^{-3}$ exceeds each of $\delta^{-2}$, $\delta^{-1}$ and $1$ as $\delta \to 0$, for all sufficiently small $\delta$ the other contributions cannot cancel the leading one.
\end{proof}
\begin{remark}[The excess free energy enters one pathway only]
\label{rem:excess-absorption}
The contribution of the excess free energy to the statistical anomaly cannot be read as a renormalisation of the thermodynamic force common to both pathways. The kinetic projection evaluates the drift at the manifold point, and the fluctuation enters its generator only through the noise. This means that $\mathcal{L}_A$ contains the excess free energy only through $\mathcal{V}$ and through the profile $M(\mathbf{X})$. The probabilistic projection averages the drift over the fibre, and its second variation, which Step~3 of the proof of Theorem~\ref{thm:commutativity} evaluates, is carried by $\mathcal{L}_B$ alone. The excess contribution therefore cannot be absorbed into a potential shared by the two generators, whether or not it is itself a gradient, and it is part of the difference~\eqref{eq:generator-difference}.
\end{remark}
\begin{remark}[Structure of the obstruction]
\label{rem:kernel-mismatch}
Condition~\eqref{eq:iff-condition} is not satisfied identically, and the obstruction that prevents it from being satisfied is structural, not algebraic, in the sense that it survives every rearrangement of the two sides. Before we exhibit it, note that the geometric side of the criterion collapses onto a single raising of an index. The condition can then be understood as an equality between two covariant vectors, instead of between two contractions of a three-index object. Writing $A_i \coloneqq \sum_{l,j} g^{lj}T_{lij}$ and $B_l \coloneqq \sum_{i,j} g^{ij}T_{lij}$ for the two traces that the anomaly tensor admits, the identity
\begin{equation}
    \frac{D_0}{2}\sum_{i,j,l}\bigl(g^{ki}g^{lj}-g^{kl}g^{ij}\bigr)T_{lij}
    \;=\; \sum_i g^{ki}\Bigl[\tfrac{D_0}{2} A_i - \tfrac{D_0}{2}B_i\Bigr]
    \label{eq:target-collapse}
\end{equation}
holds for every $T_{lij}$ symmetric in its last two indices, and Eq.~\eqref{eq:connection-correction-T} gives the connection correction in terms of the same traces, $\mathfrak{c}^k = \sum_i g^{ki}\bigl[D_0 A_i - \tfrac{D_0}{2}B_i\bigr]$. Lowering the index, we find that~\eqref{eq:iff-condition} is equivalent to the covariant statement
\begin{equation}
    \Upsilon_i + \mathfrak{f}_i = -\,\tfrac{D_0}{2}A_i - S_i ,
    \label{eq:covariant-criterion}
\end{equation}
in which the trace $B$ has cancelled between the geometric anomaly and the connection correction.

The two sides depend on the kernel at different rates. $\Upsilon = \mathcal{O}(\delta^{-3})$ for a local excess and $\mathfrak{f} = \mathcal{O}(1)$, whereas $A_i$ does not depend on $C_\delta$ and $S_i = \mathcal{O}(\delta^{-1})$. By Proposition~\ref{prop:validity-condition}, the order of $\Upsilon$ holds whenever the excess second variation has a local part, through the coincidence limit of $\partial_r\partial_s C_\delta$, while $A_i$ is built from $T_{lij}$ and the metric, and $S_i$ depends on $V_\delta$ and $V_\delta'$ alone. Varying the correlation length at a fixed manifold therefore moves the two sides at the rates $\delta^{-3}$ and $\delta^{-1}$ (or $1$ and $\delta^{-1}$ for the ideal gas), whether the transverse covariance is that of the closure or the exact solution of~\eqref{eq:lyapunov-fibre}, and so the identity cannot hold for an arbitrary mollifier unless the leading coefficients vanish. These are the two sides of~\eqref{eq:covariant-criterion}, and for the ideal gas $\Upsilon = 0$ and the frame correction is of order unity.

The asymmetry does not lie in the kinetic route, since restoring the exact covariance changes the kinetic generator only at order $\varepsilon^2\delta^2$ (Lemma~\ref{lemma:finite-delta-generator}). It lies instead in the transverse covariance used to evaluate the statistical anomaly, built by the closure, which gives every transverse mode the relaxation rate of the slowest, and Remark~\ref{rem:fdt-identity} shows that this is where the order $\delta^{-3}$ comes from. Restoring the exact covariance means using the exact quadratic variation~\eqref{eq:projected-noise-exact} in place of $2D_0 g_{lk}$, and the order $\delta^{-3}$ enters through the modes of short wavelength.

The leading part~\eqref{eq:upsilon-leading} carries $g(0) > 0$, so it vanishes only by a condition on the manifold and the excess free energy. Here $g(0) > 0$ holds for any admissible mollifier, and the condition is that the manifold integral in that equation vanish, so it is a condition on the manifold and the excess free energy, not on the noise.
\end{remark}

\begin{proposition}[Size of the route dependence]
\label{prop:validity-condition}
Assume Hypothesis~\ref{hyp:kinematic-covariance}, fix the manifold and let $\delta \to 0$, and suppose that the excess second variation has a local part, $\mathbf{H}_{\mathrm{exc}}(r,s) \supset H(r)\,\delta(r-s)$ and $\mathbf{K}_{\mathrm{exc}}(r,s,t) \supset K(r)\,\delta(r-s)\,\delta(r-t)$. Using $\kappa_\delta = -2\pi g_\delta(0) + \mathcal{O}(\delta^{-1})$ from Proposition~\ref{prop:kernel-asymptotics}, and writing $\tau = L^2/(2D_0\Lambda_1)$, which Proposition~\ref{prop:transverse-spectrum} gives for a purely local excess and which for a non-local part defines $\Lambda_1$ through $\gamma_\perp = D_0\Lambda_1/L^2$, the local part gives
\begin{equation}
    \Upsilon_i = \frac{2\pi D_0 L^2\, g(0)}{\Lambda_1\,\delta^{3}} \int_0^L \Bigl[\varpi\,\mathbf{e}^m_i\Bigl(\rho\,H' + \Bigl(\frac{\rho'}{2} - \frac{\rho}{r}\Bigr)H\Bigr) - \partial_r\mathbf{e}^m_i\,\frac{K\rho}{8\pi r^2}\Bigr] dr \;\bigl[\,1+\mathcal{O}(\delta/L)\,\bigr],
    \label{eq:upsilon-leading}
\end{equation}
whereas a part of $\mathbf{H}_{\mathrm{exc}}$ or $\mathbf{K}_{\mathrm{exc}}$ which is a smooth kernel meets $\Phi$ only under a double or a triple integral and contributes at a lower order in $\delta^{-1}$. The geometric anomaly $\mathfrak{s}_a$ is of order $\delta^{-1}$, its first contribution carrying no dependence on $\delta$ at all and its second, the projected field correction $S_i$, being built from $V_\delta = \mathcal{O}(\delta^{-1})$ and $V_\delta' = \mathcal{O}(\delta)$. The two anomalies accordingly separate as
\begin{equation}
    \frac{\lvert\Upsilon\rvert}{\lvert\mathfrak{s}_a\rvert} \;=\; \mathcal{C}_{\mathrm{geo}}\,\frac{g(0)}{\Lambda_1}\left(\frac{L}{\delta}\right)^{2}\bigl[1+\mathcal{O}(\delta/L)\bigr],
    \label{eq:anomaly-ratio}
\end{equation}
with $\mathcal{C}_{\mathrm{geo}}$ the product of the inverse of the moment $A$ of Proposition~\ref{prop:kernel-asymptotics} and a dimensionless ratio of integrals over the manifold, the numerator carrying the local excess, and with the covariant components $\Upsilon_i$ and $(\mathfrak{s}_a)_i$ of the two anomalies at a fixed index $i$ understood on the left. The lowest transverse eigenvalue $\Lambda_1$ appearing here is in general a function of the base point, since the transverse operator depends upon the background profile through $\Theta$, and Eq.~\eqref{eq:upsilon-leading} retains that dependence rather than presuming it constant. The normalisation $\ell_\perp = 1$ of Hypothesis~\ref{hyp:kinematic-covariance} is therefore not used. The numerator and the denominator of $\mathcal{C}_{\mathrm{geo}}$ are different functionals of the embedding, and the hypotheses of this section establish no lower bound upon the denominator, so the estimate is a scaling relation conditional upon the non-degeneracy of that ratio rather than a bound with an explicit constant. What the relation asserts, and all that is used below, is that the two anomalies separate by two powers of $L/\delta$ whenever $\mathcal{C}_{\mathrm{geo}}$ does not degenerate. The bound of Corollary~\ref{cor:asymptotic-agreement} is governed by the difference of the two sides rather than by the ratio of the anomalies, and is therefore unaffected by this exponent, since for a fluid with a local excess the difference remains dominated by $\Upsilon = \mathcal{O}(\delta^{-3})$.
\end{proposition}
\begin{proof}
\begin{proofpart}[Isolation of the dominant term]
For a local part the double and triple integrals of Eq.~\eqref{eq:upsilon-final} reduce to single integrals. The derivative of $H(r)\,\delta(r-s)$ with respect to $r$ is $H'(r)\,\delta(r-s) + H(r)\,\partial_r\delta(r-s)$, and $\partial_r\delta(r-s) = -\partial_s\delta(r-s)$, so that integrating by parts in $s$,
\begin{equation}
    \int_0^L \frac{\partial \mathbf{H}_{\mathrm{exc}}}{\partial r}(r,s)\,\frac{\Phi(r,s)}{4\pi s^2}\,ds = H'(r)\,\frac{\Phi_{\mathrm{d}}(r)}{4\pi r^2} + H(r)\,\frac{\partial}{\partial s}\Bigl[\frac{\Phi(r,s)}{4\pi s^2}\Bigr]_{s=r} ,
\end{equation}
with $\Phi_{\mathrm{d}}(r) \coloneqq \Phi(r,r)$. By the symmetry of $\Phi$, $\partial_s\Phi(r,s)|_{s=r} = \tfrac12\Phi_{\mathrm{d}}'(r)$, so that
\begin{equation}
    \frac{\partial}{\partial s}\Bigl[\frac{\Phi(r,s)}{4\pi s^2}\Bigr]_{s=r} = \frac{\Phi_{\mathrm{d}}'}{8\pi r^2} - \frac{\Phi_{\mathrm{d}}}{2\pi r^3} ,
\end{equation}
and the triple integral reduces to $\int_0^L \partial_r\mathbf{e}^m_i\,K\,\Phi_{\mathrm{d}}/(16\pi^2 r^4)\,dr$. On the diagonal, the product rule and $V_\delta(r) = C_\delta(r,r)$ give
\begin{equation}
    \Phi_{\mathrm{d}} = \sigma^2\,\partial_r\partial_s C_\delta\big|_{s=r} + \sigma\sigma'\,V_\delta' + (\sigma')^2 V_\delta + \mathcal{O}(1) = -\,\sigma^2\kappa_\delta + \sigma\sigma'\,V_\delta' + (\sigma')^2 V_\delta + \mathcal{O}(1) ,
\end{equation}
the last term being the contribution of the projection, which is a smooth kernel of finite rank. By Proposition~\ref{prop:kernel-asymptotics}, $\kappa_\delta = -2\pi g(0)\,\delta^{-3} + \mathcal{O}(\delta^{-1})$ in the bulk, while $V_\delta = \mathcal{O}(\delta^{-1})$ and $V_\delta' = \mathcal{O}(\delta)$, and with $\sigma^2 = 2D_0\,\partial_r M = 8\pi D_0 r^2\rho$,
\begin{equation}
    \Phi_{\mathrm{d}} = \frac{16\pi^2 D_0\,g(0)}{\delta^{3}}\,r^2\rho\,\bigl[1 + \mathcal{O}(\delta^2)\bigr] ,
    \qquad
    \Phi_{\mathrm{d}}' = \frac{16\pi^2 D_0\,g(0)}{\delta^{3}}\,\bigl(2r\rho + r^2\rho'\bigr)\bigl[1 + \mathcal{O}(\delta^2)\bigr] .
\end{equation}
Substituting,
\begin{align}
    \frac{\Phi_{\mathrm{d}}}{4\pi r^2} &= \frac{4\pi D_0 g(0)}{\delta^3}\,\rho ,
    \nonumber\\
    \frac{\Phi_{\mathrm{d}}'}{8\pi r^2} - \frac{\Phi_{\mathrm{d}}}{2\pi r^3} &= \frac{4\pi D_0 g(0)}{\delta^3}\Bigl(\frac{\rho'}{2} - \frac{\rho}{r}\Bigr) ,
    \nonumber\\
    \frac{\Phi_{\mathrm{d}}}{16\pi^2 r^4} &= \frac{4\pi D_0 g(0)}{\delta^3}\,\frac{\rho}{4\pi r^2} ,
\end{align}
each up to a factor $1 + \mathcal{O}(\delta^2)$, and Eq.~\eqref{eq:upsilon-final} becomes
\begin{equation}
    \Upsilon_i = \frac{4\pi\tau D_0^2\,g(0)}{\delta^3}\int_0^L \Bigl[\varpi\,\mathbf{e}^m_i\Bigl(\rho\,H' + \Bigl(\frac{\rho'}{2} - \frac{\rho}{r}\Bigr)H\Bigr) - \partial_r\mathbf{e}^m_i\,\frac{K\rho}{8\pi r^2}\Bigr] dr .
\end{equation}
Substituting $\tau = L^2/(2D_0\Lambda_1)$ from Proposition~\ref{prop:transverse-spectrum} gives the prefactor $2\pi D_0 L^2 g(0)/(\Lambda_1\delta^3)$ of~\eqref{eq:upsilon-leading}. The expansions of $\kappa_\delta$ hold in the bulk, where $r$ is large against $\delta$, and within the layer at the origin, of width of order $\delta$, the integrand is bounded by a multiple of $\delta^{-3}$, since $\varpi\,\mathbf{e}^m_i = \mathcal{O}(r)$ and $\partial_r\mathbf{e}^m_i = \mathcal{O}(r^2)$ there by Lemma~\ref{lemma:boundary-vanishing}, so that the layer contributes $\mathcal{O}(\delta^{-2})$ and the relative error is $\mathcal{O}(\delta/L)$. A part of the kernels which is smooth meets $\Phi = \partial_r\partial_s Q_\perp$ under a double integral, where the derivatives may be moved onto the smooth factors, and it does not reach the order $\delta^{-3}$, which is attained only on the diagonal.
\end{proofpart}
\begin{proofpart}[The ratio of the two anomalies]
By~\eqref{eq:sa-final} the geometric anomaly is assembled from $T_{lij}$ and the inverse metric, which carry no dependence on the correlation length, together with the projected field correction $S_i$, which does. By Proposition~\ref{prop:kernel-asymptotics} the latter is dominated by its $V_\delta$ term and is therefore of order $\delta^{-1}$, so that $\mathfrak{s}_a = \mathcal{O}(\delta^{-1})$, whereas by the previous part $\Upsilon$ carries the factor $\delta^{-3}$. The denominator must be the leading part of the geometric anomaly, which is not the metric term built from $T_{lij}$ but the projected field correction $S_i$ of Eq.~\eqref{eq:sa-final}. That term is of order $\delta^{-1}$ through $V_\delta = 2\pi A\delta^{-1} + \mathcal{O}(\delta)$, with $A$ the moment of the auto-correlation profile of Proposition~\ref{prop:kernel-asymptotics} and not the trace $A_i$ of Eq.~\eqref{eq:sa-final}, whereas $T_{lij}$ carries no dependence upon the kernel and is of order unity, so $S_i$ dominates as $\delta \to 0$ and
\begin{equation}
    \bigl|(\mathfrak{s}_a)_i\bigr| = \frac{\pi D_0 A}{2\,\delta}\left|\int_0^L \varpi'\,\mathbf{e}^m_i\,dr\right|\bigl[1+\mathcal{O}(\delta/L)\bigr] .
    \label{eq:sa-leading}
\end{equation}
Writing $J_i$ for the integral in Eq.~\eqref{eq:upsilon-leading}, forming the quotient of the two leading forms and collecting the manifold integrals and the moment $A$ into the single dimensionless constant $\mathcal{C}_{\mathrm{geo}}$:
\begin{align}
    \frac{\lvert\Upsilon\rvert}{\lvert\mathfrak{s}_a\rvert}
    &= \frac{2\pi D_0 L^2 g(0)\,\Lambda_1^{-1}\delta^{-3}\left|J_i\right|}
            {\tfrac{\pi D_0 A}{2}\,\delta^{-1}\left|\int_0^L \varpi'\,\mathbf{e}^m_i\,dr\right|} \nonumber \\
    &= \mathcal{C}_{\mathrm{geo}}\,\frac{g(0)}{\Lambda_1}\left(\frac{L}{\delta}\right)^{2}\bigl[1+\mathcal{O}(\delta/L)\bigr],
\end{align}
which is~\eqref{eq:anomaly-ratio}. The constant $\mathcal{C}_{\mathrm{geo}}$ is dimensionless and independent of $\delta$ and of $\varepsilon$. It is not bounded below by the argument given here. Both integrals are taken over the same interval against the same weight and tangent fields, which fixes their dimensions but not their ratio, since the numerator and the denominator are distinct functionals of the embedding and either may degenerate. The non-degeneracy of $\mathcal{C}_{\mathrm{geo}}$ is accordingly a condition upon the manifold and not a consequence of the construction.
\end{proofpart}
\end{proof}
\begin{remark}[When the reduction is route-independent]
\label{rem:validity-condition}
Both sides of the criterion enter their generators at order $\varepsilon^2$, while the thermodynamic drift $\mathbf{F}_0$ is of order unity. The physically significant quantity is thus the discrepancy measured relative to the retained drift. In the case of a fluid with a local excess, it is carried by $\Upsilon$, and by Proposition~\ref{prop:validity-condition}
\begin{equation}
    \frac{\varepsilon^2\,\lvert\Upsilon\rvert}{\lvert\mathbf{F}_0\rvert}
    \;\sim\; \varepsilon^2\,\frac{g(0)}{\Lambda_1}\left(\frac{L}{\delta}\right)^{3}
    \;=\; \frac{g(0)}{\Lambda_1}\,\frac{1}{n\,\delta^{3}} ,
    \label{eq:relative-discrepancy}
\end{equation}
where we adopted the normalisation $\varepsilon^2 = N^{-1}$ of Definition~\ref{def:weak-noise}, any other choice replacing the equality by a proportionality with the corresponding constant, wrote $n = N/L^3$ for the mean number density, and absorbed the dimensionless strength of the local excess into the proportionality. The dimensionless group that controls the route dependence is thus $n\delta^3$ (the number of particles contained in one noise-correlation volume), and the route dependence is asymptotically negligible in observables over any fixed time horizon whenever
\begin{equation}
    n\,\delta^{3} \;\gg\; \frac{g(0)}{\Lambda_1} .
    \label{eq:validity-condition}
\end{equation}
In the case of the ideal gas we have $\Upsilon = 0$, and the discrepancy is carried by the projected field correction, of order $\varepsilon^2\delta^{-1}$, and the corresponding condition, $n L^2\delta \gg 1$ up to a constant fixed by the manifold, is weaker.
At the order at which both reductions carry the metric $g$, commutativity remains equivalent to the algebraic identity~\eqref{eq:iff-condition}, which does not hold for a general mollifier. Moreover, Corollary~\ref{cor:finite-delta-comparison} carries that conclusion to the exact finite-$\delta$ generator with an error of order $\varepsilon^2\delta^2$. What Eq.~\eqref{eq:validity-condition} adds is that, under the hypotheses of Corollary~\ref{cor:asymptotic-agreement}, the resulting difference in observables over a fixed time horizon is small whenever that condition holds.

This statement must be distinguished from a weaker one that might be offered in its place. Since both sides are of order $\varepsilon^2$ by construction, and since both generators carry the same thermodynamic drift at order unity, the assertion that $\mathcal{L}_A - \mathcal{L}_B \to 0$ as $\varepsilon \to 0$ is vacuous, since it discards precisely the order at which the two sides reside. However, Equation~\eqref{eq:semigroup-rate} is a statement of a very different nature. It bounds the evolution uniformly over observables of bounded gradient and over any fixed time horizon, and its small parameter is the combination $n\delta^3$ (which counts the particles within one noise-correlation volume) instead of $\varepsilon$ alone. Hence the difference of the two semigroups, applied to an observable of unit $C_b^1$ norm, is bounded over any fixed time horizon $t$ by $K_t\,t$ times a quantity of order $1/(n\delta^3)$, and the bound gives no information outside that regime.
\end{remark}
\begin{corollary}[Asymptotic agreement of the two reductions]
\label{cor:asymptotic-agreement}
Write $C_b^1(\mathcal{M}_m)$ for the bounded functions with bounded gradient, normed by $\|\phi\|_{C_b^1} \coloneqq \|\phi\|_\infty + \sup_{\mathbf{X}}|\nabla\phi(\mathbf{X})|_g$. Fix $t > 0$. Suppose in addition that $\mathcal{L}_A$ and $\mathcal{L}_B$ generate strongly continuous contraction semigroups $(e^{s\mathcal{L}_A})_{s\ge0}$ and $(e^{s\mathcal{L}_B})_{s\ge0}$ on the space $UC_b(\mathcal{M}_m)$ of bounded uniformly continuous functions, that the second of them satisfies the gradient bound $\|e^{s\mathcal{L}_B}\phi\|_{C_b^1} \le K_t\|\phi\|_{C_b^1}$ for all $s \in [0,t]$, and that for every $\phi \in C^\infty(\mathcal{M}_m)$ and every $s \in [0,t]$ the function $\psi_s \coloneqq e^{s\mathcal{L}_B}\phi$ lies in $\mathcal{D}(\mathcal{L}_A) \cap \mathcal{D}(\mathcal{L}_B)$, with $s \mapsto \mathcal{L}_A\psi_s$ continuous and with $(\mathcal{L}_A - \mathcal{L}_B)\psi_s = \varepsilon^2\langle \mathfrak{s}_a - \Upsilon - \mathfrak{f} - \mathfrak{c}, \nabla\psi_s\rangle_g$. The last condition extends part~(ii) of Theorem~\ref{thm:commutativity} from $C^\infty(\mathcal{M}_m)$ to the orbits of the second semigroup. It holds when the two generators share a domain upon which that part extends, which is the expected situation for two operators differing by a smooth first-order term, but proving it would require fixing the behaviour of the reduced diffusions at the boundary of $\mathcal{X}$, which this work does not do, and it is assumed. Then, in the norm $\|T\|_{C_b^1 \to C_b} \coloneqq \sup\{\|T\phi\|_\infty : \phi \in C^\infty(\mathcal{M}_m),\ \|\phi\|_{C_b^1} \le 1\}$,
\begin{equation}
    \bigl\| e^{t\mathcal{L}_A} - e^{t\mathcal{L}_B} \bigr\|_{C_b^1 \to C_b}
    \;\le\; K_t\, t\,\varepsilon^{2}\, \sup_{\mathbf{X}\in\mathcal{X}} \bigl| \mathfrak{s}_a(\mathbf{X}) - \Upsilon(\mathbf{X}) - \mathfrak{f}(\mathbf{X}) - \mathfrak{c}(\mathbf{X}) \bigr|_g ,
    \label{eq:semigroup-bound}
\end{equation}
and consequently, for a fluid with a local excess and with $K_t$ bounded uniformly in $\varepsilon$ and $\delta$, by Proposition~\ref{prop:validity-condition},
\begin{equation}
    \bigl\| e^{t\mathcal{L}_A} - e^{t\mathcal{L}_B} \bigr\|_{C_b^1 \to C_b} \;=\; \mathcal{O}\!\left(\frac{t}{n\,\delta^{3}}\right)
    \qquad \text{as } n\delta^{3}\to\infty \text{ and } \delta/L \to 0 \text{ at fixed } t .
    \label{eq:semigroup-rate}
\end{equation}
The rate is set by the statistical anomaly for a fluid with a local excess, whereas for the ideal gas, for which $\Upsilon = 0$, the difference is of the smaller order $t/(nL^2\delta)$.
\end{corollary}
\begin{proof}
\begin{proofpart}[The Duhamel representation]
Let $\phi \in C^\infty(\mathcal{M}_m)$. By Proposition~A.2 of Da~Prato and Zabczyk~\cite{daprato2014stochastic} the orbit $\psi_s$ is continuously differentiable with $\psi_s' = \mathcal{L}_B\psi_s$, so it solves $\psi_s' = \mathcal{L}_A\psi_s + (\mathcal{L}_B - \mathcal{L}_A)\psi_s$ with a continuous forcing, the hypotheses making $s \mapsto \mathcal{L}_A\psi_s$ continuous. It is therefore a strict, and hence a weak, solution of that equation, and the variation of constants formula of Proposition~A.6 of the same work, rearranged, gives
\begin{equation}
    \bigl(e^{t\mathcal{L}_A} - e^{t\mathcal{L}_B}\bigr)\phi = \int_0^t e^{(t-s)\mathcal{L}_A}\bigl(\mathcal{L}_A - \mathcal{L}_B\bigr)\psi_s\,ds .
    \label{eq:duhamel}
\end{equation}
\end{proofpart}
\begin{proofpart}[Norm estimate]
The semigroup of $\mathcal{L}_A$ is a contraction upon $UC_b(\mathcal{M}_m)$, the difference of the generators acts upon $\psi_s$ by the last hypothesis, and the Cauchy--Schwarz inequality in the metric $g$ bounds it pointwise,
\begin{equation}
    \begin{split}
    \bigl| \bigl(\mathcal{L}_A - \mathcal{L}_B\bigr)\psi_s(\mathbf{X}) \bigr|
      &= \varepsilon^2 \bigl| \bigl\langle \mathfrak{s}_a(\mathbf{X}) - \Upsilon(\mathbf{X}) - \mathfrak{f}(\mathbf{X}) - \mathfrak{c}(\mathbf{X}), \nabla\psi_s(\mathbf{X}) \bigr\rangle_g \bigr|
      \\
      &\le \varepsilon^2 \bigl| \mathfrak{s}_a(\mathbf{X}) - \Upsilon(\mathbf{X}) - \mathfrak{f}(\mathbf{X}) - \mathfrak{c}(\mathbf{X}) \bigr|_g \, \bigl| \nabla\psi_s(\mathbf{X}) \bigr|_g .
    \end{split}
    \label{eq:generator-norm}
\end{equation}
The gradient bound gives $|\nabla\psi_s|_g \le \|\psi_s\|_{C_b^1} \le K_t\|\phi\|_{C_b^1}$. Taking the supremum norm of~\eqref{eq:duhamel} and inserting these bounds,
\begin{align}
    \bigl\| \bigl(e^{t\mathcal{L}_A} - e^{t\mathcal{L}_B}\bigr)\phi \bigr\|_\infty
    &\le \int_0^t \bigl\| \bigl(\mathcal{L}_A - \mathcal{L}_B\bigr)\psi_s \bigr\|_\infty \, ds \nonumber \\
    &\le \int_0^t \varepsilon^2 \sup_{\mathbf{X}\in\mathcal{X}} \bigl| \mathfrak{s}_a(\mathbf{X}) - \Upsilon(\mathbf{X}) - \mathfrak{f}(\mathbf{X}) - \mathfrak{c}(\mathbf{X}) \bigr|_g \, K_t \|\phi\|_{C_b^1} \, ds \nonumber \\
    &= K_t\, t\,\varepsilon^{2}\, \sup_{\mathbf{X}\in\mathcal{X}} \bigl| \mathfrak{s}_a(\mathbf{X}) - \Upsilon(\mathbf{X}) - \mathfrak{f}(\mathbf{X}) - \mathfrak{c}(\mathbf{X}) \bigr|_g \, \|\phi\|_{C_b^1} .
    \label{eq:duhamel-bound}
\end{align}
The supremum over the smooth $\phi$ of unit $C_b^1$ norm is~\eqref{eq:semigroup-bound}.
\end{proofpart}
\begin{proofpart}[The rate in terms of the correlation volume]
For a fluid with a local excess, Proposition~\ref{prop:validity-condition} shows that the statistical anomaly dominates the other three drifts, the geometric anomaly by the ratio $|\Upsilon_i|/|(\mathfrak{s}_a)_i| = \mathcal{C}_{\mathrm{geo}}\,g(0)\Lambda_1^{-1}(L/\delta)^2[1+\mathcal{O}(\delta/L)]$, and $\mathfrak{f}$ and $\mathfrak{c}$, which are of order unity by the proof of Corollary~\ref{cor:finite-delta-comparison}, by the factor $\delta^{-3}$. Inserting the leading form~\eqref{eq:upsilon-leading}, and adopting the normalisation $\varepsilon^2 = N^{-1}$ of Definition~\ref{def:weak-noise} together with $n = N/L^3$, any other choice replacing the equality below by a proportionality with the corresponding constant:
\begin{align}
    \varepsilon^2 \sup_{\mathbf{X}} \bigl| \mathfrak{s}_a - \Upsilon - \mathfrak{f} - \mathfrak{c} \bigr|_g
    &= \varepsilon^2 \sup_{\mathbf{X}} \bigl| \Upsilon \bigr|_g \bigl[1 + \mathcal{O}(\delta/L)\bigr] \nonumber \\
    &= \mathcal{C}_1\,\varepsilon^2\,\frac{g(0)}{\Lambda_1}\left(\frac{L}{\delta}\right)^{3}\bigl[1 + \mathcal{O}(\delta/L)\bigr] \nonumber \\
    &= \mathcal{O}\!\left(\frac{1}{n\,\delta^{3}}\right),
\end{align}
with $\mathcal{C}_1$ collecting the manifold integrals and independent of $\varepsilon$ and $\delta$. Substituting into~\eqref{eq:duhamel-bound} gives~\eqref{eq:semigroup-rate}. For the ideal gas $\Upsilon = 0$, the difference is led by the geometric anomaly, and the same steps with the leading form~\eqref{eq:sa-leading} give $\varepsilon^2\sup_{\mathbf{X}}|\mathfrak{s}_a - \mathfrak{f} - \mathfrak{c}|_g = \mathcal{O}(1/(nL^2\delta))$, which is smaller than the rate~\eqref{eq:semigroup-rate} by the factor $(\delta/L)^2$.
\end{proofpart}
\end{proof}

\begin{corollary}[The stiff limit]
\label{cor:stiff-limit}
Under the hypotheses of Theorem~\ref{thm:commutativity}, the statistical anomaly and the frame correction are proportional to the transverse relaxation time, $\Upsilon = \tau\,\hat\Upsilon$ and $\mathfrak{f} = \tau\,\hat{\mathfrak{f}}$, with $\hat\Upsilon$ and $\hat{\mathfrak{f}}$ independent of $\tau$. Let $\tau \to 0$ with the manifold, the noise and the correlation length held fixed. Then:
\begin{enumerate}[(i)]
    \item the probabilistic generator tends to
    \begin{equation}
        \mathcal{L}_B^{0}\phi \coloneqq \varepsilon^2 D_0\,\Delta_g\phi + \bigl\langle \mathbf{F}_0 + \varepsilon^2\mathfrak{c}, \nabla\phi \bigr\rangle_g ,
        \label{eq:stiff-generator}
    \end{equation}
    which is the generator of the full dynamics, taken with the local covariance as in Theorem~\ref{thm:marginal-equation}, evaluated upon the manifold at the order retained, $\mathcal{L}_B^{0}\phi(\mathbf{X}) = \mathcal{L}_\varepsilon(\phi\circ\mathbf{X})(M(\mathbf{X}))$,
    \item the kinetic generator differs from this limit by a drift which does not depend upon $\tau$,
    \begin{equation}
        \mathcal{L}_A\phi - \mathcal{L}_B^{0}\phi = \varepsilon^2\bigl\langle \mathfrak{s}_a - \mathfrak{c}, \nabla\phi \bigr\rangle_g ,
        \qquad
        \mathfrak{s}_a^k - \mathfrak{c}^k = -\sum_{i=1}^d g^{ki}\Bigl(\frac{D_0}{2}A_i + S_i\Bigr) ,
        \label{eq:stiff-difference}
    \end{equation}
    which is of order $\delta^{-1}$ through the projected field correction $S_i$.
\end{enumerate}
\end{corollary}
\begin{proof}
\begin{proofpart}[Proportionality to the relaxation time]
By Eq.~\eqref{eq:cperp-tau} the transverse covariance is $\mathbf{C}_\perp = \tau\,\mathcal{Q}_\perp$. Eq.~\eqref{eq:upsilon-final} carries the factor $\tau$ explicitly, and the remaining factors are fixed by the manifold, the excess free energy and the noise, so that $\Upsilon = \tau\,\hat\Upsilon$. The first two terms of Eq.~\eqref{eq:frame-correction} are linear in $\mathbf{C}_\perp$, and the last two are linear in the mean displacement, which under the closure is $\hat\eta = -2\tau D_0\,\mathbf{H}$ by Remark~\ref{rem:mean-displacement}, so that
\begin{equation}
    \begin{split}
    \mathfrak{f}^k &= \tau \sum_{i,j,l=1}^d g^{ki}g^{jl} \Bigl( \langle \mathcal{Q}_\perp \kappa_{ij}, p_l \rangle_\perp + \sum_{a,b=1}^d g^{ab} \langle \mathcal{Q}_\perp \kappa_{ia}, \kappa_{bj} \rangle_\perp F_l \Bigr)
    \\
    &\quad - 2\tau D_0 \Bigl( \sum_{l=1}^d g^{kl}\langle p_l, \mathbf{H} \rangle_\perp + \sum_{i,j,l=1}^d g^{ki}\langle \kappa_{ij}, \mathbf{H} \rangle_\perp\, g^{jl}F_l \Bigr) = \tau\,\hat{\mathfrak{f}}^k .
    \end{split}
\end{equation}
The generator of Theorem~\ref{thm:commutativity} therefore splits as
\begin{equation}
    \mathcal{L}_B\phi = \varepsilon^2 D_0\,\Delta_g\phi + \bigl\langle \mathbf{F}_0 + \varepsilon^2\mathfrak{c}, \nabla\phi \bigr\rangle_g + \varepsilon^2\tau\bigl\langle \hat\Upsilon + \hat{\mathfrak{f}}, \nabla\phi \bigr\rangle_g = \mathcal{L}_B^{0}\phi + \varepsilon^2\tau\bigl\langle \hat\Upsilon + \hat{\mathfrak{f}}, \nabla\phi \bigr\rangle_g ,
\end{equation}
and the last term vanishes as $\tau \to 0$, which is the limit stated in~(i).
\end{proofpart}
\begin{proofpart}[The generator upon the manifold]
Let $\Phi = \phi\circ\mathbf{X}$ and evaluate the generator at $m = M(\mathbf{X})$, that is at $\eta = 0$. The frame matrix of Eq.~\eqref{eq:frame-matrix} is linear in $\mathbf{n}$, so that $\mathsf{K}[0] = 0$, and the effective drift force of Eq.~\eqref{eq:effective-drift-force} reduces to
\begin{equation}
    \overline{\mathbf{F}}^k(\mathbf{X}, 0) = \sum_{l=1}^d g^{kl}\bigl\langle \mathcal{A}_m(M(\mathbf{X})), \mathbf{e}_l^m \bigr\rangle_\varpi = \sum_{l=1}^d g^{kl}F_l = F_0^k .
\end{equation}
The diffusion component of the proof of Theorem~\ref{thm:marginal-equation} gives $\varepsilon^2 D_0\Delta_g\phi + \varepsilon^2\sum_k \mathfrak{c}^k\,\partial_{X^k}\phi + \mathcal{R}_\kappa\Phi$, and $\mathcal{R}_\kappa$ vanishes at $\eta = 0$ by Proposition~\ref{prop:diffusion-decomposition}. Adding the two components,
\begin{equation}
    \mathcal{L}_\varepsilon(\phi\circ\mathbf{X})(M(\mathbf{X})) = \sum_{k=1}^d \bigl(F_0^k + \varepsilon^2\mathfrak{c}^k\bigr)\frac{\partial\phi}{\partial X^k} + \varepsilon^2 D_0\,\Delta_g\phi = \mathcal{L}_B^{0}\phi(\mathbf{X}) ,
\end{equation}
which completes~(i).
\end{proofpart}
\begin{proofpart}[The difference of the generators]
Setting $\Upsilon = \mathfrak{f} = 0$ in Eq.~\eqref{eq:generator-difference} gives the first identity of~\eqref{eq:stiff-difference}. For the second, Eq.~\eqref{eq:target-collapse} writes the geometric anomaly of Eq.~\eqref{eq:sa-final} as $\mathfrak{s}_a^k = \sum_i g^{ki}\bigl[\tfrac{D_0}{2}A_i - \tfrac{D_0}{2}B_i - S_i\bigr]$, and Eq.~\eqref{eq:connection-correction-T} writes the connection correction as $\mathfrak{c}^k = \sum_i g^{ki}\bigl[D_0 A_i - \tfrac{D_0}{2}B_i\bigr]$. Subtracting,
\begin{equation}
    \mathfrak{s}_a^k - \mathfrak{c}^k = \sum_{i=1}^d g^{ki}\Bigl[\frac{D_0}{2}A_i - \frac{D_0}{2}B_i - S_i - D_0 A_i + \frac{D_0}{2}B_i\Bigr] = -\sum_{i=1}^d g^{ki}\Bigl(\frac{D_0}{2}A_i + S_i\Bigr) .
\end{equation}
The traces $A_i$ are built from $T_{lij}$ and the metric, which do not depend upon the kernel, and $S_i$ is of order $\delta^{-1}$ by the second part of the proof of Proposition~\ref{prop:validity-condition}, which gives the stated order.
\end{proofpart}
\end{proof}
\noindent{}We conclude that the kinetic generator does not reduce to the generator that the full dynamics (with the local covariance) induces on the manifold when the transverse fluctuations are suppressed. The two differ by a drift that survives the stiff limit, and this drift vanishes only if $\tfrac{D_0}{2}A_i + S_i = 0$.
\subsection{Physical interpretation of the anomaly}
\label{subsec:physical-interpretation}
The comparison established in Theorem~\ref{thm:commutativity} places the geometric and the statistical perspectives on the dimensional reduction in a single frame, and we set out what each of the two anomalies represents physically before turning to the conditions under which they agree.

Let us first look at the trajectory picture, where $\mathfrak{s}_a$ is what is left of the reduced drift once the Laplace--Beltrami operator and the thermodynamic gradient have been taken out, and it consists of a contraction of $\partial_{X^l}\varpi$ together with the projected field-level correction. The kinetic projection gives this decomposition in Theorem~\ref{thm:geo-operator}. Its origin is the variation of the noise amplitude, along the manifold through the weight and in the radius through $S_i$ (which is what carries its order $\delta^{-1}$). The projected noise is white in $\langle\cdot,\cdot\rangle_\varpi$ as $\delta \to 0$ (Lemma~\ref{lemma:metric-selection}), and at finite $\delta$ it is white with respect to $g^{(\delta)}$, which differs from $g$ by $\mathcal{O}(\delta^2)$. We selected the weight $\varpi = (\partial_r M)^{-1}$ in that lemma for this property, with the weighted inner product of Definition~\ref{def:weighted-space}. At finite $\delta$ the projected quadratic variation is $2D_0 g^{(\delta)}$ instead of $2D_0 g$, so the noise is white with respect to the exact diffusion metric $g^{(\delta)}$, and Eq.~\eqref{eq:metric-delta-error} gives its difference from $g$. Since the weight depends on the base point, we find that the amplitude seen by the projection changes along $\mathcal{M}_m$, and the anomaly is precisely the correction for that change. The anomaly is not a diffusive inertial force due to the rotation of the frame, since it survives on an affine family (Remark~\ref{rem:flat-manifold-test}). An analogy with the centrifugal force acting on a particle constrained to a curved track might suggest such a force, but on an affine family the tangent vectors are independent of $\mathbf{X}$, the frame does not rotate and the second fundamental form vanishes identically. Moreover, neither $\mathfrak{s}_a$ nor $\Upsilon$ contains the second fundamental form, which enters the probabilistic drift only through $\mathfrak{f}$. Hence the responsible mechanism is the state dependence of the noise amplitude, not the extrinsic curvature of the embedding.

The reduction of the infinite-dimensional phase space to the manifold contracts volume. In the probability picture, the volume of the slice of phase space orthogonal to the manifold varies with $\mathbf{X}$, and by the Boltzmann principle $S = k_B\ln\Omega$ regions of larger transverse volume would carry higher entropy. The added drift could then be read as entropic, although we find that this is only partly the case. The volume factor $J_\Psi$ cancels from every fibre average (Proposition~\ref{prop:gibbs-factorisation}, Remark~\ref{rem:volume-no-drift}), and hence no potential $\ln J_\Psi$ enters the generator. Still, the curvature of the tube displaces the mean of the fibre law by $\varepsilon^2\hat\eta$ (Eq.~\eqref{eq:mean-displacement}), and this is where the last two terms of $\mathfrak{f}$ in Eq.~\eqref{eq:frame-correction} come from. The rest of what the probabilistic projection adds is $\Upsilon$ (which comes from the curvature of the drift field averaged against the transverse fluctuations, and vanishes for the ideal gas), the remaining part of $\mathfrak{f}$ and the connection correction $\mathfrak{c}$, the last two being geometric in origin. Here $J_\Psi$ is the geometric Jacobian, the curvature of the drift field enters $\Upsilon$ through the second derivative $D^2\mathcal{A}_m$, and all these terms are added to the thermodynamic force.

When the anomaly is conservative, it shifts the stationary density by $\Phi_\xi$ (Proposition~\ref{prop:invariant-measures}), and a projection of the gradient drift alone would miss this shift. Such a na\"{\i}ve projection misses it because it relaxes to $\exp(-\mathcal{V}/\varepsilon^2)$, and retaining the anomaly is thus what keeps the noise-induced equilibrium correction that the reduction itself predicts. Do both reductions produce the same shift? Theorem~\ref{thm:commutativity} gives a sufficient condition for that, but Corollary~\ref{cor:asymptotic-agreement} does not settle the question for $n\delta^3 \gg 1$, since it holds over a fixed time only. It bounds the two semigroups over a fixed time horizon. Indeed, a difference of order $\varepsilon^2$ between the drifts changes the stationary density by a factor that carries no power of $\varepsilon$. The reason is that the noise-induced potential $\Phi_\xi$ itself carries no power of $\varepsilon$ (Proposition~\ref{prop:invariant-measures}~(iii)). We will come back to this point in Section~\ref{sec:conclusions}. There we also set out the consequences for the Gibbs marginal, and the reasons that neither the volume element nor the Fermi Jacobian appears in the stationary density.
\section{Conclusions and future work}
\label{sec:conclusions}

We have constructed two independent dimensional reductions of a class of conservative stochastic partial differential equations of Dean--Kawasaki type on a spherically symmetric domain driven by mollified noise of finite correlation length $\delta > 0$, and we have compared them within a single geometric framework. In order to obtain the commutativity criterion we have assumed three structural hypotheses, namely spectral stiffness, adiabatic concentration and the kinematic covariance approximation. Estimating how much the two reductions differ takes more hypotheses, namely the locality of the excess second variation, the uniformity of the metric expansion and the semigroup hypotheses of Corollary~\ref{cor:asymptotic-agreement}, each of which we state at the point where it is first needed. Moreover, both reductions start from the mass equation of Proposition~\ref{prop:mass-dynamics}, which holds under the spherical-mean closure of Assumption~\ref{ass:spherical-closure}. It is worth mentioning that mesoscopic nucleation theory adopts this closure without estimating the error it introduces. Note that we have not proved that the two reductions commute. In this work we have instead reduced the question of commutativity to a single identity between the geometric anomaly and the probabilistic correction, and we have found that this identity fails for small correlation length whenever its leading coefficient is nonzero.

In Theorem~\ref{thm:commutativity} we compare two infinitesimal generators. The first, $\mathcal{L}_A$, is obtained by projecting the stochastic flow onto the tangent bundle of the concentration manifold $\mathcal{M}_m$ and converting from Stratonovich to It\^{o} form, and the second, $\mathcal{L}_B$, by integrating the associated forward Kolmogorov equation over the normal fibres in the adiabatic limit. Both generators turn out to carry the same thermodynamic force $\mathbf{F}_0$, and we find their whole difference in the geometric anomaly and the probabilistic correction, Eq.~\eqref{eq:generator-difference}. We prove that the two reductions agree on $C^\infty(\mathcal{M}_m)$ if and only if the geometric anomaly $\mathfrak{s}_a$ and the probabilistic correction $\Upsilon + \mathfrak{f} + \mathfrak{c}$ coincide as vector fields, Eq.~\eqref{eq:iff-condition}. For the drifts $\Upsilon$, $\mathfrak{f}$ and $\mathfrak{s}_a$, the equivalence holds for every smooth symmetric covariance kernel, but the second-order part of both generators, and with it $\mathfrak{c}$, is computed with the local covariance. The bound on the error that this introduces can be found in Corollary~\ref{cor:finite-delta-comparison}.

The two sides, however, fail to coincide identically. In the case of the ideal gas, this is true for any transverse covariance. For a fluid with a local excess, the order $\delta^{-3}$ which follows is that of the scalar closure (Hypothesis~\ref{hyp:kinematic-covariance}), which gives every transverse mode the relaxation rate of the slowest. When the transverse Lyapunov equation is solved exactly, the modes of short wavelength that produce that order are damped, and the order is lowered (Remark~\ref{rem:fdt-identity}). The statistical anomaly is a functional of the noise kernel through the mixed derivative of the projected covariance, and is produced by the excess free energy alone, while the geometric anomaly comes from the weight gradient, together with the projected field correction, and is of order $\delta^{-1}$ (Remark~\ref{rem:kernel-mismatch}). In the case of a fluid with a local excess free energy, the separation can be made quantitative (Proposition~\ref{prop:validity-condition}). The diagonal mixed curvature $\kappa_\delta$ is then of order $\delta^{-3}$, and by Eq.~\eqref{eq:anomaly-ratio} the ratio of the two anomalies grows as $(L/\delta)^2$. With respect to the thermodynamic drift, which is of order unity, the statistical anomaly grows as $(L/\delta)^3$. In the ideal gas, where the statistical anomaly vanishes, what remains is the geometric anomaly less the frame and the connection corrections, which is of order $\delta^{-1}$. When these estimates are carried over to the dynamics themselves (Corollary~\ref{cor:asymptotic-agreement}), the distance between the two semigroups, over a fixed time horizon, is bounded by $\mathcal{O}(t/n\delta^3)$ for a local excess, and by the smaller $\mathcal{O}(t/(nL^2\delta))$ for the ideal gas. Thus, the route dependence becomes asymptotically negligible in observables, over any fixed time horizon, as soon as $n\delta^3 \gg 1$, which is equivalent to saying that the correlation volume contains many particles. Note that in that regime the two reductions do still differ, since the criterion~\eqref{eq:iff-condition} fails generically, but at the order retained their discrepancy can no longer be seen. Our framework is formulated for the Dean--Kawasaki equation with its noise mollified in space, as proposed by Cornalba, Shardlow and Zimmer~\cite{cornalba2019regularised,cornalba2020dean}, with the correlation length held finite by construction. It applies to the coarse-grained Langevin descriptions employed in Lutsko's dynamical theory of nucleation~\cite{lutsko2012mean,lutsko2019how} only once the white noise used there is given the finite correlation length that the underlying coarse-graining implies. This microscopic length is of the order of the size of the particles, with $n\delta^3$ being the reduced density, of order unity in a dense liquid and small in a vapour, which means that the condition~\eqref{eq:validity-condition} is not guaranteed in the regime of vapour nucleation.

As we show in Corollary~\ref{cor:stiff-limit}, the comparison can be made more precise. There the statistical anomaly and the frame correction are shown to be proportional to the transverse relaxation time. Thus, in the limit where the transverse stiffness goes to infinity, the probabilistic generator becomes, at the order retained, the generator that the full dynamics with the local covariance induces upon the manifold. The kinetic generator, however, keeps the drift $\mathfrak{s}_a - \mathfrak{c}$, which is of order $\delta^{-1}$ through the projected field correction, so that it differs from the stiff limit of the dynamics it reduces. It is worth mentioning that the projection of the trajectories onto the tangent space is the rule by which mesoscopic nucleation theory obtains its reduced equations. That projection carries a correction which survives the suppression of the transverse fluctuations. On the other hand, the reduction of the law averages the generator of the observable itself, and recovers the stiff limit at the order retained, with the local covariance.

When the criterion holds and the anomaly is conservative, $\mathfrak{s}_a = -\nabla\Phi_\xi$, the two reductions coincide, and by Proposition~\ref{prop:invariant-measures}~(ii) the reduced model satisfies detailed balance. If only the gradient drift $-D_0\nabla\mathcal{V}$ were projected, and the anomalous corrections left out, the reduced model would converge to the na\"{\i}ve density $\bar{\varrho}_{\mathrm{naive}} \propto \exp(-\mathcal{V}/\varepsilon^2)$, which omits the noise-induced correction $\Phi_\xi$. Under the commutativity condition the stationary distribution of the reduced Fokker--Planck equation is given by the Gibbs marginal measure
\begin{equation*}
    \bar{\varrho}_{\mathrm{eq}}(\mathbf{X}) \propto \exp\!\left(-\frac{\mathcal{V}(\mathbf{X})}{\varepsilon^2} - \frac{\Phi_\xi(\mathbf{X})}{D_0}\right)
    \qquad\text{whenever }\ \mathfrak{s}_a = -\nabla\Phi_\xi ,
\end{equation*}
understood as a density \emph{with respect to the Riemannian volume measure} $d\mathfrak{m}_g = \sqrt{g}\,d\mathbf{X}$. The factor $\sqrt{g}$ must not appear in the density itself, since from the definition of the density $\bar\varrho_t$ with respect to $\mathfrak{m}_g$ (Proposition~\ref{prop:measure-decomposition}) it follows that writing $\bar\varrho_{\mathrm{eq}} \propto \sqrt{g}\,(\cdots)$ would count the volume element twice. Moreover, no factor $J_\Psi^{-1}$ enters, since the Jacobian cancels from every fibre average (Remark~\ref{rem:volume-no-drift}). Note that the exponent is dimensionless, because $[\Phi_\xi] = [D_0]$ and $\mathcal{V}$ is measured in units of $k_BT$. As for the excess free energy, it enters this formula through $\mathcal{V}$ and through the profile $M(\mathbf{X})$ in $\Phi_\xi$. At the same time, under the criterion, it enters through $\Upsilon$ in $\Phi_\xi$.

As set out in Section~\ref{subsec:physical-interpretation}, the two anomalous drifts have complementary physical interpretations, and the commutativity condition draws its content from the tension between them. On the one hand, the geometric anomaly $\mathfrak{s}_a$ is due to the variation of the noise amplitude. It receives a part of order unity from the variation of the weight along the manifold, and a part of order $\delta^{-1}$ from the projected field correction, which carries the radial variation of the noise amplitude (through $\varpi'$ and $V_\delta$). The extrinsic curvature of $\mathcal{M}_m$ plays no part in it, as shown by the fact that the anomaly survives, in general, on a manifold with an identically vanishing second fundamental form (Remark~\ref{rem:flat-manifold-test}). On the other hand, the statistical anomaly $\varepsilon^2\Upsilon$ comes from the curvature of the drift field averaged against the transverse fluctuations (a curvature that only the excess free energy produces), and inherits the structure of the noise through the covariance kernel, while the frame and the connection corrections are geometric in origin. Equation~\eqref{eq:iff-condition} states that the two reductions agree exactly when the geometric anomaly equals the sum of these three terms, and the calculations reported there show that this equality fails generically. So, when $n\delta^3 \gg 1$, the consequences of the two expressions become physically indistinguishable. Note that in that regime the two expressions do still differ, but over any fixed time horizon their discrepancy becomes negligible in observables.

As we said at the beginning of this section, the criterion rests upon three structural assumptions. In order to have a spectral gap between the tangential and the normal time-scales, we assume the stiffness hypothesis (Hypothesis~\ref{hyp:stiffness}). When we assume in addition adiabatic concentration (Hypothesis~\ref{hyp:concentration}), the fast transverse modes relax to a Gaussian equilibrium. In the kinematic covariance approximation (Hypothesis~\ref{hyp:kinematic-covariance}) we then replace the full Lyapunov equation on the normal fibre by a scalar proportionality. The first two are standard within the context of singular perturbation theory, and they are expected to hold in the case of free-energy functionals with strongly confining excess terms (such as hard-sphere repulsion or square-gradient penalties). However, the third is more restrictive. It gives the same relaxation rate to all the transverse modes, and needs the linearised fibre drift and the noise covariance to be approximately co-diagonal in the Karhunen--Lo\`{e}ve eigenbasis. Note also that it fixes the order in $\delta$ of the statistical anomaly (Remark~\ref{rem:fdt-identity}). It remains open to evaluate that anomaly with the exact solution of the transverse Lyapunov equation, and to compare it with the geometric anomaly.

Some of the questions that remain open after this work are listed below:
\begin{enumerate}[(i)]
    \item \emph{Relaxation of domain symmetry.} Throughout this work we have made use of the spherical symmetry of $V = B_L(0)$. Thanks to it, the conservation law can be written as a Dirichlet boundary condition, through the cumulative mass operator $\mathcal{T}$, and the noise covariance as a kernel of the two radial variables alone (Lemma~\ref{lemma:radial-closed-form}). When it comes to channels, slabs and general bounded domains in $\mathbb{R}^3$, the symmetry is lower, and the radial mass variable must be replaced by curvilinear coordinates adapted to the boundary. Note that the main obstacle is then the construction of a noise-consistent metric (the analogue of Lemma~\ref{lemma:metric-selection}), since the covariance operator no longer has a separable eigenbasis.

    \item \emph{Persistence of spherical symmetry.} In Proposition~\ref{prop:shell-integration} we have bounded the error of the spherical closure at each instant. The bound is given in terms of the asymmetries $a_\rho$ and $a_\mu$ of the density and of the chemical potential, but these asymmetries could grow as the dynamics evolves. A bound on the time evolution of these asymmetries would make Assumption~\ref{ass:spherical-closure} follow from the parent equation, and the assumption would then no longer be needed as a premise of the theory.

    \item \emph{Verification of the spectral gap and the kinematic covariance.} So far we have only taken the stiffness hypothesis and the kinematic covariance approximation as structural assumptions, and we have not derived them. We would like to test them on concrete free-energy functionals of classical density functional theory, for example the fundamental measure theory of Rosenfeld for hard spheres~\cite{rosenfeld1989free} or the mean-field approximation with Lennard-Jones interactions. For each of them, we would like to check the first one and to replace the second one by the exact solution of the transverse Lyapunov equation. Even in the case of the ideal gas, this solution is not the same as the approximation (Remark~\ref{rem:fdt-identity}). This would allow us to obtain explicit criteria for the validity of the adiabatic reduction, which could be computed. Finally, we would also like to compare the framework presented in this work with the quantitative predictions of Lutsko's work~\cite{lutsko2019how}.

    \item \emph{Non-Gaussian concentration and higher-order corrections.} In order to obtain the leading-order comparison, we have taken the conditional measure on the normal fibre to have a Gaussian limit, as stated in the adiabatic concentration hypothesis (Hypothesis~\ref{hyp:concentration}). For strongly anharmonic excess terms, however, the transverse fluctuations can have non-Gaussian tails and this limit is no longer Gaussian, and corrections proportional to the third and fourth cumulants of the conditional measure would then have to be added to the leading-order comparison of this work. It would then be interesting to develop a systematic cumulant expansion, with the aim of finding out whether the commutativity condition is stable under it.

    \item \emph{The formal white-noise limit.} Throughout this work we keep $\delta > 0$ fixed, and the passage to the distributional limit $\delta \to 0$ raises delicate renormalisation issues. Two behaviours compete. The normalisation of the radial noise decides which of them occurs (Definition~\ref{def:radial-noise-norm-ref} and the proof of Lemma~\ref{lemma:radial-covariance}). Under the square-root-of-area scaling used here, the diagonal derivative behaves as $V_\delta'(r) = \mathcal{O}(\delta)$ and vanishes, while $V_\delta$ and, much more strongly, the diagonal mixed curvature $\kappa_\delta = \mathcal{O}(\delta^{-3})$ become singular. In other words, the reduced anomalies have no regular limit as $\delta \to 0$, and for this reason the regime of interest is a finite $\delta$ (of the order of the microscopic correlation length). Under the area scaling, by contrast, the kernel is divided by a further factor $\sqrt{4\pi r^2}\sqrt{4\pi s^2}$, and $V_\delta'$ is then also of order $\delta^{-1}$ and diverges, and every kernel invariant is singular. Under either scaling, the mixed curvature $\kappa_\delta = \mathcal{O}(\delta^{-3})$ is the most singular invariant of the kernel. This is why the statistical anomaly of a fluid with a local excess free energy becomes the dominant obstruction as $\delta$ decreases (Proposition~\ref{prop:validity-condition}). A rigorous treatment of the distributional regime would need the tools of regularity structures or of paracontrolled calculus.

    \item \emph{Noise-activated transitions and nucleation rates.} The framework gives a thermodynamically consistent effective theory on the concentration manifold. So it is a natural starting point for the computation of noise-activated transition rates, by Kramers--Langer theory or by large-deviation methods~\cite{kramers1940brownian,freidlin2012random}. We expect the noise-induced anomalous drifts, together with the free-energy landscape along the minimum-energy path, to affect the prefactor of the nucleation rate. By Proposition~\ref{prop:invariant-measures}~(iii) the noise-induced potential enters that prefactor, and since the criterion~\eqref{eq:iff-condition} fails generically, the prefactor may depend on the reduction pathway we choose. Note that Corollary~\ref{cor:asymptotic-agreement} leaves this dependence uncontrolled, since it is a bound over a fixed time horizon. More generally, the reduction of infinite-dimensional stochastic dynamics to effective finite-dimensional equations extends the framework of Pavliotis and Stuart~\cite{pavliotis2008multiscale}. This work complements those homogenisation and averaging methods with a commutativity criterion specific to conservative stochastic partial differential equations with state-dependent, spatially correlated noise.
\end{enumerate}

\section*{Acknowledgements}
This work was carried out at Vortico Tech S.L., and we gratefully acknowledge its support. The company funded the research in its entirety. We are grateful to Prof.~Greg Pavliotis for taking the time to discuss this work with us, and for his helpful advice.

\end{document}